\documentclass{amsart}

\usepackage{mathrsfs}
\usepackage{enumerate}
\usepackage[matrix,arrow,graph,curve,frame]{xy}

\theoremstyle{plain}      
    \newtheorem{theorem}{Theorem}[section]
    \newtheorem{proposition}[theorem]{Proposition}
    \newtheorem{lemma}[theorem]{Lemma}
    \newtheorem{corollary}[theorem]{Corollary}
\theoremstyle{definition}
    \newtheorem{definition}[theorem]{Definition}
    
\theoremstyle{remark}

\allowdisplaybreaks[1]
\newcommand{\boxwidth}{78ex}

\newcommand{\C}{\ensuremath{\mathscr{C}}}
\newcommand{\Q}{\ensuremath{\mathcal{Q}}}
\newcommand{\V}{\ensuremath{\mathscr{V}}}
\newcommand{\Comod}{\ensuremath{\mathbf{Comod}}}
\newcommand{\Bicomod}{\ensuremath{\mathbf{Bicomod}}}
\renewcommand{\o}{\ensuremath{\circ}}
\newcommand{\ox}{\ensuremath{\otimes}}
\newcommand{\oxC}{\ensuremath{\ox_C}}

\newcommand{\ob}{\ensuremath{\mathrm{ob}}}
\newcommand{\ra}{\ensuremath{\xymatrix@1@C=16pt{\ar[r]&}}}
\newcommand{\dra}{\ensuremath{\xymatrix@1@C=20pt{\ar[r]&}}}
\newcommand{\mra}{\ensuremath{\xymatrix@1@C=20pt{\ar[r] |<(0.4){\object@{|}}&}}}

\begin{document}
\title{Weak Hopf monoids in braided monoidal categories}
\author{Craig Pastro}
\author{Ross Street}
\thanks{The first author gratefully acknowledges support of an international
Macquarie University Research Scholarship while the second gratefully
acknowledges support of the Australian Research Council Discovery Grant
DP0771252.}
\address{Department of Mathematics \\
         Macquarie University \\
         New South Wales 2109 Australia}
\email{\{craig,street\}@maths.mq.edu.au}
\date{\today}

\begin{abstract}
We develop the theory of weak bimonoids in braided monoidal categories and
show them to be quantum categories in a certain sense. Weak Hopf monoids are
shown to be quantum groupoids. Each separable Frobenius monoid $R$ leads to
a weak Hopf monoid $R \ox R$.
\end{abstract}
\maketitle

\tableofcontents

%=========================================================================%
\section{Introduction}
%=========================================================================%

Weak Hopf algebras were introduced by B\"ohm, Nill, and Szlach\'anyi
in a series of papers~\cite{BS,N,Sz,BNS}. They are generalizations of Hopf
algebras and were proposed as an alternative to weak quasi-Hopf algebras. A
weak bialgebra is both an associative algebra and a coassociative coalgebra,
but instead of requiring that the multiplication and unit morphism are
coalgebra morphisms (or equivalently that the comultiplication and the
counit are algebra morphisms) other ``weakened'' axioms are imposed. The
multiplication is still required
to be comultiplicative (equivalently, the comultiplication is still required
to be multiplicative), but the counit is no longer required to be an algebra
morphism and the unit is no longer required to be a coalgebra morphism.
Instead, these requirements are replaced by weakened versions (see
equations~(v) and~(w) below). As the name suggests, any bialgebra satisfies
these weakened axioms and is therefore a weak bialgebra.

Given a weak bialgebra $A$ one may define source and target morphisms $s,t:A
\ra A$ whose images $s(A)$ and $t(A)$ are called the ``source and target
(counital) subalgebras''. It has been shown by Nill~\cite{N} that Hayashi's
face algebras~\cite{H} are special cases of weak bialgebras for which the,
say, target subalgebra is commutative.

A weak Hopf algebra is a weak bialgebra $H$ equipped with an antipode
$\nu:H \ra H$ satisfying the axioms\footnote{There may be some discrepancy
with what we call the source and target morphisms and what exists in the
literature. This arises from our convention of taking multiplication in the
groupoid algebra to be $f \cdot g = g \o f$ (whenever $g \o f$ is defined).}
\[
    \mu(\nu \ox 1)\delta = t, \qquad
    \mu(1 \ox \nu)\delta = s, \qquad \text{and} \qquad 
    \mu_3(\nu \ox 1 \ox \nu)\delta_3 = \nu,
\]
where $\mu_3 = \mu(\mu \ox 1)$ and $\delta_3 = (\delta \ox 1)\delta$.
Again, any Hopf algebra satisfies these weakened axioms and so is a weak
Hopf algebra. Also in~\cite{N} Nill has shown that the (finite dimensional)
generalized Kac algebras of Yamanouchi~\cite{Y} are examples of weak Hopf
algebras with involutive antipode. Weak Hopf algebras have also been called
``quantum groupoids''~\cite{NV} and in this paper this is \emph{not} what we
mean by quantum groupoid.

Perhaps the simplest example of weak bialgebras and weak Hopf algebras are,
respectively, category algebras and groupoid algebras. Suppose that $k$ is
a field and let $\C$ be a category with set of object $\C_0$ and set of
morphism $\C_1$. The \emph{category algebra} $k[\C]$ is the vector space
$k[\C_1]$ over $k$ with basis $\C_1$. Elements are formal linear
combinations of the elements of $\C_1$ with coefficients in $k$, i.e., 
\[
    \alpha f + \beta g + \cdots
\]
with $\alpha,\beta \in k$ and $f,g \in \C_1$. An associative multiplication
on $k[\C]$ is defined by
\[
    \mu(f,g) = f \cdot g =
    \begin{cases}
        g \o f & \text{if $g \o f$ exists} \\
        0 & \text{otherwise}
    \end{cases}
\]
and extended by linearity to $k[\C]$. This algebra does not have a unit
unless $\C_0$ is finite, in which case the unit is
\[
    \eta(1) = e = \sum_{A \in \ob \C} 1_A,
\]
making $k[\C]$ into a unital algebra; all algebras (monoids) considered in
this paper will be unital. A comultiplication and counit may be defined on
$k[\C]$ as
\begin{align*}
    \delta(f) &= f \ox f \\
    \epsilon(f) &= 1
\end{align*}
making $k[\C]$ into a coalgebra. Note that $k[\C]$ equipped with this
algebra and coalgebra structure will not satisfy any of the following
usual bialgebra axioms:
\[
    \epsilon\mu = \epsilon \ox \epsilon \quad\qquad
    \delta\eta = \eta \ox \eta \quad\qquad
    \epsilon\eta = 1_k.
\]
The one bialgebra axiom that does hold is
$\delta\mu = (\mu \ox \mu)(1 \ox c \ox 1)(\delta \ox \delta)$. Equipped with
this algebra and coalgebra structure $k[\C]$ does, however, satisfy the
axioms of a weak bialgebra. Furthermore, if $\C$ is a groupoid, then $k[\C]$,
which is then called the \emph{groupoid algebra}, is an example of a weak
Hopf algebra with antipode $\nu:k[\C] \ra k[\C]$ defined by
\[
    \nu(f) = f^{-1}.
\]
and extended by linearity. If $f:A \ra B \in \C$, the source and target
morphisms $s,t:k[\C] \ra k[\C]$ are given by
\[
    s(f) = 1_A \qquad \text{and} \qquad t(f) = 1_B,
\]
as one would expect.
    
In this paper we define weak bialgebras and weak Hopf algebras in a braided
monoidal category $\V$, where prefer to call them ``weak bimonoids'' and
``weak Hopf monoids''. To define a weak bimonoid in $\V$ the only difference
from the definition given by B\"ohm, Nill, and Szlach\'anyi~\cite{BNS} is
that a choice of ``crossing'' must be made in the axioms. Our definition is
not as general as the one given by J. N. Alonso \'Alvarez, J. M.
Fern\'andez Vilaboa, and R. Gonz\'alez Rodr\'iguez in~\cite{AFG1,AFG2}, but,
in the case that
their weak Yang-Baxter operator $t_{A,A}$ is the braiding $c_{A,A}$ and their
idempotent $\nabla_{A \ox A} = 1_{A \ox A}$, then our choices of crossings
are the same. Our difference in defining weak bimonoids occurs in the choice
of source and target morphisms. We have chosen $s:A \ra A$ and $t:A \ra A$
so that:
\begin{enumerate}
\item the ``globular'' identities $ts = s$ and $st = t$ hold;
\item the source subcomonoid and target subcomonoid coincide (up to
isomorphism), and is denoted by $C$;
\item $s:A \ra C^\o$ and $t:A \ra C$ are comonoid morphisms.
\end{enumerate}
These properties of the source and target morphisms are essential for our
view of quantum categories. These are $s = \bar{\Pi}^L_A$ and $t = \Pi^R_A$
in the notation of~\cite{AFG1,AFG2} and $s = \epsilon'_{s}$ and $t=\epsilon_s$
in the notation of~\cite{Sch}, with the appropriate choice of crossings.

We choose to work in the Cauchy completion $\Q\V$ of $\V$. The category
$\Q\V$ is also called the ``completion under idempotents'' of $\V$ or the
``Karoubi envelope'' of $\V$. This is done rather than assume that
idempotents split in $\V$. Suppose that $A$ is a weak bimonoid in $\Q\V$. In
this case we find $C$ by splitting either the source or target morphism. As
in~\cite[Prop.~4.2]{Sch}, $C$ is a separable Frobenius monoid in $\Q\V$,
meaning that $(C,\mu,\eta,\delta,\epsilon)$ is a Frobenius monoid with
$\mu\delta = 1_C$.

It turns out that our definition of weak Hopf monoid is (in the symmetric
case) the same as what is proposed in~\cite{BNS}, and in the braided case
in~\cite{AFG1,AFG2}. A weak bimonoid $H$ is a weak Hopf monoid if it is
equipped with an antipode $\nu:H \ra H$ satisfying
\[
    \mu(\nu \ox 1)\delta = t, \qquad
    \mu(1 \ox \nu)\delta = r, \qquad \text{and} \qquad 
    \mu_3(\nu \ox 1 \ox \nu)\delta_3 = \nu,
\]
where $r = \nu s$. This $r:H \ra H$ here turns out to be the ``usual''
source morphism; $\Pi^L_H$ in the notation of~\cite{AFG1,AFG2}. Ignoring
crossings $r$ is $\epsilon_t$ in the notation of~\cite{Sch} and our $r$ and
$t$ correspond respectively to $\sqcap^L$ and $\sqcap^R$ in the notation
of~\cite{BNS}; the morphism $s$ does not appear in~\cite{BNS}. Usually, in
the second axiom above, $\mu(1 \ox nu)\delta = r$, the right-hand
side is equal to the chosen source map $s$ of the weak bimonoid $H$. The
reason that this $r$ does not work as a source morphism for us is that it
does not satisfy all three requirements for the source morphism mentioned
above. This choice of $r$ allows us to show that any Frobenius monoid in
$\V$ yields a weak Hopf monoid $R \ox R$ with bijective antipode (cf. the
example in the Appendix of~\cite{BNS}).

There are a number of generalizations of bialgebras and Hopf algebras to their
``many object'' versions. For example, Sweedler's generalized
bialgebras~\cite{Sw}, which were later generalized by Takeuchi to
$\times_R$-bialgebras~\cite{T}, the quantum groupoids of Lu~\cite{L} and
Xu~\cite{X}, Schauenburg's $\times_R$-Hopf algebras~\cite{SH}, the
bialgebroids and Hopf
algebroids of B\"ohm and Szlach\'anyi~\cite{BSb}, the earlier mentioned 
face algebras~\cite{H} and generalized Kac algebras~\cite{Y}, and, the ones
of interest in this paper, the quantum categories and quantum groupoids of
Day and Street~\cite{DS}. It has been shown by Brzezi\'nski and Militaru
that the quantum groupoids of Lu and Xu are equivalent to Takeuchi's
$\times_R$-bialgebras~\cite[Thm.~3.1]{BM}. Schauenburg has shown
in~\cite{Sch2} that face algebras are an example of $\times_R$-bialgebras
for which $R$ is commutative and separable. In ~\cite[Thm.~5.1]{Sch}
Schauenburg has shown that weak bialgebras are also examples of
$\times_R$-bialgebras for which $R$ is separable Frobenius (there called
Frobenius-separable). Schauenburg also shows in~\cite[Thm.~6.1]{Sch} that a
weak Hopf algebra may be characterized as a weak bialgebra $H$ for which a
certain canonical map $H \ox_C H \ra \mu(\delta(\eta(1)),H \ox H)$ is a
bijection. As a corollary he shows that a weak Hopf algebra is a
$\times_R$-Hopf algebra.

Quantum groupoids were introduced in~\cite{DS}. They first introduce
quantum categories. A quantum category in $\V$ consists of two comonoids
$A$ and $C$ in $\V$, with $A$ playing the role of the object-of-morphisms
and $C$ the object-of-objects. There are source and target morphisms
$s,t:A \ra C$, a ``composition'' morphism $\mu:A \ox_C A \ra A$, and a
``unit'' morphism $\eta:C \ra A$ all in $\V$. This data must satisfy a
number of axioms. Indeed, ordinary categories are examples of quantum
categories. Motivated by the duality found in $*$-autonomous
categories~\cite{B}, they then define a quantum groupoid to be a quantum
category equipped with a generalized antipode coming from a $*$-autonomous
structure.  

In this paper we show that weak bimonoids are examples of quantum categories
for which the object-of-objects $C$ is a separable Frobenius monoid, and
that weak Hopf monoids with invertible antipode are quantum groupoids.

An outline of this paper is as follows:

In \S\ref{sec-wb} we provide the definition of weak bimonoid $A$ in a
braided monoidal category $\V$ and define the source and target morphisms.
We then move to the Cauchy completion $\Q\V$ and prove the three required
properties of our source and target morphisms mentioned above. In this
section we also prove that $C$, the object-of-objects of $A$, is a separable
Frobenius monoid.

Weak Hopf monoids in braided monoidal categories are introduced in
\S\ref{sec-whm}.

In \S\ref{sec-monoidal} we describe a monoidal structure on the categories
$\Bicomod(C)$ of $C$-bicomodules in $\V$, and $\Comod(A)$ of right
$A$-comodules in $\V$, such that the underlying functor
\[
    U:\Comod(A) \dra \Bicomod(C)
\]
is strong monoidal. If $H$ is a weak Hopf monoid, then we are able to show
that the category $\Comod_f(H)$, consisting of the dualizable objects of
$\Comod(H)$, is left autonomous.

In \S\ref{sec-frobeg} we prove that any separable Frobenius monoid $R$ in a
braided monoidal category $\V$ yields an example of a weak Hopf monoid
$R \ox R$ with invertible antipode in $\V$.

The definitions of quantum categories and quantum groupoids are recalled in
\S\ref{sec-quangroup}, and in \S\ref{sec-whqg} we show that any weak
bimonoid is a quantum category and any weak Hopf monoid with invertible
antipode is a quantum groupoid.

This paper depends heavily on of the string diagrams in braided monoidal
categories of Joyal and Street~\cite{JS:braided}, which were shown to be
rigorous in~\cite{JS:geometry}. The reader unfamiliar with string diagrams
may first want to read Appendix~\ref{sec-prelim} where we review some
preliminary concepts using these diagrams.

We would like to thank J. N. Alonso \'Alvarez, J. M. Fern\'andez Vilaboa,
and R. Gonz\'alez Rodr\'iguez for sending us copies of their
preprints~\cite{AFG1,AFG2}.

%=========================================================================%
\section{Weak bimonoids}\label{sec-wb}
%=========================================================================%

A weak bialgebra~\cite{BS,N,Sz,BNS} is a generalization of a bialgebra with
weakened axioms. These weakened axioms replace the three axioms that say
that the unit is a coalgebra morphism and the counit is an algebra morphism.
With the appropriate choices of under and over crossings the definition of
a weak bialgebra carries over rather straightforwardly into braided
monoidal categories, where we prefer to call it a ``weak bimonoid''.

%=========================================================================%
\subsection{Weak bimonoids}
%=========================================================================%

Suppose that $\V = (\V,\ox,I,c)$ is a braided monoidal category.

\begin{definition}
A \emph{weak bimonoid} $A = (A,\mu,\eta,\delta,\epsilon)$ in $\V$ is an
object $A \in \V$ equipped with the structure of a monoid $(A,\mu,\eta)$
and a comonoid $(A,\delta,\epsilon)$ satisfying the following equations.
\begin{equation*}\tag{b}\label{b}
    \vcenter{\xy %0;/r.16pc/:
        (3,5)="tl";
        (-3,5)="tr";
        (0,2)="mt";
        (0,-2)="mb";
        (3,-5)="bl";
        (-3,-5)="br";
        "tl";"mt" **\dir{-};
        "tr";"mt" **\dir{-};
        "mt";"mb" **\dir{-};
        "mb";"bl" **\dir{-};
        "mb";"br" **\dir{-};
    \endxy}
    ~=~
    \vcenter{\xy
        (-3,5)="1";
        (3,5)="2";
        (-3,2)="3";
        (3,2)="4";
        (-3,-2)="5";
        (3,-2)="6";
        (-3,-5)="7";
        (3,-5)="8";
        "1";"7" **\dir{-};
        "2";"8" **\dir{-};
        "3";"6" **\dir{-}?(0.5)*{\hole}="x";
        "4";"x" **\dir{-};
        "x";"5" **\dir{-};
    \endxy}
\end{equation*}

% (v)

\begin{equation*}\tag{v}\label{v}
    \vcenter{\xy
        (-4,4)="l";
        (0,4)="c";
        (4,4)="r";
        (0,-1)="m";
        (0,-5)*[o]=<5pt>[Fo]{~}="b";
        "l";"m" **\dir{-};
        "c";"m" **\dir{-};
        "r";"m" **\dir{-};
        "m";"b" **\dir{-};
    \endxy}
    ~=~
    \vcenter{\xy
        (-4,4)="bl";
        (0,4)="bc";
        (4,4)="br";
        (-4,-5)*[o]=<5pt>[Fo]{~}="tl";
        (-4,-2)="ml";
        (0,1)="mc";
        (4,-2)="mr";
        (4,-5)*[o]=<5pt>[Fo]{~}="tr";
        "bl";"tl" **\dir{-};
        "br";"tr" **\dir{-};
        "bc";"mc" **\dir{-};
        "mc";"ml" **\dir{-};
        "mc";"mr" **\dir{-};
    \endxy}
    ~=~
    \vcenter{\xy
        (-4,4)="bl";
        (0,4)="bc";
        (4,4)="br";
        (-4,-5)*[o]=<5pt>[Fo]{~}="tl";
        (-4,-2)="ml";
        (0,2)="mc";
        (4,-2)="mr";
        (4,-5)*[o]=<5pt>[Fo]{~}="tr";
        "bl";"tl" **\dir{-};
        "br";"tr" **\dir{-};
        "bc";"mc" **\dir{-};
        "mc";"ml" **\crv{(4,0)}?(0.6)*{\hole}="x";
        "mc";"x" **\crv{(-3,0.5)};
        "mr";"x" **\crv{(1,-1.5)};
    \endxy}
\end{equation*}

% (w)

\begin{equation*}\tag{w}\label{w}
    \vcenter{\xy
        (-4,-4)="l";
        (0,-4)="c";
        (4,-4)="r";
        (0,1)="m";
        (0,5)*[o]=<5pt>[Fo]{~}="b";
        "l";"m" **\dir{-};
        "c";"m" **\dir{-};
        "r";"m" **\dir{-};
        "m";"b" **\dir{-};
    \endxy}
    ~=~
    \vcenter{\xy
        (-4,-4)="bl";
        (0,-4)="bc";
        (4,-4)="br";
        (-4,5)*[o]=<5pt>[Fo]{~}="tl";
        (-4,2)="ml";
        (0,-1)="mc";
        (4,2)="mr";
        (4,5)*[o]=<5pt>[Fo]{~}="tr";
        "bl";"tl" **\dir{-};
        "br";"tr" **\dir{-};
        "bc";"mc" **\dir{-};
        "mc";"ml" **\dir{-};
        "mc";"mr" **\dir{-};
    \endxy}
    ~=~
    \vcenter{\xy
        (4,-4)="bl";
        (0,-4)="bc";
        (-4,-4)="br";
        (4,5)*[o]=<5pt>[Fo]{~}="tl";
        (4,2)="ml";
        (0,-2)="mc";
        (-4,2)="mr";
        (-4,5)*[o]=<5pt>[Fo]{~}="tr";
        "bl";"tl" **\dir{-};
        "br";"tr" **\dir{-};
        "bc";"mc" **\dir{-};
        "mc";"ml" **\crv{(-4,0)}?(0.6)*{\hole}="x";
        "mc";"x" **\crv{( 3,-0.5)};
        "mr";"x" **\crv{(-1,1.5)};
    \endxy}
\end{equation*}
\end{definition}

Suppose $A$ and $B$ are weak bimonoids in $\V$. A \emph{morphism of weak
bimonoids} $f:A \ra B$ is a morphism $f:A \ra B$ in $\V$ which is both a
monoid morphism and a comonoid morphism.

Let $A$ be a weak bimonoid and define the \emph{source} and \emph{target}
morphisms $s,t:A \ra A$ of $A$ as follows:
\[
    s =~
    \vcenter{\xy
        ( 3, 6)="0";
        ( 0, 4)*[o]=<5pt>[Fo]{~}="1";
        ( 0, 2)="2";
        ( 0,-2)="3";
        ( 0,-4)*[o]=<5pt>[Fo]{~}="4";
        (-3, -6)="5";
        "1";"4" **\dir{-};
        "0";"3" **\crv{(3,2)};
        "5";"2" **\crv{(-3,-2)};
    \endxy}
    \quad\qquad\qquad
    t =~
    \vcenter{\xy
        (-4, 6)="0";
        ( 0, 4)*[o]=<5pt>[Fo]{~}="1";
        ( 0, 2)="2";
        ( 0,-2)="3";
        ( 0,-4)*[o]=<5pt>[Fo]{~}="4";
        (-4, -6)="5";
        "1";"4" **\dir{-};
        "0";"3" **\crv{(-4,0)}?(0.6)*{\hole}="x";
        "5";"x" **\crv{(-4,-2)};
        "2";"x" **\crv{(-2,1)};
    \endxy}~.
\]

Notice that $s:A \ra A$ is invariant under rotation by $\pi$, while $t:A \ra
A$ is invariant under horizontal reflection and the inverse braiding.
Importantly, under either of these transformations
\begin{itemize}
\item (m) and (c) are interchanged\footnote{The (m) and (c) here refer to
the monoid and comonoid identities found in Appendix~\ref{sec-prelim}.},
\item (b) is invariant, and
\item (v) and (w) are interchanged.
\end{itemize}

Note that these are not the ``usual'' source and target morphisms. They
were chosen, as mentioned in the introduction, precisely because of the
need for them to satisfy the following three properties:
\begin{enumerate}
\item the ``globular'' identities $ts = s$ and $st = t$ hold;
\item the source subcomonoid and target subcomonoid coincide (up to
isomorphism), and is denoted by $C$;
\item $s:A \ra C^\o$ and $t:A \ra C$ are comonoid morphisms.
\end{enumerate}
These properties will be proved in this section. Note that we will run into
the usual source morphism (which we call $r$) in the definition of weak Hopf
monoids (Definition~\ref{def-whm}).

\begin{figure}
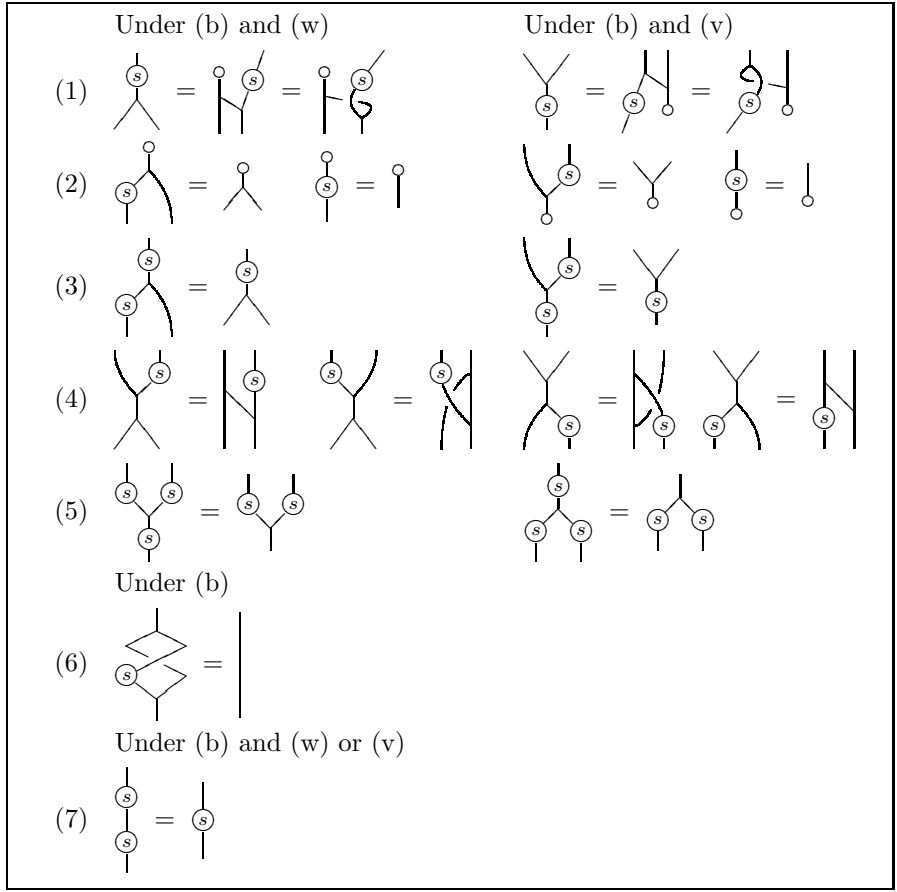
 % PROPERTIES OF s
\begin{center}
\framebox[\boxwidth]{
\begin{tabular}{rll}
& Under (b) and (w) & Under (b) and (v) \\

(1) & $
    \vcenter{\xy
        (0,5)="t";
        (0,2)*[o]=<9pt>[Fo]{\scriptstyle s}="f";
        (0,-1)="m";
        (-3,-5)="l";
        ( 3,-5)="r";
        "t";"f" **\dir{-};
        "f";"m" **\dir{-};
        "m";"l" **\dir{-};
        "m";"r" **\dir{-};
    \endxy}
    ~=~ 
    \vcenter{\xy
        (0,7)*{~};
        (0,-6)*{~};
        (-3,3)*[o]=<5pt>[Fo]{~}="1";
        (-3,0)="2";
        (-3,-5)="3";
        (3,6)="4";
        (1.5,2)*[o]=<9pt>[Fo]{\scriptstyle s}="5";
        (0,-2)="6";
        (0,-5)="7";
        "1";"3" **\dir{-};
        "2";"6" **\dir{-};
        "6";"7" **\dir{-};
        "5";"6" **\dir{-};
        "4";"5" **\dir{-};
    \endxy}
    ~=~
    \vcenter{\xy
        (0,8)*{~};
        (0,-5)*{~};
        (-5,4)*[o]=<5pt>[Fo]{~}="1";
        (-5,1)="2";
        (-5,-4)="3";
        (3,7)="4";
        (0,3)*[o]=<9pt>[Fo]{\scriptstyle s}="5";
        (0,-2)="6";
        (0,-4)="7";
        "1";"3" **\dir{-};
        "6";"7" **\dir{-};
        "5";"6" **\crv{(-3,0)}?(0.5)*{\hole}="x";
        "2";"x" **\dir{-};
        "6";"x" **\crv{(3,0)};
        "4";"5" **\dir{-};
    \endxy}

$ & $

    \vcenter{\xy
        (0,-5)="t";
        (0,-2)*[o]=<9pt>[Fo]{\scriptstyle s}="f";
        (0,1)="m";
        (-3,5)="l";
        ( 3,5)="r";
        "t";"f" **\dir{-};
        "f";"m" **\dir{-};
        "m";"l" **\dir{-};
        "m";"r" **\dir{-};
    \endxy}
    ~=~ 
    \vcenter{\xy
        (0,-7)*{~};
        (0, 6)*{~};
        (3,-3)*[o]=<5pt>[Fo]{~}="1";
        (3,0)="2";
        (3,5)="3";
        (-3,-6)="4";
        (-1.5,-2)*[o]=<9pt>[Fo]{\scriptstyle s}="5";
        (0,2)="6";
        (0,5)="7";
        "1";"3" **\dir{-};
        "2";"6" **\dir{-};
        "6";"7" **\dir{-};
        "5";"6" **\dir{-};
        "4";"5" **\dir{-};
    \endxy}
~=~
    \vcenter{\xy
        (0,-8)*{~};
        (0, 5)*{~};
        (5,-4)*[o]=<5pt>[Fo]{~}="1";
        (5,-1)="2";
        (5, 4)="3";
        (-3,-7)="4";
        (0,-3)*[o]=<9pt>[Fo]{\scriptstyle s}="5";
        (0, 2)="6";
        (0, 4)="7";
        "1";"3" **\dir{-};
        "6";"7" **\dir{-};
        "5";"6" **\crv{(3,0)}?(0.5)*{\hole}="x";
        "2";"x" **\dir{-};
        "6";"x" **\crv{(-3,0)};
        "4";"5" **\dir{-};
    \endxy}
$ \\

(2) & $
\vcenter{\xy
    (0,4)*{~};
    (0,-8)*{~};
    (0,3)*[o]=<5pt>[Fo]{~}="t";
    (0,0)="m";
    (-3,-3)*[o]=<9pt>[Fo]{\scriptstyle s}="g";
    (-3,-7)="l";
    ( 3,-7)="r";
    "t";"m" **\dir{-};
    "r";"m" **\crv{(3,-3)};
    "m";"g" **\dir{-};
    "g";"l" **\dir{-};
\endxy}
~=~
\vcenter{\xy
    (0,2.5)*[o]=<5pt>[Fo]{~}="1";
    (0,0)="2";
    (-2.5,-3)="3";
    ( 2.5,-3)="4";
    "1";"2" **\dir{-};
    "2";"3" **\dir{-};
    "2";"4" **\dir{-};
\endxy}

\qquad

\vcenter{\xy
    (0,-4.5)="t";
    (0,0)*[o]=<9pt>[Fo]{\scriptstyle s}="g";
    (0,4)*[o]=<5pt>[Fo]{}="b";
    "t";"g" **\dir{-};
    "g";"b" **\dir{-};
\endxy}
~=~
\vcenter{\xy
    (0,2)*[o]=<5pt>[Fo]{}="t";
    (0,-3)="b";
    "t";"b" **\dir{-};
\endxy}

$ & $

\vcenter{\xy
    (0,-4)*[o]=<5pt>[Fo]{~}="t";
    (0,-1)="m";
    (3,2)*[o]=<9pt>[Fo]{\scriptstyle s}="g";
    (3,6)="l";
    (-3,6)="r";
    "t";"m" **\dir{-};
    "r";"m" **\crv{(-3,2)};
    "m";"g" **\dir{-};
    "g";"l" **\dir{-};
\endxy}
~=~
\vcenter{\xy
    (0,-2.5)*[o]=<5pt>[Fo]{~}="1";
    (0,0)="2";
    (-2.5,3)="3";
    ( 2.5,3)="4";
    "1";"2" **\dir{-};
    "2";"3" **\dir{-};
    "2";"4" **\dir{-};
\endxy}

\qquad

\vcenter{\xy
    (0,4)="t";
    (0,0)*[o]=<9pt>[Fo]{\scriptstyle s}="g";
    (0,-4.5)*[o]=<5pt>[Fo]{}="b";
    "t";"g" **\dir{-};
    "g";"b" **\dir{-};
\endxy}
~=~
\vcenter{\xy
    (0,2.5)="t";
    (0,-2.5)*[o]=<5pt>[Fo]{}="b";
    "t";"b" **\dir{-};
\endxy}
$ \\

(3) & $
\vcenter{\xy
    (0,7)*{~};
    (0,-8)*{~};
    (0,6)="t";
    (0,3)*[o]=<9pt>[Fo]{\scriptstyle s}="f";
    (0,-0)="m";
    (-3,-3)*[o]=<9pt>[Fo]{\scriptstyle s}="g";
    (-3,-7)="l";
    ( 3,-7)="r";
    "t";"f" **\dir{-};
    "f";"m" **\dir{-};
    "r";"m" **\crv{(3,-3)};
    "m";"g" **\dir{-};
    "g";"l" **\dir{-};
\endxy}
~=~
\vcenter{\xy
    (0,5)="t";
    (0,2)*[o]=<9pt>[Fo]{\scriptstyle s}="f";
    (0,-1)="m";
    (-3,-5)="l";
    ( 3,-5)="r";
    "t";"f" **\dir{-};
    "f";"m" **\dir{-};
    "m";"l" **\dir{-};
    "m";"r" **\dir{-};
\endxy}

$ & $

\vcenter{\xy
    (0,-6)="t";
    (0,-3)*[o]=<9pt>[Fo]{\scriptstyle s}="f";
    (0,0)="m";
    (3,3)*[o]=<9pt>[Fo]{\scriptstyle s}="g";
    ( 3,7)="l";
    (-3,7)="r";
    "t";"f" **\dir{-};
    "f";"m" **\dir{-};
    "r";"m" **\crv{(-3,3)};
    "m";"g" **\dir{-};
    "g";"l" **\dir{-};
\endxy}
~=~ 
\vcenter{\xy
    (0,-5)="t";
    (0,-2)*[o]=<9pt>[Fo]{\scriptstyle s}="f";
    (0,1)="m";
    (-3,5)="l";
    ( 3,5)="r";
    "t";"f" **\dir{-};
    "f";"m" **\dir{-};
    "m";"l" **\dir{-};
    "m";"r" **\dir{-};
\endxy}
$ \\

(4) & $
\vcenter{\xy
    (0,9)*{~};
    (0,-6)*{~};
    (-3,8)="1";
    ( 3,8)="2";
    ( 3,5)*[o]=<9pt>[Fo]{\scriptstyle s}="3";
    ( 0,2)="4";
    ( 0,-1)="5";
    (-3,-5)="6";
    ( 3,-5)="7";
    "1";"4" **\crv{(-3,5)};
    "2";"3" **\dir{-};
    "3";"4" **\dir{-};
    "4";"5" **\dir{-};
    "5";"6" **\dir{-};
    "5";"7" **\dir{-};
\endxy}
~=~
\vcenter{\xy
    (-2,8)="1";
    ( 2,8)="2";
    ( 2,4)*[o]=<9pt>[Fo]{\scriptstyle s}="3";
    (-2,3)="4";
    ( 2,-1)="5";
    (-2,-5)="6";
    ( 2,-5)="7";
    "1";"6" **\dir{-};
    "2";"3" **\dir{-};
    "3";"7" **\dir{-};
    "4";"5" **\dir{-};
\endxy}

\qquad

\vcenter{\xy
    (0,9)*{~};
    (0,-6)*{~};
    ( 3,8)="1";
    (-3,8)="2";
    (-3,5)*[o]=<9pt>[Fo]{\scriptstyle s}="3";
    ( 0,2)="4";
    ( 0,-1)="5";
    ( 3,-5)="6";
    (-3,-5)="7";
    "1";"4" **\crv{(3,5)};
    "2";"3" **\dir{-};
    "3";"4" **\dir{-};
    "4";"5" **\dir{-};
    "5";"6" **\dir{-};
    "5";"7" **\dir{-};
\endxy}
~=~
\vcenter{\xy
    (-2,8)="1";
    ( 2,8)="2";
    ( 2,5)="3";
    ( 2,-2)="4";
    (-2,-5)="5";
    ( 2,-5)="6";
    (-2,5)*[o]=<9pt>[Fo]{\scriptstyle s}="7";
    "2";"6" **\dir{-};
    "1";"7" **\dir{-};
    "7";"4" **\crv{(-2,2)}?(0.5)*{\hole}="x";
    "3";"x" **\crv{(1,5)};
    "5";"x" **\crv{(-2,-2)};
\endxy}

\quad$ & $

\vcenter{\xy
    (0,-9)*{~};
    (0,6)*{~};
    (-3,-8)="1";
    ( 3,-8)="2";
    ( 3,-5)*[o]=<9pt>[Fo]{\scriptstyle s}="3";
    ( 0,-2)="4";
    ( 0,1)="5";
    (-3,5)="6";
    ( 3,5)="7";
    "1";"4" **\crv{(-3,-5)};
    "2";"3" **\dir{-};
    "3";"4" **\dir{-};
    "4";"5" **\dir{-};
    "5";"6" **\dir{-};
    "5";"7" **\dir{-};
\endxy}
~=~
\vcenter{\xy
    ( 2,-8)="1";
    (-2,-8)="2";
    (-2,-5)="3";
    (-2,2)="4";
    ( 2,5)="5";
    (-2,5)="6";
    ( 2,-5)*[o]=<9pt>[Fo]{\scriptstyle s}="7";
    "2";"6" **\dir{-};
    "1";"7" **\dir{-};
    "7";"4" **\crv{(2,-2)}?(0.5)*{\hole}="x";
    "3";"x" **\crv{(-1,-5)};
    "5";"x" **\crv{(2,2)};
\endxy}

\quad

\vcenter{\xy
    (0,-9)*{~};
    (0,6)*{~};
    ( 3,-8)="1";
    (-3,-8)="2";
    (-3,-5)*[o]=<9pt>[Fo]{\scriptstyle s}="3";
    ( 0,-2)="4";
    ( 0,1)="5";
    ( 3,5)="6";
    (-3,5)="7";
    "1";"4" **\crv{(3,-5)};
    "2";"3" **\dir{-};
    "3";"4" **\dir{-};
    "4";"5" **\dir{-};
    "5";"6" **\dir{-};
    "5";"7" **\dir{-};
\endxy}
~=~
\vcenter{\xy
    (-2,8)="1";
    ( 2,8)="2";
    (-2,4)="3";
    ( 2,0)="4";
    (-2,-1)*[o]=<9pt>[Fo]{\scriptstyle s}="5";
    (-2,-5)="6";
    ( 2,-5)="7";
    "1";"5" **\dir{-};
    "5";"6" **\dir{-};
    "2";"7" **\dir{-};
    "3";"4" **\dir{-};
\endxy}
$ \\

(5) & $
\vcenter{\xy
    (0,8)*{~};
    (0,-7)*{~};
    ( 0,-6)="1";
    ( 0,-3)*[o]=<9pt>[Fo]{\scriptstyle s}="2";
    ( 0,-0)="3";
    (-3, 3)*[o]=<9pt>[Fo]{\scriptstyle s}="4";
    ( 3, 3)*[o]=<9pt>[Fo]{\scriptstyle s}="5";
    (-3, 7)="6";
    ( 3, 7)="7";
    "1";"2" **\dir{-};
    "2";"3" **\dir{-};
    "3";"4" **\dir{-};
    "4";"6" **\dir{-};
    "3";"5" **\dir{-};
    "5";"7" **\dir{-};
\endxy}
~=~
\vcenter{\xy
    ( 0,-3)="2";
    ( 0, 0)="3";
    (-3, 3)*[o]=<9pt>[Fo]{\scriptstyle s}="4";
    ( 3, 3)*[o]=<9pt>[Fo]{\scriptstyle s}="5";
    (-3, 7)="6";
    ( 3, 7)="7";
    "2";"3" **\dir{-};
    "3";"4" **\dir{-};
    "4";"6" **\dir{-};
    "3";"5" **\dir{-};
    "5";"7" **\dir{-};
\endxy}

$ & $

\vcenter{\xy
    ( 0, 6)="1";
    ( 0, 3)*[o]=<9pt>[Fo]{\scriptstyle s}="2";
    ( 0, 0)="3";
    (-3,-3)*[o]=<9pt>[Fo]{\scriptstyle s}="4";
    ( 3,-3)*[o]=<9pt>[Fo]{\scriptstyle s}="5";
    (-3,-7)="6";
    ( 3,-7)="7";
    "1";"2" **\dir{-};
    "2";"3" **\dir{-};
    "3";"4" **\dir{-};
    "4";"6" **\dir{-};
    "3";"5" **\dir{-};
    "5";"7" **\dir{-};
\endxy}
~=~
\vcenter{\xy
    ( 0, 3)="2";
    ( 0, 0)="3";
    (-3,-3)*[o]=<9pt>[Fo]{\scriptstyle s}="4";
    ( 3,-3)*[o]=<9pt>[Fo]{\scriptstyle s}="5";
    (-3,-7)="6";
    ( 3,-7)="7";
    "2";"3" **\dir{-};
    "3";"4" **\dir{-};
    "4";"6" **\dir{-};
    "3";"5" **\dir{-};
    "5";"7" **\dir{-};
\endxy}
$ \\

& Under (b) \\

(6) & $
    \vcenter{\xy
        (0,9)*{~};
        (0,-8)*{~};
        (0,8)="1";
        (0,5)="2";
        (-4,3)="3";
        ( 4,3)="4";
        (-4,-1)*[o]=<9pt>[Fo]{\scriptstyle s}="5";
        ( 4,-1)="6";
        ( 0,-4)="7";
        ( 0,-7)="8";
        "1";"2" **\dir{-}; 
        "2";"3" **\dir{-}; 
        "2";"4" **\dir{-}; 
        "4";"5" **\dir{-}; 
        {\ar@{-} |\hole "3";"6"}; 
        "5";"7" **\dir{-}; 
        "6";"7" **\dir{-}; 
        "8";"7" **\dir{-}; 
    \endxy}
    ~=~
    \vcenter{\xy
        (0,7)="1";
        (0,-7)="2";
        **\dir{-}; 
    \endxy}
$ \\

& Under (b) and (w) or (v) \\

(7) & $
\vcenter{\xy
    (0,8)*{~};
    (0,-8)*{~};
    (0,7)="t";
    (0,3)*[o]=<9pt>[Fo]{\scriptstyle s}="m1";
    (0,-3)*[o]=<9pt>[Fo]{\scriptstyle s}="m2";
    (0,-7)="b";
    "t";"m1" **\dir{-};
    "m1";"m2" **\dir{-};
    "m2";"b" **\dir{-};
\endxy}
~=~
\vcenter{\xy
    (0,5)="t";
    (0,0)*[o]=<9pt>[Fo]{\scriptstyle s}="m";
    (0,-5)="b";
    "t";"m" **\dir{-};
    "m";"b" **\dir{-};
\endxy}
$ \\
\end{tabular}}
\end{center}
\caption{Properties of $s$} \label{prop-s}
\end{figure}

\begin{figure}
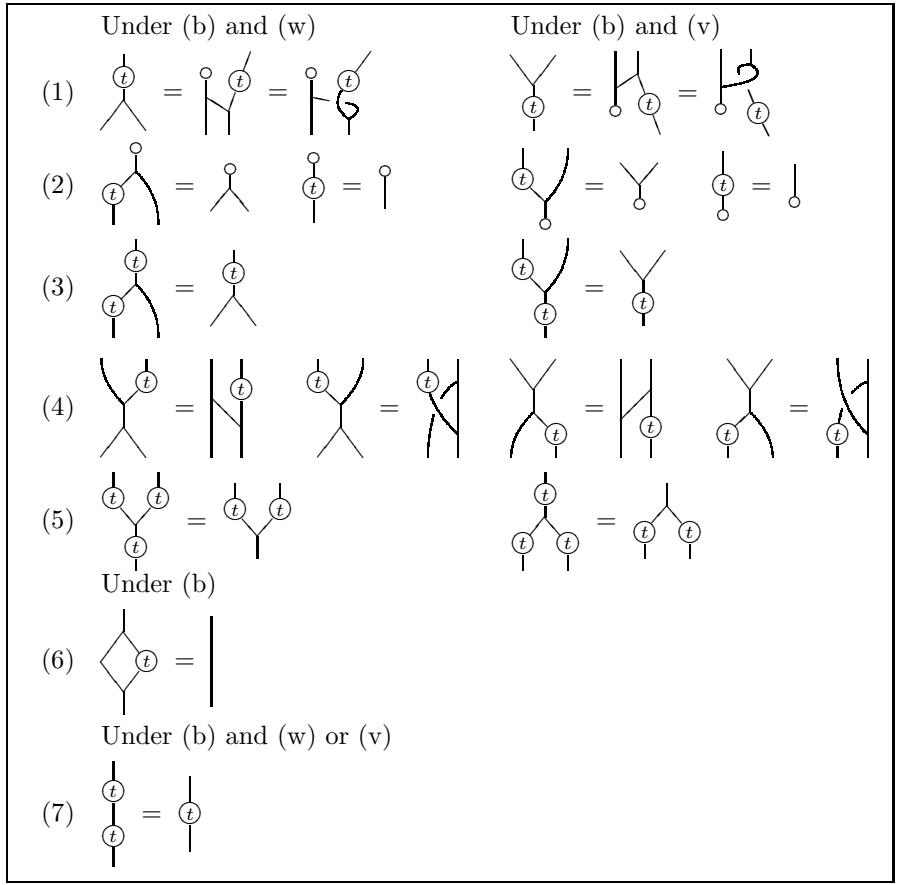
 % PROPERTIES OF t
\begin{center}
\framebox[\boxwidth]{
\begin{tabular}{rll}
& Under (b) and (w) & Under (b) and (v) \\

(1) & $
    \vcenter{\xy
        (0,5)="t";
        (0,2)*[o]=<9pt>[Fo]{\scriptstyle t}="f";
        (0,-1)="m";
        (-3,-5)="l";
        ( 3,-5)="r";
        "t";"f" **\dir{-};
        "f";"m" **\dir{-};
        "m";"l" **\dir{-};
        "m";"r" **\dir{-};
    \endxy}
    ~=~ 
    \vcenter{\xy
        (0,7)*{~};
        (0,-6)*{~};
        (-3,3)*[o]=<5pt>[Fo]{~}="1";
        (-3,0)="2";
        (-3,-5)="3";
        (3,6)="4";
        (1.5,2)*[o]=<9pt>[Fo]{\scriptstyle t}="5";
        (0,-2)="6";
        (0,-5)="7";
        "1";"3" **\dir{-};
        "2";"6" **\dir{-};
        "6";"7" **\dir{-};
        "5";"6" **\dir{-};
        "4";"5" **\dir{-};
    \endxy}
~=~
    \vcenter{\xy
        (0,8)*{~};
        (0,-5)*{~};
        (-5,4)*[o]=<5pt>[Fo]{~}="1";
        (-5,1)="2";
        (-5,-4)="3";
        (3,7)="4";
        (0,3)*[o]=<9pt>[Fo]{\scriptstyle t}="5";
        (0,-2)="6";
        (0,-4)="7";
        "1";"3" **\dir{-};
        "6";"7" **\dir{-};
        "5";"6" **\crv{(-3,0)}?(0.5)*{\hole}="x";
        "2";"x" **\dir{-};
        "6";"x" **\crv{(3,0)};
        "4";"5" **\dir{-};
    \endxy}

$ & $

\vcenter{\xy
    (0,-5)="t";
    (0,-2)*[o]=<9pt>[Fo]{\scriptstyle t}="f";
    (0,1)="m";
    (-3,5)="l";
    ( 3,5)="r";
    "t";"f" **\dir{-};
    "f";"m" **\dir{-};
    "m";"l" **\dir{-};
    "m";"r" **\dir{-};
\endxy}
~=~ 
\vcenter{\xy
    (0,-7)*{~};
    (0, 6)*{~};
    (-3,-3)*[o]=<5pt>[Fo]{~}="1";
    (-3,0)="2";
    (-3,5)="3";
    (3,-6)="4";
    (1.5,-2)*[o]=<9pt>[Fo]{\scriptstyle t}="5";
    (0,2)="6";
    (0,5)="7";
    "1";"3" **\dir{-};
    "2";"6" **\dir{-};
    "6";"7" **\dir{-};
    "5";"6" **\dir{-};
    "4";"5" **\dir{-};
\endxy}
~=~
\vcenter{\xy
    (-5,-4)*[o]=<5pt>[Fo]{~}="1";
    (-5,-1)="2";
    (-5,4)="3";
    (0,-4.5)*[o]=<9pt>[Fo]{\scriptstyle t}="5";
    (-1,2)="6";
    (-1,4)="7";
    (1.5,-7.5)="4";
    "1";"3" **\dir{-};
    "2";"6" **\crv{(2,0)}?(0.3)*{\hole}="x";
    "6";"7" **\dir{-};
    "5";"x" **\dir{-};
    "6";"x" **\crv{(-4,2)};
    "4";"5" **\dir{-};
\endxy}
$ \\

(2) & $
\vcenter{\xy
    (0,5)*{~};
    (0,-7)*{~};
    (0,4)*[o]=<5pt>[Fo]{~}="t";
    (0,1)="m";
    (-3,-2)*[o]=<9pt>[Fo]{\scriptstyle t}="g";
    (-3,-6)="l";
    ( 3,-6)="r";
    "t";"m" **\dir{-};
    "r";"m" **\crv{(3,-2)};
    "m";"g" **\dir{-};
    "g";"l" **\dir{-};
\endxy}
~=~
\vcenter{\xy
    (0,2.5)*[o]=<5pt>[Fo]{~}="1";
    (0,0)="2";
    (-2.5,-3)="3";
    ( 2.5,-3)="4";
    "1";"2" **\dir{-};
    "2";"3" **\dir{-};
    "2";"4" **\dir{-};
\endxy}
\qquad
\vcenter{\xy
    (0,-4.5)="t";
    (0,0)*[o]=<9pt>[Fo]{\scriptstyle t}="g";
    (0,4)*[o]=<5pt>[Fo]{}="b";
    "t";"g" **\dir{-};
    "g";"b" **\dir{-};
\endxy}
~=~
\vcenter{\xy 
    (0,2)*[o]=<5pt>[Fo]{}="t";
    (0,-3)="b";
    "t";"b" **\dir{-};
\endxy}

$ & $

\vcenter{\xy
    (0,-5)*{~};
    (0, 7)*{~};
    (0,-4)*[o]=<5pt>[Fo]{~}="t";
    (0,-1)="m";
    (-3,2)*[o]=<9pt>[Fo]{\scriptstyle t}="g";
    (-3,6)="l";
    ( 3,6)="r";
    "t";"m" **\dir{-};
    "r";"m" **\crv{(3,2)};
    "m";"g" **\dir{-};
    "g";"l" **\dir{-};
\endxy}
~=~
\vcenter{\xy
    (0,-2.5)*[o]=<5pt>[Fo]{~}="1";
    (0,0)="2";
    (-2.5,3)="3";
    ( 2.5,3)="4";
    "1";"2" **\dir{-};
    "2";"3" **\dir{-};
    "2";"4" **\dir{-};
\endxy}
\qquad
\vcenter{\xy
    (0,4.5)="t";
    (0,0)*[o]=<9pt>[Fo]{\scriptstyle t}="g";
    (0,-4)*[o]=<5pt>[Fo]{}="b";
    "t";"g" **\dir{-};
    "g";"b" **\dir{-};
\endxy}
~=~
\vcenter{\xy
    (0,3)="t";
    (0,-2)*[o]=<5pt>[Fo]{}="b";
    "t";"b" **\dir{-};
\endxy}
$ \\

(3) & $
\vcenter{\xy
    (0,7)*{~};
    (0,-8)*{~};
    (0,6)="t";
    (0,3)*[o]=<9pt>[Fo]{\scriptstyle t}="f";
    (0,-0)="m";
    (-3,-3)*[o]=<9pt>[Fo]{\scriptstyle t}="g";
    (-3,-7)="l";
    ( 3,-7)="r";
    "t";"f" **\dir{-};
    "f";"m" **\dir{-};
    "r";"m" **\crv{(3,-3)};
    "m";"g" **\dir{-};
    "g";"l" **\dir{-};
\endxy}
~=~ 
\vcenter{\xy
    (0,5)="t";
    (0,2)*[o]=<9pt>[Fo]{\scriptstyle t}="f";
    (0,-1)="m";
    (-3,-5)="l";
    ( 3,-5)="r";
    "t";"f" **\dir{-};
    "f";"m" **\dir{-};
    "m";"l" **\dir{-};
    "m";"r" **\dir{-};
\endxy}

$ & $

\vcenter{\xy
    (0,-6)="t";
    (0,-3)*[o]=<9pt>[Fo]{\scriptstyle t}="f";
    (0,0)="m";
    (-3,3)*[o]=<9pt>[Fo]{\scriptstyle t}="g";
    (-3,7)="l";
    ( 3,7)="r";
    "t";"f" **\dir{-};
    "f";"m" **\dir{-};
    "r";"m" **\crv{(3,3)};
    "m";"g" **\dir{-};
    "g";"l" **\dir{-};
\endxy}
~=~
\vcenter{\xy
    (0,-5)="t";
    (0,-2)*[o]=<9pt>[Fo]{\scriptstyle t}="f";
    (0,1)="m";
    (-3,5)="l";
    ( 3,5)="r";
    "t";"f" **\dir{-};
    "f";"m" **\dir{-};
    "m";"l" **\dir{-};
    "m";"r" **\dir{-};
\endxy}
$ \\

(4) & $
\vcenter{\xy
    (0,10)*{~};
    (0,-7)*{~};
    (-3,8)="1";
    ( 3,8)="2";
    ( 3,5)*[o]=<9pt>[Fo]{\scriptstyle t}="3";
    ( 0,2)="4";
    ( 0,-1)="5";
    (-3,-5)="6";
    ( 3,-5)="7";
    "1";"4" **\crv{(-3,5)};
    "2";"3" **\dir{-};
    "3";"4" **\dir{-};
    "4";"5" **\dir{-};
    "5";"6" **\dir{-};
    "5";"7" **\dir{-};
\endxy}
~=~
\vcenter{\xy
    (0,10)*{~};
    (0,-7)*{~};
    (-2,8)="1";
    ( 2,8)="2";
    ( 2,4)*[o]=<9pt>[Fo]{\scriptstyle t}="3";
    (-2,3)="4";
    ( 2,-1)="5";
    (-2,-5)="6";
    ( 2,-5)="7";
    "1";"6" **\dir{-};
    "2";"3" **\dir{-};
    "3";"7" **\dir{-};
    "4";"5" **\dir{-};
\endxy}
\qquad
\vcenter{\xy
    ( 3,8)="1";
    (-3,8)="2";
    (-3,5)*[o]=<9pt>[Fo]{\scriptstyle t}="3";
    ( 0,2)="4";
    ( 0,-1)="5";
    ( 3,-5)="6";
    (-3,-5)="7";
    "1";"4" **\crv{(3,5)};
    "2";"3" **\dir{-};
    "3";"4" **\dir{-};
    "4";"5" **\dir{-};
    "5";"6" **\dir{-};
    "5";"7" **\dir{-};
\endxy}
~=~
\vcenter{\xy
    (0,9)*{~};
    (0,-6)*{~};
    (-2,8)="1";
    ( 2,8)="2";
    ( 2,5)="3";
    ( 2,-2)="4";
    (-2,-5)="5";
    ( 2,-5)="6";
    (-2,5)*[o]=<9pt>[Fo]{\scriptstyle t}="7";
    "2";"6" **\dir{-};
    "1";"7" **\dir{-};
    "7";"4" **\crv{(-2,2)}?(0.5)*{\hole}="x";
    "3";"x" **\crv{(1,5)};
    "5";"x" **\crv{(-2,-2)};
\endxy}

\quad$ & $

\vcenter{\xy
    (-3,-8)="1";
    ( 3,-8)="2";
    ( 3,-5)*[o]=<9pt>[Fo]{\scriptstyle t}="3";
    ( 0,-2)="4";
    ( 0, 1)="5";
    (-3, 5)="6";
    ( 3, 5)="7";
    "1";"4" **\crv{(-3,-5)};
    "2";"3" **\dir{-};
    "3";"4" **\dir{-};
    "4";"5" **\dir{-};
    "5";"6" **\dir{-};
    "5";"7" **\dir{-};
\endxy}
~=~
\vcenter{\xy
    (0,-10)*{~};
    (0,7)*{~};
    (-2,-8)="1";
    ( 2,-8)="2";
    ( 2,-4)*[o]=<9pt>[Fo]{\scriptstyle t}="3";
    (-2,-3)="4";
    ( 2,1)="5";
    (-2,5)="6";
    ( 2,5)="7";
    "1";"6" **\dir{-};
    "2";"3" **\dir{-};
    "3";"7" **\dir{-};
    "4";"5" **\dir{-};
\endxy}

\qquad

\vcenter{\xy
    ( 3,-8)="1";
    (-3,-8)="2";
    (-3,-5)*[o]=<9pt>[Fo]{\scriptstyle t}="3";
    ( 0,-2)="4";
    ( 0, 1)="5";
    ( 3, 5)="6";
    (-3, 5)="7";
    "1";"4" **\crv{(3,-5)};
    "2";"3" **\dir{-};
    "3";"4" **\dir{-};
    "4";"5" **\dir{-};
    "5";"6" **\dir{-};
    "5";"7" **\dir{-};
\endxy}
~=~
\vcenter{\xy
    (0,9)*{~};
    (0,-6)*{~};
    (-2,8)="1";
    ( 2,8)="2";
    ( 2,5)="3";
    ( 2,-2)="4";
    (-2,-5)="5";
    ( 2,-5)="6";
    (-2,-2)*[o]=<9pt>[Fo]{\scriptstyle t}="7";
    "2";"6" **\dir{-};
    "1";"4" **\crv{(-2,2)}?(0.5)*{\hole}="x";
    "3";"x" **\crv{(1,5)};
    "5";"7" **\dir{-};
    "7";"x" **\crv{(-2,0)};
\endxy}
$ \\

(5) & $
\vcenter{\xy
    ( 0,-6)="1";
    ( 0,-3)*[o]=<9pt>[Fo]{\scriptstyle t}="2";
    ( 0,-0)="3";
    (-3, 3.5)*[o]=<9pt>[Fo]{\scriptstyle t}="4";
    ( 3, 3.5)*[o]=<9pt>[Fo]{\scriptstyle t}="5";
    (-3, 7)="6";
    ( 3, 7)="7";
    "1";"2" **\dir{-};
    "2";"3" **\dir{-};
    "3";"4" **\dir{-};
    "4";"6" **\dir{-};
    "3";"5" **\dir{-};
    "5";"7" **\dir{-};
\endxy}
~=~
\vcenter{\xy
    ( 0,-3)="2";
    ( 0, 0)="3";
    (-3, 3.5)*[o]=<9pt>[Fo]{\scriptstyle t}="4";
    ( 3, 3.5)*[o]=<9pt>[Fo]{\scriptstyle t}="5";
    (-3, 7)="6";
    ( 3, 7)="7";
    "2";"3" **\dir{-};
    "3";"4" **\dir{-};
    "4";"6" **\dir{-};
    "3";"5" **\dir{-};
    "5";"7" **\dir{-};
\endxy}

$ & $

\vcenter{\xy
    ( 0, 6)="1";
    ( 0, 3)*[o]=<9pt>[Fo]{\scriptstyle t}="2";
    ( 0, 0)="3";
    (-3,-3.5)*[o]=<9pt>[Fo]{\scriptstyle t}="4";
    ( 3,-3.5)*[o]=<9pt>[Fo]{\scriptstyle t}="5";
    (-3,-7)="6";
    ( 3,-7)="7";
    "1";"2" **\dir{-};
    "2";"3" **\dir{-};
    "3";"4" **\dir{-};
    "4";"6" **\dir{-};
    "3";"5" **\dir{-};
    "5";"7" **\dir{-};
\endxy}
~=~
\vcenter{\xy
    ( 0, 3)="2";
    ( 0, 0)="3";
    (-3,-3.5)*[o]=<9pt>[Fo]{\scriptstyle t}="4";
    ( 3,-3.5)*[o]=<9pt>[Fo]{\scriptstyle t}="5";
    (-3,-7)="6";
    ( 3,-7)="7";
    "2";"3" **\dir{-};
    "3";"4" **\dir{-};
    "4";"6" **\dir{-};
    "3";"5" **\dir{-};
    "5";"7" **\dir{-};
\endxy}
$ \\

& Under (b) \\

(6) & $
\vcenter{\xy
    (0,8)*{~};
    (0,-8)*{~};
    (0,7)="1";
    (0,4)="2";
    (3,0)*[o]=<9pt>[Fo]{\scriptstyle t}="3";
    (-3,0)="4";
    (0,-4)="5";
    (0,-7)="6";
    "1";"2" **\dir{-};
    "3";"2" **\dir{-};
    "3";"5" **\dir{-};
    "4";"2" **\dir{-};
    "4";"5" **\dir{-};
    "5";"6" **\dir{-};
\endxy}
~=~
\vcenter{\xy
    (0, 6)="1";
    (0,-6)="2";
    "1";"2" **\dir{-};
\endxy}
$ \\

& Under (b) and (w) or (v) \\

(7) & $
\vcenter{\xy
    (0,8)*{~};
    (0,-8)*{~};
    (0,7)="t";
    (0,3)*[o]=<9pt>[Fo]{\scriptstyle t}="m1";
    (0,-3)*[o]=<9pt>[Fo]{\scriptstyle t}="m2";
    (0,-7)="b";
    "t";"m1" **\dir{-};
    "m1";"m2" **\dir{-};
    "m2";"b" **\dir{-};
\endxy}
~=~
\vcenter{\xy
    (0,5)="t";
    (0,0)*[o]=<9pt>[Fo]{\scriptstyle t}="m";
    (0,-5)="b";
    "t";"m" **\dir{-};
    "m";"b" **\dir{-};
\endxy}$

\end{tabular}}
\end{center}
\caption{Properties of $t$} \label{prop-t}
\end{figure}

A table of properties of the source morphism $s$ is given in
Figure~\ref{prop-s} and table of properties of the target morphism $t$ in
Figure~\ref{prop-t}. Properties involving the interaction of $s$ and $t$
are given in Figure~\ref{prop-st}. Proofs of these properties may be found
in Appendix~\ref{sec-proofst}.

\begin{figure}
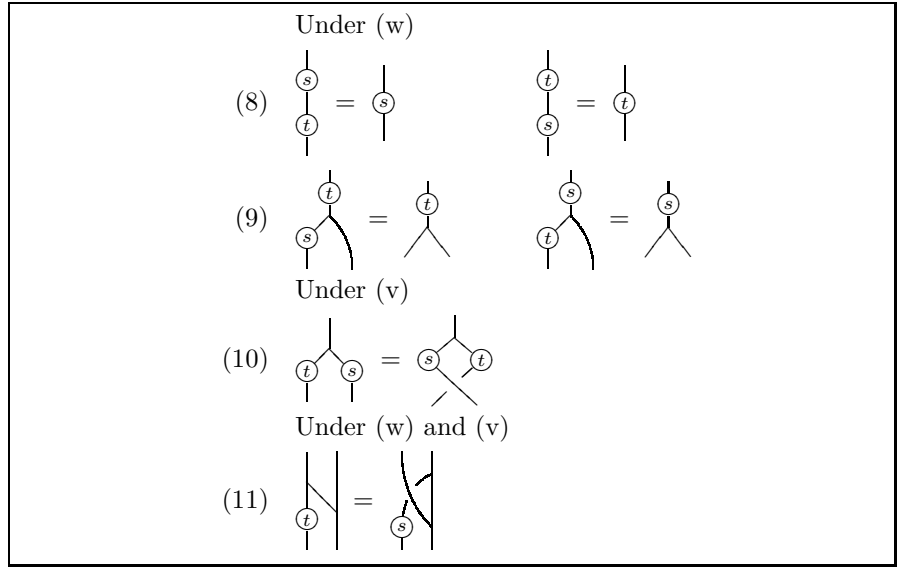
 % Interactions of s and t
\begin{center}
\framebox[\boxwidth]{
\begin{tabular}{rll}
& Under (w) 
\\
(8) & $
\vcenter{\xy
    (0,8)*{~};
    (0,-8)*{~};
    (0,7)="t";
    (0,3)*[o]=<9pt>[Fo]{\scriptstyle s}="m1";
    (0,-3)*[o]=<9pt>[Fo]{\scriptstyle t}="m2";
    (0,-7)="b";
    "t";"m1" **\dir{-};
    "m1";"m2" **\dir{-};
    "m2";"b" **\dir{-};
\endxy}
~=~
\vcenter{\xy
    (0,5)="t";
    (0,0)*[o]=<9pt>[Fo]{\scriptstyle s}="m";
    (0,-5)="b";
    "t";"m" **\dir{-};
    "m";"b" **\dir{-};
\endxy}

$ & $

\vcenter{\xy
    (0,7)="t";
    (0,3)*[o]=<9pt>[Fo]{\scriptstyle t}="m1";
    (0,-3)*[o]=<9pt>[Fo]{\scriptstyle s}="m2";
    (0,-7)="b";
    "t";"m1" **\dir{-};
    "m1";"m2" **\dir{-};
    "m2";"b" **\dir{-};
\endxy}
~=~
\vcenter{\xy
    (0,5)="t";
    (0,0)*[o]=<9pt>[Fo]{\scriptstyle t}="m";
    (0,-5)="b";
    "t";"m" **\dir{-};
    "m";"b" **\dir{-};
\endxy}
$ \\

(9) & $
\vcenter{\xy
    (0,7)*{~};
    (0,-8)*{~};
    (0,6)="t";
    (0,3)*[o]=<9pt>[Fo]{\scriptstyle t}="f";
    (0,-0)="m";
    (-3,-3)*[o]=<9pt>[Fo]{\scriptstyle s}="g";
    (-3,-7)="l";
    ( 3,-7)="r";
    "t";"f" **\dir{-};
    "f";"m" **\dir{-};
    "r";"m" **\crv{(3,-3)};
    "m";"g" **\dir{-};
    "g";"l" **\dir{-};
\endxy}
~=~ 
\vcenter{\xy
    (0,5)="t";
    (0,2)*[o]=<9pt>[Fo]{\scriptstyle t}="f";
    (0,-1)="m";
    (-3,-5)="l";
    ( 3,-5)="r";
    "t";"f" **\dir{-};
    "f";"m" **\dir{-};
    "m";"l" **\dir{-};
    "m";"r" **\dir{-};
\endxy}

$ & $

\vcenter{\xy
    (0,6)="t";
    (0,3)*[o]=<9pt>[Fo]{\scriptstyle s}="f";
    (0,-0)="m";
    (-3,-3)*[o]=<9pt>[Fo]{\scriptstyle t}="g";
    (-3,-7)="l";
    ( 3,-7)="r";
    "t";"f" **\dir{-};
    "f";"m" **\dir{-};
    "r";"m" **\crv{(3,-3)};
    "m";"g" **\dir{-};
    "g";"l" **\dir{-};
\endxy}
~=~
\vcenter{\xy
    (0,5)="t";
    (0,2)*[o]=<9pt>[Fo]{\scriptstyle s}="f";
    (0,-1)="m";
    (-3,-5)="l";
    ( 3,-5)="r";
    "t";"f" **\dir{-};
    "f";"m" **\dir{-};
    "m";"l" **\dir{-};
    "m";"r" **\dir{-};
\endxy}
$ \\
& Under (v) \\

(10) & $
\vcenter{\xy
    (0,5)*{~};
    (0,-8)*{~};
    ( 0, 4)="2";
    ( 0, 0)="3";
    (-3,-3)*[o]=<9pt>[Fo]{\scriptstyle t}="4";
    ( 3,-3)*[o]=<9pt>[Fo]{\scriptstyle s}="5";
    (-3,-7)="6";
    ( 3,-7)="7";
    "2";"3" **\dir{-};
    "3";"4" **\dir{-};
    "4";"6" **\dir{-};
    "3";"5" **\dir{-};
    "5";"7" **\dir{-};
\endxy}
~=~
\vcenter{\xy
    (0,5)*{~};
    (0,-9)*{~};
    ( 0, 4)="2";
    ( 0, 1)="3";
    (-3.5,-2)*[o]=<9pt>[Fo]{\scriptstyle s}="4";
    ( 3.5,-2)*[o]=<9pt>[Fo]{\scriptstyle t}="5";
    (-3,-8)="6";
    ( 3,-8)="7";
    "2";"3" **\dir{-};
    "3";"4" **\dir{-};
    "3";"5" **\dir{-};
    "4";"7" **\dir{-}?(0.53)*{\hole}="x";
    "5";"x" **\dir{-};
    "x";"6" **\dir{-};
\endxy}
$ \\

& Under (w) and (v) \\
(11) & $
\vcenter{\xy
    (0,9)*{~};
    (0,-6)*{~};
    (-2,8)="1";
    ( 2,8)="2";
    (-2,4)="3";
    ( 2,0)="4";
    (-2,-1)*[o]=<9pt>[Fo]{\scriptstyle t}="5";
    (-2,-5)="6";
    ( 2,-5)="7";
    "1";"5" **\dir{-};
    "5";"6" **\dir{-};
    "2";"7" **\dir{-};
    "3";"4" **\dir{-};
\endxy}
~=~
\vcenter{\xy
    (0,9)*{~};
    (0,-6)*{~};
    (-2,8)="1";
    ( 2,8)="2";
    ( 2,5)="3";
    ( 2,-2)="4";
    (-2,-5)="5";
    ( 2,-5)="6";
    (-2,-2)*[o]=<9pt>[Fo]{\scriptstyle s}="7";
    "2";"6" **\dir{-};
    "1";"4" **\crv{(-2,2)}?(0.5)*{\hole}="x";
    "3";"x" **\crv{(1,5)};
    "5";"7" **\dir{-};
    "7";"x" **\crv{(-2,0)};
\endxy}
$
\end{tabular}}
\end{center}
\caption{Interactions of $s$ and $t$} \label{prop-st}
\end{figure}

In the sequel $A=(A,\mu,\eta,\delta,\epsilon)$ will always denote a weak
bimonoid and $s,t:A \ra A$ the source and target morphisms.

We see from property (7) in Figures~\ref{prop-s} and~\ref{prop-t}
respectively that both $s$ and $t$ are idempotents. In the following we will
work in the Cauchy completion (= completion under idempotents = Karoubi
envelope) $\Q\V$ of $\V$. We do this rather than assume that idempotents
split in $\V$.

%=========================================================================%
\subsection{Cauchy completion}
%=========================================================================%

Given a category $\V$, its \emph{Cauchy completion} $\Q\V$ is the category
whose objects are pairs $(X,e)$ with $X \in \V$ and $e:X \ra X \in \V$ an
idempotent. A morphism $(X,e) \ra (X',e')$ in $\Q\V$ is a morphism
$f:X \ra X' \in \V$ such that $e'fe = f$. Note that the identity morphism
of $(X,e)$ is $e$ itself.

The point of working in the Cauchy completion is that every idempotent
$f:(X,e) \ra (X,e)$ in $\Q\V$ has a splitting, viz.,
\[
    \xygraph{{(X,e)}="1"
        [r(2.2)]{(X,e)}="2"
        [d(1.2)l(1.1)]{(X,f).}="3"
        "1":"2" ^-f 
        "1":"3" _-f 
        "3":"2" _-f}
\]
If $\V$ is a monoidal category, then $\Q\V$ is a monoidal category via
\[
    (X,e) \ox (X',e') = (X \ox X',e \ox e').
\]

The category $\V$ may be fully embedded in $\Q\V$ by sending $X \in \V$ to
$(X,1) \in \Q\V$ and $f:X \ra Y \in \V$ to $f:(X,1) \ra (Y,1)$, which is
obviously a morphism in $\Q\V$. When working in $\Q\V$ we will often
identify an object $X \in \V$ with $(X,1) \in \Q\V$.

%=========================================================================%
\subsection{Properties of the source and target morphisms}
%=========================================================================%

Let $A = (A,1)$ be a weak bimonoid in $\Q\V$. From the definition of the Cauchy
completion the result of splitting the source morphism $s$ is $(A,s)$, and
similarly, the result of splitting the target morphism $t$ is $(A,t)$. The
following proposition shows that these two objects are isomorphic.

\begin{proposition}
The idempotent $t:(A,1) \ra (A,1)$ has the following two splittings.
\[
\xygraph{{(A,1)}="1"
    [r(2.2)]{(A,1)}="2"
    [d(1.2)l(1.1)]{(A,t)}="3"
    "1":"2" ^-t 
    "1":"3" _-t 
    "3":"2" _-t}
\qquad
\xygraph{{(A,1)}="1"
    [r(2.2)]{(A,1)}="2"
    [d(1.2)l(1.1)]{(A,s)}="3"
    "1":"2" ^-t 
    "1":"3" _-t 
    "3":"2" _-s}
\]
In this case $s:(A,s) \ra (A,t)$ and $t:(A,t) \ra (A,s)$ are inverse
morphisms, and hence $(A,t) \cong (A,s)$.
\end{proposition}

\begin{proof}
This result follows from the identities $ts=s$ and $st=t$ (property (8) in
Figure~\ref{prop-st}).
\end{proof}

We will denote this object by $C = (A,t)$ and call it the
\emph{object-of-objects} of $A$. In the next propositions we will show that
$C$ is a comonoid, and furthermore, that it is a separable Frobenius monoid,
similar to what was done in~\cite{Sch} (there called Frobenius-separable).

\begin{proposition} \label{prop-Ccomon}
The object $C = (A,t)$ equipped with
\begin{align*}
\delta &= \big(\xygraph{{C}
    :[r(1.4)] {C \ox C} ^-\delta
    :[r(1.8)] {C \ox C} ^-{t \ox t}}\big) \\
\epsilon &= \xygraph{{C} :[r(1)] {I} ^-\epsilon}
\end{align*}
is a comonoid in $\Q\V$, and if furthermore equipped with
\begin{align*}
\mu &= \big(\xygraph{{C \ox C}
        :[r(1.8)] {C \ox C} ^-{t \ox t}
        :[r(1.4)] {C} ^-\mu}\big) \\
\eta &= \xygraph{{I} :[r(1)] {C} ^-\eta}
\end{align*}
then $C$ is a separable Frobenius monoid in $\Q\V$ (see
Definition~\ref{def-sepfrob}).
\end{proposition}

\begin{proof}
We first observe that $(t \ox t)\delta: C \ra C \ox C$ and $\epsilon:C \ra
I$ are in $\Q\V$ which follows from~(5) and~(2) respectively.

The comonoid identities are given as
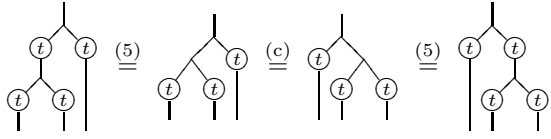
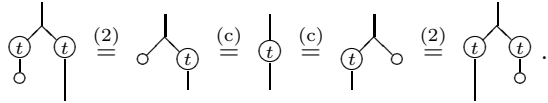
\[
    \vcenter{\xy
        ( 0, 9)="1";
        ( 0, 6)="2";
        (-3, 3)*[o]=<9pt>[Fo]{\scriptstyle t}="3";
        ( 3, 3)*[o]=<9pt>[Fo]{\scriptstyle t}="4";
        (-3,-1)="5";
        (-6,-4)*[o]=<9pt>[Fo]{\scriptstyle t}="6";
        ( 0,-4)*[o]=<9pt>[Fo]{\scriptstyle t}="7";
        (-6,-8)="8";
        ( 0,-8)="9";
        ( 3,-8)="10";
        "1";"2" **\dir{-};
        "2";"3" **\dir{-};
        "2";"4" **\dir{-};
        "3";"5" **\dir{-};
        "5";"6" **\dir{-};
        "5";"7" **\dir{-};
        "6";"8" **\dir{-};
        "7";"9" **\dir{-};
        "4";"10" **\dir{-};
    \endxy}
    ~\overset{\text{(5)}}{=}~
    \vcenter{\xy
        ( 0, 9)="1";
        ( 0, 6)="2";
        ( 3, 3)*[o]=<9pt>[Fo]{\scriptstyle t}="4";
        (-3, 3)="5";
        (-6,-1)*[o]=<9pt>[Fo]{\scriptstyle t}="6";
        ( 0,-1)*[o]=<9pt>[Fo]{\scriptstyle t}="7";
        (-6,-5)="8";
        ( 0,-5)="9";
        ( 3,-5)="10";
        "1";"2" **\dir{-};
        "2";"5" **\dir{-};
        "2";"4" **\dir{-};
        "5";"6" **\dir{-};
        "5";"7" **\dir{-};
        "6";"8" **\dir{-};
        "7";"9" **\dir{-};
        "4";"10" **\dir{-};
    \endxy}
    ~\overset{\text{(c)}}{=}~
    \vcenter{\xy
        ( 0, 9)="1";
        ( 0, 6)="2";
        (-3, 3)*[o]=<9pt>[Fo]{\scriptstyle t}="4";
        ( 3, 3)="5";
        ( 6,-1)*[o]=<9pt>[Fo]{\scriptstyle t}="6";
        ( 0,-1)*[o]=<9pt>[Fo]{\scriptstyle t}="7";
        ( 6,-5)="8";
        ( 0,-5)="9";
        (-3,-5)="10";
        "1";"2" **\dir{-};
        "2";"5" **\dir{-};
        "2";"4" **\dir{-};
        "5";"6" **\dir{-};
        "5";"7" **\dir{-};
        "6";"8" **\dir{-};
        "7";"9" **\dir{-};
        "4";"10" **\dir{-};
    \endxy}
    ~\overset{\text{(5)}}{=}~
    \vcenter{\xy
        ( 0, 9)="1";
        ( 0, 6)="2";
        ( 3, 3)*[o]=<9pt>[Fo]{\scriptstyle t}="3";
        (-3, 3)*[o]=<9pt>[Fo]{\scriptstyle t}="4";
        ( 3,-1)="5";
        ( 6,-4)*[o]=<9pt>[Fo]{\scriptstyle t}="6";
        ( 0,-4)*[o]=<9pt>[Fo]{\scriptstyle t}="7";
        ( 6,-8)="8";
        ( 0,-8)="9";
        (-3,-8)="10";
        "1";"2" **\dir{-};
        "2";"3" **\dir{-};
        "2";"4" **\dir{-};
        "3";"5" **\dir{-};
        "5";"6" **\dir{-};
        "5";"7" **\dir{-};
        "6";"8" **\dir{-};
        "7";"9" **\dir{-};
        "4";"10" **\dir{-};
    \endxy}
\]
and
\[
    \vcenter{\xy
        ( 0, 3)="2";
        ( 0, 0)="3";
        (-3,-3)*[o]=<9pt>[Fo]{\scriptstyle t}="4";
        ( 3,-3)*[o]=<9pt>[Fo]{\scriptstyle t}="5";
        (-3,-7)*[o]=<5pt>[Fo]{~}="6";
        ( 3,-10)="7";
        "2";"3" **\dir{-};
        "3";"4" **\dir{-};
        "4";"6" **\dir{-};
        "3";"5" **\dir{-};
        "7";"5" **\dir{-};
    \endxy}
    ~\overset{\text{(2)}}{=}~
    \vcenter{\xy
        ( 0, 3)="2";
        ( 0, 0)="3";
        (-3,-3)*[o]=<5pt>[Fo]{~}="4";
        ( 3,-3)*[o]=<9pt>[Fo]{\scriptstyle t}="5";
        ( 3,-7)="7";
        "2";"3" **\dir{-};
        "3";"4" **\dir{-};
        "3";"5" **\dir{-};
        "7";"5" **\dir{-};
    \endxy}
    ~\overset{\text{(c)}}{=}~
    \vcenter{\xy
        (0,5)="1";
        (0,0)*[o]=<9pt>[Fo]{\scriptstyle t}="2";
        (0,-5)="3";
        "1";"2" **\dir{-};
        "2";"3" **\dir{-};
    \endxy}
    ~\overset{\text{(c)}}{=}~
    \vcenter{\xy
        ( 0, 3)="2";
        ( 0, 0)="3";
        ( 3,-3)*[o]=<5pt>[Fo]{~}="4";
        (-3,-3)*[o]=<9pt>[Fo]{\scriptstyle t}="5";
        (-3,-7)="7";
        "2";"3" **\dir{-};
        "3";"4" **\dir{-};
        "3";"5" **\dir{-};
        "7";"5" **\dir{-};
        \endxy}
    ~\overset{\text{(2)}}{=}~
    \vcenter{\xy
        ( 0, 3)="2";
        ( 0, 0)="3";
        ( 3,-3)*[o]=<9pt>[Fo]{\scriptstyle t}="4";
        (-3,-3)*[o]=<9pt>[Fo]{\scriptstyle t}="5";
        ( 3,-7)*[o]=<5pt>[Fo]{~}="6";
        (-3,-10)="7";
        "2";"3" **\dir{-};
        "3";"4" **\dir{-};
        "4";"6" **\dir{-};
        "3";"5" **\dir{-};
        "7";"5" **\dir{-};
    \endxy}~.
\]

To see that $C$ is a separable Frobenius monoid we first observe that $\mu$
and $\eta$ are morphisms in $\Q\V$ from (5) and (2), and the monoid
identities are dual to the comonoid identities. The following calculation
proves that the Frobenius condition holds.
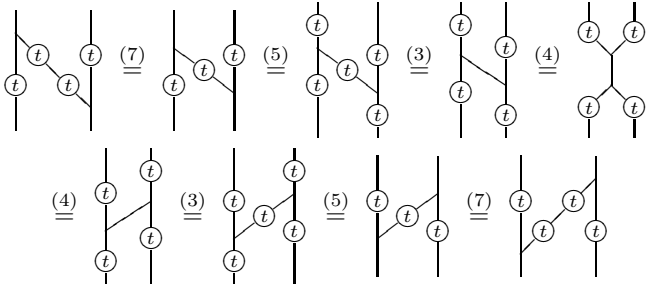
\begin{gather*}
   \vcenter{\xy
        (-5, 8)="a";
        (-5,-2)*[o]=<9pt>[Fo]{\scriptstyle t}="b";
        (-5,-8)="c";
        ( 5, 8)="d";
        ( 5, 2)*[o]=<9pt>[Fo]{\scriptstyle t}="e";
        ( 5,-8)="f";
        (-5, 5)="2";
        (-2, 2)*[o]=<9pt>[Fo]{\scriptstyle t}="3";
        ( 2,-2)*[o]=<9pt>[Fo]{\scriptstyle t}="4";
        ( 5,-5)="5";
        "a";"b" **\dir{-};
        "b";"c" **\dir{-};
        "d";"e" **\dir{-};
        "e";"f" **\dir{-};
        "2";"3" **\dir{-};
        "3";"4" **\dir{-};
        "4";"5" **\dir{-};
    \endxy}
    ~\overset{\text{(7)}}{=}~
    \vcenter{\xy
        (-4, 8)="a";
        (-4,-2)*[o]=<9pt>[Fo]{\scriptstyle t}="b";
        (-4,-8)="c";
        ( 4, 8)="d";
        ( 4, 2)*[o]=<9pt>[Fo]{\scriptstyle t}="e";
        ( 4,-8)="f";
        (-4, 3)="2";
        ( 0, 0)*[o]=<9pt>[Fo]{\scriptstyle t}="3";
        ( 4,-3)="4";
        "a";"b" **\dir{-};
        "b";"c" **\dir{-};
        "d";"e" **\dir{-};
        "e";"f" **\dir{-};
        "2";"3" **\dir{-};
        "3";"4" **\dir{-};
    \endxy}
    ~\overset{\text{(5)}}{=}~
    \vcenter{\xy
        (-4, 9)="a";
        (-4, 6)*[o]=<9pt>[Fo]{\scriptstyle t}="x";
        (-4,-2)*[o]=<9pt>[Fo]{\scriptstyle t}="b";
        (-4,-9)="c";
        ( 4, 9)="d";
        ( 4, 2)*[o]=<9pt>[Fo]{\scriptstyle t}="e";
        ( 4,-6)*[o]=<9pt>[Fo]{\scriptstyle t}="y";
        ( 4,-9)="f";
        (-4, 3)="2";
        ( 0, 0)*[o]=<9pt>[Fo]{\scriptstyle t}="3";
        ( 4,-3)="4";
        "a";"x" **\dir{-};
        "x";"b" **\dir{-};
        "b";"c" **\dir{-};
        "d";"e" **\dir{-};
        "e";"y" **\dir{-};
        "y";"f" **\dir{-};
        "2";"3" **\dir{-};
        "3";"4" **\dir{-};
    \endxy}
    ~\overset{\text{(3)}}{=}~
    \vcenter{\xy
        (-3, 9)="a";
        (-3, 6)*[o]=<9pt>[Fo]{\scriptstyle t}="x";
        (-3,-3)*[o]=<9pt>[Fo]{\scriptstyle t}="b";
        (-3,-9)="c";
        ( 3, 9)="d";
        ( 3, 3)*[o]=<9pt>[Fo]{\scriptstyle t}="e";
        ( 3,-6)*[o]=<9pt>[Fo]{\scriptstyle t}="y";
        ( 3,-9)="f";
        (-3, 2)="2";
        ( 3,-2)="4";
        "a";"x" **\dir{-};
        "x";"b" **\dir{-};
        "b";"c" **\dir{-};
        "d";"e" **\dir{-};
        "e";"y" **\dir{-};
        "y";"f" **\dir{-};
        "2";"4" **\dir{-};
    \endxy}
    ~\overset{\text{(4)}}{=}~
    \vcenter{\xy
        (-3, 9)="a";
        (-3, 5)*[o]=<9pt>[Fo]{\scriptstyle t}="b";
        (-3,-5)*[o]=<9pt>[Fo]{\scriptstyle t}="c";
        (-3,-9)="d";
        ( 3, 9)="w";
        ( 3, 5)*[o]=<9pt>[Fo]{\scriptstyle t}="x";
        ( 3,-5)*[o]=<9pt>[Fo]{\scriptstyle t}="y";
        ( 3,-9)="z";
        ( 0, 2)="1";
        ( 0,-2)="2";
        "a";"b" **\dir{-};
        "c";"d" **\dir{-};
        "w";"x" **\dir{-};
        "y";"z" **\dir{-};
        "b";"1" **\dir{-};
        "x";"1" **\dir{-};
        "c";"2" **\dir{-};
        "y";"2" **\dir{-};
        "1";"2" **\dir{-};
    \endxy}
\\
    ~\overset{\text{(4)}}{=}~
    \vcenter{\xy
        ( 3, 9)="a";
        ( 3, 6)*[o]=<9pt>[Fo]{\scriptstyle t}="x";
        ( 3,-3)*[o]=<9pt>[Fo]{\scriptstyle t}="b";
        ( 3,-9)="c";
        (-3, 9)="d";
        (-3, 3)*[o]=<9pt>[Fo]{\scriptstyle t}="e";
        (-3,-6)*[o]=<9pt>[Fo]{\scriptstyle t}="y";
        (-3,-9)="f";
        ( 3, 2)="2";
        (-3,-2)="4";
        "a";"x" **\dir{-};
        "x";"b" **\dir{-};
        "b";"c" **\dir{-};
        "d";"e" **\dir{-};
        "e";"y" **\dir{-};
        "y";"f" **\dir{-};
        "2";"4" **\dir{-};
    \endxy}
    ~\overset{\text{(3)}}{=}~
    \vcenter{\xy
        ( 4, 9)="a";
        ( 4, 6)*[o]=<9pt>[Fo]{\scriptstyle t}="x";
        ( 4,-2)*[o]=<9pt>[Fo]{\scriptstyle t}="b";
        ( 4,-9)="c";
        (-4, 9)="d";
        (-4, 2)*[o]=<9pt>[Fo]{\scriptstyle t}="e";
        (-4,-6)*[o]=<9pt>[Fo]{\scriptstyle t}="y";
        (-4,-9)="f";
        ( 4, 3)="2";
        ( 0, 0)*[o]=<9pt>[Fo]{\scriptstyle t}="3";
        (-4,-3)="4";
        "a";"x" **\dir{-};
        "x";"b" **\dir{-};
        "b";"c" **\dir{-};
        "d";"e" **\dir{-};
        "e";"y" **\dir{-};
        "y";"f" **\dir{-};
        "2";"3" **\dir{-};
        "3";"4" **\dir{-};
    \endxy}
    ~\overset{\text{(5)}}{=}~
    \vcenter{\xy
        ( 4, 8)="a";
        ( 4,-2)*[o]=<9pt>[Fo]{\scriptstyle t}="b";
        ( 4,-8)="c";
        (-4, 8)="d";
        (-4, 2)*[o]=<9pt>[Fo]{\scriptstyle t}="e";
        (-4,-8)="f";
        ( 4, 3)="2";
        ( 0, 0)*[o]=<9pt>[Fo]{\scriptstyle t}="3";
        (-4,-3)="4";
        "a";"b" **\dir{-};
        "b";"c" **\dir{-};
        "d";"e" **\dir{-};
        "e";"f" **\dir{-};
        "2";"3" **\dir{-};
        "3";"4" **\dir{-};
    \endxy}
    ~\overset{\text{(7)}}{=}~
    \vcenter{\xy
        ( 5, 8)="a";
        ( 5,-2)*[o]=<9pt>[Fo]{\scriptstyle t}="b";
        ( 5,-8)="c";
        (-5, 8)="d";
        (-5, 2)*[o]=<9pt>[Fo]{\scriptstyle t}="e";
        (-5,-8)="f";
        ( 5, 5)="2";
        ( 2, 2)*[o]=<9pt>[Fo]{\scriptstyle t}="3";
        (-2,-2)*[o]=<9pt>[Fo]{\scriptstyle t}="4";
        (-5,-5)="5";
        "a";"b" **\dir{-};
        "b";"c" **\dir{-};
        "d";"e" **\dir{-};
        "e";"f" **\dir{-};
        "2";"3" **\dir{-};
        "3";"4" **\dir{-};
        "4";"5" **\dir{-};
    \endxy}
\end{gather*}
Finally, that this is a separable Frobenius monoid follows from
\[
    \mu\delta =~
    \vcenter{\xy
        (0,9)="1";
        (0,6)="2";
        (-3,3)*[o]=<9pt>[Fo]{\scriptstyle t}="a";
        (3,3)*[o]=<9pt>[Fo]{\scriptstyle t}="b";
        (-3,-3)*[o]=<9pt>[Fo]{\scriptstyle t}="x";
        (3,-3)*[o]=<9pt>[Fo]{\scriptstyle t}="y";
        (0,-6)="3";
        (0,-9)="4";
        "1";"2" **\dir{-};
        "2";"a" **\dir{-};
        "2";"b" **\dir{-};
        "a";"x" **\dir{-};
        "b";"y" **\dir{-};
        "3";"x" **\dir{-};
        "3";"y" **\dir{-};
        "3";"4" **\dir{-};
    \endxy}
    ~\overset{\text{(7)}}{=}~
    \vcenter{\xy
        (0,7)="1";
        (0,4)="2";
        (-3,0)*[o]=<9pt>[Fo]{\scriptstyle t}="a";
        (3,0)*[o]=<9pt>[Fo]{\scriptstyle t}="b";
        (0,-4)="3";
        (0,-7)="4";
        "1";"2" **\dir{-};
        "2";"a" **\dir{-};
        "2";"b" **\dir{-};
        "3";"a" **\dir{-};
        "3";"b" **\dir{-};
        "3";"4" **\dir{-};
    \endxy}
    ~\overset{\text{(5)}}{=}~
    \vcenter{\xy
        (0,11)="1";
        (0,7)*[o]=<9pt>[Fo]{\scriptstyle t}="x";
        (0,4)="2";
        (-3,0)*[o]=<9pt>[Fo]{\scriptstyle t}="a";
        (3,0)*[o]=<9pt>[Fo]{\scriptstyle t}="b";
        (0,-4)="3";
        (0,-7)="4";
        "1";"x" **\dir{-};
        "x";"2" **\dir{-};
        "2";"a" **\dir{-};
        "2";"b" **\dir{-};
        "3";"a" **\dir{-};
        "3";"b" **\dir{-};
        "3";"4" **\dir{-};
    \endxy}
    ~\overset{\text{(3)}}{=}~
    \vcenter{\xy
        (0,11)="1";
        (0,7)*[o]=<9pt>[Fo]{\scriptstyle t}="x";
        (0,4)="2";
        (-3,0)="a";
        (3,0)*[o]=<9pt>[Fo]{\scriptstyle t}="b";
        (0,-4)="3";
        (0,-7)="4";
        "1";"x" **\dir{-};
        "x";"2" **\dir{-};
        "2";"a" **\dir{-};
        "2";"b" **\dir{-};
        "3";"a" **\dir{-};
        "3";"b" **\dir{-};
        "3";"4" **\dir{-};
    \endxy}
    ~\overset{\text{(6)}}{=}~
    \vcenter{\xy
        (0,7)="1";
        (0,0)*[o]=<9pt>[Fo]{\scriptstyle t}="x";
        (0,-7)="2";
        "1";"x" **\dir{-};
        "x";"2" **\dir{-};
    \endxy}
    ~= 1_C.
\]
\end{proof}

\begin{corollary}
Every morphism of weak bimonoids induces an isomorphism on the
``objects-of-objects''. That is, if $(A,1)$ and $(B,1)$ are weak bimonoids,
and $f:(A,1) \ra (B,1)$ is a morphism of weak bimonoids, then the induced
morphism $tft:(A,t) \ra (B,t)$ is an isomorphism. 
\end{corollary}

\begin{proof}
Note that if $f:A \ra B$ is a morphism of weak bimonoids then $ft = tf$ and
$fs = st$. The corollary now follows from Proposition~\ref{prop-Ccomon} and
Proposition~\ref{prop-frobiso}.
\end{proof}

\begin{proposition}\label{prop-stcm}
If we write $C^\circ$ for the comonoid $C$ with the ``opposite''
comultiplication defined via
\[
    \xygraph{{C}
        :[r(1.4)]{C \ox C} ^-\delta
        :[r(1.8)]{C \ox C} ^-{t \ox t}
        :[r(1.6)]{C \ox C} ^-c}
    =~
    \vcenter{\xy
        ( 0, 5)="2";
        ( 0, 1)="3";
        (-3.5,-2)*[o]=<9pt>[Fo]{\scriptstyle t}="4";
        ( 3.5,-2)*[o]=<9pt>[Fo]{\scriptstyle t}="5";
        (-3,-8)="6";
        ( 3,-8)="7";
        "2";"3" **\dir{-};
        "3";"4" **\dir{-};
        "3";"5" **\dir{-};
        "4";"7" **\dir{-}?(0.53)*{\hole}="x";
        "5";"x" **\dir{-};
        "x";"6" **\dir{-};
    \endxy}
\]
then $s:A \ra C^\circ$ and $t:A \ra C$ are comonoid morphisms. That is, the
diagrams
\[
    \vcenter{\xygraph{{A}="1"
        [r(2)] {C}="2"
        [d(1.2)] {C \ox C}="3"
        [l(2)] {A \ox A}="4"
        "1":"2" ^-s
        "2":"3" ^-{c(t \ox t)\delta}
        "1":"4" _-\delta
        "4":"3" ^-{s \ox s}}}
    \qquad\qquad
    \vcenter{\xygraph{{A}="1"
        [r(2)] {C}="2"
        [d(1.2)] {C \ox C}="3"
        [l(2)] {A \ox A}="4"
        "1":"2" ^-t
        "2":"3" ^-{(t \ox t)\delta}
        "1":"4" _-\delta
        "4":"3" ^-{t \ox t}}}
\]
commute.
\end{proposition}

\begin{proof}
The second diagram expresses
\[
    \vcenter{\xy
        ( 0, 6)="1";
        ( 0, 3)*[o]=<9pt>[Fo]{\scriptstyle t}="2";
        ( 0, 0)="3";
        (-3,-3)*[o]=<9pt>[Fo]{\scriptstyle t}="4";
        ( 3,-3)*[o]=<9pt>[Fo]{\scriptstyle t}="5";
        (-3,-7)="6";
        ( 3,-7)="7";
        "1";"2" **\dir{-};
        "2";"3" **\dir{-};
        "3";"4" **\dir{-};
        "4";"6" **\dir{-};
        "3";"5" **\dir{-};
        "5";"7" **\dir{-};
    \endxy}
    ~=~
    \vcenter{\xy
        ( 0, 3)="2";
        ( 0, 0)="3";
        (-3,-3)*[o]=<9pt>[Fo]{\scriptstyle t}="4";
        ( 3,-3)*[o]=<9pt>[Fo]{\scriptstyle t}="5";
        (-3,-7)="6";
        ( 3,-7)="7";
        "2";"3" **\dir{-};
        "3";"4" **\dir{-};
        "4";"6" **\dir{-};
        "3";"5" **\dir{-};
        "5";"7" **\dir{-};
    \endxy}
\]
which is exactly (5), and the following calculation
\[
    \vcenter{\xy
        ( 0, 3)="2";
        ( 0, 0)="3";
        (-3,-3)*[o]=<9pt>[Fo]{\scriptstyle s}="4";
        ( 3,-3)*[o]=<9pt>[Fo]{\scriptstyle s}="5";
        (-3,-7)="6";
        ( 3,-7)="7";
        "2";"3" **\dir{-};
        "3";"4" **\dir{-};
        "4";"6" **\dir{-};
        "3";"5" **\dir{-};
        "5";"7" **\dir{-};
    \endxy}
    ~\overset{\text{(5)}}{=}~
    \vcenter{\xy
        ( 0, 6)="1";
        ( 0, 3)*[o]=<9pt>[Fo]{\scriptstyle s}="2";
        ( 0, 0)="3";
        (-3,-3)*[o]=<9pt>[Fo]{\scriptstyle s}="4";
        ( 3,-3)*[o]=<9pt>[Fo]{\scriptstyle s}="5";
        (-3,-7)="6";
        ( 3,-7)="7";
        "1";"2" **\dir{-};
        "2";"3" **\dir{-};
        "3";"4" **\dir{-};
        "4";"6" **\dir{-};
        "3";"5" **\dir{-};
        "5";"7" **\dir{-};
    \endxy}
    ~\overset{\text{(8)}}{=}~
    \vcenter{\xy
        ( 0, 6)="1";
        ( 0, 3)*[o]=<9pt>[Fo]{\scriptstyle s}="2";
        ( 0, 0)="3";
        (-3,-3)*[o]=<9pt>[Fo]{\scriptstyle s}="4";
        ( 3,-3)*[o]=<9pt>[Fo]{\scriptstyle s}="5";
        (-3,-8)*[o]=<9pt>[Fo]{\scriptstyle t}="6";
        (-3,-12)="7";
        ( 3,-12)="8";
        "1";"2" **\dir{-};
        "2";"3" **\dir{-};
        "3";"4" **\dir{-};
        "3";"5" **\dir{-};
        "6";"4" **\dir{-};
        "7";"6" **\dir{-};
        "8";"5" **\dir{-};
    \endxy}
    ~\overset{\text{(3)}}{=}~
    \vcenter{\xy
        ( 0, 6)="1";
        ( 0, 3)*[o]=<9pt>[Fo]{\scriptstyle s}="2";
        ( 0, 0)="3";
        (-3,-3)*[o]=<9pt>[Fo]{\scriptstyle t}="4";
        ( 3,-3)*[o]=<9pt>[Fo]{\scriptstyle s}="5";
        (-3,-7)="6";
        ( 3,-7)="7";
        "1";"2" **\dir{-};
        "2";"3" **\dir{-};
        "3";"4" **\dir{-};
        "4";"6" **\dir{-};
        "3";"5" **\dir{-};
        "5";"7" **\dir{-};
    \endxy}
    ~\overset{\text{(10)}}{=}~
    \vcenter{\xy
        ( 0, 7)="1";
        ( 0, 4)*[o]=<9pt>[Fo]{\scriptstyle s}="2";
        ( 0, 1)="3";
        (-3.5,-2)*[o]=<9pt>[Fo]{\scriptstyle s}="4";
        ( 3.5,-2)*[o]=<9pt>[Fo]{\scriptstyle t}="5";
        (-3,-8)="6";
        ( 3,-8)="7";
        "1";"2" **\dir{-};
        "2";"3" **\dir{-};
        "3";"4" **\dir{-};
        "3";"5" **\dir{-};
        "4";"7" **\dir{-}?(0.53)*{\hole}="x";
        "5";"x" **\dir{-};
        "x";"6" **\dir{-};
    \endxy}
    ~\overset{\text{(9)}}{=}~
    \vcenter{\xy
        ( 0, 7)="1";
        ( 0, 4)*[o]=<9pt>[Fo]{\scriptstyle s}="2";
        ( 0, 1)="3";
        (-3.5,-2)*[o]=<9pt>[Fo]{\scriptstyle t}="4";
        ( 3.5,-2)*[o]=<9pt>[Fo]{\scriptstyle t}="5";
        (-3.5,-7)*[o]=<9pt>[Fo]{\scriptstyle s}="a";
        ( 3.5,-7)="b";
        (-3,-12)="6";
        ( 3,-12)="7";
        "1";"2" **\dir{-};
        "2";"3" **\dir{-};
        "3";"4" **\dir{-};
        "3";"5" **\dir{-};
        "4";"a" **\dir{-};
        "5";"b" **\dir{-};
        "a";"7" **\dir{-}?(0.53)*{\hole}="x";
        "b";"x" **\dir{-};
        "x";"6" **\dir{-};
    \endxy}
    ~\overset{\text{(8)}}{=}~
    \vcenter{\xy
        ( 0, 7)="1";
        ( 0, 4)*[o]=<9pt>[Fo]{\scriptstyle s}="2";
        ( 0, 1)="3";
        (-3.5,-2)*[o]=<9pt>[Fo]{\scriptstyle t}="4";
        ( 3.5,-2)*[o]=<9pt>[Fo]{\scriptstyle t}="5";
        (-3,-8)="6";
        ( 3,-8)="7";
        "1";"2" **\dir{-};
        "2";"3" **\dir{-};
        "3";"4" **\dir{-};
        "3";"5" **\dir{-};
        "4";"7" **\dir{-}?(0.53)*{\hole}="x";
        "5";"x" **\dir{-};
        "x";"6" **\dir{-};
    \endxy}
\]
shows that the first diagram commutes.
\end{proof}

%=========================================================================%
\section{Weak Hopf monoids} \label{sec-whm}
%=========================================================================%

In this section we introduce weak Hopf monoids. Usually in the literature, a
weak Hopf monoid is a weak bimonoid $H$ equipped with an antipode $\nu:H \ra
H$ satisfying the three axioms
\[
  \nu * 1 = t, \qquad 1*\nu = s, \qquad \text{and} \qquad \nu * 1 * \nu = \nu,
\]
where $f * g = \mu(f \ox g)\delta$ is the convolution product. Our
definition is slightly different as, instead of choosing our source morphism
in the second axiom, we replace it with
\[
    1*\nu = r,
\]
where $r$ is defined below. This turns out to be the usual definition of
weak Hopf monoids as found in the literature; in the symmetric case
see~\cite{BNS}, and in the braided case see~\cite{AFG1,AFG2}.

%=========================================================================%
\subsection{The endomorphism $r$ and weak Hopf monoids}
%=========================================================================%

Define an endomorphism $r:A \ra A$ by rotating the target morphism $t:A \ra
A$ by $\pi$, i.e.,
\[
    r =~
    \vcenter{\xy
        ( 5,0)*{~};
        ( 4,-6)="0";
        ( 0,-4)*[o]=<5pt>[Fo]{~}="1";
        ( 0,-2)="2";
        ( 0, 2)="3";
        ( 0, 4)*[o]=<5pt>[Fo]{~}="4";
        ( 4, 6)="5";
        "1";"4" **\dir{-};
        "0";"3" **\crv{( 4,0)}?(0.6)*{\hole}="x";
        "5";"x" **\crv{( 4,2)};
        "2";"x" **\crv{( 2,-1)};
    \endxy}~.
\]
Since $r$ is just $t$ rotated by $\pi$, all the identities for $t$ in
Figure~\ref{prop-t} rotated by $\pi$ hold for $r$. We list some additional
identities of $r$ interacting with $s$ and $t$.

\begin{equation}\tag{12}
    \vcenter{\xy
        (0,7)="1";
        (0,3)*[o]=<9pt>[Fo]{\scriptstyle r}="2";
        (0,-3)*[o]=<9pt>[Fo]{\scriptstyle s}="3";
        (0,-7)="4";
        "1";"2" **\dir{-};
        "2";"3" **\dir{-};
        "3";"4" **\dir{-};
    \endxy}
    ~=~
    \vcenter{\xy
        (0,5)="1";
        (0,0)*[o]=<9pt>[Fo]{\scriptstyle s}="2";
        (0,-5)="3";
        "1";"2" **\dir{-};
        "2";"3" **\dir{-};
    \endxy}
    \qquad\qquad
    \vcenter{\xy
        (0,7)="1";
        (0,3)*[o]=<9pt>[Fo]{\scriptstyle s}="2";
        (0,-3)*[o]=<9pt>[Fo]{\scriptstyle r}="3";
        (0,-7)="4";
        "1";"2" **\dir{-};
        "2";"3" **\dir{-};
        "3";"4" **\dir{-};
    \endxy}
    ~=~
    \vcenter{\xy
        (0,5)="1";
        (0,0)*[o]=<9pt>[Fo]{\scriptstyle r}="2";
        (0,-5)="3";
        "1";"2" **\dir{-};
        "2";"3" **\dir{-};
    \endxy}
\end{equation}

\begin{equation}\tag{13}
    \vcenter{\xy
        ( 0, 4)="2";
        ( 0, 0)="3";
        (-3,-3)*[o]=<9pt>[Fo]{\scriptstyle t}="4";
        ( 3,-3)*[o]=<9pt>[Fo]{\scriptstyle r}="5";
        (-3,-7)="6";
        ( 3,-7)="7";
        "2";"3" **\dir{-};
        "3";"4" **\dir{-};
        "4";"6" **\dir{-};
        "3";"5" **\dir{-};
        "5";"7" **\dir{-};
    \endxy}
    ~=~
    \vcenter{\xy
        ( 0, 4)="2";
        ( 0, 1)="3";
        (-3.5,-2)*[o]=<9pt>[Fo]{\scriptstyle r}="4";
        ( 3.5,-2)*[o]=<9pt>[Fo]{\scriptstyle t}="5";
        (-3,-8)="6";
        ( 3,-8)="7";
        "2";"3" **\dir{-};
        "3";"4" **\dir{-};
        "3";"5" **\dir{-};
        "4";"7" **\dir{-}?(0.53)*{\hole}="x";
        "5";"x" **\dir{-};
        "x";"6" **\dir{-};
    \endxy}
    \qquad\qquad
    \vcenter{\xy
        ( 0,-4)="2";
        ( 0, 0)="3";
        (-3, 3)*[o]=<9pt>[Fo]{\scriptstyle t}="4";
        ( 3, 3)*[o]=<9pt>[Fo]{\scriptstyle r}="5";
        (-3, 7)="6";
        ( 3, 7)="7";
        "2";"3" **\dir{-};
        "3";"4" **\dir{-};
        "4";"6" **\dir{-};
        "3";"5" **\dir{-};
        "5";"7" **\dir{-};
    \endxy}
    ~=~
    \vcenter{\xy
        (-3, 6)="1";
        ( 3, 6)="2";
        (-3, 0)*[o]=<8pt>[Fo]{\scriptstyle r}="3";
        ( 3, 0)*[o]=<8pt>[Fo]{\scriptstyle t}="4";
        ( 0,-2)="5";
        ( 0,-6)="6";
        "1";"4" **\dir{-};
        "3";"5" **\dir{-};
        "4";"5" **\dir{-};
        "5";"6" **\dir{-};
        \ar@{-} "2";"3" |<(0.45)\hole
    \endxy}
\end{equation}

\begin{equation}\tag{14}
    \vcenter{\xy
        (0,-6)="t";
        (0,-3)*[o]=<9pt>[Fo]{\scriptstyle s}="f";
        (0,0)="m";
        ( 3,3)*[o]=<9pt>[Fo]{\scriptstyle r}="g";
        ( 3,7)="l";
        (-3,7)="r";
        "t";"f" **\dir{-};
        "f";"m" **\dir{-};
        "r";"m" **\crv{(-3,3)};
        "m";"g" **\dir{-};
        "g";"l" **\dir{-};
    \endxy}
    ~=~ 
    \vcenter{\xy
        (0,-5)="t";
        (0,-2)*[o]=<9pt>[Fo]{\scriptstyle s}="f";
        (0,1)="m";
        ( 3,5)="l";
        (-3,5)="r";
        "t";"f" **\dir{-};
        "f";"m" **\dir{-};
        "m";"l" **\dir{-};
        "m";"r" **\dir{-};
    \endxy}
    \qquad\qquad
    \vcenter{\xy
        (0,-6)="t";
        (0,-3)*[o]=<9pt>[Fo]{\scriptstyle r}="f";
        (0,0)="m";
        ( 3,3)*[o]=<9pt>[Fo]{\scriptstyle s}="g";
        ( 3,7)="l";
        (-3,7)="r";
        "t";"f" **\dir{-};
        "f";"m" **\dir{-};
        "r";"m" **\crv{(-3,3)};
        "m";"g" **\dir{-};
        "g";"l" **\dir{-};
    \endxy}
    ~=~ 
    \vcenter{\xy
        (0,-5)="t";
        (0,-2)*[o]=<9pt>[Fo]{\scriptstyle r}="f";
        (0,1)="m";
        ( 3,5)="l";
        (-3,5)="r";
        "t";"f" **\dir{-};
        "f";"m" **\dir{-};
        "m";"l" **\dir{-};
        "m";"r" **\dir{-};
    \endxy}
\end{equation}
The proofs of these properties may also be found in Appendix~\ref{sec-proofst}.

\begin{definition}\label{def-whm}
A weak bimonoid $H$ is called a \emph{weak Hopf monoid} if it is equipped
with an endomorphism $\nu:H \ra H$, called the \emph{antipode}, satisfying
\[
    \vcenter{\xy
        (0,7)="t";
        (0,4)="t1";
        (-4,0)*[o]=<9pt>[Fo]{\scriptstyle \nu}="m";
        (0,-4)="b1";
        (0,-7)="b";
        "t";"b" **\dir{-};
        "t1";"m" **\dir{-};
        "b1";"m" **\dir{-};
    \endxy}
    ~=~
    \vcenter{\xy
        (0,7)="t";
        (0,0)*[o]=<9pt>[Fo]{\scriptstyle t}="m";
        (0,-7)="b";
        "t";"m" **\dir{-};
        "m";"b" **\dir{-};
    \endxy}
    \qquad, \qquad
    \vcenter{\xy
        (0,7)="t";
        (0,4)="t1";
        (4,0)*[o]=<9pt>[Fo]{\scriptstyle \nu}="m";
        (0,-4)="b1";
        (0,-7)="b";
        "t";"b" **\dir{-};
        "t1";"m" **\dir{-};
        "b1";"m" **\dir{-};
    \endxy}
    ~=~
    \vcenter{\xy
        (0,7)="t";
        (0,0)*[o]=<9pt>[Fo]{\scriptstyle r}="m";
        (0,-7)="b";
        "t";"m" **\dir{-};
        "m";"b" **\dir{-};
    \endxy}
    \qquad, \qquad
    \vcenter{\xy
        (0,7)="t";
        (0,4)="t1";
        (4,0)*[o]=<9pt>[Fo]{\scriptstyle \nu}="m1";
        (-4,0)*[o]=<9pt>[Fo]{\scriptstyle \nu}="m2";
        (0,-4)="b1";
        (0,-7)="b";
        "t";"b" **\dir{-};
        "t1";"m1" **\dir{-};
        "b1";"m1" **\dir{-};
        "t1";"m2" **\dir{-};
        "b1";"m2" **\dir{-};
    \endxy}
    ~=~
    \vcenter{\xy
        (0,7)="t";
        (0,0)*[o]=<9pt>[Fo]{\scriptstyle \nu}="m";
        (0,-7)="b";
        "t";"m" **\dir{-};
        "m";"b" **\dir{-};
    \endxy}~.
\]
\end{definition}

The axioms of a weak Hopf monoid immediately imply the following identities
\[
    \vcenter{\xy
        (0,6)="1";
        (0,0)*[o]=<9pt>[Fo]{\scriptstyle \nu}="2";
        (0,-6)="3";
        "1";"2" **\dir{-};
        "2";"3" **\dir{-};
    \endxy}
    ~=~
    \vcenter{\xy
        (0,7)="t";
        (0,4)="t1";
        (-4,0)*[o]=<9pt>[Fo]{\scriptstyle t}="m2";
        (4,0)*[o]=<9pt>[Fo]{\scriptstyle \nu}="m1";
        (0,-4)="b1";
        (0,-7)="b";
        "t";"t1" **\dir{-};
        "b1";"b" **\dir{-};
        "t1";"m1" **\dir{-};
        "b1";"m1" **\dir{-};
        "t1";"m2" **\dir{-};
        "b1";"m2" **\dir{-};
    \endxy}
    ~=~
    \vcenter{\xy
        (0,7)="t";
        (0,4)="t1";
        (-4,0)*[o]=<9pt>[Fo]{\scriptstyle \nu}="m2";
        (4,0)*[o]=<9pt>[Fo]{\scriptstyle r}="m1";
        (0,-4)="b1";
        (0,-7)="b";
        "t";"t1" **\dir{-};
        "b1";"b" **\dir{-};
        "t1";"m1" **\dir{-};
        "b1";"m1" **\dir{-};
        "t1";"m2" **\dir{-};
        "b1";"m2" **\dir{-};
    \endxy}~.
\]
The antipode is unique since if $\nu'$ is another
\begin{align*}
    \nu' = \nu' * 1 * \nu' = t * \nu' = \nu * 1 * \nu' = \nu * r = \nu * 1 *
\nu = \nu.
\end{align*}

If $H$ and $K$ are weak Hopf monoids in $\V$, then a \emph{morphism of weak
Hopf monoids} $f:H \ra K$ is a morphism $f:H \ra K$ in $\V$ which is a
monoid and comonoid morphism that also preserves the antipode,
i.e., $f \nu = \nu f$.

We list some properties of the antipode $\nu:H \ra H$.

\begin{proposition}\label{lem-nu-s}
\begin{equation}\tag{15}
    \vcenter{\xy
        (0,7)="1";
        (0,3)*[o]=<9pt>[Fo]{\scriptstyle s}="2";
        (0,-3)*[o]=<9pt>[Fo]{\scriptstyle \nu}="3";
        (0,-7)="4";
        "1";"2" **\dir{-};
        "2";"3" **\dir{-};
        "3";"4" **\dir{-};
    \endxy}
    ~=~
    \vcenter{\xy
        (0,5)="1";
        (0,0)*[o]=<9pt>[Fo]{\scriptstyle r}="2";
        (0,-5)="3";
        "1";"2" **\dir{-};
        "2";"3" **\dir{-};
    \endxy}
\end{equation}

\begin{equation}\tag{16}
    \vcenter{\xy
        (0,7)="t";
        (0,3)*[o]=<9pt>[Fo]{\scriptstyle \nu}="m1";
        (0,-3)*[o]=<9pt>[Fo]{\scriptstyle t}="m2";
        (0,-7)="b";
        "t";"m1" **\dir{-};
        "m1";"m2" **\dir{-};
        "m2";"b" **\dir{-};
    \endxy}
    ~=~
    \vcenter{\xy
        (0,7)="t";
        (0,3)*[o]=<9pt>[Fo]{\scriptstyle r}="m1";
        (0,-3)*[o]=<9pt>[Fo]{\scriptstyle \nu}="m2";
        (0,-7)="b";
        "t";"m1" **\dir{-};
        "m1";"m2" **\dir{-};
        "m2";"b" **\dir{-};
    \endxy}
    ~=~
    \vcenter{\xy
        (0,7)="t";
        (0,3)*[o]=<9pt>[Fo]{\scriptstyle r}="m1";
        (0,-3)*[o]=<9pt>[Fo]{\scriptstyle t}="m2";
        (0,-7)="b";
        "t";"m1" **\dir{-};
        "m1";"m2" **\dir{-};
        "m2";"b" **\dir{-};
    \endxy}
\qquad\qquad
    \vcenter{\xy
        (0,7)="t";
        (0,3)*[o]=<9pt>[Fo]{\scriptstyle \nu}="m1";
        (0,-3)*[o]=<9pt>[Fo]{\scriptstyle r}="m2";
        (0,-7)="b";
        "t";"m1" **\dir{-};
        "m1";"m2" **\dir{-};
        "m2";"b" **\dir{-};
    \endxy}
    ~=~
    \vcenter{\xy
        (0,7)="t";
        (0,3)*[o]=<9pt>[Fo]{\scriptstyle t}="m1";
        (0,-3)*[o]=<9pt>[Fo]{\scriptstyle \nu}="m2";
        (0,-7)="b";
        "t";"m1" **\dir{-};
        "m1";"m2" **\dir{-};
        "m2";"b" **\dir{-};
    \endxy}
    ~=~
    \vcenter{\xy
        (0,7)="t";
        (0,3)*[o]=<9pt>[Fo]{\scriptstyle t}="m1";
        (0,-3)*[o]=<9pt>[Fo]{\scriptstyle r}="m2";
        (0,-7)="b";
        "t";"m1" **\dir{-};
        "m1";"m2" **\dir{-};
        "m2";"b" **\dir{-};
    \endxy}
\end{equation}

\begin{equation}\tag{17}
\begin{split}
    & \vcenter{\xy
        (0,4)="t";
        (0,0)*[o]=<9pt>[Fo]{\scriptstyle \nu}="m";
        (0,-4.5)*[o]=<5pt>[Fo]{~}="b";
        "t";"m" **\dir{-};
        "m";"b" **\dir{-};
    \endxy}
    ~=~
    \vcenter{\xy
        (0,3)="t";
        (0,-2)*[o]=<5pt>[Fo]{~}="b";
        "t";"b" **\dir{-};
    \endxy}
\qquad\qquad
    \vcenter{\xy
        (0,7)="t";
        (0,3)*[o]=<9pt>[Fo]{\scriptstyle \nu}="f";
        (0,-1)="m";
        (-3,-5)="l";
        ( 3,-5)="r";
        "t";"f" **\dir{-};
        "f";"m" **\dir{-};
        "m";"l" **\dir{-};
        "m";"r" **\dir{-};
    \endxy}
    ~=~
    \vcenter{\xy
        ( 0,-9)*{~}; ( 0,5)*{~};
        ( 0, 4)="2";
        ( 0, 1)="3";
        (-3.5,-2)*[o]=<9pt>[Fo]{\scriptstyle \nu}="4";
        ( 3.5,-2)*[o]=<9pt>[Fo]{\scriptstyle \nu}="5";
        (-3,-8)="6";
        ( 3,-8)="7";
        "2";"3" **\dir{-};
        "3";"4" **\dir{-};
        "3";"5" **\dir{-};
        "4";"7" **\dir{-}?(0.53)*{\hole}="x";
        "5";"x" **\dir{-};
        "x";"6" **\dir{-};
    \endxy}
\\
    & \vcenter{\xy
        (0,4.5)*[o]=<5pt>[Fo]{~}="1";
        (0,0)*[o]=<9pt>[Fo]{\scriptstyle \nu}="2";
        (0,-4)="3";
        "1";"2" **\dir{-};
        "2";"3" **\dir{-};
    \endxy}
    ~=~
    \xy 
        (0,2)*[o]=<5pt>[Fo]{}="t";
        (0,-3)="b";
        "t";"b" **\dir{-};
    \endxy
\qquad\qquad
    \vcenter{\xy
        (0,-7)="t";
        (0,-3)*[o]=<9pt>[Fo]{\scriptstyle \nu}="f";
        (0,1)="m";
        (-3,5)="l";
        ( 3,5)="r";
        "t";"f" **\dir{-};
        "f";"m" **\dir{-};
        "m";"l" **\dir{-};
        "m";"r" **\dir{-};
    \endxy}
    ~=~
    \vcenter{\xy
        ( 0,9)*{~}; ( 0,-5)*{~};
        ( 0,-4)="2";
        ( 0,-1)="3";
        ( 3.5,2)*[o]=<9pt>[Fo]{\scriptstyle \nu}="4";
        (-3.5,2)*[o]=<9pt>[Fo]{\scriptstyle \nu}="5";
        ( 3,8)="6";
        (-3,8)="7";
        "2";"3" **\dir{-};
        "3";"4" **\dir{-};
        "3";"5" **\dir{-};
        "4";"7" **\dir{-}?(0.53)*{\hole}="x";
        "5";"x" **\dir{-};
        "x";"6" **\dir{-};
    \endxy}
\end{split}
\end{equation}
\end{proposition}

The last identity (17) states that $\nu:A \ra A$ is both an anti-comonoid
morphism and an anti-monoid morphism.

\begin{proof}
The calculation
\[
    \vcenter{\xy
        (0,7)="1";
        (0,3)*[o]=<9pt>[Fo]{\scriptstyle s}="2";
        (0,-3)*[o]=<9pt>[Fo]{\scriptstyle \nu}="3";
        (0,-7)="4";
        "1";"2" **\dir{-};
        "2";"3" **\dir{-};
        "3";"4" **\dir{-};
    \endxy}
    ~\overset{\text{($\nu$)}}{=}~
    \vcenter{\xy
        ( 0, 8)="1";
        ( 0, 4)*[o]=<9pt>[Fo]{\scriptstyle s}="2";
        ( 0, 0)="3";
        (-4,-4)*[o]=<9pt>[Fo]{\scriptstyle t}="4";
        ( 4,-4)*[o]=<9pt>[Fo]{\scriptstyle \nu}="5";
        ( 0,-8)="6";
        ( 0,-11)="7";
        "1";"2" **\dir{-};
        "2";"3" **\dir{-};
        "4";"3" **\dir{-};
        "4";"6" **\dir{-};
        "5";"3" **\dir{-};
        "5";"6" **\dir{-};
        "7";"6" **\dir{-};
    \endxy}
    ~\overset{\text{(9)}}{=}~
    \vcenter{\xy
        ( 0, 8)="1";
        ( 0, 4)*[o]=<9pt>[Fo]{\scriptstyle s}="2";
        ( 0, 0)="3";
        ( 4,-4)*[o]=<9pt>[Fo]{\scriptstyle \nu}="5";
        ( 0,-8)="6";
        ( 0,-11)="7";
        "1";"2" **\dir{-};
        "2";"7" **\dir{-};
        "5";"3" **\dir{-};
        "5";"6" **\dir{-};
    \endxy}
    ~\overset{\text{($\nu$)}}{=}~
    \vcenter{\xy
        (0,7)="1";
        (0,3)*[o]=<9pt>[Fo]{\scriptstyle s}="2";
        (0,-3)*[o]=<9pt>[Fo]{\scriptstyle r}="3";
        (0,-7)="4";
        "1";"2" **\dir{-};
        "2";"3" **\dir{-};
        "3";"4" **\dir{-};
    \endxy}
    ~\overset{\text{(12)}}{=}~
    \vcenter{\xy
        (0,5)="1";
        (0,0)*[o]=<9pt>[Fo]{\scriptstyle r}="2";
        (0,-5)="3";
        "1";"2" **\dir{-};
        "2";"3" **\dir{-};
    \endxy}
\]
verifies the identity (15), and the following verifies the first identity of
(16):
\[
    \vcenter{\xy
        (0,7)="1";
        (0,3)*[o]=<9pt>[Fo]{\scriptstyle \nu}="2";
        (0,-3)*[o]=<9pt>[Fo]{\scriptstyle t}="3";
        (0,-7)="4";
        "1";"2" **\dir{-};
        "2";"3" **\dir{-};
        "3";"4" **\dir{-};
    \endxy}
    ~\overset{\text{($\nu$)}}{=}~
    \vcenter{\xy
        (0,9)="1";
        (0,6)="2";
        (-4,2)*[o]=<9pt>[Fo]{\scriptstyle t}="3";
        (4,2)*[o]=<9pt>[Fo]{\scriptstyle \nu}="4";
        (0,-2)="5";
        (0,-6)*[o]=<9pt>[Fo]{\scriptstyle t}="6";
        (0,-10)="7";
        "1";"2" **\dir{-};
        "2";"3" **\dir{-};
        "2";"4" **\dir{-};
        "5";"3" **\dir{-};
        "5";"4" **\dir{-};
        "5";"6" **\dir{-};
        "6";"7" **\dir{-};
    \endxy}
    ~\overset{\text{(3)}}{=}~
    \vcenter{\xy
        (0,9)="1";
        (0,6)="2";
        (4,2)*[o]=<9pt>[Fo]{\scriptstyle \nu}="4";
        (0,-2)="5";
        (0,-6)*[o]=<9pt>[Fo]{\scriptstyle t}="6";
        (0,-10)="7";
        "1";"6" **\dir{-};
        "2";"4" **\dir{-};
        "5";"4" **\dir{-};
        "6";"7" **\dir{-};
    \endxy}
    ~\overset{\text{($\nu$)}}{=}~
    \vcenter{\xy
        (0,7)="1";
        (0,3)*[o]=<9pt>[Fo]{\scriptstyle r}="2";
        (0,-3)*[o]=<9pt>[Fo]{\scriptstyle t}="3";
        (0,-7)="4";
        "1";"2" **\dir{-};
        "2";"3" **\dir{-};
        "3";"4" **\dir{-};
    \endxy}
    ~\overset{\text{($\nu$)}}{=}~
    \vcenter{\xy
        ( 0,10)="1";
        ( 0, 6)*[o]=<9pt>[Fo]{\scriptstyle r}="2";
        ( 0, 2)="3";
        (-4,-2)*[o]=<9pt>[Fo]{\scriptstyle \nu}="4";
        ( 0,-6)="6";
        ( 0,-9)="7";
        "1";"2" **\dir{-};
        "2";"7" **\dir{-};
        "3";"4" **\dir{-};
        "6";"4" **\dir{-};
    \endxy}
    ~\overset{\text{(3)}}{=}~
    \vcenter{\xy
        ( 0,10)="1";
        ( 0, 6)*[o]=<9pt>[Fo]{\scriptstyle r}="2";
        ( 0, 2)="3";
        (-4,-2)*[o]=<9pt>[Fo]{\scriptstyle \nu}="4";
        ( 4,-2)*[o]=<9pt>[Fo]{\scriptstyle r}="5";
        ( 0,-6)="6";
        ( 0,-9)="7";
        "1";"2" **\dir{-};
        "2";"3" **\dir{-};
        "3";"4" **\dir{-};
        "3";"5" **\dir{-};
        "6";"4" **\dir{-};
        "6";"5" **\dir{-};
        "6";"7" **\dir{-};
    \endxy}
    ~\overset{\text{($\nu$)}}{=}~
    \vcenter{\xy
        (0,7)="1";
        (0,3)*[o]=<9pt>[Fo]{\scriptstyle r}="2";
        (0,-3)*[o]=<9pt>[Fo]{\scriptstyle \nu}="3";
        (0,-7)="4";
        "1";"2" **\dir{-};
        "2";"3" **\dir{-};
        "3";"4" **\dir{-};
    \endxy}~.
\]
The second identity of (16) follows from a similar calculation.

To prove (17) we will only prove that $\nu$ is an anti-comonoid morphism.
That $\nu$ is an anti-monoid morphism follows by rotating all the diagrams
used to prove this statement by $\pi$.

The proof of the counit property is easy enough:
\[
    \vcenter{\xy
        (0,4)="1";
        (0,0)*[o]=<9pt>[Fo]{\scriptstyle \nu}="2";
        (0,-4.5)*[o]=<5pt>[Fo]{~}="3";
        "1";"2" **\dir{-};
        "2";"3" **\dir{-};
    \endxy}
    ~\overset{\text{($\nu$)}}{=}~
    \vcenter{\xy
        (0,-6)*[o]=<5pt>[Fo]{~}="1";
        (0,-3)="2";
        (-3,0)*[o]=<9pt>[Fo]{\scriptstyle t}="3";
        (3,0)*[o]=<9pt>[Fo]{\scriptstyle \nu}="4";
        (0,3)="5";
        (0,6)="6";
        "1";"2" **\dir{-};
        "2";"3" **\dir{-};
        "2";"4" **\dir{-};
        "3";"5" **\dir{-};
        "4";"5" **\dir{-};
        "5";"6" **\dir{-};
    \endxy}
    ~\overset{\text{(2)}}{=}~
    \vcenter{\xy
        (0,-6)*[o]=<5pt>[Fo]{~}="1";
        (0,-3)="2";
        (3,0)*[o]=<9pt>[Fo]{\scriptstyle \nu}="4";
        (0,3)="5";
        (0,6)="6";
        "1";"6" **\dir{-};
        "2";"4" **\dir{-};
        "4";"5" **\dir{-};
    \endxy}
    ~\overset{\text{($\nu$)}}{=}~
    \vcenter{\xy
        (0,4)="1";
        (0,0)*[o]=<9pt>[Fo]{\scriptstyle r}="2";
        (0,-4.5)*[o]=<5pt>[Fo]{}="3";
        "1";"2" **\dir{-};
        "2";"3" **\dir{-};
    \endxy}
    ~\overset{\text{(2)}}{=}~
    \vcenter{\xy
        (0,3)="1";
        (0,-2)*[o]=<5pt>[Fo]{}="2"; 
        **\dir{-};
    \endxy}~.
\]
The following calculation proves that the antipode is anti-comultiplicative.
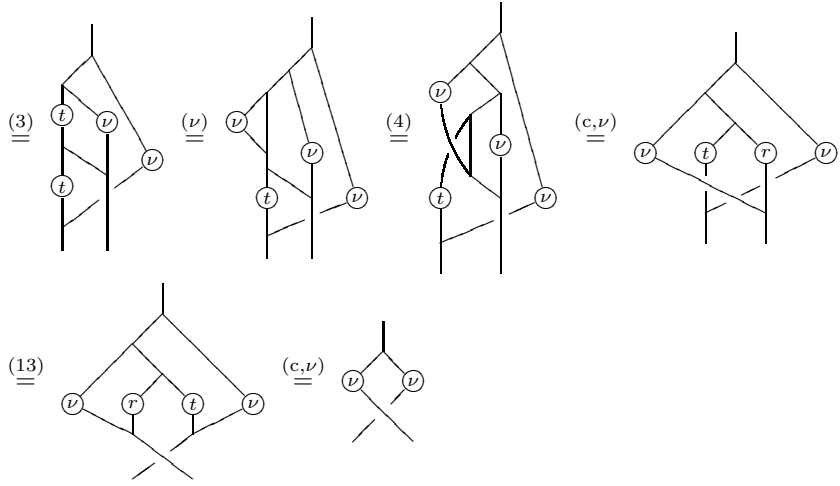
\begin{align*}
    \vcenter{\xy
        ( 0, 8)="1";
        ( 0, 4)*[o]=<9pt>[Fo]{\scriptstyle \nu}="2";
        ( 0, 0)="3";
        (-3,-3)="4";
        ( 3,-3)="5";
        "1";"2" **\dir{-};
        "2";"3" **\dir{-};
        "3";"4" **\dir{-};
        "3";"5" **\dir{-};
    \endxy}
    &~\overset{\text{($\nu$)}}{=}~
    \vcenter{\xy
        ( 0,12)="1";
        ( 0, 8)="2";
        (-4, 4)*[o]=<9pt>[Fo]{\scriptstyle \nu}="3";
        ( 4, 4)*[o]=<9pt>[Fo]{\scriptstyle r}="4";
        ( 0, 0)="5";
        ( 0,-4)="6";
        (-3,-7)="7";
        ( 3,-7)="8";
        "1";"2" **\dir{-};
        "2";"3" **\dir{-};
        "2";"4" **\dir{-};
        "5";"3" **\dir{-};
        "5";"4" **\dir{-};
        "5";"6" **\dir{-};
        "6";"7" **\dir{-};
        "6";"8" **\dir{-};
    \endxy}
    ~\overset{\text{(b)}}{=}~
    \vcenter{\xy
        ( 0,12)="1";
        ( 0, 8)="2";
        (-4, 4)*[o]=<9pt>[Fo]{\scriptstyle \nu}="3";
        ( 4, 4)*[o]=<9pt>[Fo]{\scriptstyle r}="4";
        (-4, 1)="5";
        ( 4, 1)="6";
        (-4,-3)="7";
        ( 4,-3)="8";
        (-4,-6)="9";
        ( 4,-6)="10";
        "1";"2" **\dir{-};
        "2";"3" **\dir{-};
        "2";"4" **\dir{-};
        "3";"9" **\dir{-};
        "4";"10" **\dir{-};
        "5";"8" **\dir{-};
        \ar@{-} "6";"7" |\hole
    \endxy}
    ~\overset{\text{(3)}}{=}~
    \vcenter{\xy
        ( 0,12)="1";
        ( 0, 8)="2";
        (-4, 4)*[o]=<9pt>[Fo]{\scriptstyle \nu}="3";
        ( 4, 4)*[o]=<9pt>[Fo]{\scriptstyle r}="4";
        (-4, 1)="5";
        ( 4, 1)="6";
        (-4,-7)="7";
        ( 4,-7)="8";
        (-4,-10)="9";
        ( 4,-10)="10";
        ( 4,-3)*[o]=<9pt>[Fo]{\scriptstyle r}="x";
        "1";"2" **\dir{-};
        "2";"3" **\dir{-};
        "2";"4" **\dir{-};
        "3";"9" **\dir{-};
        "4";"x" **\dir{-};
        "x";"10" **\dir{-};
        "5";"8" **\dir{-};
        \ar@{-} "6";"7" |\hole
    \endxy}
    ~\overset{\text{($\nu$)}}{=}~
    \vcenter{\xy
        ( 0,12)="1";
        ( 0, 8)="2";
        (-4, 4)*[o]=<9pt>[Fo]{\scriptstyle \nu}="3";
        ( 4, 4)="4";
        ( 8, 0)*[o]=<9pt>[Fo]{\scriptstyle \nu}="5";
        ( 4,-4)="6";
        (-4,-6)="7";
        ( 4,-6)="8";
        ( 4,-10)*[o]=<9pt>[Fo]{\scriptstyle r}="9";
        (-4,-14)="10";
        ( 4,-14)="11";
        (-4,-17)="12";
        ( 4,-17)="13";
        "1";"2" **\dir{-};
        "2";"3" **\dir{-};
        "2";"4" **\dir{-};
        "5";"6" **\dir{-};
        "4";"5" **\dir{-};
        "4";"9" **\dir{-};
        "9";"13" **\dir{-};
        "3";"12" **\dir{-};
        "7";"11" **\dir{-};
        \ar@{-} "8";"10" |\hole
    \endxy}
    ~\overset{\text{(4)}}{=}~
    \vcenter{\xy
        ( 0,12)="1";
        ( 0, 8)="2";
        (-4, 4)*[o]=<9pt>[Fo]{\scriptstyle \nu}="3";
        ( 4, 4)="4";
        ( 0, 0)="5";
        ( 6, 0)*[o]=<9pt>[Fo]{\scriptstyle \nu}="6";
        ( 0,-8)="8";
        (-4,-10)="7";
        ( 6,-10)*[o]=<9pt>[Fo]{\scriptstyle r}="9";
        (-4,-14)="10";
        ( 6,-14)="11";
        (-4,-17)="12";
        ( 6,-17)="13";
        "1";"2" **\dir{-};
        "2";"3" **\dir{-};
        "2";"4" **\dir{-};
        "4";"5" **\dir{-};
        "4";"6" **\dir{-};
        "5";"8" **\dir{-};
        "9";"13" **\dir{-};
        "3";"12" **\dir{-};
        "7";"11" **\dir{-};
        \ar@{-} "8";"10" |\hole
        \ar@{-}@/^3pt/ "5";"9"
        \ar@{-}@/^3pt/ "6";"8" |\hole
    \endxy}
\\
    &~\overset{\text{($\nu$)}}{=}~
    \vcenter{\xy
        ( 0,12)="1";
        ( 0, 8)="2";
        (-4, 4)*[o]=<9pt>[Fo]{\scriptstyle \nu}="3";
        ( 4, 4)="4";
        ( 0, 0)="5";
        ( 6, 0)*[o]=<9pt>[Fo]{\scriptstyle \nu}="6";
        ( 0,-8)="8";
        (-4,-10)="7";
        ( 6,-8)="x";
        (10,-11)*[o]=<9pt>[Fo]{\scriptstyle \nu}="z";
        ( 6,-14)="y";
        (-4,-16)="10";
        ( 6,-16)="11";
        (-4,-19)="12";
        ( 6,-19)="13";
        "1";"2" **\dir{-};
        "2";"3" **\dir{-};
        "2";"4" **\dir{-};
        "4";"5" **\dir{-};
        "4";"6" **\dir{-};
        "5";"8" **\dir{-};
        "x";"13" **\dir{-};
        "x";"z" **\dir{-};
        "y";"z" **\dir{-};
        "3";"12" **\dir{-};
        "7";"11" **\dir{-};
        \ar@{-} "8";"10" |<(0.43)\hole
        \ar@{-}@/^3pt/ "5";"x"
        \ar@{-}@/^3pt/ "6";"8" |<(0.35)\hole
    \endxy}
    ~\overset{\text{(c)}}{=}~
    \vcenter{\xy
        ( 0,16)="1";
        ( 0,12)="2";
        (-4, 8)="3";
        (-8, 4)*[o]=<9pt>[Fo]{\scriptstyle \nu}="4";
        ( 0, 4)="5";
        (-8, 0)="6";
        ( 0,0)="7";
        ( 4,-2)*[o]=<9pt>[Fo]{\scriptstyle \nu}="8";
        ( 9,-5)*[o]=<9pt>[Fo]{\scriptstyle \nu}="9";
        (-8,-4)="10";
        ( 0,-4)="11";
        ( 0,-8)="12";
        (-8,-17)="13";
        (-8,-20)="14";
        ( 0,-20)="15";
        "1";"2" **\dir{-};
        "2";"3" **\dir{-};
        "3";"4" **\dir{-};
        "3";"5" **\dir{-};
        "5";"8" **\dir{-};
        "8";"12" **\dir{-};
        "6";"11" **\dir{-};
        "2";"9" **\dir{-};
        "4";"14" **\dir{-};
        "5";"15" **\dir{-};
        \ar@{-} "7";"10" |\hole
        \ar@{-} "9";"13" |<(0.5)\hole
    \endxy}
    ~\overset{\text{(b)}}{=}~
    \vcenter{\xy
        ( 0,16)="1";
        ( 0,12)="2";
        (-4, 8)="3";
        (-8, 4)*[o]=<9pt>[Fo]{\scriptstyle \nu}="4";
        ( 0, 4)="5";
        (-4, 0)="6";
        ( 4,-4)*[o]=<9pt>[Fo]{\scriptstyle \nu}="7";
        (-4,-4)="8";
        (-8,-8)="9";
        ( 0,-8)="10";
        (9,-6)*[o]=<9pt>[Fo]{\scriptstyle \nu}="13";
        (-8,-20)="11";
        ( 0,-20)="12";
        (-8,-17)="14";
        "1";"2" **\dir{-};
        "2";"3" **\dir{-};
        "3";"4" **\dir{-};
        "3";"5" **\dir{-};
        "4";"6" **\dir{-};
        "5";"6" **\dir{-};
        "5";"7" **\dir{-};
        "6";"8" **\dir{-};
        "7";"10" **\dir{-};
        "8";"9" **\dir{-};
        "8";"10" **\dir{-};
        "9";"11" **\dir{-};
        "10";"12" **\dir{-};
        "2";"13" **\dir{-};
        \ar@{-} "13";"14" |<(0.5)\hole
    \endxy}
    ~\overset{\text{(c,$\nu$)}}{=}~
    \vcenter{\xy
        ( 0,16)="1";
        ( 0,12)="2";
        (-4, 8)="3";
        (-4, 4)*[o]=<9pt>[Fo]{\scriptstyle t}="4";
        ( 2, 3)*[o]=<9pt>[Fo]{\scriptstyle \nu}="5";
        (-4,0)="8";
        ( 2,-4)="10";
        ( 8,-2)*[o]=<9pt>[Fo]{\scriptstyle \nu}="13";
        (-4,-14)="11";
        ( 2,-14)="12";
        (-4,-11)="14";
        "1";"2" **\dir{-};
        "2";"3" **\dir{-};
        "3";"4" **\dir{-};
        "3";"5" **\dir{-};
        "4";"11" **\dir{-};
        "5";"10" **\dir{-};
        "8";"10" **\dir{-};
        "10";"12" **\dir{-};
        "2";"13" **\dir{-};
        \ar@{-} "13";"14" |<(0.44)\hole
    \endxy}
\\
    &~\overset{\text{(3)}}{=}~
    \vcenter{\xy
        ( 0,16)="1";
        ( 0,12)="2";
        (-4, 8)="3";
        (-4, 4)*[o]=<9pt>[Fo]{\scriptstyle t}="4";
        (-4,-5.5)*[o]=<9pt>[Fo]{\scriptstyle t}="x";
        ( 2, 3)*[o]=<9pt>[Fo]{\scriptstyle \nu}="5";
        (-4,0)="8";
        ( 2,-4)="10";
        ( 8,-2)*[o]=<9pt>[Fo]{\scriptstyle \nu}="13";
        (-4,-14)="11";
        ( 2,-14)="12";
        (-4,-11)="14";
        "1";"2" **\dir{-};
        "2";"3" **\dir{-};
        "3";"4" **\dir{-};
        "3";"5" **\dir{-};
        "4";"x" **\dir{-};
        "x";"11" **\dir{-};
        "5";"10" **\dir{-};
        "8";"10" **\dir{-};
        "10";"12" **\dir{-};
        "2";"13" **\dir{-};
        \ar@{-} "13";"14" |<(0.44)\hole
    \endxy}
    ~\overset{\text{($\nu$)}}{=}~
    \vcenter{\xy
        ( 0,16)="1";
        ( 0,12)="2";
        (-3,9)="3";
        (-6, 6)="4";
        (-10,2)*[o]=<9pt>[Fo]{\scriptstyle \nu}="5";
        (-6,-2)="6";
        ( 0,-2)*[o]=<9pt>[Fo]{\scriptstyle \nu}="7";
        (-6,-4)="8";
        (-6,-8)*[o]=<9pt>[Fo]{\scriptstyle t}="9";
        (0,-8)="10";
        (-6,-16)="11";
        (0,-16)="12";
        ( 6,-8)*[o]=<9pt>[Fo]{\scriptstyle \nu}="14";
        (-6,-13)="15";
        "1";"2" **\dir{-};
        "2";"3" **\dir{-};
        "3";"4" **\dir{-};
        "4";"5" **\dir{-};
        "5";"6" **\dir{-};
        "4";"8" **\dir{-};
        "8";"9" **\dir{-};
        "8";"10" **\dir{-};
        "7";"10" **\dir{-};
        "3";"7" **\dir{-};
        "9";"11" **\dir{-};
        "10";"12" **\dir{-};
        "2";"14" **\dir{-};
        \ar@{-} "14";"15" |<(0.42)\hole
    \endxy}
    ~\overset{\text{(4)}}{=}~
    \vcenter{\xy
        ( 0,16)="1";
        ( 0,12)="2";
        (-4, 8)="3";
        (-8, 4)*[o]=<9pt>[Fo]{\scriptstyle \nu}="4";
        ( 0, 4)="5";
        (-4, 1)="6";
        ( 0,-3)*[o]=<9pt>[Fo]{\scriptstyle \nu}="7";
        ( 6,-10)*[o]=<9pt>[Fo]{\scriptstyle \nu}="8";
        (-4,-7)="9";
        (-8,-10)*[o]=<9pt>[Fo]{\scriptstyle t}="10";
        ( 0,-10)="11";
        (-8,-16)="12";
        (-8,-20)="13";
        ( 0,-20)="14";
        "1";"2" **\dir{-};
        "2";"3" **\dir{-};
        "3";"4" **\dir{-};
        "3";"5" **\dir{-};
        "5";"6" **\dir{-};
        "5";"7" **\dir{-};
        "7";"14" **\dir{-};
        "6";"9" **\dir{-};
        "9";"11" **\dir{-};
        "2";"8" **\dir{-};
        "10";"13" **\dir{-};
        \ar@{-}@/_3pt/ "4";"9"
        \ar@{-}@/_3pt/ "6";"10" |<(0.4)\hole
        \ar@{-} "8";"12" |<(0.35)\hole
    \endxy}
    ~\overset{\text{(c,$\nu$)}}{=}~
    \vcenter{\xy
        ( 0,16)="1";
        ( 0,12)="2";
        (-4, 8)="3";
        ( 0, 4)="4";
        (-12, 0)*[o]=<9pt>[Fo]{\scriptstyle \nu}="5";
        (-4, 0)*[o]=<9pt>[Fo]{\scriptstyle t}="6";
        ( 4, 0)*[o]=<9pt>[Fo]{\scriptstyle r}="7";
        (12, 0)*[o]=<9pt>[Fo]{\scriptstyle \nu}="8";
        (-4,-8)="9";
        ( 4,-8)="10";
        (-4,-12)="11";
        ( 4,-12)="12";
        "1";"2" **\dir{-};
        "2";"3" **\dir{-};
        "2";"8" **\dir{-};
        "3";"4" **\dir{-};
        "3";"5" **\dir{-};
        "4";"6" **\dir{-};
        "4";"7" **\dir{-};
        "12";"7" **\dir{-};
        "5";"10" **\dir{-};
        \ar@{-} "6";"11" |<(0.24)\hole
        \ar@{-} "8";"9" |<(0.43)\hole |<(0.74)\hole
    \endxy}
\\
    &~\overset{\text{(13)}}{=}~
    \vcenter{\xy
        ( 0,16)="1";
        ( 0,12)="2";
        (-4, 8)="3";
        ( 0, 4)="4";
        (-12, 0)*[o]=<9pt>[Fo]{\scriptstyle \nu}="5";
        (-4, 0)*[o]=<9pt>[Fo]{\scriptstyle r}="6";
        ( 4, 0)*[o]=<9pt>[Fo]{\scriptstyle t}="7";
        (12, 0)*[o]=<9pt>[Fo]{\scriptstyle \nu}="8";
        (-4,-4)="9";
        ( 4,-4)="10";
        (-4,-10)="11";
        ( 4,-10)="12";
        "1";"2" **\dir{-};
        "2";"3" **\dir{-};
        "2";"8" **\dir{-};
        "3";"4" **\dir{-};
        "3";"5" **\dir{-};
        "4";"6" **\dir{-};
        "4";"7" **\dir{-};
        "5";"9" **\dir{-};
        "6";"9" **\dir{-};
        "7";"10" **\dir{-};
        "8";"10" **\dir{-};
        "9";"12" **\dir{-};
        \ar@{-} "10";"11" |\hole
    \endxy}
    ~\overset{\text{(c,$\nu$)}}{=}~
    \vcenter{\xy
        ( 0,8)="1";
        ( 0,4)="2";
        (-4, 0)*[o]=<9pt>[Fo]{\scriptstyle \nu}="3";
        ( 4, 0)*[o]=<9pt>[Fo]{\scriptstyle \nu}="4";
        (-4,-8)="5";
        ( 4,-8)="6";
        "1";"2" **\dir{-};
        "2";"3" **\dir{-};
        "2";"4" **\dir{-};
        "3";"6" **\dir{-};
        \ar@{-} "4";"5" |\hole
    \endxy}
\end{align*}
\end{proof}

%=========================================================================%
\section{The monoidal category of $A$-comodules}\label{sec-monoidal}
%=========================================================================%

Suppose $A = (A,1)$ is a weak bimonoid in $\Q\V$ and let $C = (A,t)$. In
this section we describe a monoidal structure on the categories
$\Bicomod(C)$ of $C$-bicomodules in $\Q\V$, and $\Comod(A)$ of right
$A$-comodules in $\Q\V$ such that the underlying functor
\[
    U:\Comod(A) \dra \Bicomod(C)
\]
is strong monoidal. If $A$ is a weak Hopf monoid then we show that
$\Comod_f(A)$, the subcategory consisting of the dualizable objects, is left
autonomous.

This section is fairly standard in the $\V = \mathbf{Vect}$ case
(see~\cite{BSc}, \cite{N}, or~\cite{NV} for example) and carries
over rather straightforwardly to the general braided $\V$ case
(cf.~\cite{DMS}).

%=========================================================================%
\subsection{The monoidal structure on $C$-bicomodules}\label{sec-monC}
%=========================================================================%

Suppose, for this section, that $C \in \V$ is just a comonoid, and that
$M \in \V$ is a $C$-bicomodule with coaction
\[
    \gamma:M \dra C \ox M \ox C.
\]
A left $C$-coaction and a right $C$-coaction are obtained from $\gamma$ by
involving the counit $\epsilon$:
\begin{align*}
    \gamma_l &=
    \big(\xygraph{{M}
        :[r(1.7)] {C \ox M \ox C} ^-\gamma
        :[r(2.5)] {C \ox M} ^-{1 \ox 1 \ox \epsilon}}\big)
    \\
    \gamma_r &=
    \big(\xygraph{{M}
        :[r(1.7)] {C \ox M \ox C} ^-\gamma
        :[r(2.5)] {M \ox C} ^-{\epsilon \ox 1 \ox 1}}\big).
\end{align*}

Suppose now that $N$ is another $C$-bicomodule. The tensor product of
$M$ and $N$ over $C$ is defined to be the equalizer
\[
    \xygraph{{M \oxC N}="1"
        [r(1.8)] {M \ox N}="2"
        [r(2.4)] {M \ox C \ox N}="3"
        "1":"2" ^-\iota
        "2":@<3pt>"3" ^-{\gamma_r \ox 1}
        "2":@<-3pt>"3" _-{1 \ox \gamma_l}}.
\]
Obviously the morphism
\[
    \xygraph{{M \oxC N}
        :[r(1.8)] {M \ox N} ^-\iota
        :[r(2.6)] {C \ox M \ox N \ox C} ^-{\gamma_l \ox \gamma_r}}
\]
equalizes the two morphisms $C \ox M \ox N \ox C \dra C \ox M \ox C \ox N
\ox C$ and so induces a morphism
\[
    \gamma : M \oxC N \dra C \ox M \oxC N \ox C,
\]
which is the coaction on $M \oxC N$.

That this defines a monoidal structure on the category $\Bicomod(C)$ with
tensor product $\ox_C$ and unit $C$ is standard.

%=========================================================================%
\subsection{The tensor product of $A$-comodules}
%=========================================================================%

Let $A=(A,1)$ be a weak bimonoid in $\Q\V$ and let $C=(A,t)$. The monoidal
structure on the category of right $A$-comodules will be $\ox_C$, the tensor
product over $C$, with unit $C$.

Suppose that $M$ is a right $A$-comodule. We know that $s:A \ra C^\o$ and
$t:A \ra C$ are comonoid morphisms and that property (10) holds, where
recall that property (10) expresses the commutativity of the following
diagram.
\[
    \xygraph{{A}="s"
        :[u(0.7)r(1.2)] {A \ox A} ^-{\delta}
        :[r(1.8)] {C \ox C} ^-{s \ox t}
        :[d(1.4)] {C \ox C}="t" ^-c
     "s":[d(0.7)r(1.2)] {A \ox A} _-{\delta}
        :"t" ^-{t \ox s}}
\]
Therefore, $M$ may be made into a $C$-bicomodule via
\[
    \gamma =
    \big(\xygraph{{M}
        :[r(1.3)] {M \ox A} ^-\gamma
        :[r(1.9)] {M \ox A \ox A} ^-{1 \ox \delta}
        :[r(2.4)] {A \ox M \ox A} ^-{c^{-1} \ox 1}
        :[r(2.4)] {C \ox M \ox C} ^-{s \ox 1 \ox t}}\big),
\]
which is
\[
    \gamma =~
    \vcenter{\xy
        (-2,8)*{\scriptstyle M};
        (4,3)*{\scriptstyle A};
        (-5,-9)*{\scriptstyle C};
        (0,-9)*{\scriptstyle M};
        (5,-9)*{\scriptstyle C};
        (0,7)="1";
        (0,4)="2";
        (3,1)="3";
        (-5,-3)*[o]=<9pt>[Fo]{\scriptstyle s}="4";
        (5,-3)*[o]=<9pt>[Fo]{\scriptstyle t}="5";
        (-5,-7)="6";
        ( 0,-7)="7";
        ( 5,-7)="8";
        "2";"3" **\dir{-};
        "3";"4" **\dir{-};
        "3";"5" **\dir{-};
        "4";"6" **\dir{-};
        "5";"8" **\dir{-};
        \ar@{-} |<(0.45){\hole} "1";"7" ;
    \endxy}
\]
in strings. The left and right $C$-coactions are
\[
    \gamma_l =~
    \vcenter{\xy
        (0,7)="1";
        (0,4)="2";
        (3,2)="3";
        (-5,-1)*[o]=<9pt>[Fo]{\scriptstyle s}="4";
        (-5,-6)="6";
        ( 0,-6)="7";
        "2";"3" **\dir{-};
        "3";"4" **\dir{-};
        "4";"6" **\dir{-};
        \ar@{-} |<(0.5){\hole} "1";"7";
    \endxy}
    \qquad \text{and} \qquad
    \gamma_r =~
    \vcenter{\xy
        (0,7)="1";
        (0,3)="2";
        (4,0)*[o]=<9pt>[Fo]{\scriptstyle t}="3";
        (0,-5)="7";
        (4,-5)="5";
        "1";"7" **\dir{-};
        "2";"3" **\dir{-};
        "3";"5" **\dir{-};
    \endxy}~.
\]

The tensor product of two $A$-comodules $M$ and $N$ over $C$ then may be
defined as in \S\ref{sec-monC}. We derive an explicit description of
$M \oxC N$. Before doing so we will need the following definition.

\begin{definition}
Let $f,g:X \ra Y$ be a parallel pair in $\V$. This pair is called
\emph{cosplit} when there is an arrow $d:Y \ra X$ such that 
\[
    df = 1_X \qquad \text{and} \qquad fdg = gdg.
\]
It is not hard to see that, in this case, $dg:X \ra X$ is an idempotent and
a splitting of $dg$, i.e.,
\[
    \xygraph{{X}="1"
        [r(2)]{X}="2"
        [d(1)l(1)]{Q}="3"
        "1":"2" ^-{dg}
        "1":"3" _-x 
        "3":"2" _-y}
    \qquad\qquad
    \xygraph{{Q}="1"
        [r(2)]{Q}="2"
        [d(1)l(1)]{X}="3"
        "1":"2" ^-1
        "1":"3" _-y 
        "3":"2" _-x}
\]
provides an absolute equalizer $(Q,y)$ for $f$ and $g$.
\end{definition}

Now suppose $M$ and $N$ are $A$-comodules. Two morphisms $M \ox N \ra
M \ox C \ox N$ are given as  
\[
    \gamma_r \ox 1 =~
    \vcenter{\xy
        (-4,8)*{\scriptstyle M};
        ( 4,8)*{\scriptstyle N};
        ( -2,4)*{\scriptstyle A};
        (-4,-7)*{\scriptstyle M};
        ( 0,-7)*{\scriptstyle C};
        ( 4,-7)*{\scriptstyle N};
        (-4,6)="1";
        ( 4,6)="2";
        (-4,3)="3";
        ( 0,0)*[o]=<9pt>[Fo]{\scriptstyle t}="4";
        (-4,-5)="5";
        ( 0,-5)="6";
        ( 4,-5)="7";
        "1";"5" **\dir{-};
        "2";"7" **\dir{-};
        "3";"4" **\dir{-};
        "4";"6" **\dir{-};
    \endxy}
    \qquad \text{and} \qquad
    1 \ox \gamma_l =~
    \vcenter{\xy
        (-8,11)*{\scriptstyle M};
        (0,11)*{\scriptstyle N};
        (4,7)*{\scriptstyle A};
        (-8,-6)*{\scriptstyle M};
        (-4,-6)*{\scriptstyle C};
        (0,-6)*{\scriptstyle N};
        (-8, 9)="a";
        (-8,-4)="b";
        (0, 9)="1";
        (0, 6)="2";
        (3, 4)="x";
        (-4,0)*[o]=<9pt>[Fo]{\scriptstyle s}="3";
        (-4,-4)="4";
        (0,-4)="5";
        "a";"b" **\dir{-};
        "2";"x" **\dir{-};
        "x";"3" **\dir{-};
        "3";"4" **\dir{-};
        \ar@{-} |<(0.49){\hole} "1";"5";
    \endxy}~.
\]

\begin{proposition}\label{comod-cosplit}
The pair $\gamma_r \ox 1$ and $1 \ox \gamma_l$ are cosplit by
\[
    d =~
    \vcenter{\xy
        (-4,8)*{\scriptstyle M};
        ( 0,8)*{\scriptstyle C};
        ( 4,8)*{\scriptstyle N};
        (-4,6)="1";
        ( 0,6)="a";
        ( 4,6)="2";
        ( 0,2)*[o]=<9pt>[Fo]{\scriptstyle t}="3";
        (-4,1)="4";
        (0,-2)="5";
        ( 0,-5)*[o]=<5pt>[Fo]{~}="6";
        (-4,-8)="7";
        ( 4,-8)="8";
        (-4,-10)*{\scriptstyle M};
            ( 4,-10)*{\scriptstyle N};
        "1";"7" **\dir{-};
        "2";"8" **\dir{-};
        "a";"3" **\dir{-};
        "3";"6" **\dir{-};
        "4";"5" **\dir{-};
    \endxy}~.
\]
\end{proposition}

\begin{proof}
That $d$ is a morphism in $\Q\V$ follows immediately as $t$ is idempotent.
The calculation
\[
    d(\gamma_r \ox 1) =
    \vcenter{\xy
        (-4,13)="1";
        ( 4,13)="2";
        (-4,10)="a";
        ( 0,7)*[o]=<9pt>[Fo]{\scriptstyle t}="b";
        ( 0,2)*[o]=<9pt>[Fo]{\scriptstyle t}="3";
        (-4,1)="4";
        (0,-2)="5";
        ( 0,-5)*[o]=<5pt>[Fo]{~}="6";
        (-4,-8)="7";
        ( 4,-8)="8";
        "1";"7" **\dir{-};
        "2";"8" **\dir{-};
        "a";"b" **\dir{-};
        "b";"3" **\dir{-};
        "3";"6" **\dir{-};
        "4";"5" **\dir{-};
    \endxy}
    ~\overset{\text{(7)}}{=}~
    \vcenter{\xy
        (-4,8)="1";
        ( 4,8)="2";
        ( 0,2)*[o]=<9pt>[Fo]{\scriptstyle t}="3";
        (-4,5)="a";
        (-4,1)="4";
        (0,-2)="5";
        ( 0,-5)*[o]=<5pt>[Fo]{~}="6";
        (-4,-8)="7";
        ( 4,-8)="8";
        "1";"7" **\dir{-};
        "2";"8" **\dir{-};
        "a";"3" **\dir{-};
        "3";"6" **\dir{-};
        "4";"5" **\dir{-};
    \endxy}
    ~\overset{\text{(c)}}{=}~
    \vcenter{\xy
        (-4,8)="1";
        ( 5,8)="2";
        (-4,5)="3";
        ( 0,3)="4";
        ( 2,0)*[o]=<9pt>[Fo]{\scriptstyle t}="5";
        (-2,0)="6";
        ( 0,-3)="7";
        ( 0,-6)*[o]=<5pt>[Fo]{~}="8";
        (-4,-9)="9";
        ( 5,-9)="10";
        "1";"9" **\dir{-};
        "2";"10" **\dir{-};
        "3";"4" **\dir{-};
        "4";"5" **\dir{-};
        "4";"6" **\dir{-};
        "5";"7" **\dir{-};
        "6";"7" **\dir{-};
        "8";"7" **\dir{-};
    \endxy}
    ~\overset{\text{(6)}}{=}~
    \vcenter{\xy
        (-3,5)="1";
        ( 3,5)="2";
        (-3,3)="3";
        (0,0)*[o]=<5pt>[Fo]{~}="4";
        (-3,-5)="5";
        ( 3,-5)="6";
        "1";"5" **\dir{-};
        "2";"6" **\dir{-};
        "3";"4" **\dir{-};
    \endxy}
    ~\overset{\text{(c)}}{=}~
    \vcenter{\xy
        (-2,4)="1";
        ( 2,4)="2";
        (-2,-4)="3";
        ( 2,-4)="4";
        "1";"3" **\dir{-};
        "2";"4" **\dir{-};
    \endxy}
    ~= 1_{M \ox N}
    \]
shows that $d(\gamma_r \ox 1)=1$ and the identity
$(\gamma_r \ox 1)d(1 \ox \gamma_l) = (1 \ox \gamma_l)d(1 \ox \gamma_l)$
follows from:

\begin{multline*}
(\gamma_r \ox 1)d(1 \ox \gamma_l) =~
    \vcenter{\xy
        (0,13)*{~}; (0,-17)*{~};
        (-4,12)="1";
        ( 4,12)="2";
        ( 4,9)="3";
        ( 7,7)="4";
        ( 0,4)*[o]=<9pt>[Fo]{\scriptstyle s}="5";
        ( 0,-1)*[o]=<9pt>[Fo]{\scriptstyle t}="6";
        (-4,-2)="7";
        ( 0,-5)="8";
        ( 0,-8)*[o]=<5pt>[Fo]{~}="9";
        (-4,-9)="10";
        ( 0,-12)*[o]=<9pt>[Fo]{\scriptstyle t}="11";
        (-4,-16)="12";
        ( 0,-16)="13";
        ( 4,-16)="14";
        "2";"3" **\dir{-};
        "3";"4" **\dir{-};
        "4";"5" **\dir{-};
        "5";"6" **\dir{-};
        "7";"8" **\dir{-};
        "6";"9" **\dir{-};
        "10";"11" **\dir{-};
        "11";"13" **\dir{-};
        "1";"12" **\dir{-};
        \ar@{-} |<(0.77){\hole} "2";"14";
    \endxy}
    \overset{\text{(8)}}{=}~
    \vcenter{\xy
        (-4,8)="1";
        ( 4,8)="2";
        ( 4,5)="3";
        ( 7,3)="4";
        ( 0,0)*[o]=<9pt>[Fo]{\scriptstyle s}="5";
        (-4,-2)="7";
        ( 0,-5)="8";
        ( 0,-8)*[o]=<5pt>[Fo]{~}="9";
        (-4,-9)="10";
        ( 0,-12)*[o]=<9pt>[Fo]{\scriptstyle t}="11";
        (-4,-16)="12";
        ( 0,-16)="13";
        ( 4,-16)="14";
        "2";"3" **\dir{-};
        "3";"4" **\dir{-};
        "4";"5" **\dir{-};
        "7";"8" **\dir{-};
        "5";"9" **\dir{-};
        "10";"11" **\dir{-};
        "11";"13" **\dir{-};
        "1";"12" **\dir{-};
        \ar@{-} |<(0.74){\hole} "2";"14";
    \endxy}
    \overset{\text{(2)}}{=}~
    \vcenter{\xy
        (-4,3)="1";
        ( 4,3)="2";
        ( 4,0)="3";
        ( 7,-2)="4";
        (-4,-2)="7";
        ( 0,-5)="8";
        ( 0,-8)*[o]=<5pt>[Fo]{~}="9";
        (-4,-9)="10";
        ( 0,-12)*[o]=<9pt>[Fo]{\scriptstyle t}="11";
        (-4,-16)="12";
        ( 0,-16)="13";
        ( 4,-16)="14";
        "2";"3" **\dir{-};
        "3";"4" **\dir{-};
        "4";"8" **\dir{-};
        "8";"9" **\dir{-};
        "7";"8" **\dir{-};
        "10";"11" **\dir{-};
        "11";"13" **\dir{-};
        "1";"12" **\dir{-};
        \ar@{-} |<(0.67){\hole} "2";"14";
    \endxy}
    ~\overset{\text{(c)}}{=}~
    \vcenter{\xy
        (-5,9)="1";
        ( 5,9)="2";
        (-5,6)="3";
        ( 5,6)="4";
        (-1,4)="5";
        ( 7,4)="6";
        ( 2,2)="7";
        ( 2,-1)*[o]=<5pt>[Fo]{~}="9";
        (-2,-1)*[o]=<9pt>[Fo]{\scriptstyle t}="8";
        (-5,-6)="10";
        (-2,-6)="11";
        ( 5,-6)="12";
        "1";"10" **\dir{-};
        "3";"5" **\dir{-};
        "5";"8" **\dir{-};
        "7";"9" **\dir{-};
        "4";"6" **\dir{-};
        "5";"7" **\dir{-};
        "6";"7" **\dir{-};
        "8";"11" **\dir{-};
        \ar@{-} |<(0.6){\hole} "2";"12";
    \endxy}
    ~\overset{\text{(12)}}{=}~
    \vcenter{\xy
        (-5,11)="1";
        ( 5,11)="2";
        ( 5,8)="3";
        ( 8,6)="4";
        ( 3,3)="5";
        ( 0,1)*[o]=<9pt>[Fo]{\scriptstyle s}="6";
        (-5,-1)="7";
        ( 3,-2)="8";
        (-3,-4)="9";
        (-3,-7)*[o]=<5pt>[Fo]{~}="10";
        (-5,-10)="11";
        ( 0,-10)="12";
        ( 5,-10)="13";
        "1";"11" **\dir{-};
        "7";"9" **\dir{-};
        "9";"10" **\dir{-};
        "3";"4" **\dir{-};
        "4";"5" **\dir{-};
        "5";"6" **\dir{-};
        "5";"8" **\dir{-};
        "8";"9" **\dir{-};
        \ar@{-} "6";"12" |<(0.3){\hole} 
        \ar@{-} "2";"13" |<(0.33){\hole}
    \endxy}
\\
    ~\overset{\text{(2)}}{=}~
    \vcenter{\xy
        (-5,11)="1";
        ( 7,11)="2";
        ( 7,8)="3";
        (10,6)="4";
        ( 3,3)="5";
        ( 0,1)*[o]=<9pt>[Fo]{\scriptstyle s}="6";
        (-5,-1)="7";
        ( 5,-1)*[o]=<9pt>[Fo]{\scriptstyle s}="8";
        (-3,-4)="9";
        (-3,-7)*[o]=<5pt>[Fo]{~}="10";
        (-5,-10)="11";
        ( 0,-10)="12";
        ( 7,-10)="13";
        "1";"11" **\dir{-};
        "7";"9" **\dir{-};
        "9";"10" **\dir{-};
        "3";"4" **\dir{-};
        "4";"5" **\dir{-};
        "5";"6" **\dir{-};
        "5";"8" **\dir{-};
        "8";"9" **\dir{-};
        \ar@{-} "6";"12" |<(0.28){\hole} 
        \ar@{-} "2";"13" |<(0.3){\hole}
    \endxy}
    ~\overset{\text{(c)}}{=}~
    \vcenter{\xy
        (-4,7)="1";
        ( 4,7)="2";
        ( 4,4)="3";
        ( 7,2)="4";
        ( 0,-1)*[o]=<9pt>[Fo]{\scriptstyle s}="5";
        (-4,-2)="7";
        ( 0,-5)="8";
        ( 0,-8)*[o]=<5pt>[Fo]{~}="9";
        ( 4,-8)="10";
        ( 7,-10)="x";
        ( 0,-13)*[o]=<9pt>[Fo]{\scriptstyle s}="11";
        (-4,-17)="12";
        ( 0,-17)="13";
        ( 4,-17)="14";
        "3";"4" **\dir{-};
        "4";"5" **\dir{-};
        "5";"9" **\dir{-};
        "7";"8" **\dir{-};
        "10";"x" **\dir{-};
        "x";"11" **\dir{-};
        "11";"13" **\dir{-};
        "1";"12" **\dir{-};
        \ar@{-} "2";"14" |<(0.27){\hole} |<(0.77){\hole}
    \endxy}
    ~\overset{\text{(8)}}{=}~
    \vcenter{\xy
        (-4,12)="1";
        ( 4,12)="2";
        ( 4,9)="3";
        ( 7,7)="4";
        ( 0,4)*[o]=<9pt>[Fo]{\scriptstyle s}="5";
        ( 0,-1)*[o]=<9pt>[Fo]{\scriptstyle t}="6";
        (-4,-2)="7";
        ( 0,-5)="8";
        ( 0,-8)*[o]=<5pt>[Fo]{~}="9";
        ( 4,-8)="10";
        ( 7,-10)="x";
        ( 0,-13)*[o]=<9pt>[Fo]{\scriptstyle s}="11";
        (-4,-17)="12";
        ( 0,-17)="13";
        ( 4,-17)="14";
        "3";"4" **\dir{-};
        "4";"5" **\dir{-};
        "5";"6" **\dir{-};
        "7";"8" **\dir{-};
        "6";"9" **\dir{-};
        "10";"x" **\dir{-};
        "x";"11" **\dir{-};
        "11";"13" **\dir{-};
        "1";"12" **\dir{-};
        \ar@{-} "2";"14" |<(0.23){\hole} |<(0.81){\hole}
    \endxy}
~=(1 \ox \gamma_l)d(1 \ox \gamma_l).
\end{multline*}
\end{proof}

The idempotent $d(1 \ox \gamma_l)$ will be denoted by $m$. The following
calculation gives a simpler representation of $m$:
\begin{equation*}\label{equ-m}
    \vcenter{\xy
        (0,13)*{~}; (0,-17)*{~};
        (-4,12)="1";
        ( 4,12)="2";
        ( 4,9)="3";
        ( 7,7)="4";
        ( 0,4)*[o]=<9pt>[Fo]{\scriptstyle s}="5";
        ( 0,-1)*[o]=<9pt>[Fo]{\scriptstyle t}="6";
        (-4,-2)="7";
        ( 0,-5)="8";
        ( 0,-8)*[o]=<5pt>[Fo]{~}="9";
        (-4,-11)="10";
        ( 4,-11)="14";
        "3";"4" **\dir{-};
        "4";"5" **\dir{-};
        "5";"6" **\dir{-};
        "7";"8" **\dir{-};
        "6";"9" **\dir{-};
        "1";"10" **\dir{-};
        \ar@{-} |<(0.72){\hole} "2";"14";
    \endxy}
    ~\overset{\text{(8)}}{=}~
    \vcenter{\xy
        (-4,8)="1";
        ( 4,8)="2";
        ( 4,5)="3";
        ( 7,3)="4";
        ( 0,0)*[o]=<9pt>[Fo]{\scriptstyle s}="5";
        (-4,-2)="7";
        ( 0,-5)="8";
        ( 0,-8)*[o]=<5pt>[Fo]{~}="9";
        (-4,-11)="10";
        ( 4,-11)="14";
        "3";"4" **\dir{-};
        "4";"5" **\dir{-};
        "7";"8" **\dir{-};
        "5";"9" **\dir{-};
        "1";"10" **\dir{-};
        \ar@{-} |<(0.68){\hole} "2";"14";
    \endxy}
    ~\overset{\text{(2)}}{=}~
    \vcenter{\xy
        (-4,3)="1";
        ( 4,3)="2";
        ( 4,0)="3";
        ( 7,-2)="4";
        (-4,-2)="7";
        ( 0,-5)="8";
        ( 0,-8)*[o]=<5pt>[Fo]{~}="9";
        (-4,-11)="10";
        ( 4,-11)="14";
        "3";"4" **\dir{-};
        "4";"8" **\dir{-};
        "8";"9" **\dir{-};
        "7";"8" **\dir{-};
        "1";"10" **\dir{-};
        \ar@{-} |<(0.56){\hole} "2";"14";
    \endxy}
    ~=m.
\end{equation*}
A splitting of $m$, i.e.,
\[
\xygraph{{(M \ox N,1)}="1"
    [r(2.6)]{(M \ox N,1)}="2"
    [d(1.2)l(1.3)]{(M \ox N,m)}="3"
    "1":"2" ^-m
    "1":"3" _-m
    "3":"2" _-m}
\qquad
\xygraph{{(M \ox N,m)}="1"
    [r(2.6)]{(M \ox N,m)}="2"
    [d(1.2)l(1.3)]{(M \ox N,1)}="3"
    "1":"2" ^-{m}
    "1":"3" _-{m} 
    "3":"2" _-{m}}
\]
provides an absolute equalizer $(M \ox N,m)$ of $(\gamma_r \ox 1)$ and $(1
\ox \gamma_l)$. Thus, the tensor product of $M$ and $N$ over $C$ is
\begin{equation*}\label{equ-mnc}
    M \oxC N = (M \ox N,m).
\end{equation*}

%=========================================================================%
\subsection{The coaction on the tensor product}
%=========================================================================%

If $\Comod(A)$ is to be a monoidal category with underlying functor
$U:\Comod(A) \ra \Bicomod(C)$ strong monoidal, then the tensor product of two
$A$-comodules must also be an $A$-comodule. In this section we show that the
obvious coaction on $M \oxC N$, namely,
\[
    \gamma =~
    \vcenter{\xy
        (-3,6)="1";
        ( 3,6)="2";
        (-3,3)="3";
        ( 3,3)="4";
        ( 6,-1)="7";
        (-3,-5)="5";
        ( 3,-5)="6";
        ( 6,-5)="8";
        "1";"5" **\dir{-};
        "3";"7" **\dir{-};
        "4";"7" **\dir{-};
        "8";"7" **\dir{-};
        \ar@{-} "2";"6" |<(0.55)\hole
    \endxy}
    ~:M \oxC N \ra M \oxC N \ox A
\]
does the job.

\begin{lemma}\label{lem-delta-mor}
The coaction $\gamma:M \ox_C N \ra M \ox_C N \ox A$, as defined above,
is a morphism in $\Q\V$. That is, the following equation holds.
\[
    \vcenter{\xy
        (0,4)*{~}; (0,-18)*{~};
        (-3,3)="1";
        ( 3,3)="2";
        ( 3,0)="3";
        ( 5,-2)="4";
        (-3,-2)="7";
        ( 0,-5)="8";
        ( 0,-8)*[o]=<5pt>[Fo]{~}="9";
        (-3,-10)="10";
        ( 3,-10)="14";
        ( 6,-14)="w";
        (-3,-17)="x";
        ( 3,-17)="y";
        ( 6,-17)="z";
        "2";"3" **\dir{-};
        "3";"4" **\dir{-};
        "4";"8" **\dir{-};
        "8";"9" **\dir{-};
        "7";"8" **\dir{-};
        "1";"10" **\dir{-};
        "10";"w" **\dir{-};
        "14";"w" **\dir{-};
        "10";"x" **\dir{-};
        "w";"z" **\dir{-};
        \ar@{-} "14";"y" |<(0.4){\hole}
        \ar@{-} "2";"14" |<(0.48){\hole}
    \endxy}
    \quad = \quad
    \vcenter{\xy
        (-3,6)="1";
        ( 3,6)="2";
        (-3,3)="3";
        ( 3,3)="4";
        ( 6,-1)="7";
        (-3,-5)="5";
        ( 3,-5)="6";
        ( 6,-5)="8";
        "1";"5" **\dir{-};
        "3";"7" **\dir{-};
        "4";"7" **\dir{-};
        "8";"7" **\dir{-};
        \ar@{-} "2";"6" |<(0.55)\hole
    \endxy}
    \quad = \quad
    \vcenter{\xy
        (-3,9)="a";
        ( 3,9)="b";
        (-3,6)="c";
        ( 3,6)="d";
        ( 6,2)="w";
        ( 6,-10)="z";
        (-3,3)="1";
        ( 3,1)="2";
        ( 3,0)="3";
        ( 5,-2)="4";
        (-3,-2)="7";
        ( 0,-5)="8";
        ( 0,-8)*[o]=<5pt>[Fo]{~}="9";
        (-3,-10)="10";
        ( 3,-10)="14";
        "a";"1" **\dir{-};
        "c";"w" **\dir{-};
        "d";"w" **\dir{-};
        "w";"z" **\dir{-};
        "2";"3" **\dir{-};
        "3";"4" **\dir{-};
        "4";"8" **\dir{-};
        "8";"9" **\dir{-};
        "7";"8" **\dir{-};
        "1";"10" **\dir{-};
        \ar@{-} "b";"2" |<(0.75){\hole}
        \ar@{-} "2";"14" |<(0.4){\hole}
    \endxy}
\]
\end{lemma}

\begin{proof}
The first equality is given by
\[
    \vcenter{\xy
        (0,4)*{~}; (0,-18)*{~};
        (-3,3)="1";
        ( 3,3)="2";
        ( 3,0)="3";
        ( 5,-2)="4";
        (-3,-2)="7";
        ( 0,-5)="8";
        ( 0,-8)*[o]=<5pt>[Fo]{~}="9";
        (-3,-10)="10";
        ( 3,-10)="14";
        ( 6,-14)="w";
        (-3,-17)="x";
        ( 3,-17)="y";
        ( 6,-17)="z";
        "2";"3" **\dir{-};
        "3";"4" **\dir{-};
        "4";"8" **\dir{-};
        "8";"9" **\dir{-};
        "7";"8" **\dir{-};
        "1";"10" **\dir{-};
        "10";"w" **\dir{-};
        "14";"w" **\dir{-};
        "10";"x" **\dir{-};
        "w";"z" **\dir{-};
        \ar@{-} "14";"y" |<(0.4){\hole}
        \ar@{-} "2";"14" |<(0.48){\hole}
    \endxy}
    ~\overset{\text{(c)}}{=}~
    \vcenter{\xy
        (-3,9)="1";
        ( 3,9)="2";
        (-3,6)="3";
        ( 3,6)="4";
        ( 5,2)="5";
        (11,2)="6";
        ( 5,-2)="7";
        (11,-2)="8";
        (11,-5)*[o]=<5pt>[Fo]{~}="12";
        (-3,-9)="9";
        ( 3,-9)="10";
        ( 5,-9)="11";
        "1";"9" **\dir{-};
        "3";"5" **\dir{-};
        "4";"6" **\dir{-};
        "5";"11" **\dir{-};
        "6";"12" **\dir{-};
        "5";"8" **\dir{-};
        \ar@{-} "6";"7" |<(0.5){\hole}
        \ar@{-} "2";"10" |<(0.34){\hole}
    \endxy}
    ~\overset{\text{(b)}}{=}~
    \vcenter{\xy
        (-3,9)="1";
        ( 3,9)="2";
        (-3,6)="3";
        ( 3,6)="4";
        ( 6,2)="5";
        ( 6,-1)="6";
        ( 8,-3)*[o]=<5pt>[Fo]{~}="7";
        (-3,-6)="8";
        ( 3,-6)="9";
        ( 6,-6)="10";
        "1";"8" **\dir{-};
        "3";"5" **\dir{-};
        "4";"5" **\dir{-};
        "5";"10" **\dir{-};
        "6";"7" **\dir{-};
        \ar@{-} "2";"9" |<(0.39){\hole}
    \endxy}
    ~\overset{\text{(c)}}{=}~
    \vcenter{\xy
        (-3,6)="1";
        ( 3,6)="2";
        (-3,3)="3";
        ( 3,3)="4";
        ( 6,-1)="7";
        (-3,-5)="5";
        ( 3,-5)="6";
        ( 6,-5)="8";
        "1";"5" **\dir{-};
        "3";"7" **\dir{-};
        "4";"7" **\dir{-};
        "8";"7" **\dir{-};
        \ar@{-} "2";"6" |<(0.55)\hole
    \endxy}~,
\]
and the second by a similar calculation:
\[
    \vcenter{\xy
        (0,10)*{~}; (0,-11)*{~};
        (-3,9)="a";
        ( 3,9)="b";
        (-3,6)="c";
        ( 3,6)="d";
        ( 6,2)="w";
        ( 6,-10)="z";
        (-3,3)="1";
        ( 3,1)="2";
        ( 3,0)="3";
        ( 5,-2)="4";
        (-3,-2)="7";
        ( 0,-5)="8";
        ( 0,-8)*[o]=<5pt>[Fo]{~}="9";
        (-3,-10)="10";
        ( 3,-10)="14";
        "a";"1" **\dir{-};
        "c";"w" **\dir{-};
        "d";"w" **\dir{-};
        "w";"z" **\dir{-};
        "2";"3" **\dir{-};
        "3";"4" **\dir{-};
        "4";"8" **\dir{-};
        "8";"9" **\dir{-};
        "7";"8" **\dir{-};
        "1";"10" **\dir{-};
        \ar@{-} "b";"2" |<(0.75){\hole}
        \ar@{-} "2";"14" |<(0.4){\hole}
    \endxy}
    ~\overset{\text{(c)}}{=}~
    \vcenter{\xy
        (-3,9)="1";
        ( 3,9)="2";
        (-3,6)="3";
        ( 3,6)="4";
        ( 5,2)="5";
        (11,2)="6";
        ( 5,-2)="7";
        (11,-2)="8";
        ( 5,-5)*[o]=<5pt>[Fo]{~}="11";
        (11,-9)="12";
        (-3,-9)="9";
        ( 3,-9)="10";
        "1";"9" **\dir{-};
        "3";"5" **\dir{-};
        "4";"6" **\dir{-};
        "5";"11" **\dir{-};
        "6";"12" **\dir{-};
        "5";"8" **\dir{-};
        \ar@{-} "6";"7" |<(0.5){\hole}
        \ar@{-} "2";"10" |<(0.34){\hole}
    \endxy}
    ~\overset{\text{(b)}}{=}~
    \vcenter{\xy
        (-3,9)="1";
        ( 3,9)="2";
        (-3,6)="3";
        ( 3,6)="4";
        ( 7,2)="5";
        ( 7,-1)="6";
        ( 5,-3)*[o]=<5pt>[Fo]{~}="7";
        (-3,-6)="8";
        ( 3,-6)="9";
        ( 7,-6)="10";
        "1";"8" **\dir{-};
        "3";"5" **\dir{-};
        "4";"5" **\dir{-};
        "5";"10" **\dir{-};
        "6";"7" **\dir{-};
        \ar@{-} "2";"9" |<(0.39){\hole}
    \endxy}
    ~\overset{\text{(c)}}{=}~
    \vcenter{\xy
        (-3,6)="1";
        ( 3,6)="2";
        (-3,3)="3";
        ( 3,3)="4";
        ( 6,-1)="7";
        (-3,-5)="5";
        ( 3,-5)="6";
        ( 6,-5)="8";
        "1";"5" **\dir{-};
        "3";"7" **\dir{-};
        "4";"7" **\dir{-};
        "8";"7" **\dir{-};
        \ar@{-} "2";"6" |<(0.55)\hole
    \endxy}~.
\]
\end{proof}

\begin{proposition}
$(M \oxC N,\gamma)$ is an $A$-comodule.
\end{proposition}

\begin{proof}
Coassociativity is proved as usual,
\[
    \vcenter{\xy
        (-3,9)="1";
        ( 3,9)="2";
        (-3,6)="3";
        ( 3,6)="4";
        ( 7,2)="5";
        ( 7,-1)="6";
        (-3,-5)="7";
        ( 3,-5)="8";
        ( 5,-5)="9";
        ( 9,-5)="10";
        "1";"7" **\dir{-};
        "3";"5" **\dir{-};
        "4";"5" **\dir{-};
        "5";"6" **\dir{-};
        "6";"9" **\dir{-};
        "6";"10" **\dir{-};
        \ar@{-} "2";"8" |<(0.39){\hole}
    \endxy}
    ~\overset{\text{(b)}}{=}~
    \vcenter{\xy
        (-3,9)="1";
        ( 3,9)="2";
        (-3,6)="3";
        ( 3,6)="4";
        ( 5,2)="5";
        (11,2)="6";
        ( 5,-2)="7";
        (11,-2)="8";
        ( 5,-7)="11";
        (11,-7)="12";
        (-3,-7)="9";
        ( 3,-7)="10";
        "1";"9" **\dir{-};
        "3";"5" **\dir{-};
        "4";"6" **\dir{-};
        "5";"11" **\dir{-};
        "6";"12" **\dir{-};
        "5";"8" **\dir{-};
        \ar@{-} "6";"7" |<(0.5){\hole}
        \ar@{-} "2";"10" |<(0.34){\hole}
    \endxy}
    ~\overset{\text{(c)}}{=}~
    \vcenter{\xy
        (-3,9)="1";
        ( 3,9)="2";
        (-3,6)="3";
        ( 3,6)="4";
        (-3,1)="5";
        ( 3,1)="6";
        (-3,-7)="7";
        ( 3,-7)="8";
        ( 9,2)="9";
        ( 6,-3)="10";
        ( 6,-7)="11";
        ( 9,-7)="12";
        "1";"7" **\dir{-};
        "3";"9" **\dir{-};
        "4";"9" **\dir{-};
        "5";"10" **\dir{-};
        "6";"10" **\dir{-};
        "11";"10" **\dir{-};
        "9";"12" **\dir{-};
        \ar@{-} "2";"8" |<(0.34){\hole} |<(0.68){\hole}
    \endxy}
\]
and the counit condition as
\[
    \vcenter{\xy
        (-3,3)="1";
        ( 3,3)="2";
        (-3,0)="3";
        ( 3,0)="4";
        (-3,-9)="5";
        ( 3,-9)="6";
        ( 5,-3)="7";
        ( 5,-6)*[o]=<5pt>[Fo]{~}="8";
        "1";"5" **\dir{-};
        "3";"7" **\dir{-};
        "4";"7" **\dir{-};
        "7";"8" **\dir{-};
        \ar@{-} "2";"6" |<(0.48){\hole}
    \endxy}
    \overset{\text{(\text{L.~\ref{lem-delta-mor}})}}{=}~
    \vcenter{\xy
        (0,4)*{~}; (0,-21)*{~};
        (-3,3)="1";
        ( 3,3)="2";
        ( 3,0)="3";
        ( 5,-2)="4";
        (-3,-2)="7";
        ( 0,-5)="8";
        ( 0,-8)*[o]=<5pt>[Fo]{~}="9";
        (-3,-10)="10";
        ( 3,-10)="14";
        ( 6,-14)="w";
        (-3,-20)="x";
        ( 3,-20)="y";
        ( 6,-17)*[o]=<5pt>[Fo]{~}="z";
        "2";"3" **\dir{-};
        "3";"4" **\dir{-};
        "4";"8" **\dir{-};
        "8";"9" **\dir{-};
        "7";"8" **\dir{-};
        "1";"10" **\dir{-};
        "10";"w" **\dir{-};
        "14";"w" **\dir{-};
        "10";"x" **\dir{-};
        "w";"z" **\dir{-};
        \ar@{-} "14";"y" |<(0.3){\hole}
        \ar@{-} "2";"14" |<(0.48){\hole}
    \endxy}
    ~\overset{\text{(c)}}{=}~
    \vcenter{\xy
        (-4,9)="1";
        ( 4,9)="2";
        (-4,4)="3";
        ( 4,6)="4";
        ( 6,4)="5";
        (-2,2)="6";
        ( 2,2)="7";
        (-2,-3)="8";
        ( 2,-3)="9";
        (-2,-6)*[o]=<5pt>[Fo]{~}="10";
        ( 2,-6)*[o]=<5pt>[Fo]{~}="11";
        (-4,-9)="12";
        ( 4,-9)="13";
        "1";"12" **\dir{-};
        "3";"6" **\dir{-};
        "4";"5" **\dir{-};
        "5";"7" **\dir{-};
        "6";"10" **\dir{-};
        "6";"9" **\dir{-};
        "7";"11" **\dir{-};
        \ar@{-} "2";"13" |<(0.33){\hole}
        \ar@{-} "7";"8" |<(0.5){\hole}
    \endxy}
    ~\overset{\text{(b)}}{=}~
    \vcenter{\xy
        (-4,3)="1";
        ( 4,3)="2";
        ( 4,0)="3";
        ( 6,-2)="4";
        (-4,-2)="7";
        ( 0,-5)="8";
        ( 0,-8)="9";
        (-2,-10)*[o]=<5pt>[Fo]{~}="x";
        ( 2,-10)*[o]=<5pt>[Fo]{~}="y";
        (-4,-13)="10";
        ( 4,-13)="14";
        "2";"3" **\dir{-};
        "3";"4" **\dir{-};
        "4";"8" **\dir{-};
        "8";"9" **\dir{-};
        "9";"x" **\dir{-};
        "9";"y" **\dir{-};
        "7";"8" **\dir{-};
        "1";"10" **\dir{-};
        \ar@{-} "2";"14" |<(0.38){\hole}
    \endxy}
    ~\overset{\text{(c)}}{=}~
    \vcenter{\xy
        (-3,3)="1";
        ( 3,3)="2";
        ( 3,0)="3";
        ( 5,-2)="4";
        (-3,-2)="7";
        ( 0,-5)="8";
        ( 0,-8)*[o]=<5pt>[Fo]{~}="9";
        (-3,-10)="10";
        ( 3,-10)="14";
        "2";"3" **\dir{-};
        "3";"4" **\dir{-};
        "4";"8" **\dir{-};
        "8";"9" **\dir{-};
        "7";"8" **\dir{-};
        "1";"10" **\dir{-};
        \ar@{-} "2";"14" |<(0.48){\hole}
    \endxy}
= 1_{M \oxC N}.
\]
\end{proof}

%=========================================================================%
\subsection{$\Comod(A)$ is a monoidal category}
%=========================================================================%

We now set out to prove the claim at the beginning of this section, that
$(\Comod(A),\oxC,C)$ is a monoidal category. It will turn out that
associativity is a strict equality (if it is so in $\V$) and the unit
conditions are only up to isomorphism.

We state this as a theorem and devote the remainder of this section to its
proof.

\begin{theorem}\label{thm-comod-mon}
$\Comod(A) = (\Comod(A),\oxC,C)$ is a monoidal category.
\end{theorem}

First off note that $C$ itself is an $A$-comodule with coaction
\[
    \vcenter{\xy
        (-2,10)*{\scriptstyle C};
        (-3,-5)*{\scriptstyle C};
        ( 3,-5)*{\scriptstyle A};
        (0,9)="1";
        (0,5)*[o]=<9pt>[Fo]{\scriptstyle t}="2";
        (0,1)="3";
        (-3,-3)="4";
        ( 3,-3)="5";
        "1";"2" **\dir{-};
        "2";"3" **\dir{-};
        "3";"4" **\dir{-};
        "3";"5" **\dir{-};
    \endxy}~.
\]
Before proving this theorem it will be useful to have the following lemma.

\begin{lemma}\label{lem-useful}
The following identities hold.
\[
    \vcenter{\xy
        (-3,11)="1";
        ( 3,11)="2";
        ( 3,7)*[o]=<9pt>[Fo]{\scriptstyle t}="x";
        (-3,3)="3";
        ( 3,3)="4";
        ( 6,-1)="7";
        (-3,-5)="5";
        ( 3,-5)="6";
        ( 6,-5)="8";
        "1";"5" **\dir{-};
        "3";"7" **\dir{-};
        "4";"7" **\dir{-};
        "8";"7" **\dir{-};
        "2";"x" **\dir{-};
        \ar@{-} "x";"6" |<(0.52)\hole
    \endxy}
    ~=~
    \vcenter{\xy
        (-3,10)="1";
        ( 3,10)="2";
        (-3,6)="3";
        ( 3,6)*[o]=<9pt>[Fo]{\scriptstyle t}="4";
        (-3,2)="5";
        ( 3,2)="6";
        ( 0,0)*[o]=<9pt>[Fo]{\scriptstyle t}="7";
        (-3,-4)="8";
        ( 0,-4)="9";
        ( 3,-4)="10";
        "1";"8" **\dir{-};
        "2";"4" **\dir{-};
        "4";"10" **\dir{-};
        "3";"6" **\dir{-};
        "5";"7" **\dir{-};
        "7";"9" **\dir{-};
    \endxy}
    \qquad\qquad\qquad
    \vcenter{\xy
        (-3,11)="1";
        ( 3,11)="2";
        (-3,7)*[o]=<9pt>[Fo]{\scriptstyle t}="x";
        (-3,3)="3";
        ( 3,3)="4";
        ( 6,-1)="7";
        (-3,-5)="5";
        ( 3,-5)="6";
        ( 6,-5)="8";
        "1";"x" **\dir{-};
        "x";"5" **\dir{-};
        "3";"7" **\dir{-};
        "4";"7" **\dir{-};
        "8";"7" **\dir{-};
        \ar@{-} "2";"6" |<(0.69)\hole
    \endxy}
    ~=~
    \vcenter{\xy
        (-4,9)="1";
        ( 4,9)="2";
        (-4,5)*[o]=<9pt>[Fo]{\scriptstyle t}="3";
        ( 4,5)="4";
        (-2,1)="5";
        ( 4,1)="6";
        ( 0,-1)="7";
        (-5,-4)*[o]=<9pt>[Fo]{\scriptstyle s}="8";
        (-5,-8)="9";
        (-2,-8)="10";
        ( 4,-8)="11";
        "1";"3" **\dir{-};
        "3";"6" **\dir{-};
        "2";"11" **\dir{-};
        "5";"7" **\dir{-};
        "7";"8" **\dir{-};
        "8";"9" **\dir{-};
        \ar@{-} "4";"5" |<(0.55)\hole
        \ar@{-} "5";"10" |<(0.35)\hole
    \endxy}
\]
\end{lemma}

\begin{proof}
The first identity is proved by
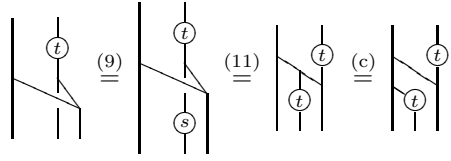
\[
    \vcenter{\xy
        (-3,11)="1";
        ( 3,11)="2";
        ( 3,7)*[o]=<9pt>[Fo]{\scriptstyle t}="x";
        (-3,3)="3";
        ( 3,3)="4";
        ( 6,-1)="7";
        (-3,-5)="5";
        ( 3,-5)="6";
        ( 6,-5)="8";
        "1";"5" **\dir{-};
        "3";"7" **\dir{-};
        "4";"7" **\dir{-};
        "8";"7" **\dir{-};
        "2";"x" **\dir{-};
        \ar@{-} "x";"6" |<(0.52)\hole
    \endxy}
    ~\overset{\text{(9)}}{=}~
    \vcenter{\xy
        (0,12)*{~};
        (0,-10)*{~};
        (-3,11)="1";
        ( 3,11)="2";
        ( 3,7)*[o]=<9pt>[Fo]{\scriptstyle t}="x";
        (-3,3)="3";
        ( 3,3)="4";
        ( 6,-1)="7";
        ( 3,-5)*[o]=<9pt>[Fo]{\scriptstyle s}="y";
        (-3,-9)="5";
        ( 3,-9)="6";
        ( 6,-9)="8";
        "1";"5" **\dir{-};
        "3";"7" **\dir{-};
        "4";"7" **\dir{-};
        "8";"7" **\dir{-};
        "2";"x" **\dir{-};
        "y";"6" **\dir{-};
        \ar@{-} "x";"y" |<(0.52)\hole
    \endxy}
    ~\overset{\text{(11)}}{=}~
    \vcenter{\xy
        (-3,10)="1";
        ( 3,10)="2";
        (-3,6)="3";
        ( 3,6)*[o]=<9pt>[Fo]{\scriptstyle t}="4";
        ( 0,4)="5";
        ( 3,2)="6";
        ( 0,0)*[o]=<9pt>[Fo]{\scriptstyle t}="7";
        (-3,-4)="8";
        ( 0,-4)="9";
        ( 3,-4)="10";
        "1";"8" **\dir{-};
        "2";"4" **\dir{-};
        "4";"10" **\dir{-};
        "3";"6" **\dir{-};
        "5";"7" **\dir{-};
        "7";"9" **\dir{-};
    \endxy}
    ~\overset{\text{(c)}}{=}~
    \vcenter{\xy
        (-3,10)="1";
        ( 3,10)="2";
        (-3,6)="3";
        ( 3,6)*[o]=<9pt>[Fo]{\scriptstyle t}="4";
        (-3,2)="5";
        ( 3,2)="6";
        ( 0,0)*[o]=<9pt>[Fo]{\scriptstyle t}="7";
        (-3,-4)="8";
        ( 0,-4)="9";
        ( 3,-4)="10";
        "1";"8" **\dir{-};
        "2";"4" **\dir{-};
        "4";"10" **\dir{-};
        "3";"6" **\dir{-};
        "5";"7" **\dir{-};
        "7";"9" **\dir{-};
    \endxy}
\]
and the second by
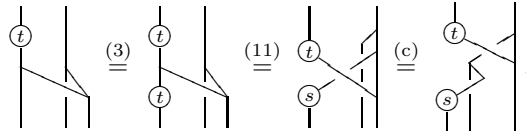
\[
    \vcenter{\xy
        (-3,11)="1";
        ( 3,11)="2";
        (-3,7)*[o]=<9pt>[Fo]{\scriptstyle t}="x";
        (-3,3)="3";
        ( 3,3)="4";
        ( 6,-1)="7";
        (-3,-5)="5";
        ( 3,-5)="6";
        ( 6,-5)="8";
        "1";"x" **\dir{-};
        "x";"5" **\dir{-};
        "3";"7" **\dir{-};
        "4";"7" **\dir{-};
        "8";"7" **\dir{-};
        \ar@{-} "2";"6" |<(0.69)\hole
    \endxy}
    ~\overset{\text{(3)}}{=}~
    \vcenter{\xy
        (-3,11)="1";
        ( 3,11)="2";
        (-3,7)*[o]=<9pt>[Fo]{\scriptstyle t}="x";
        (-3,3)="3";
        (-3,-1)*[o]=<9pt>[Fo]{\scriptstyle t}="y";
        ( 3,3)="4";
        ( 6,-1)="7";
        (-3,-5)="5";
        ( 3,-5)="6";
        ( 6,-5)="8";
        "1";"x" **\dir{-};
        "x";"y" **\dir{-};
        "y";"5" **\dir{-};
        "3";"7" **\dir{-};
        "4";"7" **\dir{-};
        "8";"7" **\dir{-};
        \ar@{-} "2";"6" |<(0.69)\hole
    \endxy}
    ~\overset{\text{(11)}}{=}~
    \vcenter{\xy
        (-4,12)="1";
        ( 5,12)="2";
        (-4,6)*[o]=<9pt>[Fo]{\scriptstyle t}="3";
        ( 5,6)="4";
        (-4,0)*[o]=<9pt>[Fo]{\scriptstyle s}="5";
        ( 5,0)="6";
        (-4,-4)="7";
        ( 5,-4)="8";
        ( 5, 9)="x";
        ( 3, 7)="y";
        ( 3,-4)="z";
        "1";"3" **\dir{-};
        "3";"6" **\dir{-};
        "2";"8" **\dir{-};
        "5";"7" **\dir{-};
        "x";"y" **\dir{-};
        \ar@{-} "y";"z" |<(0.21)\hole |<(0.55)\hole
        \ar@{-} "4";"5" |<(0.5)\hole
    \endxy}
    ~\overset{\text{(c)}}{=}~
    \vcenter{\xy
        (-4,9)="1";
        ( 4,9)="2";
        (-4,5)*[o]=<9pt>[Fo]{\scriptstyle t}="3";
        ( 4,5)="4";
        (-2,1)="5";
        ( 4,1)="6";
        ( 0,-1)="7";
        (-5,-4)*[o]=<9pt>[Fo]{\scriptstyle s}="8";
        (-5,-8)="9";
        (-2,-8)="10";
        ( 4,-8)="11";
        "1";"3" **\dir{-};
        "3";"6" **\dir{-};
        "2";"11" **\dir{-};
        "5";"7" **\dir{-};
        "7";"8" **\dir{-};
        "8";"9" **\dir{-};
        \ar@{-} "4";"5" |<(0.55)\hole
        \ar@{-} "5";"10" |<(0.35)\hole
    \endxy}~.
\]
\end{proof}

\begin{proof}[Proof of Theorem~\ref{thm-comod-mon}]
Consider $(M \oxC N) \oxC P$ and $M \oxC (N \oxC P)$ in $\Q\V$. The former
is $(M \ox N \ox P,u)$ and the latter $(M \ox N \ox P,v)$ where
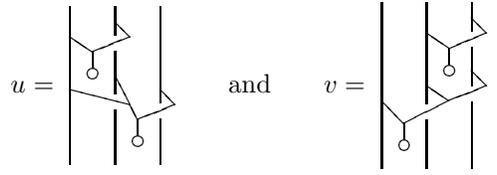
\[
    u =~
    \vcenter{\xy
        (-6,9)="1";
        ( 0,9)="2";
        ( 6,9)="3";
        (-6,5)="4";
        ( 0,7)="5";
        ( 2,5)="6";
        (-3,3)="7";
        (-3,0)*[o]=<5pt>[Fo]{~}="8";
        (-6,-2)="9";
        ( 0, 0)="10";
        ( 6,-2)="11";
        ( 8,-4)="12";
        ( 2,-4)="13";
        ( 3,-6)="14";
        ( 3,-9)*[o]=<5pt>[Fo]{~}="15";
        (-6,-12)="16";
        ( 0,-12)="17";
        ( 6,-12)="18";
        "1";"16" **\dir{-};
        "4";"7" **\dir{-};
        "7";"8" **\dir{-};
        "5";"6" **\dir{-};
        "6";"7" **\dir{-};
        "9";"13" **\dir{-};
        "10";"13" **\dir{-};
        "13";"14" **\dir{-};
        "14";"15" **\dir{-};
        "11";"12" **\dir{-};
        "12";"14" **\dir{-};
        \ar@{-} "2";"17" |<(0.23){\hole} |<(0.6){\hole}
        \ar@{-} "3";"18" |<(0.66){\hole}
    \endxy}
    \qquad \text{and} \qquad
    v =~
    \vcenter{\xy
        (0,10)*{~}; (0,-14)*{~};
        (-6,9)="1";
        ( 0,9)="2";
        ( 6,9)="3";
        ( 0,5)="4";
        ( 6,7)="5";
        ( 8,5)="6";
        ( 3,3)="7";
        ( 3,0)*[o]=<5pt>[Fo]{~}="8";
        ( 0,-2)="9";
        ( 6, 0)="10";
        ( 8,-2)="11";
        ( 3,-4)="12";
        (-6,-4)="13";
        (-3,-7)="14";
        (-3,-10)*[o]=<5pt>[Fo]{~}="15";
        (-6,-13)="16";
        ( 0,-13)="17";
        ( 6,-13)="18";
        "1";"16" **\dir{-};
        "4";"7" **\dir{-};
        "7";"8" **\dir{-};
        "5";"6" **\dir{-};
        "6";"7" **\dir{-};
        "9";"12" **\dir{-};
        "10";"11" **\dir{-};
        "11";"12" **\dir{-};
        "14";"12" **\dir{-};
        "13";"14" **\dir{-};
        "14";"15" **\dir{-};
        \ar@{-} "2";"17" |<(0.66){\hole}
        \ar@{-} "3";"18" |<(0.22){\hole} |<(0.55){\hole}
    \endxy}~.
\]
Since, by Lemma~\ref{lem-delta-mor}, $\gamma$ is a morphism in $\Q\V$, both
$u$ and $v$ may be rewritten as
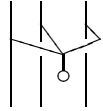
\[
    \vcenter{\xy
        (0,11)*{~}; (0,-5)*{~};
        (-4,10)="1";
        ( 0,10)="2";
        ( 6,10)="3";
        (-4,5)="4";
        ( 0,7)="5";
        ( 6,7)="6";
        ( 8,5)="7";
        ( 3,3)="8";
        ( 3,0)*[o]=<5pt>[Fo]{~}="9";
        (-4,-4)="10";
        ( 0,-4)="11";
        ( 6,-4)="12";
        "1";"10" **\dir{-};
        "4";"8" **\dir{-};
        "5";"8" **\dir{-};
        "6";"7" **\dir{-};
        "7";"8" **\dir{-};
        "8";"9" **\dir{-};
        \ar@{-} "2";"11" |<(0.45){\hole}
        \ar@{-} "3";"12" |<(0.42){\hole}
    \endxy}
\]
proving the (strict) equality $(M \oxC N) \oxC P = M \oxC (N \oxC P)$ in
$\Q\V$ (since we are writing as if $\V$ were strict).

It remains to prove $M \oxC C \cong M \cong C \oxC M$. By definition
\[
    M \oxC C = (M \ox C,
    ~\vcenter{\xy
        (-3,7)="1";
        ( 3,7)="2";
        ( 3,3)*[o]=<9pt>[Fo]{\scriptstyle t}="x";
        ( 3,0)="3";
        ( 5,-2)="4";
        (-3,-2)="7";
        ( 0,-4)="8";
        ( 0,-7)*[o]=<5pt>[Fo]{~}="9";
        (-3,-10)="10";
        ( 3,-10)="14";
        "2";"x" **\dir{-};
        "3";"4" **\dir{-};
        "4";"8" **\dir{-};
        "8";"9" **\dir{-};
        "7";"8" **\dir{-};
        "1";"10" **\dir{-};
        \ar@{-} "x";"14" |<(0.4){\hole}
    \endxy}~)
\qquad \text{and} \qquad
    C \oxC M = (C \ox M,
    ~\vcenter{\xy
        (-3,6)="1";
        ( 3,6)="2";
        (-3,2)*[o]=<9pt>[Fo]{\scriptstyle t}="x";
        ( 3,0)="3";
        ( 5,-2)="4";
        (-3,-2)="7";
        ( 0,-4)="8";
        ( 0,-7)*[o]=<5pt>[Fo]{~}="9";
        (-3,-10)="10";
        ( 3,-10)="14";
        "3";"4" **\dir{-};
        "4";"8" **\dir{-};
        "8";"9" **\dir{-};
        "7";"8" **\dir{-};
        "1";"x" **\dir{-};
        "x";"10" **\dir{-};
        \ar@{-} "2";"14" |<(0.57){\hole}
    \endxy}~)~.
\]
We will show that the morphisms
\[
    \vcenter{\xy
        (0,8)*{~}; (0,-8)*{~};
        (-2, 7)="1";
        ( 2, 7)="2";
        ( 2, 3)*[o]=<9pt>[Fo]{\scriptstyle t}="3";
        (-2, 2)="4";
        ( 2,-1)="5";
        ( 2,-4)*[o]=<5pt>[Fo]{~}="6";
        (-2,-7)="7";
        "1";"7" **\dir{-};
        "2";"3" **\dir{-};
        "3";"6" **\dir{-};
        "4";"5" **\dir{-};
    \endxy}~:
    M \oxC C \dra M
    \qquad \text{and} \qquad
    \vcenter{\xy
        (-2,11)="1";
        (-2, 8)="2";
        ( 2, 5)*[o]=<9pt>[Fo]{\scriptstyle t}="3";
        (-2, 0)="4";
        ( 2, 0)="5";
        "1";"4" **\dir{-};
        "2";"3" **\dir{-};
        "3";"5" **\dir{-};
    \endxy}~:
    M \dra M \oxC C
\]
will establish the isomorphism $M \oxC C \cong M$, and 
\[
    \vcenter{\xy
        (0,8)*{~}; (0,-8)*{~};
        (-2, 7)="1";
        ( 2, 7)="2";
        (-2, 3)*[o]=<9pt>[Fo]{\scriptstyle t}="3";
        ( 2, 4)="4";
        ( 4, 2)="5";
        (-2,-1)="6";
        (-2,-4)*[o]=<5pt>[Fo]{~}="7";
        ( 2,-7)="8";
        "1";"3" **\dir{-};
        "3";"7" **\dir{-};
        "4";"5" **\dir{-};
        "5";"6" **\dir{-};
        \ar@{-} "2";"8" |<(0.45){\hole}
    \endxy}~:
    C \oxC M \dra M
    \qquad \text{and} \qquad
    \vcenter{\xy
        ( 4, 7)="1";
        ( 4, 4)="2";
        ( 6, 2)="3";
        ( 0,-1)*[o]=<9pt>[Fo]{\scriptstyle s}="4";
        ( 0,-5)="5";
        ( 4,-5)="6";
        "2";"3" **\dir{-};
        "3";"4" **\dir{-};
        "4";"5" **\dir{-};
        \ar@{-} "1";"6" |<(0.5){\hole}
    \endxy}~:
    M \dra C \oxC M
\]
the isomorphism $M \cong C \oxC M$. These morphisms are easily seen to be in
$\Q\V$, and the fact that they are mutually inverse pairs is given in one
direction by Lemma~\ref{lem-useful}, and in the other by an easy string
calculation making use of the identity (6).

It now remains to show that these four morphisms are $A$-comodules morphisms,
i.e., that they are in $\Comod(A)$. Note that $M \oxC C$ and $C \oxC M$ are
$A$-comodules via the coactions
\[
    \vcenter{\xy
        (0,12)*{~};(0,-6)*{~};
        (-3,11)="1";
        ( 3,11)="2";
        ( 3,7)*[o]=<9pt>[Fo]{\scriptstyle t}="x";
        (-3,3)="3";
        ( 3,3)="4";
        ( 6,-1)="7";
        (-3,-5)="5";
        ( 3,-5)="6";
        ( 6,-5)="8";
        "1";"5" **\dir{-};
        "3";"7" **\dir{-};
        "4";"7" **\dir{-};
        "8";"7" **\dir{-};
        "2";"x" **\dir{-};
        \ar@{-} "x";"6" |<(0.52)\hole
    \endxy}
    \quad\qquad \text{and} \quad\qquad
    \vcenter{\xy
        (-3,11)="1";
        ( 3,11)="2";
        (-3,7)*[o]=<9pt>[Fo]{\scriptstyle t}="x";
        (-3,3)="3";
        ( 3,3)="4";
        ( 6,-1)="7";
        (-3,-5)="5";
        ( 3,-5)="6";
        ( 6,-5)="8";
        "1";"x" **\dir{-};
        "x";"5" **\dir{-};
        "3";"7" **\dir{-};
        "4";"7" **\dir{-};
        "8";"7" **\dir{-};
        \ar@{-} "2";"6" |<(0.69)\hole
    \endxy}
\]
respectively. We then have:
\begin{itemize}

\item $~\vcenter{\xy
        (-2, 7)="1";
        ( 2, 7)="2";
        ( 2, 3)*[o]=<9pt>[Fo]{\scriptstyle t}="3";
        (-2, 2)="4";
        ( 2,-1)="5";
        ( 2,-4)*[o]=<5pt>[Fo]{~}="6";
        (-2,-7)="7";
        "1";"7" **\dir{-};
        "2";"3" **\dir{-};
        "3";"6" **\dir{-};
        "4";"5" **\dir{-};
    \endxy}~:
M \oxC C \dra M$ is an $A$-comodule morphism as:
\[
    \vcenter{\xy
        (-4,12)="1";
        ( 3,12)="2";
        ( 3, 8)*[o]=<9pt>[Fo]{\scriptstyle t}="3";
        (-4, 4)="4";
        ( 3, 4)="5";
        ( 5, 0)="6";
        ( 0,-3)*[o]=<9pt>[Fo]{\scriptstyle t}="8";
        (-4,-5)="7";
        (-2,-7)="9";
        (-2,-10)*[o]=<5pt>[Fo]{~}="10";
        (-4,-13)="11";
        ( 5,-13)="12";
        "4";"6" **\dir{-};
        "5";"6" **\dir{-};
        "2";"3" **\dir{-};
        "3";"5" **\dir{-};
        "8";"9" **\dir{-};
        "10";"9" **\dir{-};
        "7";"9" **\dir{-};
        "1";"11" **\dir{-};
        "6";"12" **\dir{-};
        \ar@{-} "5";"8" |<(0.4)\hole
    \endxy}
    \overset{\text{(L.~\ref{lem-useful})}}{=}
    \vcenter{\xy
        (-3, 9)="1";
        ( 3, 9)="2";
        (-3, 4)="3";
        ( 3, 5)*[o]=<9pt>[Fo]{\scriptstyle t}="4";
        (-3, 0)="5";
        ( 3, 1)="6";
        ( 0,-2)*[o]=<9pt>[Fo]{\scriptstyle t}="7";
        ( 0,-7)*[o]=<9pt>[Fo]{\scriptstyle t}="x";
        (-3,-9)="8";
        ( 0,-11)="9";
        ( 0,-14)*[o]=<5pt>[Fo]{~}="10";
        (-3,-17)="11";
        ( 3,-17)="12";
        "1";"11" **\dir{-};
        "2";"4" **\dir{-};
        "4";"12" **\dir{-};
        "3";"6" **\dir{-};
        "5";"7" **\dir{-};
        "8";"9" **\dir{-};
        "7";"x" **\dir{-};
        "x";"10" **\dir{-};
    \endxy}
    \overset{\text{(7)}}{=}~
    \vcenter{\xy
        (-3, 9)="1";
        ( 3, 9)="2";
        (-3, 4)="3";
        ( 3, 5)*[o]=<9pt>[Fo]{\scriptstyle t}="4";
        (-3, 0)="5";
        ( 3, 1)="6";
        ( 0,-2)*[o]=<9pt>[Fo]{\scriptstyle t}="7";
        (-3,-4)="8";
        ( 0,-6)="9";
        ( 0,-9)*[o]=<5pt>[Fo]{~}="10";
        (-3,-12)="11";
        ( 3,-12)="12";
        "1";"11" **\dir{-};
        "2";"4" **\dir{-};
        "4";"12" **\dir{-};
        "3";"6" **\dir{-};
        "5";"7" **\dir{-};
        "8";"9" **\dir{-};
        "7";"10" **\dir{-};
    \endxy}
    ~\overset{\text{(c,6)}}{=}~
    \vcenter{\xy
        (-2, 7)="1";
        ( 3, 7)="2";
        ( 3, 3)*[o]=<9pt>[Fo]{\scriptstyle t}="3";
        (-2, 2)="4";
        ( 3,-1)="5";
        ( 3,-6)="6";
        (-2,-6)="7";
        (-2,-1)="x";
        (0.5,-2.5)*[o]=<5pt>[Fo]{~}="y";
        "1";"7" **\dir{-};
        "2";"3" **\dir{-};
        "3";"6" **\dir{-};
        "4";"5" **\dir{-};
        "x";"y" **\dir{-};
    \endxy}
    ~\overset{\text{(c)}}{=}~
    \vcenter{\xy
        (-2, 7)="1";
        ( 2, 7)="2";
        ( 2, 3)*[o]=<9pt>[Fo]{\scriptstyle t}="3";
        (-2, 2)="4";
        ( 2,-1)="5";
        ( 2,-4)="6";
        (-2,-4)="7";
        "1";"7" **\dir{-};
        "2";"3" **\dir{-};
        "3";"6" **\dir{-};
        "4";"5" **\dir{-};
    \endxy}
    ~\overset{\text{(c)}}{=}~
    \vcenter{\xy
        (-2, 7)="1";
        ( 3, 7)="2";
        ( 3, 3)*[o]=<9pt>[Fo]{\scriptstyle t}="3";
        (-2, 2)="4";
        ( 3,-1)="5";
        ( 3,-3)="6";
        ( 1,-8)="x";
        ( 5,-5)*[o]=<5pt>[Fo]{~}="y";
        (-2,-8)="7";
        "1";"7" **\dir{-};
        "2";"3" **\dir{-};
        "3";"6" **\dir{-};
        "4";"5" **\dir{-};
        "6";"x" **\dir{-};
        "6";"y" **\dir{-};
    \endxy}
    ~\overset{\text{(4)}}{=}~
    \vcenter{\xy
        (-2, 7)="1";
        ( 3, 7)="2";
        ( 3, 3)*[o]=<9pt>[Fo]{\scriptstyle t}="3";
        (-2, 2)="4";
        ( 3,-1)="5";
        ( 3,-4)*[o]=<5pt>[Fo]{~}="6";
        (-2,-7)="7";
        ( 0.5,0.5)="x";
        ( 0.5,-7)="y";
        "1";"7" **\dir{-};
        "2";"3" **\dir{-};
        "3";"6" **\dir{-};
        "4";"5" **\dir{-};
        "x";"y" **\dir{-};
    \endxy}
    ~\overset{\text{(c)}}{=}~
    \vcenter{\xy
        (-2, 7)="1";
        ( 2, 7)="2";
        ( 2, 3)*[o]=<9pt>[Fo]{\scriptstyle t}="3";
        (-2, 2)="4";
        ( 2,-1)="5";
        ( 2,-4)*[o]=<5pt>[Fo]{~}="6";
        (-2,-8)="7";
        (-2,-4)="x";
        ( 2,-8)="y";
        "1";"7" **\dir{-};
        "2";"3" **\dir{-};
        "3";"6" **\dir{-};
        "4";"5" **\dir{-};
        "x";"y" **\dir{-};
    \endxy}~.
\]

\item $~\vcenter{\xy
        (-2,11)="1";
        (-2, 8)="2";
        ( 2, 5)*[o]=<9pt>[Fo]{\scriptstyle t}="3";
        (-2, 0)="4";
        ( 2, 0)="5";
        "1";"4" **\dir{-};
        "2";"3" **\dir{-};
        "3";"5" **\dir{-};
    \endxy}~:
    M \dra M \oxC C$ is an $A$-comodule morphism as:
\[
    \vcenter{\xy
        (-2,11)="1";
        (-2, 8)="2";
        ( 2, 6)*[o]=<9pt>[Fo]{\scriptstyle t}="3";
        ( 2, 1)*[o]=<9pt>[Fo]{\scriptstyle t}="4";
        (-2,-3)="5";
        ( 2,-2)="6";
        ( 5,-6)="7";
        (-2,-9)="8";
        ( 2,-9)="9";
        ( 5,-9)="10";
        "1";"8" **\dir{-};
        "2";"3" **\dir{-};
        "3";"4" **\dir{-};
        "5";"7" **\dir{-};
        "6";"7" **\dir{-};
        "10";"7" **\dir{-};
        \ar@{-} "4";"9" |<(0.55)\hole
    \endxy}
    ~\overset{\text{(7)}}{=}~
    \vcenter{\xy
        (-2, 7)="1";
        (-2, 4)="2";
        ( 2, 2)*[o]=<9pt>[Fo]{\scriptstyle t}="3";
        (-2,-3)="5";
        ( 2,-2)="6";
        ( 5,-6)="7";
        (-2,-9)="8";
        ( 2,-9)="9";
        ( 5,-9)="10";
        "1";"8" **\dir{-};
        "2";"3" **\dir{-};
        "5";"7" **\dir{-};
        "6";"7" **\dir{-};
        "10";"7" **\dir{-};
        \ar@{-} "3";"9" |<(0.55)\hole
    \endxy}
    \overset{\text{(L.~\ref{lem-useful})}}{=}
    \vcenter{\xy
        (-3, 9)="1";
        (-3, 6)="2";
        (-3, 2)="3";
        (-3,-2)="4";
        ( 3, 3)*[o]=<9pt>[Fo]{\scriptstyle t}="6";
        ( 3, -1)="7";
        ( 0,-3.5)*[o]=<9pt>[Fo]{\scriptstyle t}="9";
        (-3,-8)="5";
        ( 0,-8)="10";
        ( 3,-8)="8";
        "1";"5" **\dir{-};
        "6";"8" **\dir{-};
        "2";"6" **\dir{-};
        "3";"7" **\dir{-};
        "4";"9" **\dir{-};
        "10";"9" **\dir{-};
    \endxy}
    ~\overset{\text{(c)}}{=}~
    \vcenter{\xy
        (-3, 9)="1";
        (-3, 6)="2";
        (0,4.5)="3";
        (0,0.5)="x";
        (-3,-2)="4";
        ( 3, 3)*[o]=<9pt>[Fo]{\scriptstyle t}="6";
        ( 3, -1)="7";
        ( 0,-3.5)*[o]=<9pt>[Fo]{\scriptstyle t}="9";
        (-3,-8)="5";
        ( 0,-8)="10";
        ( 3,-8)="8";
        "1";"5" **\dir{-};
        "6";"8" **\dir{-};
        "2";"6" **\dir{-};
        "3";"x" **\dir{-};
        "x";"7" **\dir{-};
        "4";"9" **\dir{-};
        "10";"9" **\dir{-};
    \endxy}
    ~\overset{\text{(6)}}{=}~
    \vcenter{\xy
        (-3, 9)="1";
        (-3, 6)="2";
        (-3, 2)="3";
        (-3,-4)="4";
        ( 0, 0)*[o]=<9pt>[Fo]{\scriptstyle t}="5";
        ( 0,-4)="6";
        ( 3,-4)="7";
        "1";"4" **\dir{-};
        "3";"5" **\dir{-};
        "5";"6" **\dir{-};
        "7";"2" **\crv{(3,5)};
    \endxy}~.
\]

\item $~\vcenter{\xy
        (-2, 7)="1";
        ( 2, 7)="2";
        (-2, 3)*[o]=<9pt>[Fo]{\scriptstyle t}="3";
        ( 2, 4)="4";
        ( 4, 2)="5";
        (-2,-1)="6";
        (-2,-4)*[o]=<5pt>[Fo]{~}="7";
        ( 2,-7)="8";
        "1";"3" **\dir{-};
        "3";"7" **\dir{-};
        "4";"5" **\dir{-};
        "5";"6" **\dir{-};
        \ar@{-} "2";"8" |<(0.45){\hole}
    \endxy}~:
C \oxC M \dra M$ is an $A$-comodule morphism as:
\begin{align*}
    \vcenter{\xy
        (-3, 9)="1";
        ( 3, 9)="2";
        (-3, 5)*[o]=<9pt>[Fo]{\scriptstyle t}="3";
        (-3, 1)="4";
        ( 3, 3)="5";
        ( 6,-1)="6";
        (-3,-3)*[o]=<9pt>[Fo]{\scriptstyle t}="7";
        ( 2,-3)="8";
        ( 4,-5)="13";
        (-1,-7)="9";
        (-1,-10)*[o]=<5pt>[Fo]{~}="10";
        ( 2,-13)="11";
        ( 6,-13)="12";
        "1";"3" **\dir{-};
        "3";"7" **\dir{-};
        "7";"9" **\dir{-};
        "9";"10" **\dir{-};
        "2";"5" **\dir{-};
        "4";"6" **\dir{-};
        "5";"6" **\dir{-};
        "6";"12" **\dir{-};
        "8";"13" **\dir{-};
        "13";"9" **\dir{-};
        \ar@{-} "5";"8" |<(0.5){\hole}
        \ar@{-} "8";"11" |<(0.3){\hole}
    \endxy}
    &\overset{\text{(L.~\ref{lem-useful})}}{=}
    \vcenter{\xy
        (0,10)*{~}; (0,-20)*{~};
        (-3, 9)="1";
        ( 3, 9)="2";
        (-3, 5)*[o]=<9pt>[Fo]{\scriptstyle t}="3";
        ( 3, 6)="4";
        ( 6, 2)="5";
        ( 1,-1)="6";
        ( 3,-3)="7";
        (-2,-5)*[o]=<9pt>[Fo]{\scriptstyle s}="8";
        (-2,-10)*[o]=<9pt>[Fo]{\scriptstyle t}="9";
        ( 1,-9)="10";
        ( 3,-11)="11";
        (-2,-13)="12";
        (-2,-16)*[o]=<5pt>[Fo]{~}="13";
        ( 1,-19)="14";
        ( 6,-19)="15";
        "1";"3" **\dir{-};
        "3";"5" **\dir{-};
        "2";"4" **\dir{-};
        "4";"5" **\dir{-};
        "5";"15" **\dir{-};
        "6";"7" **\dir{-};
        "7";"8" **\dir{-};
        "8";"9" **\dir{-};
        "13";"9" **\dir{-};
        "10";"11" **\dir{-};
        "12";"11" **\dir{-};
        \ar@{-} "4";"6" |<(0.4){\hole}
        \ar@{-} "6";"14" |<(0.16){\hole} |<(0.6){\hole}
    \endxy}
    ~\overset{\text{(8)}}{=}~
    \vcenter{\xy
        (-3, 9)="1";
        ( 3, 9)="2";
        (-3, 5)*[o]=<9pt>[Fo]{\scriptstyle t}="3";
        ( 3, 6)="4";
        ( 6, 2)="5";
        ( 1,-1)="6";
        ( 3,-3)="7";
        (-2,-5)*[o]=<9pt>[Fo]{\scriptstyle s}="8";
        ( 1,-6)="10";
        ( 3,-8)="11";
        (-2,-10)="12";
        (-2,-13)*[o]=<5pt>[Fo]{~}="13";
        ( 1,-16)="14";
        ( 6,-16)="15";
        "1";"3" **\dir{-};
        "3";"5" **\dir{-};
        "2";"4" **\dir{-};
        "4";"5" **\dir{-};
        "5";"15" **\dir{-};
        "6";"7" **\dir{-};
        "7";"8" **\dir{-};
        "8";"13" **\dir{-};
        "10";"11" **\dir{-};
        "12";"11" **\dir{-};
        \ar@{-} "4";"6" |<(0.4){\hole}
        \ar@{-} "6";"14" |<(0.17){\hole} |<(0.53){\hole}
    \endxy}
    ~\overset{\text{(c,6)}}{=}~
    \vcenter{\xy
        (-2, 9)="1";
        ( 3, 9)="2";
        (-2, 5)*[o]=<9pt>[Fo]{\scriptstyle t}="3";
        ( 3, 6)="4";
        ( 6, 2)="5";
        ( 1, 0)="6";
        ( 3,-2)*[o]=<5pt>[Fo]{~}="7";
        ( 1,-5)="8";
        ( 6,-5)="9";
        "1";"3" **\dir{-};
        "3";"5" **\dir{-};
        "5";"9" **\dir{-};
        "2";"4" **\dir{-};
        "4";"5" **\dir{-};
        "6";"7" **\dir{-};
        "6";"8" **\dir{-};
        \ar@{-} "4";"6" |<(0.5){\hole}
    \endxy}
    \overset{\text{(c)}}{=}~
    \vcenter{\xy
        (-1, 9)="1";
        ( 3, 9)="2";
        (-1, 5)*[o]=<9pt>[Fo]{\scriptstyle t}="3";
        ( 3, 6)="4";
        ( 6, 1)="5";
        ( 3,-3)="8";
        ( 6,-3)="9";
        "1";"3" **\dir{-};
        "3";"5" **\dir{-};
        "5";"9" **\dir{-};
        "4";"5" **\dir{-};
        \ar@{-} "2";"8" |<(0.54){\hole}
    \endxy}
\\
    &\quad \overset{\text{(c)}}{=}~
    \vcenter{\xy
        (0,10)*{~}; (0,-5)*{~};
        (-1, 9)="1";
        ( 3, 9)="2";
        (-1, 5)*[o]=<9pt>[Fo]{\scriptstyle t}="3";
        ( 3, 6)="4";
        ( 6, 1)="5";
        ( 3,-7)="8";
        ( 6,-7)="9";
        ( 6,-2)="x";
        ( 8,-4)*[o]=<5pt>[Fo]{~}="y";
        "1";"3" **\dir{-};
        "3";"5" **\dir{-};
        "5";"9" **\dir{-};
        "4";"5" **\dir{-};
        "x";"y" **\dir{-};
        \ar@{-} "2";"8" |<(0.4){\hole}
    \endxy}
    ~\overset{\text{(4)}}{=}~
    \vcenter{\xy
        (-2, 9)="1";
        ( 2, 9)="2";
        (-2, 4)*[o]=<9pt>[Fo]{\scriptstyle t}="3";
        ( 2, 6)="4";
        ( 5, 4)="5";
        ( 8, 0)="6";
        ( 2,-6)="7";
        ( 5,-6)="8";
        ( 8,-3)*[o]=<5pt>[Fo]{~}="9";
        "1";"3" **\dir{-};
        "3";"6" **\dir{-};
        "6";"9" **\dir{-};
        "4";"5" **\dir{-};
        "5";"6" **\dir{-};
        \ar@{-} "2";"7" |<(0.45){\hole}
        \ar@{-} "5";"8" |<(0.3){\hole}
    \endxy}
    ~\overset{\text{(c)}}{=}~
    \vcenter{\xy
        (-3, 9)="1";
        ( 3, 9)="2";
        (-3, 5)*[o]=<9pt>[Fo]{\scriptstyle t}="3";
        ( 3, 6)="4";
        ( 6, 4)="5";
        ( 0, 2)="6";
        ( 0,-1)*[o]=<5pt>[Fo]{~}="7";
        ( 3,-2)="8";
        ( 1,-6)="9";
        ( 5,-6)="10";
        "1";"3" **\dir{-};
        "3";"6" **\dir{-};
        "6";"7" **\dir{-};
        "4";"5" **\dir{-};
        "5";"6" **\dir{-};
        "8";"9" **\dir{-};
        "8";"10" **\dir{-};
        \ar@{-} "2";"8" |<(0.6){\hole}
    \endxy}~.
\end{align*}

\item $~\vcenter{\xy
        ( 4, 7)="1";
        ( 4, 4)="2";
        ( 6, 2)="3";
        ( 0,-1)*[o]=<9pt>[Fo]{\scriptstyle s}="4";
        ( 0,-5)="5";
        ( 4,-5)="6";
        "2";"3" **\dir{-};
        "3";"4" **\dir{-};
        "4";"5" **\dir{-};
        \ar@{-} "1";"6" |<(0.5){\hole}
    \endxy}~:
M \dra C \oxC M$ is an $A$-comodule morphism as:
\[
    \vcenter{\xy
        ( 4, 7)="1";
        ( 4, 4)="2";
        ( 6, 2)="3";
        ( 0,-1)*[o]=<9pt>[Fo]{\scriptstyle s}="4";
        ( 0,-6)*[o]=<9pt>[Fo]{\scriptstyle t}="5";
        ( 0,-9)="6";
        ( 4,-7)="7";
        ( 7,-11)="8";
        ( 0,-14)="9";
        ( 4,-14)="10";
        ( 7,-14)="11";
        "2";"3" **\dir{-};
        "3";"4" **\dir{-};
        "4";"5" **\dir{-};
        "5";"9" **\dir{-};
        "6";"8" **\dir{-};
        "7";"8" **\dir{-};
        "8";"11" **\dir{-};
        \ar@{-} "1";"10" |<(0.3){\hole} |<(0.83){\hole}
    \endxy}
    ~\overset{\text{(L.~\ref{lem-useful})}}{=}~
    \vcenter{\xy
        (0,8)*{~}; (0,-19)*{~};
        ( 4, 7)="1";
        ( 4, 4)="2";
        ( 6, 2)="3";
        ( 0,-1)*[o]=<9pt>[Fo]{\scriptstyle s}="4";
        ( 0,-6)*[o]=<9pt>[Fo]{\scriptstyle t}="5";
        ( 4,-4)="7";
        ( 7,-8)="8";
        ( 4,-9)="w";
        ( 6,-11)="x";
        ( 0,-14)*[o]=<9pt>[Fo]{\scriptstyle s}="y";
        ( 0,-18)="9";
        ( 4,-18)="10";
        ( 7,-18)="11";
        "2";"3" **\dir{-};
        "3";"4" **\dir{-};
        "4";"5" **\dir{-};
        "5";"8" **\dir{-};
        "7";"8" **\dir{-};
        "8";"11" **\dir{-};
        "w";"x" **\dir{-};
        "y";"x" **\dir{-};
        "y";"9" **\dir{-};
        \ar@{-} "1";"10" |<(0.25){\hole} |<(0.56){\hole} |<(0.77){\hole}
    \endxy}
    ~\overset{\text{(8)}}{=}~
    \vcenter{\xy
        ( 4, 4)="1";
        ( 4, 1)="2";
        ( 6,-1)="3";
        ( 0,-4)*[o]=<9pt>[Fo]{\scriptstyle s}="4";
        ( 4,-4)="7";
        ( 7,-8)="8";
        ( 4,-9)="w";
        ( 6,-11)="x";
        ( 0,-14)*[o]=<9pt>[Fo]{\scriptstyle s}="y";
        ( 0,-18)="9";
        ( 4,-18)="10";
        ( 7,-18)="11";
        "2";"3" **\dir{-};
        "3";"4" **\dir{-};
        "4";"8" **\dir{-};
        "7";"8" **\dir{-};
        "8";"11" **\dir{-};
        "w";"x" **\dir{-};
        "y";"x" **\dir{-};
        "y";"9" **\dir{-};
        \ar@{-} "1";"10" |<(0.26){\hole} |<(0.48){\hole} |<(0.73){\hole}
    \endxy}
    ~\overset{\text{(c,6)}}{=}
    \vcenter{\xy
        ( 4,-2)="6";
        ( 4,-5)="7";
        ( 7,-8)="8";
        ( 4,-9)="w";
        ( 6,-11)="x";
        ( 0,-14)*[o]=<9pt>[Fo]{\scriptstyle s}="y";
        ( 0,-18)="9";
        ( 4,-18)="10";
        ( 7,-18)="11";
        "7";"8" **\dir{-};
        "8";"11" **\dir{-};
        "w";"x" **\dir{-};
        "y";"x" **\dir{-};
        "y";"9" **\dir{-};
        \ar@{-} "6";"10" |<(0.63){\hole}
    \endxy}~.
\]
\end{itemize}
Thus, $M \oxC C \cong M \cong C \oxC M$ in $\Q\V$.
\end{proof}

Thus, $\Comod(A) = (\Comod(A),\oxC,C)$ is a monoidal category.

%=========================================================================%
\subsection{The forgetful functor from $A$-comodules to $C$-bicomodules}
\label{sec-forget}
%=========================================================================%

There is a forgetful functor $U:\Comod(A) \ra \Bicomod(C)$ which assigns to
each $A$-comodule $M$ a $C$-bicomodule $UM$ which is $M$ itself with coaction
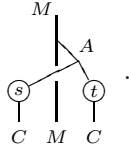
\[
    \vcenter{\xy
        (-2,8)*{\scriptstyle M};
        (4,3)*{\scriptstyle A};
        (-5,-9)*{\scriptstyle C};
        (0,-9)*{\scriptstyle M};
        (5,-9)*{\scriptstyle C};
        (0,7)="1";
        (0,4)="2";
        (3,1)="3";
        (-5,-3)*[o]=<9pt>[Fo]{\scriptstyle s}="4";
        (5,-3)*[o]=<9pt>[Fo]{\scriptstyle t}="5";
        (-5,-7)="6";
        ( 0,-7)="7";
        ( 5,-7)="8";
        "2";"3" **\dir{-};
        "3";"4" **\dir{-};
        "3";"5" **\dir{-};
        "4";"6" **\dir{-};
        "5";"8" **\dir{-};
        \ar@{-} |<(0.45){\hole} "1";"7" ;
    \endxy}\ \ .
\]
A morphism of $A$-comodules $f:M \ra N$ is automatically a morphism of the
underlying $C$-bicomodules $f:UM \ra UN$.

\begin{proposition}
The forgetful functor $U:\Comod(A) \ra \Bicomod(C)$ is strong monoidal.
\end{proposition}

\begin{proof}
We must establish the $C$-bicomodule isomorphisms
\[
    C \cong UC \qquad \text{and} \qquad UM \oxC UN \cong U(M \oxC N).
\]
The first is obvious. To establish the second isomorphism we observe that
the object $UM \oxC UN$ is $(M \oxC N,m)$ with coaction
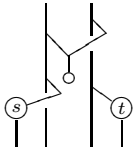
\[
    \vcenter{\xy
        (0,10)*{~}; (0,-11)*{~};
        (-3,9)="1";
        ( 3,9)="2";
        (-3,5)="3";
        ( 3,7)="4";
        ( 5,5)="5";
        ( 0,2)="6";
        ( 0,-1)*[o]=<5pt>[Fo]{~}="7";
        (-3,-1)="8";
        (-1,-3)="x";
        ( 3,-2)="9";
        (-7,-5)*[o]=<9pt>[Fo]{\scriptstyle s}="10";
        (7,-5)*[o]=<9pt>[Fo]{\scriptstyle t}="11";
        (-7,-10)="12";
        (-3,-10)="13";
        ( 3,-10)="14";
        ( 7,-10)="15";
        "3";"6" **\dir{-};
        "4";"5" **\dir{-};
        "5";"6" **\dir{-};
        "6";"7" **\dir{-};
        "8";"x" **\dir{-};
        "x";"10" **\dir{-};
        "9";"11" **\dir{-};
        "10";"12" **\dir{-};
        "11";"15" **\dir{-};
        \ar@{-} "1";"13" |<(0.68){\hole}
        \ar@{-} "2";"14" |<(0.28){\hole}
    \endxy}
\]
and $U(M \oxC N)$ is also $(M \oxC N,m)$ but with coaction
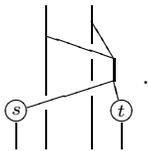
\[
    \vcenter{\xy
        (0,10)*{~}; (0,-11)*{~};
        (-3,9)="1";
        ( 3,9)="2";
        (-3,5)="3";
        ( 3,7)="4";
        ( 6,2)="5";
        ( 6,-1)="6";
        (-7,-5)*[o]=<9pt>[Fo]{\scriptstyle s}="7";
        (7,-5)*[o]=<9pt>[Fo]{\scriptstyle t}="8";
        (-7,-10)="9";
        (-3,-10)="10";
        ( 3,-10)="11";
        ( 7,-10)="12";
        "3";"5" **\dir{-};
        "4";"5" **\dir{-};
        "5";"6" **\dir{-};
        "6";"7" **\dir{-};
        "6";"8" **\dir{-};
        "7";"9" **\dir{-};
        "8";"12" **\dir{-};
        \ar@{-} "1";"10" |<(0.68){\hole}
        \ar@{-} "2";"11" |<(0.32){\hole} |<(0.58){\hole}
    \endxy}~.
\]
The following calculation shows that these two coactions are the same, and
hence the isomorphism $U(M \oxC N) \cong UM \oxC UN$.
\[
    \vcenter{\xy
        (-3,9)="1";
        ( 3,9)="2";
        (-3,5)="3";
        ( 3,7)="4";
        ( 5,5)="5";
        ( 0,2)="6";
        ( 0,-1)*[o]=<5pt>[Fo]{~}="7";
        (-3,-1)="8";
        (-1,-3)="x";
        ( 3,-2)="9";
        (-7,-5)*[o]=<9pt>[Fo]{\scriptstyle s}="10";
        (7,-5)*[o]=<9pt>[Fo]{\scriptstyle t}="11";
        (-7,-10)="12";
        (-3,-10)="13";
        ( 3,-10)="14";
        ( 7,-10)="15";
        "3";"6" **\dir{-};
        "4";"5" **\dir{-};
        "5";"6" **\dir{-};
        "6";"7" **\dir{-};
        "8";"x" **\dir{-};
        "x";"10" **\dir{-};
        "9";"11" **\dir{-};
        "10";"12" **\dir{-};
        "11";"15" **\dir{-};
        \ar@{-} "1";"13" |<(0.68){\hole}
        \ar@{-} "2";"14" |<(0.28){\hole}
    \endxy}
    \overset{\text{(2)}}{=}
    \vcenter{\xy
        (0,11)*{~}; (0,-11)*{~};
        (-3,10)="1";
        ( 6,10)="2";
        (-3,6)="3";
        ( 6,8)="4";
        ( 8,6)="5";
        ( 3,3.5)*[o]=<9pt>[Fo]{\scriptstyle s}="y";
        ( 0,2)="6";
        ( 0,-1)*[o]=<5pt>[Fo]{~}="7";
        (-3,-1)="8";
        (-1,-3)="x";
        ( 6,-2)="9";
        (-7,-5)*[o]=<9pt>[Fo]{\scriptstyle s}="10";
        (10,-5)*[o]=<9pt>[Fo]{\scriptstyle t}="11";
        (-7,-10)="12";
        (-3,-10)="13";
        ( 6,-10)="14";
        (10,-10)="15";
        "3";"6" **\dir{-};
        "4";"5" **\dir{-};
        "5";"y" **\dir{-};
        "y";"6" **\dir{-};
        "6";"7" **\dir{-};
        "8";"x" **\dir{-};
        "x";"10" **\dir{-};
        "9";"11" **\dir{-};
        "10";"12" **\dir{-};
        "11";"15" **\dir{-};
        \ar@{-} "1";"13" |<(0.68){\hole}
        \ar@{-} "2";"14" |<(0.25){\hole}
    \endxy}
    \overset{\text{(c,10)}}{=}
    \vcenter{\xy
        (-3,10)="1";
        ( 6,10)="2";
        (-3,6)="3";
        ( 6,8)="4";
        ( 8,6)="5";
        ( 3,3.5)*[o]=<9pt>[Fo]{\scriptstyle s}="y";
        ( 0,2)="6";
        ( 0,-1)*[o]=<5pt>[Fo]{~}="7";
        (-3,-1)="8";
        (-1,-3)="x";
        (-7,-5)*[o]=<9pt>[Fo]{\scriptstyle s}="10";
        (10, 0)*[o]=<9pt>[Fo]{\scriptstyle t}="11";
        (-7,-10)="12";
        (-3,-10)="13";
        ( 6,-10)="14";
        (10,-10)="15";
        "3";"6" **\dir{-};
        "4";"5" **\dir{-};
        "5";"y" **\dir{-};
        "y";"6" **\dir{-};
        "6";"7" **\dir{-};
        "8";"x" **\dir{-};
        "x";"10" **\dir{-};
        "10";"12" **\dir{-};
        "5";"11" **\dir{-};
        "11";"15" **\dir{-};
        \ar@{-} "1";"13" |<(0.68){\hole}
        \ar@{-} "2";"14" |<(0.25){\hole}
    \endxy}
    ~\overset{\text{(2)}}{=}
    \vcenter{\xy
        (-3,8)="1";
        ( 3,8)="2";
        (-3,4)="3";
        ( 3,6)="4";
        ( 5,4)="5";
        ( 0,2)="6";
        ( 0,-1)*[o]=<5pt>[Fo]{~}="7";
        (-3,-1)="8";
        (-1,-3)="x";
        (-7,-5)*[o]=<9pt>[Fo]{\scriptstyle s}="10";
        ( 7,-2)*[o]=<9pt>[Fo]{\scriptstyle t}="11";
        (-7,-10)="12";
        (-3,-10)="13";
        ( 3,-10)="14";
        ( 7,-10)="15";
        "3";"6" **\dir{-};
        "4";"5" **\dir{-};
        "5";"6" **\dir{-};
        "6";"7" **\dir{-};
        "8";"x" **\dir{-};
        "x";"10" **\dir{-};
        "10";"12" **\dir{-};
        "5";"11" **\dir{-};
        "11";"15" **\dir{-};
        \ar@{-} "1";"13" |<(0.65){\hole}
        \ar@{-} "2";"14" |<(0.27){\hole}
    \endxy}
\]
\[
    \overset{\text{(4)}}{=}
    \vcenter{\xy
        (-3,9)="1";
        ( 3,9)="2";
        (-3,5)="3";
        ( 3,7)="4";
        ( 7,3)="5";
        ( 7,0)="6";
        ( 5,-2)*[o]=<5pt>[Fo]{~}="7";
        ( 9,-3)*[o]=<9pt>[Fo]{\scriptstyle t}="8";
        (-3,0)="9";
        (-1,-2)="10";
        (-6,-4)*[o]=<9pt>[Fo]{\scriptstyle s}="11";
        (-6,-8)="12";
        (-3,-8)="13";
        ( 3,-8)="14";
        ( 9,-8)="15";
        "3";"5" **\dir{-};
        "4";"5" **\dir{-};
        "5";"6" **\dir{-};
        "6";"7" **\dir{-};
        "6";"8" **\dir{-};
        "8";"15" **\dir{-};
        "9";"10" **\dir{-};
        "10";"11" **\dir{-};
        "11";"12" **\dir{-};
        \ar@{-} "1";"13" |<(0.7){\hole}
        \ar@{-} "2";"14" |<(0.31){\hole}
    \endxy}
    ~\overset{\text{(c)}}{=}
    \vcenter{\xy
        (0,9)*{~}; (0,-13)*{~};
        (-3,8)="1";
        ( 3,8)="2";
        (-3,4)="3";
        ( 3,6)="4";
        ( 6,1)="5";
        (11,2)="6";
        (11,-2)="7";
        ( 6,-4)*[o]=<9pt>[Fo]{\scriptstyle s}="8";
        (11,-7)*[o]=<9pt>[Fo]{\scriptstyle t}="9";
        (-6,-8)="11";
        (-6,-12)="15";
        (-3,-12)="12";
        ( 3,-12)="13";
        (11,-12)="14";
        "3";"5" **\dir{-};
        "4";"6" **\dir{-};
        "6";"9" **\dir{-};
        "5";"7" **\dir{-};
        "5";"8" **\dir{-};
        "8";"11" **\dir{-};
        "11";"15" **\dir{-};
        "9";"14" **\dir{-};
        \ar@{-} "1";"12" |<(0.75){\hole}
        \ar@{-} "2";"13" |<(0.3){\hole} |<(0.65){\hole}
    \endxy}
    ~\overset{\text{(4)}}{=}
    \vcenter{\xy
        (-3,9)="1";
        ( 3,9)="2";
        (-3,4)="3";
        ( 3,6)="4";
        ( 6,2)="5";
        ( 6,-1)="6";
        (-6,-5)*[o]=<9pt>[Fo]{\scriptstyle s}="7";
        ( 8,-5)*[o]=<9pt>[Fo]{\scriptstyle t}="8";
        (-6,-9)="9";
        (-3,-9)="10";
        ( 3,-9)="11";
        ( 8,-9)="12";
        "3";"5" **\dir{-};
        "4";"5" **\dir{-};
        "5";"6" **\dir{-};
        "6";"7" **\dir{-};
        "6";"8" **\dir{-};
        "7";"9" **\dir{-};
        "8";"12" **\dir{-};
        \ar@{-} "1";"10" |<(0.72){\hole}
        \ar@{-} "2";"11" |<(0.36){\hole} |<(0.61){\hole}
    \endxy}
\]
\end{proof}

This may seem to be a strict equality, but as tensor products are really
only defined up to isomorphism we prefer ``strong''.

%=========================================================================%
\subsection{$\Comod_f(H)$ is left autonomous}
%=========================================================================%

Let $\V_f$ denote the subcategory of $\V$ consisting of the objects with a
left dual (since $\V$ is braided left duals are right duals). There is a
forgetful functor $U_l:\Comod(H) \ra \V$ defined as the composite of the two
forgetful functors $\Comod(H) \ra \Bicomod(C)$ and $\Bicomod(C) \ra \V$.
This composite $U_l:\Comod(H) \ra \V$ is sometimes called the \emph{long
forgetful functor}, as opposed to the \emph{short forgetful functor}
$U:\Comod(H) \ra \Bicomod(C)$.

Let us say an object $M \in \Comod(H)$ is \emph{dualizable} if $U_l M$ has
a left dual in $\V$, i.e., $U_l M \in \V_f$. Denote by $\Comod_f(H)$ the
subcategory of $\Comod(H)$ consisting of the dualizable objects.

The goal of this section is to prove the following proposition.

\begin{proposition}
If $H$ is a weak Hopf monoid then the category $\Comod_f(H)$ is left
autonomous (= left rigid = left compact).
\end{proposition}

Suppose $M \in \Comod_f(H)$ has a left dual $M^*$ in $\V$. Using the antipode
of $H$ a coaction on $M^*$ is defined as
\[
    \vcenter{\xy
        (-5, 7)*{\scriptstyle M^*};
        ( 3,-11)*{\scriptstyle M^*};
        ( 6,-11)*{\scriptstyle A};
        (-3, 6)="1";
        (-3, 0)="2";
        ( 0, 0)="3";
        ( 3, 0)="4";
        ( 6,-4)*[o]=<9pt>[Fo]{\scriptstyle \nu}="5";
        ( 3,-9)="6";
        ( 6,-9)="7";
        "1";"2" **\dir{-};
        "5";"7" **\dir{-};
        "3";"5" **\dir{-};
        "2";"3" **\crv{(-3,-3)&(0,-3)};
        "3";"4" **\crv{(0,3)&(3,3)};
        \ar@{-} "4";"6" |<(0.25)\hole
    \endxy}~.
\]
By (17) it is easy to see that this defines a comodule structure on $M^*$.
We claim that $M^*$ is the left dual of $M$ in $\Comod_f(H)$. Define
morphisms $e:M^* \oxC M \ra C$ and $n:C \ra M \oxC M^*$ via
\[
    e=~\vcenter{\xy
        (-3, 9)="1";
        ( 0, 9)="2";
        (-3, 6)="3";
        ( 0, 6)="4";
        ( 3, 3)*[o]=<9pt>[Fo]{\scriptstyle t}="5";
        ( 3,-2)="6";
        "1";"3" **\dir{-};
        "2";"4" **\dir{-};
        "3";"4" **\crv{(-3,3)&(0,3)};
        "4";"5" **\dir{-};
        "5";"6" **\dir{-};
    \endxy}
\qquad \text{and} \qquad
    n=~\vcenter{\xy
        (0,7)*{~}; (0,-11)*{~};
        ( 9, 6)="1";
        (-3, 0)="2";
        ( 0, 0)="3";
        ( 3,-3)*[o]=<9pt>[Fo]{\scriptstyle r}="4";
        ( 9,-1)*[o]=<9pt>[Fo]{\scriptstyle t}="5";
        ( 6,-4.5)="6";
        ( 6,-7.5)*[o]=<5pt>[Fo]{~}="7";
        (-3,-10)="8";
        ( 0,-10)="9";
        "2";"8" **\dir{-};
        "2";"4" **\dir{-};
        "2";"3" **\crv{(-3,3)&(0,3)};
        "4";"6" **\dir{-};
        "1";"5" **\dir{-};
        "5";"6" **\dir{-};
        "6";"7" **\dir{-};
        \ar@{-} "3";"9" |<(0.15)\hole
    \endxy}~.
\]

\begin{proposition}\label{prop-leftdual}
Suppose $M \in \Comod_f(H)$ with underlying left dual $M^*$. Then $M^*$
with evaluation and coevaluation morphisms $e$ and $n$ respectively is the
left dual of $M$ in $\Comod_f(H)$. That is, $\Comod_f(H)$ is left autonomous.
\end{proposition}

\begin{proof}
Let $M$, $M^*$, $e$, and $n$ be as above. We will first show that $e$ and
$n$ are comodule morphisms, and secondly that they satisfy the triangle
identities.

The following calculation shows that $e$ is a comodule morphism.
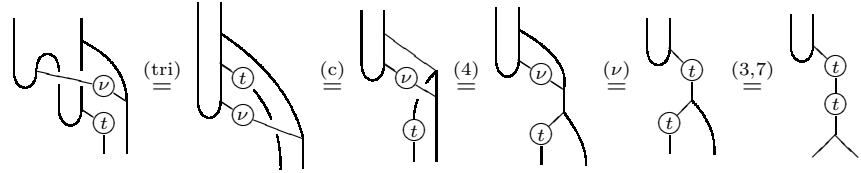
\[
    \vcenter{\xy
        (-3, 7)="1";
        ( 6, 7)="2";
        (-3, 0)="3";
        ( 0, 0)="4";
        ( 3, 0)="5";
        ( 6, 4)="6";
        ( 9,-2)*[o]=<9pt>[Fo]{\scriptstyle \nu}="7";
        (12,-4)="8";
        ( 3,-5)="9";
        ( 6,-5)="10";
        ( 9,-7)*[o]=<9pt>[Fo]{\scriptstyle t}="11";
        ( 9,-11)="12";
        (12,-11)="13";
        "1";"3" **\dir{-};
        "4";"7" **\dir{-};
        "7";"8" **\dir{-};
        "8";"13" **\dir{-};
        "10";"11" **\dir{-};
        "11";"12" **\dir{-};
        "8";"6" **\crv{(12,1)};
        "3";"4" **\crv{(-3,-3)&(0,-3)};
        "4";"5" **\crv{(0,3)&(3,3)};
        "9";"10" **\crv{(3,-8)&(6,-8)};
        \ar@{-} "5";"9" |<(0.15)\hole
        \ar@{-} "2";"10" |<(0.67)\hole
    \endxy}
    ~\overset{\text{(tri)}}{=}~
    \vcenter{\xy
        (-3,10)="1";
        ( 0,10)="2";
        ( 0, 6)="3";
        ( 0, 2)="4";
        (-3,-3)="5";
        ( 0,-3)="6";
        ( 3, 0)*[o]=<9pt>[Fo]{\scriptstyle t}="7";
        ( 3,-5)*[o]=<9pt>[Fo]{\scriptstyle \nu}="8";
        (11,-8)="9";
        ( 8,-12)="10";
        (11,-12)="11";
        "1";"5" **\dir{-};
        "2";"6" **\dir{-};
        "4";"7" **\dir{-};
        "6";"8" **\dir{-};
        "8";"9" **\dir{-};
        "11";"9" **\dir{-};
        "5";"6" **\crv{(-3,-6)&(0,-6)};
        "9";"3" **\crv{(9,0)};
        \ar@{-}@/_4pt/ "10";"7" |<(0.42)\hole
    \endxy}
    ~\overset{\text{(c)}}{=}~
    \vcenter{\xy
        (-3,9)="1";
        ( 0,9)="2";
        ( 0,6)="3";
        (-3,2)="4";
        ( 0,2)="5";
        ( 3,0)*[o]=<9pt>[Fo]{\scriptstyle \nu}="6";
        ( 7,1)="7";
        ( 7,-2.5)="8";
        ( 4,-7)*[o]=<9pt>[Fo]{\scriptstyle t}="9";
        ( 4,-11)="10";
        ( 7,-11)="11";
        "1";"4" **\dir{-};
        "2";"5" **\dir{-};
        "4";"5" **\crv{(-3,-1)&(0,-1)};
        "3";"7" **\dir{-};
        "5";"6" **\dir{-};
        "6";"8" **\dir{-};
        "11";"7" **\dir{-};
        "9";"10" **\dir{-};
        \ar@{-}@/^3pt/ "9";"7" |<(0.53)\hole
    \endxy}
    ~\overset{\text{(4)}}{=}~
    \vcenter{\xy
        (-3,9)="1";
        ( 0,9)="2";
        ( 0,6)="3";
        (-3,2)="4";
        ( 0,2)="5";
        ( 3,0)*[o]=<9pt>[Fo]{\scriptstyle \nu}="6";
        ( 6,-2)="7";
        ( 6,-5)="8";
        ( 3,-8)*[o]=<9pt>[Fo]{\scriptstyle t}="9";
        ( 3,-12)="10";
        ( 9,-12)="11";
        "1";"4" **\dir{-};
        "2";"5" **\dir{-};
        "4";"5" **\crv{(-3,-1)&(0,-1)};
        "7";"3" **\crv{(7,2)};
        "5";"6" **\dir{-};
        "6";"7" **\dir{-};
        "8";"7" **\dir{-};
        "8";"9" **\dir{-};
        "10";"9" **\dir{-};
        "11";"8" **\crv{(9,-9)};
    \endxy}
    ~\overset{\text{($\nu$)}}{=}~
    \vcenter{\xy
        (-3,6)="1";
        ( 0,6)="2";
        (-3,2)="3";
        ( 0,2)="4";
        ( 3,-1)*[o]=<9pt>[Fo]{\scriptstyle t}="5";
        ( 3,-5)="6";
        ( 0,-8)*[o]=<9pt>[Fo]{\scriptstyle t}="7";
        ( 0,-12)="8";
        ( 6,-12)="9";
        "1";"3" **\dir{-};
        "2";"4" **\dir{-};
        "3";"4" **\crv{(-3,-1)&(0,-1)};
        "4";"5" **\dir{-};
        "5";"6" **\dir{-};
        "6";"7" **\dir{-};
        "8";"7" **\dir{-};
        "9";"6" **\crv{(6,-8)};
    \endxy}
    ~\overset{\text{(3,7)}}{=}~
    \vcenter{\xy
        (-3,6)="1";
        ( 0,6)="2";
        (-3,2)="3";
        ( 0,2)="4";
        ( 3,-1)*[o]=<9pt>[Fo]{\scriptstyle t}="5";
        ( 3,-6)*[o]=<9pt>[Fo]{\scriptstyle t}="6";
        ( 3,-10)="7";
        ( 0,-13)="8";
        ( 6,-13)="9";
        "1";"3" **\dir{-};
        "2";"4" **\dir{-};
        "3";"4" **\crv{(-3,-1)&(0,-1)};
        "4";"5" **\dir{-};
        "5";"6" **\dir{-};
        "6";"7" **\dir{-};
        "8";"7" **\dir{-};
        "9";"7" **\dir{-};
    \endxy}
\]

To show that $n$ is a comodule morphism we must establish the equality
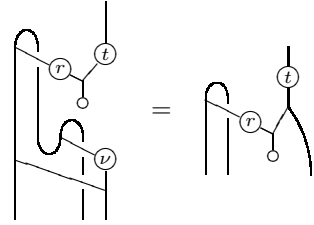
\[
    \vcenter{\xy
        (0,7)*{~}; (0,-24)*{~};
        ( 9, 6)="1";
        (-3, 0)="2";
        ( 0, 0)="3";
        ( 3,-3)*[o]=<9pt>[Fo]{\scriptstyle r}="4";
        ( 9,-1)*[o]=<9pt>[Fo]{\scriptstyle t}="5";
        ( 6,-4.5)="6";
        ( 6,-7.5)*[o]=<5pt>[Fo]{~}="7";
        ( 0,-12)="8";
        ( 3,-12)="9";
        ( 6,-12)="10";
        (-3,-15)="11";
        ( 9,-15)*[o]=<9pt>[Fo]{\scriptstyle \nu}="12";
        ( 9,-19)="13";
        (-3,-23)="14";
        ( 6,-23)="15";
        ( 9,-23)="16";
        "2";"14" **\dir{-};
        "2";"4" **\dir{-};
        "2";"3" **\crv{(-3,3)&(0,3)};
        "4";"6" **\dir{-};
        "1";"5" **\dir{-};
        "5";"6" **\dir{-};
        "6";"7" **\dir{-};
        "8";"9" **\crv{(0,-15)&(3,-15)};
        "9";"10" **\crv{(3,-9)&(6,-9)};
        "9";"12" **\dir{-};
        "11";"13" **\dir{-};
        "12";"16" **\dir{-};
        \ar@{-} "3";"8" |<(0.14)\hole
        \ar@{-} "10";"15" |<(0.15)\hole |<(0.54)\hole
    \endxy}
    \quad = \quad
    \vcenter{\xy
        ( 8, 7)="1";
        (-3, 0)="2";
        ( 0, 0)="3";
        ( 3,-3)*[o]=<9pt>[Fo]{\scriptstyle r}="4";
        ( 8, 3)*[o]=<9pt>[Fo]{\scriptstyle t}="x";
        ( 8,-1)="5";
        ( 6,-4.5)="6";
        ( 6,-7.5)*[o]=<5pt>[Fo]{~}="7";
        (-3,-10)="8";
        ( 0,-10)="9";
        (11,-10)="y";
        "2";"8" **\dir{-};
        "2";"4" **\dir{-};
        "2";"3" **\crv{(-3,3)&(0,3)};
        "4";"6" **\dir{-};
        "1";"x" **\dir{-};
        "x";"5" **\dir{-};
        "5";"6" **\dir{-};
        "6";"7" **\dir{-};
        "y";"5" **\crv{(11,-5)};
        \ar@{-} "3";"9" |<(0.15)\hole
    \endxy}
\]
which is proved by the following calculation.
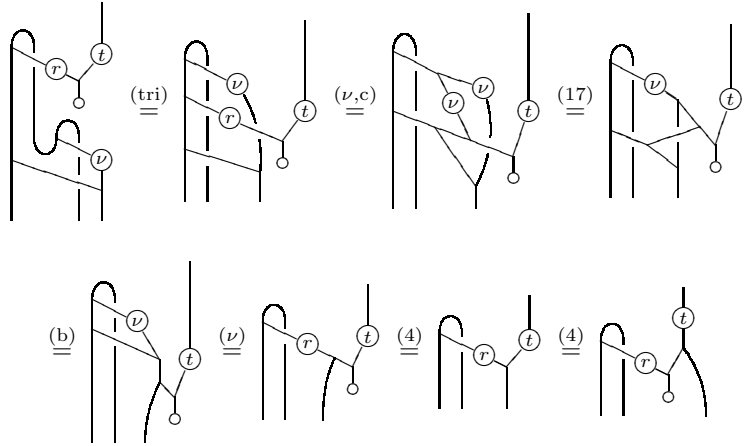
\[
    \vcenter{\xy
        ( 9, 6)="1";
        (-3, 0)="2";
        ( 0, 0)="3";
        ( 3,-3)*[o]=<9pt>[Fo]{\scriptstyle r}="4";
        ( 9,-1)*[o]=<9pt>[Fo]{\scriptstyle t}="5";
        ( 6,-4.5)="6";
        ( 6,-7.5)*[o]=<5pt>[Fo]{~}="7";
        ( 0,-12)="8";
        ( 3,-12)="9";
        ( 6,-12)="10";
        (-3,-15)="11";
        ( 9,-15)*[o]=<9pt>[Fo]{\scriptstyle \nu}="12";
        ( 9,-19)="13";
        (-3,-23)="14";
        ( 6,-23)="15";
        ( 9,-23)="16";
        "2";"14" **\dir{-};
        "2";"4" **\dir{-};
        "2";"3" **\crv{(-3,3)&(0,3)};
        "4";"6" **\dir{-};
        "1";"5" **\dir{-};
        "5";"6" **\dir{-};
        "6";"7" **\dir{-};
        "8";"9" **\crv{(0,-15)&(3,-15)};
        "9";"10" **\crv{(3,-9)&(6,-9)};
        "9";"12" **\dir{-};
        "11";"13" **\dir{-};
        "12";"16" **\dir{-};
        \ar@{-} "3";"8" |<(0.14)\hole
        \ar@{-} "10";"15" |<(0.15)\hole |<(0.54)\hole
    \endxy}
    ~\overset{\text{(tri)}}{=}~
    \vcenter{\xy
        (13, 10)="1";
        (-3, 5)="2";
        ( 0, 5)="3";
        ( 4, 1.5)*[o]=<9pt>[Fo]{\scriptstyle \nu}="4";
        (-3, 0)="5";
        ( 3,-3)*[o]=<9pt>[Fo]{\scriptstyle r}="6";
        (13,-2)*[o]=<9pt>[Fo]{\scriptstyle t}="7";
        (10,-6)="8";
        (-3,-7)="9";
        (10,-9)*[o]=<5pt>[Fo]{~}="10";
        ( 7,-10)="11";
        (-3,-14)="12";
        ( 0,-14)="13";
        ( 7,-14)="14";
        "2";"12" **\dir{-};
        "2";"4" **\dir{-};
        "2";"3" **\crv{(-3,8)&(0,8)};
        "5";"6" **\dir{-};
        "6";"8" **\dir{-};
        "8";"10" **\dir{-};
        "7";"1" **\dir{-};
        "7";"8" **\dir{-};
        "9";"11" **\dir{-};
        "14";"11" **\dir{-};
        \ar@{-} "3";"13" |<(0.1)\hole |<(0.35)\hole |<(0.7)\hole
        \ar@{-}@/^3pt/ "4";"11" |<(0.5)\hole
    \endxy}
    ~\overset{\text{($\nu$,c)}}{=}~
    \vcenter{\xy
        (15, 10)="1";
        (-3, 5)="2";
        ( 0, 5)="3";
        ( 3, 2)="4";
        ( 5,-2)*[o]=<9pt>[Fo]{\scriptstyle \nu}="5";
        ( 9, 0)*[o]=<9pt>[Fo]{\scriptstyle \nu}="6";
        (-3,-3)="7";
        (15,-3)*[o]=<9pt>[Fo]{\scriptstyle t}="10";
        (13,-9)="11";
        ( 8,-13)="12";
        (13,-12)*[o]=<5pt>[Fo]{~}="13";
        (-3,-16)="14";
        ( 0,-16)="15";
        ( 8,-16)="16";
        "2";"14" **\dir{-};
        "2";"3" **\crv{(-3,8)&(0,8)};
        "2";"4" **\dir{-};
        "5";"4" **\dir{-};
        "6";"4" **\dir{-};
        "7";"11" **\dir{-}?(0.35)="8" ?(0.65)="9";
        "5";"9" **\dir{-};
        "11";"13" **\dir{-};
        "1";"10" **\dir{-};
        "10";"11" **\dir{-};
        "12";"8" **\dir{-};
        "12";"16" **\dir{-};
        \ar@{-} "3";"15" |<(0.08)\hole |<(0.43)\hole
        \ar@{-}@/^4pt/ "6";"12" |<(0.52)\hole
    \endxy}
    ~\overset{\text{(17)}}{=}~
    \vcenter{\xy
        (13, 10)="1";
        (-3, 5)="2";
        ( 0, 5)="3";
        ( 3, 2)*[o]=<9pt>[Fo]{\scriptstyle \nu}="4";
        ( 6, 0)="5";
        (-3,-4)="6";
        ( 2,-6)="7";
        (13, 0)*[o]=<9pt>[Fo]{\scriptstyle t}="9";
        (11,-6)="10";
        ( 6,-9)="11";
        (11,-9)*[o]=<5pt>[Fo]{~}="12";
        (-3,-13)="13";
        ( 0,-13)="14";
        ( 6,-13)="15";
        "2";"13" **\dir{-};
        "2";"3" **\crv{(-3,8)&(0,8)};
        "2";"4" **\dir{-};
        "5";"4" **\dir{-};
        "5";"10" **\dir{-}?(0.6)="8";
        "1";"9" **\dir{-};
        "9";"10" **\dir{-};
        "10";"12" **\dir{-};
        "6";"7" **\dir{-};
        "7";"8" **\dir{-};
        "11";"7" **\dir{-};
        \ar@{-} "3";"14" |<(0.08)\hole |<(0.57)\hole
        \ar@{-} "5";"15" |<(0.35)\hole
    \endxy}
\]
\[
    \overset{\text{(b)}}{=}~
    \vcenter{\xy
        (10, 5)="1";
        (-3, 0)="2";
        ( 0, 0)="3";
        ( 3,-3)*[o]=<9pt>[Fo]{\scriptstyle \nu}="4";
        (-3,-4)="5";
        ( 6,-8)="6";
        (10,-8)*[o]=<9pt>[Fo]{\scriptstyle t}="7";
        ( 6,-11)="8";
        ( 8,-13)="9";
        ( 8,-16)*[o]=<5pt>[Fo]{~}="10";
        (-3,-19)="11";
        ( 0,-19)="12";
        ( 4,-19)="13";
        "2";"11" **\dir{-};
        "2";"4" **\dir{-};
        "2";"3" **\crv{(-3,3)&(0,3)};
        "4";"6" **\dir{-};
        "5";"6" **\dir{-};
        "8";"6" **\dir{-};
        "8";"9" **\dir{-};
        "10";"9" **\dir{-};
        "1";"7" **\dir{-};
        "9";"7" **\dir{-};
        "13";"8" **\crv{(4,-15)};
        \ar@{-} "3";"12" |<(0.08)\hole |<(0.28)\hole
    \endxy}
    ~\overset{\text{($\nu$)}}{=}~
    \vcenter{\xy
        (11, 5)="1";
        (-3, 0)="2";
        ( 0, 0)="3";
        ( 3,-3)*[o]=<9pt>[Fo]{\scriptstyle r}="4";
        (11,-2)*[o]=<9pt>[Fo]{\scriptstyle t}="5";
        ( 9,-6)="7";
        ( 9,-9)*[o]=<5pt>[Fo]{~}="x";
        (-3,-13)="8";
        ( 0,-13)="9";
        ( 5,-13)="10";
        "2";"8" **\dir{-};
        "2";"4" **\dir{-};
        "2";"3" **\crv{(-3,3)&(0,3)};
        "4";"7" **\dir{-}?(0.6)="6";
        "x";"7" **\dir{-};
        "1";"5" **\dir{-};
        "5";"7" **\dir{-};
        "10";"6" **\crv{(5,-8)};
        \ar@{-} "3";"9" |<(0.15)\hole
    \endxy}
    ~\overset{\text{(4)}}{=}~
    \vcenter{\xy
        ( 9, 5)="1";
        (-3, 0)="2";
        ( 0, 0)="3";
        ( 3,-3)*[o]=<9pt>[Fo]{\scriptstyle r}="4";
        ( 9,-1)*[o]=<9pt>[Fo]{\scriptstyle t}="5";
        ( 6,-4.5)="6";
        ( 6,-10)="7";
        (-3,-10)="8";
        ( 0,-10)="9";
        "2";"8" **\dir{-};
        "2";"4" **\dir{-};
        "2";"3" **\crv{(-3,3)&(0,3)};
        "4";"6" **\dir{-};
        "1";"5" **\dir{-};
        "5";"6" **\dir{-};
        "6";"7" **\dir{-};
        \ar@{-} "3";"9" |<(0.15)\hole
    \endxy}
    ~\overset{\text{(4)}}{=}~
    \vcenter{\xy
        ( 8, 7)="1";
        (-3, 0)="2";
        ( 0, 0)="3";
        ( 3,-3)*[o]=<9pt>[Fo]{\scriptstyle r}="4";
        ( 8, 3)*[o]=<9pt>[Fo]{\scriptstyle t}="x";
        ( 8,-1)="5";
        ( 6,-4.5)="6";
        ( 6,-7.5)*[o]=<5pt>[Fo]{~}="7";
        (-3,-10)="8";
        ( 0,-10)="9";
        (11,-10)="y";
        "2";"8" **\dir{-};
        "2";"4" **\dir{-};
        "2";"3" **\crv{(-3,3)&(0,3)};
        "4";"6" **\dir{-};
        "1";"x" **\dir{-};
        "x";"5" **\dir{-};
        "5";"6" **\dir{-};
        "6";"7" **\dir{-};
        "y";"5" **\crv{(11,-5)};
        \ar@{-} "3";"9" |<(0.15)\hole
    \endxy}
\]

It remains to show that $e$ and $n$ satisfy the triangle identities, i.e.,
that the following composites are the identity:

\begin{enumerate}[{\upshape (i)}]
\item $\xygraph{{M \cong C \oxC M}
            :[r(2.9)] {M \oxC M^* \oxC M} ^-{n \ox 1}
            :[r(2.9)] {M \oxC C \cong M} ^-{1 \ox e}}$

\item $\xygraph{{M^* \cong M^* \oxC C}
            :[r(3.1)] {M^* \oxC M \oxC M^*} ^-{1 \ox n}
            :[r(3)] {C \oxC M^* \cong M^*} ^-{e \ox 1}}$.
\end{enumerate}
Recall that $M \cong M \oxC C$ and $M \cong C \oxC M$ via
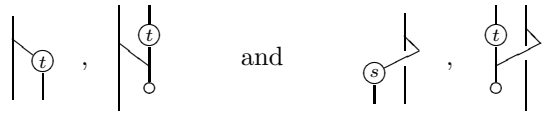
\[
    \vcenter{\xy
        (-2,11)="1";
        (-2, 8)="2";
        ( 2, 5)*[o]=<9pt>[Fo]{\scriptstyle t}="3";
        (-2, 0)="4";
        ( 2, 0)="5";
        "1";"4" **\dir{-};
        "2";"3" **\dir{-};
        "3";"5" **\dir{-};
    \endxy}
    \quad,\quad
    \vcenter{\xy
        (0,8)*{~}; (0,-8)*{~};
        (-2, 7)="1";
        ( 2, 7)="2";
        ( 2, 3)*[o]=<9pt>[Fo]{\scriptstyle t}="3";
        (-2, 2)="4";
        ( 2,-1)="5";
        ( 2,-4)*[o]=<5pt>[Fo]{~}="6";
        (-2,-7)="7";
        "1";"7" **\dir{-};
        "2";"3" **\dir{-};
        "3";"6" **\dir{-};
        "4";"5" **\dir{-};
    \endxy}
\quad\qquad \text{and} \quad\qquad
    \vcenter{\xy
        ( 4, 7)="1";
        ( 4, 4)="2";
        ( 6, 2)="3";
        ( 0,-1)*[o]=<9pt>[Fo]{\scriptstyle s}="4";
        ( 0,-5)="5";
        ( 4,-5)="6";
        "2";"3" **\dir{-};
        "3";"4" **\dir{-};
        "4";"5" **\dir{-};
        \ar@{-} "1";"6" |<(0.5){\hole}
    \endxy}
    \quad,\quad
    \vcenter{\xy
        (0,8)*{~}; (0,-8)*{~};
        (-2, 7)="1";
        ( 2, 7)="2";
        (-2, 3)*[o]=<9pt>[Fo]{\scriptstyle t}="3";
        ( 2, 4)="4";
        ( 4, 2)="5";
        (-2,-1)="6";
        (-2,-4)*[o]=<5pt>[Fo]{~}="7";
        ( 2,-7)="8";
        "1";"3" **\dir{-};
        "3";"7" **\dir{-};
        "4";"5" **\dir{-};
        "5";"6" **\dir{-};
        \ar@{-} "2";"8" |<(0.45){\hole}
    \endxy}
\]
respectively.

The following calculation proves (i):
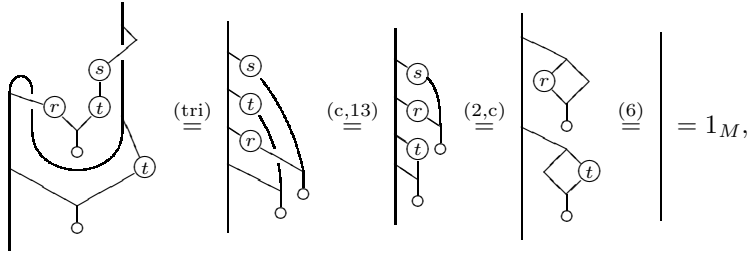
\[
    \vcenter{\xy
        (0,16)*{~}; (0,-19)*{~};
        (15,15)="1";
        (15,12)="2";
        (17,10)="x";
        (12, 6)*[o]=<9pt>[Fo]{\scriptstyle s}="3";
        ( 0, 3)="4";
        ( 3, 3)="5";
        ( 6, 1)*[o]=<9pt>[Fo]{\scriptstyle r}="6";
        (12, 1)*[o]=<9pt>[Fo]{\scriptstyle t}="7";
        ( 9,-2)="8";
        ( 9,-5)*[o]=<5pt>[Fo]{~}="9";
        ( 3,-2)="10";
        (15, 0)="y";
        (15,-2)="11";
        (18,-7)*[o]=<9pt>[Fo]{\scriptstyle t}="12";
        ( 0,-7)="13";
        ( 9,-12)="14";
        ( 9,-15)*[o]=<5pt>[Fo]{~}="15";
        ( 0,-18)="16";
        "4";"5" **\crv{(0,6)&(3,6)};
        "4";"16" **\dir{-};
        "4";"6" **\dir{-};
        "8";"6" **\dir{-};
        "8";"9" **\dir{-};
        "8";"7" **\dir{-};
        "3";"7" **\dir{-};
        "13";"14" **\dir{-};
        "15";"14" **\dir{-};
        "12";"14" **\dir{-};
        "12";"y" **\dir{-};
        "2";"x" **\dir{-};
        "3";"x" **\dir{-};
        "10";"11" **\crv{(3,-9)&(15,-9)};
        \ar@{-} "1";"11" |<(0.39)\hole
        \ar@{-} "5";"10" |<(0.23)\hole
    \endxy}
    ~\overset{\text{(tri)}}{=}~
    \vcenter{\xy
        ( 0,14)="1";
        ( 0,10)="2";
        ( 0, 5)="3";
        ( 0, 0)="4";
        ( 0,-5)="5";
        ( 0,-14)="6";
        ( 3, 8)*[o]=<9pt>[Fo]{\scriptstyle s}="7";
        ( 3, 3)*[o]=<9pt>[Fo]{\scriptstyle t}="8";
        ( 3,-2)*[o]=<9pt>[Fo]{\scriptstyle r}="9";
        (10,-6)="10";
        (10,-9)*[o]=<5pt>[Fo]{~}="11";
        ( 7,-8.5)="12";
        ( 7,-11.5)*[o]=<5pt>[Fo]{~}="13";
        "1";"6" **\dir{-};
        "2";"7" **\dir{-};
        "3";"8" **\dir{-};
        "4";"9" **\dir{-};
        "9";"10" **\dir{-};
        "10";"11" **\dir{-};
        "5";"12" **\dir{-};
        "13";"12" **\dir{-};
        \ar@{-}@/_3pt/ "12";"8" |<(0.37)\hole
        \ar@{-}@/_3pt/ "10";"7"
    \endxy}
    ~\overset{\text{(c,13)}}{=}~
    \vcenter{\xy
        ( 0,14)="1";
        ( 0,10)="2";
        ( 0, 5)="3";
        ( 0, 0)="4";
        ( 0,-4)="5";
        ( 0,-12)="6";
        ( 3, 8)*[o]=<9pt>[Fo]{\scriptstyle s}="7";
        ( 3, 3)*[o]=<9pt>[Fo]{\scriptstyle r}="8";
        ( 3,-2)*[o]=<9pt>[Fo]{\scriptstyle t}="9";
        ( 3,-6)="10";
        ( 3,-9)*[o]=<5pt>[Fo]{~}="11";
        ( 6, 1)="12";
        ( 6,-2)*[o]=<5pt>[Fo]{~}="13";
        "1";"6" **\dir{-};
        "2";"7" **\dir{-};
        "3";"8" **\dir{-};
        "4";"9" **\dir{-};
        "5";"10" **\dir{-};
        "9";"11" **\dir{-};
        "8";"12" **\dir{-};
        "13";"12" **\dir{-};
        "12";(4.2,7.5) **\crv{(6,6)};
    \endxy}
    ~\overset{\text{(2,c)}}{=}~
    \vcenter{\xy
        ( 0,14)="1";
        ( 0,11)="2";
        ( 0,-1)="3";
        ( 0,-16)="4";
        ( 6, 8)="5";
        ( 3, 5)*[o]=<9pt>[Fo]{\scriptstyle r}="6";
        ( 9, 5)="7";
        ( 6, 2)="8";
        ( 6,-1)*[o]=<5pt>[Fo]{~}="9";
        ( 6,-4)="10";
        ( 3,-7)="11";
        ( 9,-7)*[o]=<9pt>[Fo]{\scriptstyle t}="12";
        ( 6,-10)="13";
        ( 6,-13)*[o]=<5pt>[Fo]{~}="14";
        "1";"4" **\dir{-};
        "2";"5" **\dir{-};
        "3";"10" **\dir{-};
        "5";"6" **\dir{-};
        "5";"7" **\dir{-};
        "8";"6" **\dir{-};
        "8";"7" **\dir{-};
        "8";"9" **\dir{-};
        "10";"11" **\dir{-};
        "10";"12" **\dir{-};
        "13";"11" **\dir{-};
        "13";"12" **\dir{-};
        "13";"14" **\dir{-};
    \endxy}
    ~\overset{\text{(6)}}{=}~
    \vcenter{\xy
        ( 0,25); ( 0,0); **\dir{-};
    \endxy}
    ~= 1_{M},
\]
and (ii) is given by:
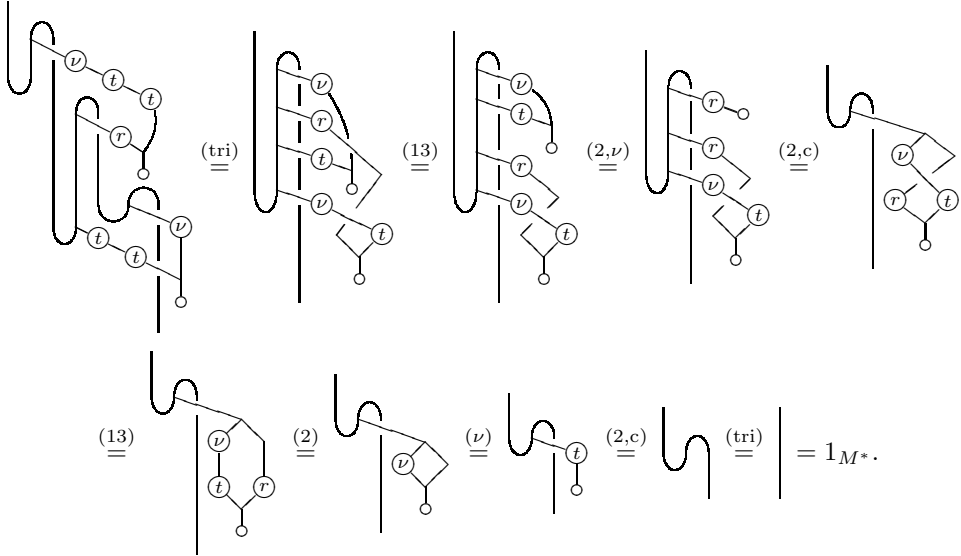
\[
    \vcenter{\xy
        (-3,30)="1"; ( 0,25)="2"; ( 3,25)="3";
        ( 6,22)*[o]=<9pt>[Fo]{\scriptstyle \nu}="4";
        (11,19.5)*[o]=<9pt>[Fo]{\scriptstyle t}="5";
        (16,17)*[o]=<9pt>[Fo]{\scriptstyle t}="6";
        (-3,20)="7"; ( 0,20)="8"; ( 6,15)="9"; ( 9,15)="10";
        (12,12)*[o]=<9pt>[Fo]{\scriptstyle r}="11";
        (15,10)="12";
        (15,7)*[o]=<5pt>[Fo]{~}="13";
        (9,3)="14"; (13,3)="15"; (17,3)="16";
        (20,0)*[o]=<9pt>[Fo]{\scriptstyle \nu}="17";
        ( 3,0)="18"; ( 6,0)="19";
        ( 9,-1.5)*[o]=<9pt>[Fo]{\scriptstyle t}="20";
        (14,-4)*[o]=<9pt>[Fo]{\scriptstyle t}="21";
        (20,-7)="22";
        (20,-10)*[o]=<5pt>[Fo]{~}="23";
        (17,-14)="24";
        "2";"3" **\crv{(0,28)&(3,28)};
        "7";"8" **\crv{(-3,17)&(0,17)};
        "9";"10" **\crv{(6,18)&(9,18)};
        "14";"15" **\crv{(9,0)&(13,0)};
        "15";"16" **\crv{(13,6)&(17,6)};
        "18";"19" **\crv{(3,-3)&(6,-3)};
        "1";"7" **\dir{-}; "2";"8" **\dir{-}; "9";"11" **\dir{-};
        "12";"11" **\dir{-}; "12";"13" **\dir{-}; "2";"4" **\dir{-};
        "4";"5" **\dir{-}; "5";"6" **\dir{-}; "9";"19" **\dir{-};
        "20";"19" **\dir{-}; "20";"21" **\dir{-}; "22";"21" **\dir{-};
        "15";"17" **\dir{-}; "23";"17" **\dir{-};
        \ar@{-}@/^3pt/ "6";"12"
        \ar@{-} "3";"18" |<(0.07)\hole
        \ar@{-} "10";"14" |<(0.13)\hole
        \ar@{-} "16";"24" |<(0.1)\hole |<(0.5)\hole
    \endxy}
    \overset{\text{(tri)}}{=}~
    \vcenter{\xy
        (-3,20)="1"; ( 0,15)="2"; ( 3,15)="3";
        ( 6,13)*[o]=<9pt>[Fo]{\scriptstyle \nu}="4";
        ( 0,10)="5";
        ( 6, 8)*[o]=<9pt>[Fo]{\scriptstyle r}="6";
        ( 0, 5)="7";
        ( 6, 3)*[o]=<9pt>[Fo]{\scriptstyle t}="8";
        (-3,-1)="9"; ( 0,-1)="10";
        ( 6,-3)*[o]=<9pt>[Fo]{\scriptstyle \nu}="11";
        (10,1.5)="12";
        (10,-1)*[o]=<5pt>[Fo]{~}="13";
        (14, 1)="14"; ( 8,-7)="15";
        (14,-7)*[o]=<9pt>[Fo]{\scriptstyle t}="16";
        (11,-10)="17";
        (11,-13)*[o]=<5pt>[Fo]{~}="18";
        ( 3,-16)="19";
        "9";"10" **\crv{(-3,-5)&(0,-5)};
        "2";"3" **\crv{(0,18)&(3,18)};
        "1";"9" **\dir{-}; "2";"10" **\dir{-}; "2";"4" **\dir{-};
        "5";"6" **\dir{-}; "7";"8" **\dir{-}; "10";"11" **\dir{-};
        "8";"12" **\dir{-}; "12";"13" **\dir{-}; "6";"14" **\dir{-};
        "11";"16" **\dir{-}; "17";"16" **\dir{-}; "17";"15" **\dir{-};
        "17";"18" **\dir{-};
        \ar@{-}@/^3pt/ "4";"12" |<(0.7)\hole
        \ar@{-} "14";"15" |<(0.69)\hole
        \ar@{-} "3";"19" |<(0.03)\hole |<(0.19)\hole |<(0.35)\hole |<(0.55)\hole
    \endxy}
    \overset{\text{(13)}}{=}~
    \vcenter{\xy
        (-3,20)="1"; ( 0,15)="2"; ( 0,11)="a"; ( 3,15)="3";
        ( 6,13)*[o]=<9pt>[Fo]{\scriptstyle \nu}="4";
        ( 6, 9)*[o]=<9pt>[Fo]{\scriptstyle t}="5";
        (10,7.5)="6";
        (10,4.5)*[o]=<5pt>[Fo]{~}="7";
        ( 0, 4)="8";
        ( 6, 2)*[o]=<9pt>[Fo]{\scriptstyle r}="9";
        (-3,-1)="z"; ( 0,-1)="10";
        ( 6,-3)*[o]=<9pt>[Fo]{\scriptstyle \nu}="11";
        (11,-2)="x";
        (12,-7)*[o]=<9pt>[Fo]{\scriptstyle t}="12";
        ( 6,-7)="y"; ( 9,-10)="13";
        ( 9,-13)*[o]=<5pt>[Fo]{~}="14";
        ( 3,-16)="15";
        "1";"z" **\dir{-}; "2";"10" **\dir{-}; "10";"11" **dir{-};
        "z";"10" **\crv{(-3,-5)&(0,-5)}; "2";"3" **\crv{(0,18)&(3,18)};
        "2";"4" **\dir{-}; "a";"5" **\dir{-}; "6";"5" **\dir{-};
        "6";"7" **\dir{-}; "8";"9" **\dir{-}; "x";"9" **\dir{-};
        "10";"11" **\dir{-}; "12";"11" **\dir{-}; "12";"13" **\dir{-};
        "y";"13" **\dir{-}; "14";"13" **\dir{-};
        \ar@{-}@/^3pt/ "4";"6"
        \ar@{-} "3";"15" |<(0.03)\hole |<(0.16)\hole |<(0.38)\hole |<(0.55)\hole
        \ar@{-} "x";"y" |\hole
    \endxy}
    \overset{\text{(2,$\nu$)}}{=}~
    \vcenter{\xy
        (-3,15)="1"; ( 0,10)="2"; ( 0,6)="a"; ( 3,10)="3";
        ( 6,8)*[o]=<9pt>[Fo]{\scriptstyle r}="4";
        (10,6.5)*[o]=<5pt>[Fo]{~}="5";
        ( 0, 4)="8";
        ( 6, 2)*[o]=<9pt>[Fo]{\scriptstyle r}="9";
        (-3,-1)="z";
        ( 0,-1)="10";
        ( 6,-3)*[o]=<9pt>[Fo]{\scriptstyle \nu}="11";
        (11,-2)="x";
        (12,-7)*[o]=<9pt>[Fo]{\scriptstyle t}="12";
        ( 6,-7)="y"; ( 9,-10)="13";
        ( 9,-13)*[o]=<5pt>[Fo]{~}="14";
        ( 3,-16)="15";
        "1";"z" **\dir{-}; "2";"10" **\dir{-}; "10";"11" **dir{-};
        "z";"10" **\crv{(-3,-5)&(0,-5)};
        "2";"3" **\crv{(0,13)&(3,13)};
        "2";"4" **\dir{-}; "4";"5" **\dir{-}; "8";"9" **\dir{-};
        "x";"9" **\dir{-}; "10";"11" **\dir{-}; "12";"11" **\dir{-};
        "12";"13" **\dir{-}; "y";"13" **\dir{-}; "14";"13" **\dir{-};
        \ar@{-} "3";"15" |<(0.03)\hole |<(0.27)\hole |<(0.46)\hole
        \ar@{-} "x";"y" |\hole
    \endxy}
    \overset{\text{(2,c)}}{=}
    \vcenter{\xy
        (-3, 6)="1";
        (-3, 0)="2";
        ( 0, 0)="3";
        ( 3, 0)="4";
        (10,-3)="5";
        ( 7,-6)*[o]=<9pt>[Fo]{\scriptstyle \nu}="6";
        (14,-6)="7";
        ( 6,-12)*[o]=<9pt>[Fo]{\scriptstyle r}="8";
        (13,-12)*[o]=<9pt>[Fo]{\scriptstyle t}="9";
        (10,-15)="10";
        (10,-18)*[o]=<5pt>[Fo]{~}="11";
        ( 3,-21)="12";
        "1";"2" **\dir{-};
        "2";"3" **\crv{(-3,-3)&(0,-3)};
        "3";"4" **\crv{(0,3)&(3,3)};
        "3";"5" **\dir{-};
        "5";"6" **\dir{-};
        "5";"7" **\dir{-};
        "6";"9" **\dir{-};
        "10";"8" **\dir{-};
        "10";"9" **\dir{-};
        "10";"11" **\dir{-};
        \ar@{-} "4";"12" |<(0.07)\hole
        \ar@{-} "7";"8" |\hole
    \endxy}
\]
\[
    ~\overset{\text{(13)}}{=}~
    \vcenter{\xy
        (-3, 6)="1";
        (-3, 0)="2";
        ( 0, 0)="3";
        ( 3, 0)="4";
        ( 9,-3)="5";
        ( 6,-6)*[o]=<9pt>[Fo]{\scriptstyle \nu}="6";
        (12,-6)="7";
        ( 6,-12)*[o]=<9pt>[Fo]{\scriptstyle t}="8";
        (12,-12)*[o]=<9pt>[Fo]{\scriptstyle r}="9";
        ( 9,-15)="10";
        ( 9,-18)*[o]=<5pt>[Fo]{~}="11";
        ( 3,-21)="12";
        "1";"2" **\dir{-};
        "2";"3" **\crv{(-3,-3)&(0,-3)};
        "3";"4" **\crv{(0,3)&(3,3)};
        "3";"5" **\dir{-};
        "5";"6" **\dir{-};
        "5";"7" **\dir{-};
        "6";"8" **\dir{-};
        "7";"9" **\dir{-};
        "10";"8" **\dir{-};
        "10";"9" **\dir{-};
        "10";"11" **\dir{-};
        \ar@{-} "4";"12" |<(0.07)\hole
    \endxy}
    ~\overset{\text{(2)}}{=}~
    \vcenter{\xy
        (-3, 6)="1";
        (-3, 0)="2";
        ( 0, 0)="3";
        ( 3, 0)="4";
        ( 9,-3)="5";
        ( 6,-6)*[o]=<9pt>[Fo]{\scriptstyle \nu}="6";
        (12,-6)="7";
        ( 9,-9)="8";
        ( 9,-12)*[o]=<5pt>[Fo]{~}="9";
        ( 3,-15)="10";
        "1";"2" **\dir{-};
        "2";"3" **\crv{(-3,-3)&(0,-3)};
        "3";"4" **\crv{(0,3)&(3,3)};
        "3";"5" **\dir{-};
        "5";"6" **\dir{-};
        "5";"7" **\dir{-};
        "8";"6" **\dir{-};
        "8";"7" **\dir{-};
        "8";"9" **\dir{-};
        \ar@{-} "4";"10" |<(0.07)\hole
    \endxy}
    ~\overset{\text{($\nu$)}}{=}~
    \vcenter{\xy
        (-3, 6)="1";
        (-3, 0)="2";
        ( 0, 0)="3";
        ( 3, 0)="4";
        ( 6,-2)*[o]=<9pt>[Fo]{\scriptstyle t}="5";
        ( 6,-7)*[o]=<5pt>[Fo]{~}="6";
        ( 3,-10)="7";
        "1";"2" **\dir{-};
        "2";"3" **\crv{(-3,-3)&(0,-3)};
        "3";"4" **\crv{(0,3)&(3,3)};
        "3";"5" **\dir{-};
        "5";"6" **\dir{-};
        \ar@{-} "4";"7" |<(0.1)\hole
    \endxy}
    ~\overset{\text{(2,c)}}{=}~
    \vcenter{\xy
        (-3, 6)="1";
        (-3, 0)="2";
        ( 0, 0)="3";
        ( 3, 0)="4";
        ( 3,-6)="5";
        "1";"2" **\dir{-};
        "2";"3" **\crv{(-3,-3)&(0,-3)};
        "3";"4" **\crv{(0,3)&(3,3)};
        "4";"5" **\dir{-};
    \endxy}
    ~\overset{\text{(tri)}}{=}~
    \vcenter{\xy
        (0,12); (0,0); **\dir{-}
    \endxy}
    ~= 1_{M^*}.
\]

This completes the proof that $M^*$ is the left dual of $M$ in
$\Comod_f(H)$, and hence that $\Comod_f(H)$ is left autonomous.
\end{proof}

%=========================================================================%
\section{Frobenius monoid example}\label{sec-frobeg}
%=========================================================================%

Let $R$ be a separable Frobenius monoid in $\V$. In this section we prove
that $R \ox R$ is an example of a weak Hopf monoid with an invertible
antipode. In the case $\V = \mathbf{Vect}$, this example is essentially the
same as in~\cite[Appendix]{BNS}.

Let $R$ be a Frobenius monoid in $\V$. Then $R \ox R$ becomes a comonoid
via
\[
    \delta =~
    \vcenter{\xy
        (0,5)*{~}; (0,-5)*{~};
        ( 0, 4)="1";
        ( 3, 4)="2";
        (-3,-4)="3";
        ( 0,-4)="4";
        ( 3,-4)="5";
        ( 6,-4)="6";
        "4";"5" **\crv{(0,-1)&(3,-1)};
        "1";"3" **\crv{(0,1)&(-3,-2)};
        "2";"6" **\crv{(3,1)&(6,-2)};
    \endxy}
    \qquad \text{and} \qquad
    \epsilon =~
    \vcenter{\xy
        (0,0);(3,0) **\crv{(0,-3)&(3,-3)};
    \endxy}~.
\]
(where, for simplicity, in this section we will adopt the simpler notation
\[
    \vcenter{\xy
        (-3,-3)="1"; 
        ( 3,-3)="2"; 
        ( 0, 0)="3";
        ( 0, 3)*[o]=<5pt>[Fo]{~}="4";
        "1";"3" **\dir{-}; 
        "2";"3" **\dir{-}; 
        "4";"3" **\dir{-}; 
    \endxy}
    ~=~
    \vcenter{\xy (0,0); (5,0) **\crv{(0,5)&(5,5)}; \endxy}
    \qquad \text{and} \qquad
    \vcenter{\xy
        (-3, 3)="1"; 
        ( 3, 3)="2"; 
        ( 0, 0)="3";
        ( 0,-3)*[o]=<5pt>[Fo]{~}="4";
        "1";"3" **\dir{-}; 
        "2";"3" **\dir{-}; 
        "4";"3" **\dir{-}; 
    \endxy}
    ~=~
    \vcenter{\xy (0,0); (5,0) **\crv{(0,-5)&(5,-5)}; \endxy}\ \ )~,
\]
and a monoid via
\[
    \mu =~
    \vcenter{\xy
        (0,11)*{~}; (0,-3)*{~};
        (-8, 10)="1"; 
        (-2.5, 10)="2"; 
        ( 2.5, 10)="3";
        ( 8, 10)="4";
        (-6.6, 4)="5";
        (-1.5, 4)="6";
        (-4, 2)="7";
        (-4,-2)="9";
        ( 4, 3)="8";
        ( 4,-2)="10";
        "1";"6" **\dir{-}; 
        "2";"8" **\dir{-}; 
        "4";"8" **\dir{-}; 
        "5";"7" **\dir{-}; 
        "6";"7" **\dir{-}; 
        "9";"7" **\dir{-}; 
        "10";"8" **\dir{-}; 
        \ar@{-} "3";"5" |<(0.3)\hole |<(0.7)\hole
    \endxy}
    \qquad \text{and} \qquad
    \eta =~
    \vcenter{\xy
        (-2,5)*[o]=<5pt>[Fo]{~}; (-2,0) **\dir{-};
        ( 2,5)*[o]=<5pt>[Fo]{~}; ( 2,0) **\dir{-};
    \endxy}\ \ .
\]
The comonoid structure is via the comonad generated by the adjunction $R
\dashv R$. The monoid structure is the usual monoid structure (viewing $R$
as a monoid) on the tensor product $R^\o \ox R$, where $R^\o$ is the
opposite monoid of $R$.

\begin{proposition}\label{prop-frob-wha}
If $R$ is separable, meaning $\mu \delta = 1_R$, then $R \ox R$ is a weak
bimonoid. An invertible antipode $\nu$ on $R \ox R$ is given by
\[
    \nu =~
    \vcenter{\xy 0;/r.15pc/:
        (0, 9)="1";
        (8,-9)="6";
        (8, 9)="2";
        (1.3,-1.5)="3";
        (0,-3)="4";
        (2,-9)="5";
        "1";"6" **\crv{(0,5)&(8,2)};
        "4";"5" **\crv{(-2,-6)&(-4,-3)&(-2,1)&(2,-3)};
        \ar@{-}@/^1pt/ "2";"3" |<(0.6)\hole
    \endxy}~.
\]
which makes $R \ox R$ into a weak Hopf monoid.
\end{proposition}

The following three sets of calculations establish respectively the
axioms (b), (v), and (w), and hence the first claim.

The axiom (b) is given by:
\[
    (\mu \ox \mu)(1 \ox c \ox 1)(\delta \ox \delta) =~
    \vcenter{\xy
        (0,11)*{~}; (0,-13)*{~};
        (-11,10)="1"; 
        (-3,10)="2"; 
        ( 2,10)="3";
        ( 9,10)="4";
        (-6, 3)="5";
        (-3, 3)="6";
        ( 3.5,4)="7";
        ( 6, 3)="8";
        (-9,-9)="9";
        (-3,-9)="10";
        ( 3,-9)="11";
        ( 9,-5)="12";
        (-9,-12)="13";
        (-3,-12)="14";
        ( 3,-12)="15";
        ( 9,-12)="16";
        (-11,-7)="w";
        (-8,-7)="x";
        ( 1,-7)="y";
        ( 5,-7)="z";
        "1";"x" **\crv{(-11,0)&(-8,-4)}; 
        "6";"z" **\dir{-}; 
        "9";"w" **\dir{-}; 
        "9";"x" **\dir{-}; 
        "11";"y" **\dir{-}; 
        "11";"z" **\dir{-}; 
        "2";"12" **\dir{-}; 
        "4";"12" **\dir{-}; 
        "9";"13" **\dir{-}; 
        "10";"14" **\dir{-}; 
        "11";"15" **\dir{-}; 
        "12";"16" **\dir{-}; 
        "5";"10" **\dir{-}; 
        "5";"6" **\crv{(-7,6)&(-4,6)}; 
        "7";"8" **\crv{(4,6)&(7,6)}; 
        \ar@{-} "7";"10" |<(0.1)\hole |<(0.45)\hole
        \ar@{-} "3";"w" |<(0.15)\hole |<(0.37)\hole |<(0.6)\hole |<(0.82)\hole
        \ar@{-} "8";"y" |<(0.32)\hole |<(0.7)\hole
    \endxy}
    ~\overset{\text{(nat)}}{=}~
    \vcenter{\xy
        (-8, 9)="1"; 
        (-3, 9)="2"; 
        ( 3, 9)="3";
        ( 8, 9)="4";
        (-7, 4)="5";
        (-1, 4)="6";
        ( 3, 4)="7";
        (-3, 2)="8";
        (-2,-5)="x";
        ( 2,-5)="y";
        (-6,-8)="9";
        (-2,-8)="10";
        ( 2,-8)="11";
        ( 6,-8)="12";
        "1";"6" **\dir{-}; 
        "2";"7" **\dir{-}; 
        "4";"7" **\dir{-}; 
        "5";"8" **\dir{-}; 
        "6";"8" **\dir{-}; 
        "9";"8" **\crv{(-6,-4)&(-3,-1)}; 
        "12";"7" **\crv{(6,-4)&(3,-1)}; 
        "x";"10" **\dir{-}; 
        "y";"11" **\dir{-}; 
        "x";"y" **\crv{(-1,-2)&(1,-2)}; 
        "x";"y" **\crv{(-6,1)&(6,1)}; 
        \ar@{-} "3";"5" |<(0.32)\hole |<(0.7)\hole
    \endxy}
    \overset{\text{(sep)}}{=}
    \vcenter{\xy
        (-8, 9)="1"; 
        (-3, 9)="2"; 
        ( 3, 9)="3";
        ( 8, 9)="4";
        (-7, 4)="5";
        (-1, 4)="6";
        ( 3, 4)="7";
        (-3, 2)="8";
        (-6,-6)="9";
        (-2,-6)="10";
        ( 2,-6)="11";
        ( 6,-6)="12";
        "1";"6" **\dir{-}; 
        "2";"7" **\dir{-}; 
        "4";"7" **\dir{-}; 
        "5";"8" **\dir{-}; 
        "6";"8" **\dir{-}; 
        "8";"9" **\crv{(-3,0)&(-6,-4)}; 
        "7";"12" **\crv{(3,0)&(6,-4)}; 
        "10";"11" **\crv{(-2,-2)&(2,-2)}; 
        \ar@{-} "3";"5" |<(0.32)\hole |<(0.7)\hole
    \endxy}
    = \delta\mu.
\]

Axiom (v) is seen from the diagrams: 

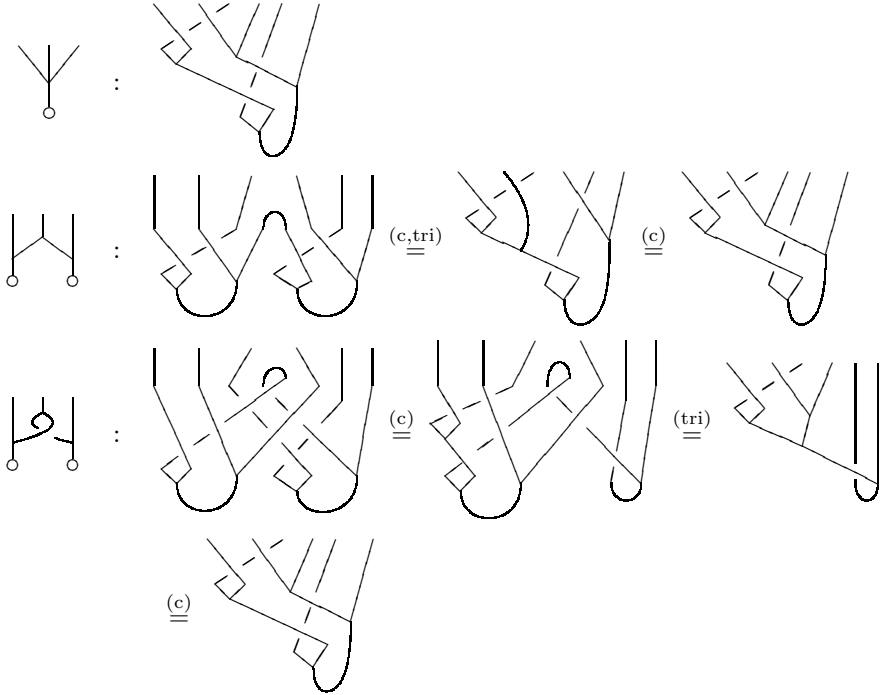
\begin{align*}
    \vcenter{\xy
        (-4,4)="l";
        (0,4)="c";
        (4,4)="r";
        (0,-1)="m";
        (0,-5)*[o]=<5pt>[Fo]{~}="b";
        "l";"m" **\dir{-};
        "c";"m" **\dir{-};
        "r";"m" **\dir{-};
        "m";"b" **\dir{-};
    \endxy}
    \quad &: \quad
    \vcenter{\xy
        (-1,12)="1"; 
        ( 5,12)="2"; 
        ( 9,12)="3";
        (13,12)="4";
        (16,12)="5";
        (21,12)="6";
        (0,6)="7";
        (4,6)="8";
        (10,5)="9";
        (2,4)="10";
        (10.5,-3)="11";
        (15,-2)="12";
        (18,1)="13";
        (13,-5)="14";
        "2";"9" **\dir{-}; 
        "4";"9" **\dir{-}; 
        "13";"9" **\dir{-}; 
        "13";"6" **\dir{-}; 
        "1";"8" **\dir{-}; 
        "7";"10" **\dir{-}; 
        "8";"10" **\dir{-}; 
        "12";"10" **\dir{-}; 
        "11";"14" **\dir{-}; 
        "12";"14" **\dir{-}; 
        "14";"13" **\crv{(13,-9)&(18,-9)&(18,0)}; 
        \ar@{-} "3";"7" |<(0.3)\hole |<(0.7)\hole
        \ar@{-} "5";"11" |<(0.54)\hole |<(0.81)\hole
    \endxy}
\\
    \vcenter{\xy
        (-4,4)="bl";
        (0,4)="bc";
        (4,4)="br";
        (-4,-5)*[o]=<5pt>[Fo]{~}="tl";
        (-4,-2)="ml";
        (0,1)="mc";
        (4,-2)="mr";
        (4,-5)*[o]=<5pt>[Fo]{~}="tr";
        "bl";"tl" **\dir{-};
        "br";"tr" **\dir{-};
        "bc";"mc" **\dir{-};
        "mc";"ml" **\dir{-};
        "mc";"mr" **\dir{-};
    \endxy}
    \quad &: \quad
    \vcenter{\xy
        (13.5,12)="x"; 
        (16.5,12)="y"; 
        (-1,12)="1"; 
        ( 5,12)="2"; 
        (10,12)="3";
        (0,6)="7";
        (4,6)="8";
        (10,5)="9";
        (2,4)="10";
        (20,12)="4";
        (24,12)="5";
        (28,12)="6";
        (15,6)="11";
        (20,5)="12";
        (26,5)="13";
        (18,4)="14";
        "x";"y" **\crv{(13.5,15)&(16.5,15)};
        "x";"9" **\dir{-}; 
        "y";"12" **\dir{-}; 
        (-1,19);"1" **\dir{-}; 
        ( 5,19);"2" **\dir{-}; 
        (12,19);"3" **\dir{-}; 
        (18,19);"4" **\dir{-}; 
        (24,19);"5" **\dir{-}; 
        (28,19);"6" **\dir{-}; 
        "1";"8" **\dir{-}; 
        "2";"9" **\dir{-}; 
        "4";"13" **\dir{-}; 
        "6";"13" **\dir{-}; 
        "7";"10" **\dir{-}; 
        "8";"10" **\dir{-}; 
        "11";"14" **\dir{-}; 
        "12";"14" **\dir{-}; 
        "10";"9" **\crv{(2,-1)&(10,-1)}; 
        "14";"13" **\crv{(18,-1)&(26,-1)}; 
        \ar@{-} "3";"7" |<(0.3)\hole |<(0.7)\hole
        \ar@{-} "5";"11" |<(0.3)\hole |<(0.7)\hole
    \endxy}
    ~\overset{\text{(c,tri)}}{=}~
    \vcenter{\xy
        (-1,12)="1"; 
        ( 5,12)="2"; 
        ( 9,12)="3";
        (13,12)="4";
        (17,12)="5";
        (21,12)="6";
        (0,6)="7";
        (4,6)="8";
        (2,4)="10";
        (10.5,-3)="11";
        (15,-2)="12";
        (19,3)="13";
        (13,-5)="14";
        "13";"4" **\dir{-}; 
        "13";"6" **\dir{-}; 
        "1";"8" **\dir{-}; 
        "7";"10" **\dir{-}; 
        "8";"10" **\dir{-}; 
        "10";"12" **\dir{-}?(0.4)="x"; 
        "2";"x" **\crv{(10,6)}; 
        "11";"14" **\dir{-}; 
        "12";"14" **\dir{-}; 
        "14";"13" **\crv{(13,-9)&(19,-9)&(19,0)}; 
        \ar@{-} "3";"7" |<(0.3)\hole |<(0.7)\hole
        \ar@{-} "5";"11" |<(0.25)\hole |<(0.81)\hole
    \endxy}
    ~\overset{\text{(c)}}{=}~
    \vcenter{\xy
        (-1,12)="1"; 
        ( 5,12)="2"; 
        ( 9,12)="3";
        (13,12)="4";
        (16,12)="5";
        (21,12)="6";
        (0,6)="7";
        (4,6)="8";
        (10,5)="9";
        (2,4)="10";
        (10.5,-3)="11";
        (15,-2)="12";
        (18,1)="13";
        (13,-5)="14";
        "2";"9" **\dir{-}; 
        "4";"9" **\dir{-}; 
        "13";"9" **\dir{-}; 
        "13";"6" **\dir{-}; 
        "1";"8" **\dir{-}; 
        "7";"10" **\dir{-}; 
        "8";"10" **\dir{-}; 
        "12";"10" **\dir{-}; 
        "11";"14" **\dir{-}; 
        "12";"14" **\dir{-}; 
        "14";"13" **\crv{(13,-9)&(18,-9)&(18,0)}; 
        \ar@{-} "3";"7" |<(0.3)\hole |<(0.7)\hole
        \ar@{-} "5";"11" |<(0.54)\hole |<(0.81)\hole
    \endxy}
\\
    \vcenter{\xy
        (-4,4)="bl";
        (0,4)="bc";
        (4,4)="br";
        (-4,-5)*[o]=<5pt>[Fo]{~}="tl";
        (-4,-2)="ml";
        (0,2)="mc";
        (4,-2)="mr";
        (4,-5)*[o]=<5pt>[Fo]{~}="tr";
        "bl";"tl" **\dir{-};
        "br";"tr" **\dir{-};
        "bc";"mc" **\dir{-};
        "mc";"ml" **\crv{(4,0)}?(0.6)*{\hole}="x";
        "mc";"x" **\crv{(-3,0.5)};
        "mr";"x" **\crv{(1,-1.5)};
    \endxy}
    \quad &: \quad
    \vcenter{\xy
        (13.5,17)="x"; 
        (16.5,18)="y"; 
        (-1,17)="1"; 
        ( 5,17)="2"; 
        ( 9,17)="3";
        (21,17)="4";
        (24,15)="5";
        (28,17)="6";
        (0,6)="7";
        (4,6)="8";
        (10,5)="9";
        (2,4)="10";
        (15,5)="11";
        (20,5)="12";
        (26,5)="13";
        (18,3)="14";
        "x";"y" **\crv{(13.5,20)&(16.5,20)};
        (-1,22);"1" **\dir{-}; 
        ( 5,22);"2" **\dir{-}; 
        (12,22);"3" **\dir{-}; 
        (18,22);"4" **\dir{-}; 
        (24,22);"5" **\dir{-}; 
        (28,22);"6" **\dir{-}; 
        "1";"8" **\dir{-}; 
        "2";"9" **\dir{-}; 
        "4";"9" **\dir{-}; 
        "6";"13" **\dir{-}; 
        "7";"10" **\dir{-}; 
        "8";"10" **\dir{-}; 
        "11";"14" **\dir{-}; 
        "12";"14" **\dir{-}; 
        "10";"9" **\crv{(2,-1)&(10,-1)}; 
        "14";"13" **\crv{(18,-1)&(26,-1)}; 
        \ar@{-} "y";"7" |<(0.55)\hole |<(0.8)\hole
        \ar@{-} "3";"12" |<(0.2)\hole |<(0.55)\hole
        \ar@{-} "x";"13" |<(0.05)\hole |<(0.35)\hole
        \ar@{-} "5";"11" |<(0.40)\hole |<(0.75)\hole
    \endxy}
    ~\overset{\text{(c)}}{=}~
    \vcenter{\xy
        (13.5,18)="x"; 
        (16.5,19)="y"; 
        (-2,13)="a";
        (2,13)="b";
        (0,11)="c";
        (-1,18)="1"; 
        ( 5,18)="2"; 
        ( 9,18)="3";
        (21,18)="4";
        (24,16)="5";
        (28,18)="6";
        (0,6)="7";
        (4,6)="8";
        (10,5)="9";
        (2,4)="10";
        (26,5)="13";
        (22,5)="14";
        "a";"c" **\dir{-}; 
        "b";"c" **\dir{-}; 
        "x";"y" **\crv{(13.5,22)&(16.5,22)};
        (-1,24);"1" **\dir{-}; 
        ( 5,24);"2" **\dir{-}; 
        (12,24);"3" **\dir{-}; 
        (18,24);"4" **\dir{-}; 
        (24,24);"5" **\dir{-}; 
        (28,24);"6" **\dir{-}; 
        "1";"b" **\dir{-}; 
        "c";"8" **\dir{-}; 
        "2";"9" **\dir{-}; 
        "4";"9" **\dir{-}; 
        "6";"13" **\dir{-}; 
        "7";"10" **\dir{-}; 
        "8";"10" **\dir{-}; 
        "10";"9" **\crv{(2,-1)&(10,-1)}; 
        "14";"13" **\crv{(22,2)&(26,2)}; 
        \ar@{-} "y";"7" |<(0.55)\hole |<(0.8)\hole
        \ar@{-} "3";"a" |<(0.3)\hole |<(0.7)\hole
        \ar@{-} "x";"13" |<(0.03)\hole |<(0.35)\hole
        \ar@{-} "5";"14" |<(0.7)\hole
    \endxy}
    ~\overset{\text{(tri)}}{=}~
    \vcenter{\xy
        (-1,12)="1"; 
        ( 5,12)="2"; 
        ( 9,12)="3";
        (13,12)="4";
        (16,12)="5";
        (19,12)="6";
        (0, 6)="7";
        (4, 6)="8";
        (2, 4)="9";
        (10,5)="10";
        (9, 1)="11";
        (16,-4)="12";
        (19,-4)="13";
        "1";"8" **\dir{-}; 
        "7";"9" **\dir{-}; 
        "8";"9" **\dir{-}; 
        "2";"10" **\dir{-}; 
        "4";"10" **\dir{-}; 
        "9";"11" **\dir{-}; 
        "10";"11" **\dir{-}; 
        "11";"13" **\dir{-}; 
        "6";"13" **\dir{-}; 
        "12";"13" **\crv{(16,-7)&(19,-7)}; 
        \ar@{-} "3";"7" |<(0.3)\hole |<(0.7)\hole
        \ar@{-} "5";"12"|<(0.9)\hole
    \endxy}
\\
    &\qquad \overset{\text{(c)}}{=}~
    \vcenter{\xy
        (-1,12)="1"; 
        ( 5,12)="2"; 
        ( 9,12)="3";
        (13,12)="4";
        (16,12)="5";
        (21,12)="6";
        (0,6)="7";
        (4,6)="8";
        (10,5)="9";
        (2,4)="10";
        (10.5,-3)="11";
        (15,-2)="12";
        (18,1)="13";
        (13,-5)="14";
        "2";"9" **\dir{-}; 
        "4";"9" **\dir{-}; 
        "13";"9" **\dir{-}; 
        "13";"6" **\dir{-}; 
        "1";"8" **\dir{-}; 
        "7";"10" **\dir{-}; 
        "8";"10" **\dir{-}; 
        "12";"10" **\dir{-}; 
        "11";"14" **\dir{-}; 
        "12";"14" **\dir{-}; 
        "14";"13" **\crv{(13,-9)&(18,-9)&(18,0)}; 
        \ar@{-} "3";"7" |<(0.3)\hole |<(0.7)\hole
        \ar@{-} "5";"11" |<(0.54)\hole |<(0.81)\hole
    \endxy}
\end{align*}

For (w), by the naturality of the braiding and the counit property of $R$
each equation in (w), i.e., 
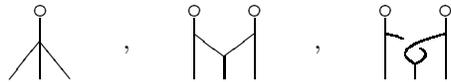
\[
    \vcenter{\xy
        (-4,-4)="l";
        (0,-4)="c";
        (4,-4)="r";
        (0,1)="m";
        (0,5)*[o]=<5pt>[Fo]{~}="b";
        "l";"m" **\dir{-};
        "c";"m" **\dir{-};
        "r";"m" **\dir{-};
        "m";"b" **\dir{-};
    \endxy}
    \qquad , \qquad
    \vcenter{\xy
        (-4,-4)="bl";
        (0,-4)="bc";
        (4,-4)="br";
        (-4,5)*[o]=<5pt>[Fo]{~}="tl";
        (-4,2)="ml";
        (0,-1)="mc";
        (4,2)="mr";
        (4,5)*[o]=<5pt>[Fo]{~}="tr";
        "bl";"tl" **\dir{-};
        "br";"tr" **\dir{-};
        "bc";"mc" **\dir{-};
        "mc";"ml" **\dir{-};
        "mc";"mr" **\dir{-};
    \endxy}
    \qquad , \qquad
    \vcenter{\xy
        (4,-4)="bl";
        (0,-4)="bc";
        (-4,-4)="br";
        (4,5)*[o]=<5pt>[Fo]{~}="tl";
        (4,2)="ml";
        (0,-2)="mc";
        (-4,2)="mr";
        (-4,5)*[o]=<5pt>[Fo]{~}="tr";
        "bl";"tl" **\dir{-};
        "br";"tr" **\dir{-};
        "bc";"mc" **\dir{-};
        "mc";"ml" **\crv{(-4,0)}?(0.6)*{\hole}="x";
        "mc";"x" **\crv{( 3,-0.5)};
        "mr";"x" **\crv{(-1,1.5)};
    \endxy}
\]
is easily seen to be equal to the following diagram
\[
    \vcenter{\xy
        (0,6)*{~}; (0,-1)*{~}; 
        ( 1,0)="1"; 
        ( 3,0)="2"; 
        ( 6,0)="3";
        ( 8,0)="4";
        (11,0)="5";
        (13,0)="6";
        ( 1,5)*[o]=<5pt>[Fo]{~}="x"; 
        (13,5)*[o]=<5pt>[Fo]{~}="y";
        "1";"x" **\dir{-}; 
        "6";"y" **\dir{-}; 
        "2";"3" **\crv{(3,3)&(6,3)};
        "4";"5" **\crv{(8,3)&(11,3)};
    \endxy}~.
\]
Thus, $R \ox R$ is a weak bimonoid. We next prove that that $R \ox R$ is a
weak Hopf monoid with invertible antipode
\[
    \nu =~
    \vcenter{\xy 0;/r.15pc/:
        (0, 9)="1";
        (8,-9)="6";
        (8, 9)="2";
        (1.3,-1.5)="3";
        (0,-3)="4";
        (2,-9)="5";
        "1";"6" **\crv{(0,5)&(8,2)};
        "4";"5" **\crv{(-2,-6)&(-4,-3)&(-2,1)&(2,-3)};
        \ar@{-}@/^1pt/ "2";"3" |<(0.6)\hole
    \endxy}~.
\]
An inverse to $\nu$ is easily seen to be given by
\[
    \nu^{-1} =~
    \vcenter{\xy 0;/r.15pc/:
        (0,9)="a";
        (-8,-9)="b";
        (-8,9)="x";
        (-1,-1.5)="y";
        (0,-3)="0";
        (-2,-9)="1";
        "a";"b" **\crv{(0,5)&(-8,2)}?(0.52)*{\hole}="h";
        "x";"h" **\crv{(-8,6)};
        "h";"y" **\dir{-};
        "0";"1" **\crv{(2,-6)&(4,-3)&(2,1)&(-2,-3)};
    \endxy}~,
\]
and so the antipode is invertible. We note that (in simplified form)
\[
    r =~
    \vcenter{\xy 0;/r.15pc/:
        (0, 9)="1";
        (8, 9)="2";
        (1.3,-1.5)="3";
        (0,-3)="4";
        (2,-9)="5";
        (8,-9)="6";
        (5,-13)="7";
        (5,-17)="8";
        (12,-6)*[o]=<5pt>[Fo]{~}="a"; 
        (12,-17)="b";
        "a";"b" **\dir{-};
        "1";"6" **\crv{(0,5)&(8,2)};
        "4";"5" **\crv{(-2,-6)&(-4,-3)&(-2,1)&(2,-3)};
        "5";"7" **\dir{-};
        "6";"7" **\dir{-};
        "8";"7" **\dir{-};
        \ar@{-}@/^1pt/ "2";"3" |<(0.6)\hole
    \endxy}
    \quad\qquad\text{and}\quad\qquad
    t =~
    \vcenter{\xy
        ( 1,5)*[o]=<5pt>[Fo]{~}="1"; 
        ( 1,0)="2";
        ( 2,9)="3"; 
        ( 8,9)="4"; 
        ( 5,6)="5";
        ( 5,0)="6";
        "1";"2" **\dir{-};
        "3";"5" **\dir{-};
        "4";"5" **\dir{-};
        "6";"5" **\dir{-};
    \endxy}~.
\]
The following calculations then prove the antipode axioms.
\[
    \mu(\nu \ox 1)\delta = ~
    \vcenter{\xy
        (0,32)*{~}; (0,-3)*{~}; 
        (-4, 31)="a"; 
        ( 4, 31)="b"; 
        (-4, 26)="c"; 
        ( 4, 26)="d"; 
        (-8, 21)="e"; 
        (-2, 21)="f"; 
        ( 2, 21)="g"; 
        ( 8, 21)="h"; 
        (-8, 16)="x"; 
        (-9.7, 15.7)="y"; 
        "f";"g" **\crv{(-2,24)&(2,24)};
        "a";"c" **\dir{-};
        "c";"e" **\dir{-};
        "b";"d" **\dir{-};
        "d";"h" **\dir{-};
        (-8, 11)="1"; 
        (-2, 11)="2"; 
        ( 2, 11)="3";
        ( 8, 11)="4";
        (-6.6, 4)="5";
        (-1.5, 4)="6";
        (-4, 2)="7";
        (-4,-2)="9";
        ( 4, 3)="8";
        ( 4,-2)="10";
        "e";"2" **\dir{-}; 
        "g";"3" **\dir{-}; 
        "h";"4" **\dir{-}; 
        "1";"6" **\dir{-}; 
        "2";"8" **\dir{-}; 
        "4";"8" **\dir{-}; 
        "5";"7" **\dir{-}; 
        "6";"7" **\dir{-}; 
        "9";"7" **\dir{-}; 
        "10";"8" **\dir{-}; 
        "y";"1" **\crv{(-12,15)&(-13,19)&(-8,19)&(-9,12)}; 
        \ar@{-} "3";"5" |<(0.29)\hole |<(0.7)\hole
        \ar@{-} "f";"x" |<(0.6)\hole
    \endxy}
    ~\overset{\text{(tri)}}{=}~
    \vcenter{\xy
        ( 0,3)="1"; 
        ( 0,0)="2";
        ( 2,11)="3"; 
        ( 8,11)="4"; 
        ( 5,8)="5";
        ( 5,0)="6";
        "1";"1" **\crv{(-8,11)&(8,11)};
        "1";"2" **\dir{-};
        "3";"5" **\dir{-};
        "4";"5" **\dir{-};
        "6";"5" **\dir{-};
    \endxy}
    ~\overset{\text{(sep)}}{=}~
    \vcenter{\xy
        ( 1,5)*[o]=<5pt>[Fo]{~}="1"; 
        ( 1,0)="2";
        ( 2,9)="3"; 
        ( 8,9)="4"; 
        ( 5,6)="5";
        ( 5,0)="6";
        "1";"2" **\dir{-};
        "3";"5" **\dir{-};
        "4";"5" **\dir{-};
        "6";"5" **\dir{-};
    \endxy}
    ~ = t
\]
\[
    \mu(1 \ox \nu)\delta = ~
    \vcenter{\xy
        (0,32)*{~}; (0,-3)*{~}; 
        (-4, 31)="a"; 
        ( 4, 31)="b"; 
        (-4, 26)="c"; 
        ( 4, 26)="d"; 
        (-8, 21)="e"; 
        (-4, 21)="f"; 
        ( 3, 21)="g"; 
        ( 8, 21)="h"; 
        ( 2, 16)="x"; 
        (0.7, 15.7)="y"; 
        "f";"g" **\crv{(-4,24)&(3,24)};
        "a";"c" **\dir{-};
        "c";"e" **\dir{-};
        "b";"d" **\dir{-};
        "d";"h" **\dir{-};
        (-8, 11)="1"; 
        (-2, 11)="2"; 
        ( 2, 11)="3";
        ( 8, 11)="4";
        (-6.6, 4)="5";
        (-1.5, 4)="6";
        (-4, 2)="7";
        (-4,-2)="9";
        ( 4, 3)="8";
        ( 4,-2)="10";
        "e";"1" **\dir{-}; 
        "f";"2" **\dir{-}; 
        "g";"4" **\dir{-}; 
        "1";"6" **\dir{-}; 
        "2";"8" **\dir{-}; 
        "4";"8" **\dir{-}; 
        "5";"7" **\dir{-}; 
        "6";"7" **\dir{-}; 
        "9";"7" **\dir{-}; 
        "10";"8" **\dir{-}; 
        "y";"3" **\crv{(-2,15)&(-3,19)&(1.7,19)&(2,12)}; 
        \ar@{-} "3";"5" |<(0.29)\hole |<(0.7)\hole
        \ar@{-} "h";"x" |<(0.6)\hole
    \endxy}
    ~\overset{\text{(nat)}}{=}~
    \vcenter{\xy 0;/r.15pc/:
        (0, 9)="1";
        (8, 9)="2";
        (1.3,-1.5)="3";
        (0,-3)="4";
        (2,-9)="5";
        (8,-9)="6";
        (5,-13)="7";
        (5,-17)="8";
        (14,-10)="a"; 
        (14,-17)="b";
        "a";"b" **\dir{-};
        "a";"a" **\crv{(2,3)&(26,3)};
        "a";"b" **\dir{-};
        "1";"6" **\crv{(0,5)&(8,2)};
        "4";"5" **\crv{(-2,-6)&(-4,-3)&(-2,1)&(2,-3)};
        "5";"7" **\dir{-};
        "6";"7" **\dir{-};
        "8";"7" **\dir{-};
        \ar@{-}@/^1pt/ "2";"3" |<(0.6)\hole
    \endxy}
    ~\overset{\text{(sep)}}{=}~
    \vcenter{\xy 0;/r.15pc/:
        (0, 9)="1";
        (8, 9)="2";
        (1.3,-1.5)="3";
        (0,-3)="4";
        (2,-9)="5";
        (8,-9)="6";
        (5,-13)="7";
        (5,-17)="8";
        (12,-6)*[o]=<5pt>[Fo]{~}="a"; 
        (12,-17)="b";
        "a";"b" **\dir{-};
        "1";"6" **\crv{(0,5)&(8,2)};
        "4";"5" **\crv{(-2,-6)&(-4,-3)&(-2,1)&(2,-3)};
        "5";"7" **\dir{-};
        "6";"7" **\dir{-};
        "8";"7" **\dir{-};
        \ar@{-}@/^1pt/ "2";"3" |<(0.6)\hole
    \endxy}
    ~= r
\]
\[
    \mu_3(\nu \ox 1 \ox \nu)\delta_3 =
    \vcenter{\xy
        (0,34)*{~}; (0,-15)*{~};
        (-8,33);(-8, 29) **\dir{-};
        ( 0,33);( 0, 29) **\dir{-};
        (-8, 29)="a";
        ( 0, 29)="b";
        (-16, 21)="k"; 
        (-11, 21)="i"; 
        (-8, 21)="j"; 
        "i";"j" **\crv{(-11,24)&(-8,24)};
        (-4, 21)="f"; 
        ( 3, 21)="g"; 
        ( 8, 21)="h"; 
        ( 2, 16)="x"; 
        (0.7, 15.7)="y"; 
        "f";"g" **\crv{(-4,24)&(3,24)};
        "a";"k" **\dir{-};
        "b";"h" **\dir{-};
        (-8, 11)="1"; 
        (-2, 11)="2"; 
        ( 2, 11)="3";
        ( 8, 11)="4";
        (-6.6, 4)="5";
        (-1.5, 4)="6";
        (-4, 2)="7";
        (-4,-2)="9";
        ( 4, 3)="8";
        ( 4,-2)="10";
        "j";"1" **\dir{-}; 
        "f";"2" **\dir{-}; 
        "g";"4" **\dir{-}; 
        "1";"6" **\dir{-}; 
        "2";"8" **\dir{-}; 
        "4";"8" **\dir{-}; 
        "5";"7" **\dir{-}; 
        "6";"7" **\dir{-}; 
        "9";"7" **\dir{-}; 
        "10";"8" **\dir{-}; 
        "y";"3" **\crv{(-2,15)&(-3,19)&(1.7,19)&(2,12)}; 
        (-17,11)="l";
        (-11,11)="m";
        (-17,-2)="c";
        (-11,-2)="d";
        (-17,16)="x1";
        (-18.7,15.7)="x2";
        "x2";"l" **\crv{(-21,15)&(-22,19)&(-17,19)&(-17,12)};
        "k";"m" **\dir{-};
        "d";"m" **\dir{-};
        "l";"c" **\dir{-};
        (0,-8);"10" **\dir{-};
        (0,-8);"d" **\dir{-};
        (0,-8);(0,-14) **\dir{-};
        (-10,-8);"c" **\dir{-};
        (-10,-8);(-13,-10) **\dir{-};
        (-16,-8);(-13,-10) **\dir{-};
        (-13,-14);(-13,-10) **\dir{-};
        \ar@{-} "9";(-16,-8) |<(0.3)\hole |<(0.7)\hole
        \ar@{-} "3";"5" |<(0.29)\hole |<(0.7)\hole
        \ar@{-} "h";"x" |<(0.6)\hole
        \ar@{-} "i";"x1" |<(0.6)\hole
    \endxy}
    ~\overset{\text{(sep)}}{=}
    \vcenter{\xy
        ( 6,15)="1";
        (12,15)="2";
        ( 6, 5)="3";
        (13, 5)="4";
        (16, 8)*[o]=<5pt>[Fo]{~}="5"; 
        ( 6,10)="6";
        (4.4,9.5)="7";
        (-2,1)="10";
        (-2,-1)="11";
        ( 8,0)="12";
        ( 6,-2)="13";
        ( 8,-7)="14";
        ( 6,-9)="15";
        ( 6,-14)="16";
        (13,-14)="17";
        "10";"12" **\crv{(0,7)&(2,7)};
        "13";"12" **\dir{-};
        "4";"17" **\dir{-};
        "1";"4" **\dir{-};
        "5";"4" **\dir{-};
        "14";"15" **\dir{-};
        "11";"14" **\crv{(-2,-4)&(-5,-4)&(-5,1)&(-2,1)};
        "7";"3" **\crv{(2,9)&(0,13)&(6,13)&(5.8,6)};
        \ar@{-} "2";"6" |<(0.6)\hole
        \ar@{-} "3";"16" |<(0.15)\hole |<(0.6)\hole
    \endxy}
    \overset{\text{(sep)}}{=}~
    \vcenter{\xy
        ( 6,15)="1";
        (12,15)="2";
        ( 6, 4)="3";
        (14, 4)="4";
        ( 9, 7)*[o]=<5pt>[Fo]{~}="x"; 
        ( 6,10)="6";
        (4.4,9.5)="7";
        "7";"3" **\crv{(2,9)&(0,13)&(6,13)&(5.8,6)};
        "1";"4" **\dir{-};
        "x";"3" **\dir{-};
        ( 6, 0);"3" **\dir{-};
        (14, 0);"4" **\dir{-};
        \ar@{-} "2";"6" |<(0.6)\hole
    \endxy}
    \overset{\text{(c)}}{=}~
    \vcenter{\xy 0;/r.15pc/:
        (0, 9)="1";
        (8,-9)="6";
        (8, 9)="2";
        (1.3,-1.5)="3";
        (0,-3)="4";
        (2,-9)="5";
        "1";"6" **\crv{(0,5)&(8,2)};
        "4";"5" **\crv{(-2,-6)&(-4,-3)&(-2,1)&(2,-3)};
        \ar@{-}@/^1pt/ "2";"3" |<(0.6)\hole
    \endxy}
    ~= \nu
\]

Thus, $R \ox R$ is a weak Hopf monoid with invertible antipode.

%=========================================================================%
\section{Quantum groupoids}\label{sec-quangroup}
%=========================================================================%

In this section we recall the quantum categories and quantum groupoids of
Day and Street~\cite{DS}. There is a succinct
definition given in~\cite[p.~216]{DS} in terms of ``basic data'' and
``Hopf basic data''. Here we give the unpacked definition of quantum
category and quantum groupoid which is essentially found
in~\cite[p.~221]{DS}; however, we do make a correction.

Our setting is a braided monoidal category $\V = (\V,\ox,I,c)$ in which the
functors
\[
    A \ox -:\V \ra \V
\]
with $A \in \V$, preserve coreflexive equalizers, i.e., equalizers of pairs
of morphisms with a common left inverse.

%=========================================================================%
\subsection{Quantum categories}
%=========================================================================%

Suppose $A$ and $C$ are comonoids in $\V$ and $s:A \ra C^\o$ and
$t:A \ra C$ are comonoid morphisms such that the diagram
\[
    \xygraph{{A}="s"
        :[u(0.7)r(1.2)] {A \ox A} ^-{\delta}
        :[r(1.8)] {C \ox C} ^-{s \ox t}
        :[d(1.4)] {C \ox C}="t" ^-c
     "s":[d(0.7)r(1.2)] {A \ox A} _-{\delta}
        :"t" ^-{t \ox s}}
\]
commutes. Then $A$ may be viewed as a $C$-bicomodule with left and right
coactions defined respectively via
\begin{align*}
    \gamma_l &= \big(\xygraph{{A}
        :[r(1.3)] {A \ox A} ^-\delta
        :[r(1.8)] {A \ox C} ^-{1 \ox s}
        :[r(1.7)] {C \ox A} ^-{c^{-1}}}\big) \\
    \gamma_r &= \big(\xygraph{{A}
        :[r(1.3)] {A \ox A} ^-\delta
        :[r(1.8)] {A \ox C} ^-{1 \ox t}}\big).
\end{align*}
Recall that the tensor product $P = A \ox_C A$ of $A$ with itself over $C$
is defined as the equalizer
\[
    \xygraph{{P}
        :[r(1.3)] {A \ox A} ^-\iota
        (:@<3pt>[r(2.3)] {A \ox C \ox A}="0" ^-{\gamma_r \ox 1}
        ):@<-3pt>"0" _-{1 \ox \gamma_l}}.
\]

The following diagrams may be seen to commute
\begin{align*}
    \xygraph{{P}
        :[r(1.3)] {A \ox A} ^-{\iota}
        :[r(2.2)] {C \ox A \ox A} ^-{\gamma_l \ox 1}
        (:@<3pt>[r(3.1)] {C \ox A \ox C \ox A}="0" ^-{1 \ox \gamma_r \ox 1}
        ):@<-3pt>"0" _-{1 \ox 1 \ox \gamma_l}}
\\
    \xygraph{{P}
        :[r(1.3)] {A \ox A} ^-{\iota}
        :[r(2.2)] {A \ox A \ox C} ^-{1 \ox \gamma_r}
        (:@<3pt>[r(3.1)] {A \ox C \ox A \ox C}="0" ^-{\gamma_r \ox 1 \ox 1}
        ):@<-3pt>"0" _-{1 \ox \gamma_l \ox 1}}
\end{align*}
and therefore induce respectively a left $C$- and right $C$-coaction on $P$.
These coactions make $P$ into a $C$-bicomodule.

The commutativity of the diagram
\[
\xygraph{{P}
    :[r(1.3)] {A \ox A} ^-{\iota}
    :[r(1.7)] {A^{\ox 4}} ^-{\delta \ox \delta}
    :[r(1.8)] {A^{\ox 4}} ^-{1 \ox c \ox 1}
    (:@<3pt>[r(3.2)] {A \ox A \ox A \ox C \ox A}="0"
                        ^-{1 \ox 1 \ox \gamma_r \ox 1}
    ):@<-3pt>"0" _-{1 \ox 1 \ox 1 \ox \gamma_l}}
\]
may be seen from
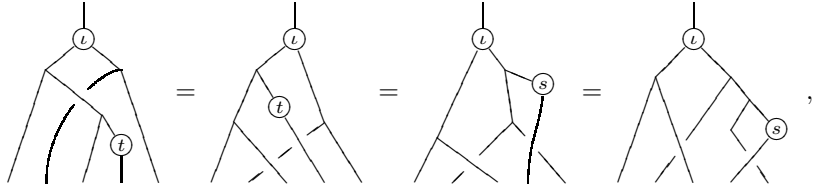
\[
    \vcenter{\xy
        (0,-13)*{~}; (0, 13)*{~};
        (-10,-12)="11";
        (-5,-12)="12";
        ( 0,-12)="13";
        ( 5,-12)="14";
        (10,-12)="15";
        (5,-7)*[o]=<9pt>[Fo]{\scriptstyle t}="2";
        (2.5,-3)="3";
        (-5,3)="41";
        (5,3)="42";
        (0,7)*[o]=<9pt>[Fo]{\scriptstyle \iota}="5";
        (0,12)="6";
        "14";"2" **\dir{-};
        "2";"3" **\dir{-};
        "13";"3" **\dir{-};
        "11";"41" **\dir{-};
        "3";"41" **\dir{-}?(0.4)*{\hole}="x"; 
        "12";"x" **\crv{(-5,-6)};
        "42";"x" **\crv{(3,3)};
        "15";"42" **\dir{-};
        "41";"5" **\dir{-};
        "42";"5" **\dir{-};
        "5";"6" **\dir{-};
    \endxy}
    ~=~
    \vcenter{\xy
        (-10,-12)="11";
        (-5,-12)="12";
        ( 0,-12)="13";
        ( 5,-12)="14";
        (10,-12)="15";
        (-7,-4)="21";
        (5,-4)="22";
        (-1,-2)*[o]=<9pt>[Fo]{\scriptstyle t}="3";
        (-4,3)="4";
        (1,7)*[o]=<9pt>[Fo]{\scriptstyle \iota}="5";
        (1,12)="6";
        "11";"21" **\dir{-};
        "13";"21" **\dir{-}?(0.3)*{\hole}="x"; % here
        "14";"3" **\dir{-}?(0.55)*{\hole}="y"; % and here
        "12";"x" **\dir{-}; "x";"y" **\dir{-}; "y";"22" **\dir{-};
        "15";"22" **\dir{-};
        "3";"4" **\dir{-};
        "21";"4" **\dir{-};
        "4";"5" **\dir{-};
        "22";"5" **\dir{-};
        "5";"6" **\dir{-};
    \endxy}
    ~=~
    \vcenter{\xy
        (-10,-12)="11";
        (-5,-12)="12";
        ( 1,-12)="13";
        ( 5,-12)="14";
        (10,-12)="15";
        (-7,-6)="21";
        (3,-4)="22";
        (7,1)*[o]=<9pt>[Fo]{\scriptstyle s}="3";
        (2,3)="4";
        (-1,7)*[o]=<9pt>[Fo]{\scriptstyle \iota}="5";
        (-1,12)="6";
        "11";"21" **\dir{-};
        "13";"21" **\dir{-}?(0.43)*{\hole}="x";
        "3";"14" **\crv{(7,-4)&(5,-6)}?(0.65)*{\hole}="y";
        "12";"x" **\dir{-}; "x";"22" **\dir{-};
        "15";"y" **\dir{-}; "y";"22" **\dir{-};
        "22";"4" **\dir{-};
        "3";"4" **\dir{-};
        "21";"5" **\dir{-};
        "4";"5" **\dir{-};
        "5";"6" **\dir{-};
    \endxy}
    ~=~
    \vcenter{\xy
        (-10,-12)="11";
        (-5,-12)="12";
        ( 0,-12)="13";
        ( 5,-12)="14";
        (10,-12)="15";
        ( 5,-5)="24";
        (11,-5)*[o]=<9pt>[Fo]{\scriptstyle s}="25";
        (7.5,-1)="3";
        (-5,2)="41";
        (5,2)="42";
        (0,7)*[o]=<9pt>[Fo]{\scriptstyle \iota}="5";
        (0,12)="6";
        "14";"25" **\dir{-}?(0.46)*{\hole}="x";
        "15";"x" **\dir{-}; "x";"24" **\dir{-};
        "24";"3" **\dir{-};
        "25";"3" **\dir{-}; 
        "11";"41" **\dir{-}; 
        "13";"41" **\dir{-}?(0.33)*{\hole}="y";
        "12";"y" **\dir{-}; "y";"42" **\dir{-};
        "3";"42" **\dir{-}; 
        "41";"5" **\dir{-};
        "42";"5" **\dir{-};
        "5";"6" **\dir{-};
    \endxy}\ \ ,
\]
and as $1 \ox 1 \ox \iota$ is the equalizer of $1 \ox 1 \ox \gamma_r \ox 1$
and $1 \ox 1 \ox 1 \ox \gamma_l$, there is a unique morphism
\[
    \delta_l:P \dra A \ox A \ox P
\]
making the diagram
\[
    \xygraph{{P}="s"
        :[r(1.3)] {A \ox A} ^-{\iota}
        :[r(2.2)] {A \ox A \ox A \ox A} ^-{\delta \ox \delta}
        :[d(1.5)] {A \ox A \ox A \ox A}="t" ^-{1 \ox c \ox 1}
     "s":[d(1.5)] {A \ox A \ox P} _-{\delta_l}
        :"t" _{1 \ox 1 \ox \iota}}
\]
commute. In strings,
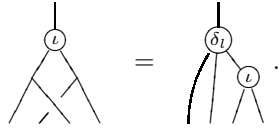
\[
    \vcenter{\xy
        (-6,-8)="11";
        (-2,-8)="12";
        ( 2,-8)="13";
        ( 6,-8)="14";
        (-3,-2)="21";
        ( 3,-2)="22";
        (0,3)*[o]=<9pt>[Fo]{\scriptstyle \iota}="3";
        (0,8)="4";
        "11";"21" **\dir{-};
        "13";"21" **\dir{-}?(0.4)*{\hole}="x";
        "12";"x" **\dir{-}; "x";"22" **\dir{-};
        "14";"22" **\dir{-};
        "21";"3" **\dir{-};
        "22";"3" **\dir{-};
        "3";"4" **\dir{-};
    \endxy}
    \quad = \quad
    \vcenter{\xy
        (-4,-8)="11";
        (-1,-8)="12";
        ( 2,-8)="13";
        ( 6,-8)="14";
        (4,-2)*[o]=<9pt>[Fo]{\scriptstyle \iota}="2";
        (0,3)*[o]=<10pt>[Fo]{\scriptstyle \delta_l}="3";
        (0,8)="4";
        "11";"3" **\crv{(-4,-2)};
        "12";"3" **\dir{-};
        "13";"2" **\dir{-};
        "14";"2" **\dir{-};
        "2";"3" **\dir{-};
        "3";"4" **\dir{-};
    \endxy}~.
\]

It is easy to see (postcompose with the monomorphism $1 \ox 1 \ox 1 \ox 1 \ox
\iota$) that the morphism $\delta_l$ is the left coaction of the comonoid
$A \ox A$ on $P$ making $P$ into a (left) $A \ox A$-comodule. This means
that the diagrams
\[
    \vcenter{\xygraph{{P}="s"
        :[d] {A \ox A \ox P} _-{\delta_l}
        :[d] {A \ox A \ox A \ox A \ox P} _-{\delta \ox \delta \ox 1}
        :[r(3.9)] {A \ox A \ox A \ox A \ox P}="t" ^-{1 \ox c \ox 1 \ox 1}
     "s":[r(3.9)] {A \ox A \ox P} ^-{\delta_l}
        :"t" ^-{1 \ox 1 \ox \delta_l}}}
    \qquad 
    \vcenter{\xygraph{{P}="s"
        :[r(1.6)] {A \ox A \ox P} ^-{\delta_l}
        :[d(1.3)] {P}="t" ^-{\epsilon \ox \epsilon \ox 1}
     "s":"t" _-1}}
\]
commute.

We are now ready to state the definition. A \emph{quantum category} in $\V$
consists of the data $\mathbf{A} = (A,C,s,t,\mu,\eta)$ where $A$, $C$, $s$,
$t$ are as above, and $\mu:P = A \ox_C A \ra A$ and $\eta:C \ra A$ are
morphisms in $\V$, called the \emph{composition morphism} and \emph{unit
morphism} respectively. This data must satisfy axioms (B1) through (B6) below.

\begin{itemize}
\item[(B1)] $(A,\mu,\eta)$ is a monoid in $\Bicomod(C)$.

\item[(B2)] The following diagram commutes.
\[
\xygraph{{P}="s"
    :[r(1.6)] {A \ox A \ox P} ^-{\delta_l}
    (:@<3pt>[r(2.3)] {C \ox P}="0" ^-{t \ox \epsilon \ox 1}
    ):@<-3pt>"0" _-{\epsilon \ox s \ox 1}
    :[r(1.8)] {C \ox A} ^-{1 \ox \mu}}
\]
\end{itemize}
Before stating (B3), we use (B2) to show that the diagram
\[
\xygraph{{P}
    :[r(1.6)] {A \ox A \ox P} ^-{\delta_l}
    :[r(2.6)] {A \ox A \ox A} ^-{1 \ox 1 \ox \mu} 
    (:@<3pt>[r(3)] {A \ox C \ox A \ox A}="0" ^-{\gamma_r \ox 1 \ox 1}
    ):@<-3pt>"0" _-{1 \ox \gamma_l \ox 1}}
\]
commutes, as seen by the calculation
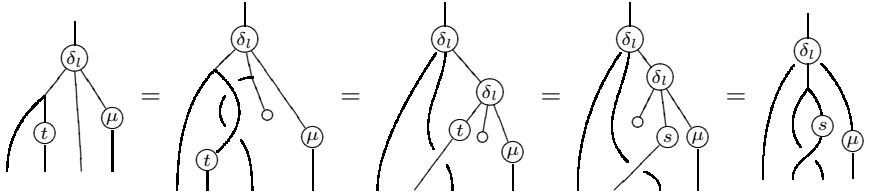
\[
    \vcenter{\xy
        (-9,-10)="11";
        (-4,-10)="12";
        ( 1,-10)="13";
        ( 5,-10)="14";
        (-4,-5)*[o]=<9pt>[Fo]{\scriptstyle t}="21";
        (5,-3)*[o]=<9pt>[Fo]{\scriptstyle \mu}="23";
        (-4,0)="31";
        (0,5)*[o]=<10pt>[Fo]{\scriptstyle \delta_l}="4";
        (0,10)="5";
        "12";"21" **\dir{-};
        "14";"23" **\dir{-};
        "11";"31" **\crv{(-9,-4)};
        "21";"31" **\dir{-};
        "13";"4" **\dir{-};
        "31";"4" **\dir{-};
        "23";"4" **\dir{-};
        "4";"5" **\dir{-};
    \endxy}
    ~=~
    \vcenter{\xy
        (-9,-15)="11";
        (-5,-15)="12";
        ( 1,-15)="13";
        ( 9,-15)="14";
        (-5,-11)*[o]=<9pt>[Fo]{\scriptstyle t}="21";
        (3,-5)*[o]=<5pt>[Fo]{~}="22";
        (9,-8)*[o]=<9pt>[Fo]{\scriptstyle \mu}="23";
        (-4,1)="31";
        (1,0)="32";
        (0,5)*[o]=<10pt>[Fo]{\scriptstyle \delta_l}="4";
        (0,10)="5";
        "21";"31" **\crv{(3,-5)}?(0.3)*{\hole}="y"?(0.8)*{\hole}="x";
        "13";"y" **\crv{(1,-9)};
        "y";"x" **\crv{(-5,-5)};
        "x";"32" **\crv{(-1,0)};
        "12";"21" **\dir{-};
        "14";"23" **\dir{-};
        "11";"31" **\crv{(-9,-4)};
        "22";"32" **\dir{-};
        "31";"4" **\dir{-};
        "32";"4" **\dir{-};
        "23";"4" **\dir{-};
        "4";"5" **\dir{-};
    \endxy}
    ~=~
    \vcenter{\xy
        (-9,-15)="11";
        (-4,-15)="12";
        ( 1,-15)="13";
        (5,-8)*[o]=<5pt>[Fo]{~}="14";
        (2,-7)*[o]=<9pt>[Fo]{\scriptstyle t}="t";
        ( 9,-15)="15";
        (9,-10)*[o]=<9pt>[Fo]{\scriptstyle \mu}="2";
        (6,-2)*[o]=<10pt>[Fo]{\scriptstyle \delta_l}="3";
        (0,5)*[o]=<10pt>[Fo]{\scriptstyle \delta_l}="4";
        (0,10)="5";
        "11";"4" **\crv{(-9,-8)};
        "12";"t" **\dir{-}?(0.55)*{\hole}="x";
        "13";"x" **\crv{(1,-12)};
        "x";"4" **\crv{(-5,-6)&(1,-0)};
        "14";"3" **\dir{-};
        "15";"2" **\dir{-};
        "2";"3" **\dir{-};
        "3";"4" **\dir{-};
        "4";"5" **\dir{-};
        "t";"3" **\dir{-};
    \endxy}
    ~=~
    \vcenter{\xy
        (-8,-15)="11";
        (-3,-15)="12";
        ( 3,-15)="13";
        ( 8,-15)="14";
        (0,-6)*[o]=<5pt>[Fo]{~}="o";
        (4,-8)*[o]=<9pt>[Fo]{\scriptstyle s}="s";
        (8,-8)*[o]=<9pt>[Fo]{\scriptstyle \mu}="d";
        (3,0)*[o]=<10pt>[Fo]{\scriptstyle \delta_l}="ml";
        (-1,5)*[o]=<10pt>[Fo]{\scriptstyle \delta_l}="m";
        (-1,10)="f";
        "12";"s" **\dir{-}?(0.4)*{\hole}="x";
        "13";"x" **\crv{(3,-13)};
        "m";"x" **\crv{(-1,-2)&(-7,-7)};
        "14";"d" **\dir{-};
        "o";"ml" **\dir{-};
        "s";"ml" **\dir{-};
        "d";"ml" **\dir{-};
        "ml";"m" **\dir{-};
        "11";"m" **\crv{(-8,-4)};
        "m";"f" **\dir{-};
    \endxy}
    ~=~
    \vcenter{\xy
        (-6,-17)="1";
        (-2,-17)="2";
        ( 2,-17)="3";
        ( 6,-17)="4";
        ( 6,-12)*[o]=<9pt>[Fo]{\scriptstyle \mu}="d";
        ( 2,-10)*[o]=<9pt>[Fo]{\scriptstyle s}="s";
        (0,-5)="m";
        (0,0)*[o]=<10pt>[Fo]{\scriptstyle \delta_l}="mu";
        (0,5)="b";
        "1";"mu" **\crv{(-6,-5)};
        "2";"s" **\crv{(-2,-14)&(2,-14)}?(0.5)*{\hole}="x";
        "3";"x" **\crv{(2,-15)};
        "x";"m" **\crv{(-4,-10)};
        "s";"m" **\crv{(2,-7)};
        "m";"mu" **\crv{-};
        "4";"d" **\dir{-};
        "d";"mu" **\crv{(6,-5)};
        "mu";"b" **\dir{-};
    \endxy}~.
\]
As $\iota \ox 1$ is the equalizer of $\gamma_r \ox 1 \ox 1$ and $1
\ox \gamma_l \ox 1$ there is a unique morphism $\delta_r:P \ra P \ox A$
making the square
\[
    \xygraph{{P}="s"
        :[r(2)] {A \ox A \ox P} ^-{\delta_l}
        :[d(1.2)] {A \ox A \ox A}="t" ^-{1 \ox 1 \ox \mu}
     "s":[d(1.2)] {P \ox A} _-{\delta_r}
        :"t" ^-{\iota \ox 1}}
\]
commute. We can now state (B3).

\begin{itemize}
\item[(B3)] The following diagram commutes.
\[
    \xygraph{{P}="s"
        :[r(1.9)] {A} ^-\mu
        :[d(1.2)] {A \ox A}="t" ^-\delta
     "s":[d(1.2)] {P \ox A} _-{\delta_r}
        :"t" ^-{\mu \ox 1}}
\]

\item[(B4)] The following diagram commutes.
\[
    \xygraph{{P}="s"
        :[r(1.5)] {A} ^-\mu
        :[d(1.2)] {I}="t" ^-\epsilon
     "s":[d(1.2)] {A \ox A} _-\iota
        :"t" ^-{\epsilon \ox \epsilon}}
\]

\item[(B5)] The following diagram commutes.
\[
    \xygraph{{C}="s"
        :[d(1.2)] {A} _-\eta
        :[r(1.2)u(0.6)] {I}="t" _-\epsilon
     "s":"t" ^-\epsilon}
\]

\item[(B6)] The following diagram commutes.
\[
    \xygraph{{C}="s"
        :[u(0.9)r(1)] {A} ^-\eta
        :[r(1.4)] {A \ox A} ^-\delta
        :[r(1.8)] {C \ox A} ^-{s \ox 1}
        :[d(0.9)r(1)] {A \ox A}="t" ^-{\eta \ox 1}
     "s":[d(0.9)r(1)] {A} _-\eta
        :[r(1.4)] {A \ox A} ^-\delta
        :[r(1.8)] {C \ox A} ^-{t \ox 1}
        :"t" _-{\eta \ox 1}
     "s":[r(2.6)] {A} ^-\eta
        :"t" ^-\delta}
\]
\end{itemize}

A consequence of these axioms is that $P$ becomes a left $A \ox A$-, right
$A$-bicomodule.

The axiom (B6) makes $C$ into a right $A$-comodule via
\[
    \xygraph{{C}
        :[r] {A} ^-\eta
        :[r(1.3)] {A \ox A} ^-\delta
        :[r(1.7)] {C \ox A} ^-{s \ox 1}}.
\]

We refer to $A$ as the \emph{object-of-arrows} and $C$ as the
\emph{object-of-objects}.

%=========================================================================%
\subsection{Quantum groupoids}\label{sec-quangpd}
%=========================================================================%

Suppose we have comonoid isomorphisms
\[
    \xymatrix{\upsilon:C^{\o \o} \ar[r]^-\cong & C}
    \qquad \text{and} \qquad
    \xymatrix{\nu:A^\o \ar[r]^-\cong & A}.
\]
Denote by $P_l$ the left $A^{\ox 3}$-comodule $P$ with coaction defined by
\[
    \xygraph{{P} 
        :[r(1.7)] {A \ox A \ox P \ox A} ^-{\delta}
        :[r(3.2)] {A \ox A \ox P \ox A} ^-{1 \ox 1 \ox 1 \ox \nu}
        :[r(3.2)] {A \ox A \ox A \ox P} ^-{1 \ox 1 \ox c_{P,A}}},
\]
and by $P_r$ the left $A^{\ox 3}$-comodule $P$ with coaction
defined by
\[
    \xygraph{{P} 
        :[r(1.7)] {A \ox A \ox P \ox A} ^-{\delta}
        :[r(3.2)] {A \ox A \ox P \ox A} ^-{1 \ox 1 \ox 1 \ox \nu^{-1}}
        :[r(3.2)] {A \ox A \ox A \ox P} ^-{c^{-1}_{A \ox A \ox P,A}}}.
\]
Furthermore, suppose that $\theta:P_l \ra P_r$ is a left
$A^{\ox 3}$-comodule isomorphism. We define a \emph{quantum groupoid} in
$\V$ to be a quantum category $\textbf{A}$ in $\V$ equipped with an
$\upsilon$, $\nu$, and $\theta$ satisfying (G1) through (G3) below.
\begin{itemize}
\item[(G1)] $s \nu = t$,
\item[(G2)] $t \nu = \upsilon s$, and
\item[(G3)] the diagram\footnote{This corrects~\cite[\S 12, p. 223]{DS}.}
\begin{equation*}
    \vcenter{\xygraph{{P}="s"
        :[r(1.6)] {C \ox C \ox C} ^-{\varsigma}
        :[r(2.4)] {C \ox C \ox C} ^-{c_{C,C \ox C}}
        :[d(1.2)] {C \ox C \ox C}="t" ^-{1 \ox 1 \ox \upsilon}
     "s":[d(1.2)] {P} _-{\theta}
        :"t" ^-\varsigma}}
\end{equation*}
commutes, where the morphism $\varsigma:P \ra C^{\ox 3}$ is defined by
taking either of the equal routes
\[
    \xygraph{{P}
        :[r(1.2)] {A \ox A} ^-{\iota}
        (:@<3pt>[r(2.1)] {A \ox C \ox A}="n" ^-{\gamma_r \ox 1}
        ):@<-3pt>"n" _-{1 \ox \gamma_l}
        :[r(2)] {C^{\ox 3}} ^-{s \ox 1 \ox t}}.
\]
\end{itemize}

%=========================================================================%
\section{Weak Hopf monoids are quantum groupoids}\label{sec-whqg}
%=========================================================================%

The goal of this section is to prove the following theorem.

\begin{theorem}\label{prop-whm-qg}
A weak bimonoid in $\Q\V$ is a quantum category in $\Q\V$ whose
object-of-objects is a separable Frobenius monoid. If the weak bimonoid is
equipped with an invertible antipode, making it a weak Hopf monoid, then the
quantum category becomes a quantum groupoid.
\end{theorem}

%========================================================================%
\subsection{Weak bimonoids are quantum categories}
%========================================================================%

Let $A=(A,1)$ be a weak bimonoid in $\Q\V$ with source morphism $s$ and
target morphism $t$ and set $C = (A,t)$. This data along with
\begin{align*}
    \mu &= \vcenter{\xy
        (0,5)*{~}; (0,-6)*{~};
        (-3, 4)="1";
        ( 3, 4)="2";
        ( 0, 0)="3";
        ( 0,-5)="4";
        "1";"3" **\dir{-};
        "2";"3" **\dir{-};
        "4";"3" **\dir{-};
    \endxy}
    ~:P \dra A
\\
    \eta &= t:C \dra A
\end{align*}
forms a quantum category in $\Q\V$. The morphisms $s$ and $t$ are obviously
in $\Q\V$, hence so is $\eta = t$, and
\[
    \vcenter{\xy
        (-3, 8)="1";
        ( 3, 8)="2";
        (-3, 3)="3";
        ( 3, 5)="4";
        ( 5, 3)="5";
        ( 0, 1)="6";
        ( 0,-2)*[o]=<5pt>[Fo]{~}="7";
        (-3,-2)="8";
        ( 3,-2)="9";
        ( 0,-6)="10";
        ( 0,-9)="11";
        "1";"8" **\dir{-};
        "3";"6" **\dir{-};
        "4";"5" **\dir{-};
        "5";"6" **\dir{-};
        "6";"7" **\dir{-};
        "8";"10" **\dir{-};
        "9";"10" **\dir{-};
        "11";"10" **\dir{-};
        \ar@{-} "2";"9" |<(0.6){\hole} 
    \endxy}
    ~\overset{\text{(nat)}}{=}~
    \vcenter{\xy
        (-3, 9)="1";
        ( 3, 9)="2";
        (-3, 6)="3";
        ( 3, 6)="4";
        (-3, 1)="5";
        ( 3, 1)="6";
        ( 3,-2)*[o]=<5pt>[Fo]{~}="7";
        (-3,-5)="8";
        "1";"8" **\dir{-};
        "2";"7" **\dir{-};
        "3";"6" **\dir{-};
        \ar@{-} "4";"5" |<(0.5){\hole} 
    \endxy}
    ~\overset{\text{(b)}}{=}~
    \vcenter{\xy
        (-3, 9)="1";
        ( 3, 9)="2";
        ( 0, 5)="3";
        ( 0, 0)="4";
        "1";"3" **\dir{-};
        "2";"3" **\dir{-};
        "4";"3" **\dir{-};
    \endxy}
\]
shows that $\mu$ is as well. Recall that $P = (A \ox A, m)$ where
\[
    m =~
    \vcenter{\xy
        (-3, 8)="1";
        ( 3, 8)="2";
        (-3, 3)="3";
        ( 3, 5)="4";
        ( 5, 3)="5";
        ( 0, 1)="6";
        ( 0,-2)*[o]=<5pt>[Fo]{~}="7";
        (-3,-4)="8";
        ( 3,-4)="9";
        "1";"8" **\dir{-};
        "3";"6" **\dir{-};
        "4";"5" **\dir{-};
        "5";"6" **\dir{-};
        "6";"7" **\dir{-};
        \ar@{-} "2";"9" |<(0.5){\hole} 
    \endxy}~.
\]
The morphisms $\delta_l:P \ra A \ox A \ox P$ and $\delta_r:P \ra P \ox A$
are given by
\[
    \delta_l =~
    \vcenter{\xy
        (0,9)*{~}; (0,-10)*{~};
        (-3, 8)="1";
        ( 3, 8)="2";
        (-3, 3)="3";
        ( 3, 5)="4";
        ( 5, 3)="5";
        ( 0, 1)="6";
        ( 0,-2)*[o]=<5pt>[Fo]{~}="7";
        (-3,-3)="8";
        ( 3,-3)="9";
        (-5,-9)="10";
        (-2,-9)="11";
        ( 2,-9)="12";
        ( 5,-9)="13";
        "1";"8" **\dir{-};
        "3";"6" **\dir{-};
        "4";"5" **\dir{-};
        "5";"6" **\dir{-};
        "6";"7" **\dir{-};
        "1";"8" **\dir{-};
        "3";"6" **\dir{-};
        "4";"5" **\dir{-};
        "5";"6" **\dir{-};
        "6";"7" **\dir{-};
        "8";"10" **\dir{-};
        "8";"12" **\dir{-};
        "9";"13" **\dir{-};
        \ar@{-} "2";"9" |<(0.55){\hole} 
        \ar@{-} "9";"11" |<(0.6){\hole} 
    \endxy}
\qquad\qquad
    \delta_r =~
    \vcenter{\xy
        (0,16)*{~}; (0,-5)*{~};
        (-3, 8)="1";
        ( 3, 5)="2";
        (-3, 3)="3";
        ( 3, 5)="4";
        ( 5, 3)="5";
        ( 0, 1)="6";
        ( 0,-2)*[o]=<5pt>[Fo]{~}="7";
        (-3,-4)="8";
        ( 3,-4)="9";
        (-3,15)="10";
        ( 6,15)="11";
        ( 6,12)="12";
        ( 6, 8)="13";
        ( 6,-4)="14";
        (-3,12)="15";
        "1";"8" **\dir{-};
        "3";"6" **\dir{-};
        "4";"5" **\dir{-};
        "5";"6" **\dir{-};
        "6";"7" **\dir{-};
        "11";"14" **\dir{-};
        "10";"1" **\dir{-};
        "15";"13" **\dir{-};
        \ar@{-} "2";"9" |<(0.35){\hole} 
        \ar@{-}@/^3pt/ "2";"12" |<(0.5){\hole} 
    \endxy}\ \ .
\]
The two calculations
\[
    \vcenter{\xy
        (0,9)*{~}; (0,-20)*{~};
        (-3, 8)="1";
        ( 3, 8)="2";
        (-3, 3)="3";
        ( 3, 5)="4";
        ( 5, 3)="5";
        ( 0, 1)="6";
        ( 0,-2)*[o]=<5pt>[Fo]{~}="7";
        (-3,-3)="8";
        ( 3,-3)="9";
        (-7,-9)="10";
        (-2,-9)="11";
        ( 1,-9)="12";
        ( 7,-9)="13";
        "1";"8" **\dir{-};
        "3";"6" **\dir{-};
        "4";"5" **\dir{-};
        "5";"6" **\dir{-};
        "6";"7" **\dir{-};
        "1";"8" **\dir{-};
        "3";"6" **\dir{-};
        "4";"5" **\dir{-};
        "5";"6" **\dir{-};
        "6";"7" **\dir{-};
        "8";"10" **\dir{-};
        "8";"12" **\dir{-};
        "9";"13" **\dir{-};
        ( 9,-11)="14";
        ( 4,-13)="15";
        ( 4,-16)*[o]=<5pt>[Fo]{~}="16";
        (-7,-19)="17";
        (-2,-19)="18";
        ( 1,-19)="19";
        ( 7,-19)="20";
        "10";"17" **\dir{-};
        "11";"18" **\dir{-};
        "12";"19" **\dir{-};
        (1,-11);"15" **\dir{-};
        "14";"13" **\dir{-};
        "14";"15" **\dir{-};
        "16";"15" **\dir{-};
        \ar@{-} "13";"20" |<(0.3){\hole} 
        \ar@{-} "2";"9" |<(0.55){\hole} 
        \ar@{-} "9";"11" |<(0.6){\hole} 
    \endxy}
    ~\overset{\text{(c)}}{=}~
    \vcenter{\xy
        (-3,15)="1";
        ( 3,15)="2";
        (-3,10)="a";
        ( 3,12)="b";
        ( 5,10)="x";
        ( 0,8)="y";
        ( 0,5)*[o]=<5pt>[Fo]{~}="z";
        "b";"x" **\dir{-};
        "x";"y" **\dir{-};
        "a";"y" **\dir{-};
        "y";"z" **\dir{-};
        (-3, 3)="3";
        ( 3, 5)="4";
        ( 5, 3)="5";
        ( 0, 1)="6";
        ( 0,-2)*[o]=<5pt>[Fo]{~}="7";
        (-3,-3)="8";
        ( 3,-3)="9";
        (-7,-9)="10";
        (-2,-9)="11";
        ( 1,-9)="12";
        ( 7,-9)="13";
        "1";"8" **\dir{-};
        "3";"6" **\dir{-};
        "4";"5" **\dir{-};
        "5";"6" **\dir{-};
        "6";"7" **\dir{-};
        "1";"8" **\dir{-};
        "3";"6" **\dir{-};
        "4";"5" **\dir{-};
        "5";"6" **\dir{-};
        "6";"7" **\dir{-};
        "8";"10" **\dir{-};
        "8";"12" **\dir{-};
        "9";"13" **\dir{-};
        \ar@{-} "2";"9" |<(0.33){\hole}  |<(0.72){\hole} 
        \ar@{-} "9";"11" |<(0.6){\hole} 
    \endxy}
    ~=~
    \vcenter{\xy
        (-3, 8)="1";
        ( 3, 8)="2";
        (-3, 3)="3";
        ( 3, 5)="4";
        ( 5, 3)="5";
        ( 0, 1)="6";
        ( 0,-2)*[o]=<5pt>[Fo]{~}="7";
        (-3,-3)="8";
        ( 3,-3)="9";
        (-7,-9)="10";
        (-2,-9)="11";
        ( 1,-9)="12";
        ( 7,-9)="13";
        "1";"8" **\dir{-};
        "3";"6" **\dir{-};
        "4";"5" **\dir{-};
        "5";"6" **\dir{-};
        "6";"7" **\dir{-};
        "1";"8" **\dir{-};
        "3";"6" **\dir{-};
        "4";"5" **\dir{-};
        "5";"6" **\dir{-};
        "6";"7" **\dir{-};
        "8";"10" **\dir{-};
        "8";"12" **\dir{-};
        "9";"13" **\dir{-};
        \ar@{-} "2";"9" |<(0.55){\hole} 
        \ar@{-} "9";"11" |<(0.6){\hole} 
    \endxy}
\]
and
\[
    \vcenter{\xy
        (-3,23)="16";
        ( 6,23)="17";
        (-3,18)="18";
        ( 6,20)="19";
        ( 8,18)="20";
        (1.5,16)="21";
        (1.5,13)*[o]=<5pt>[Fo]{~}="22";
        "18";"21" **\dir{-};
        "20";"21" **\dir{-};
        "19";"20" **\dir{-};
        "21";"22" **\dir{-};
        (-3, 8)="1";
        ( 3, 5)="2";
        (-3, 3)="3";
        ( 3, 5)="4";
        ( 5, 3)="5";
        ( 0, 1)="6";
        ( 0,-2)*[o]=<5pt>[Fo]{~}="7";
        (-3,-4)="8";
        ( 3,-4)="9";
        (-3,15)="10";
        ( 6,15)="11";
        ( 6,12)="12";
        ( 6, 8)="13";
        ( 6,-4)="14";
        (-3,12)="15";
        "16";"10" **\dir{-};
        "1";"8" **\dir{-};
        "3";"6" **\dir{-};
        "4";"5" **\dir{-};
        "5";"6" **\dir{-};
        "6";"7" **\dir{-};
        "11";"14" **\dir{-};
        "10";"1" **\dir{-};
        "15";"13" **\dir{-};
        \ar@{-} "17";"11" |<(0.72){\hole} 
        \ar@{-} "2";"9" |<(0.35){\hole} 
        \ar@{-}@/^3pt/ "2";"12" |<(0.5){\hole} 
    \endxy}
    ~\overset{\text{(c)}}{=}~
    \vcenter{\xy
        (-6,23)="1";
        ( 9,23)="2";
        (-6,20)="3";
        ( 9,20)="4";
        ( 3,17)="5";
        ( 9,17)="6";
        ( 3,12)="7";
        ( 9,12)="8";
        ( 9, 9)*[o]=<5pt>[Fo]{~}="9";
        (-6,10)="10";
        ( 0,12)="11";
        ( 2,10)="12";
        (-3, 8)="13";
        (-3, 5)*[o]=<5pt>[Fo]{~}="14";
        (-6, 2)="15";
        ( 0, 2)="16";
        ( 3, 2)="17";
        "1";"15" **\dir{-};
        "2";"9" **\dir{-};
        "3";"5" **\dir{-};
        "5";"17" **\dir{-};
        "5";"8" **\dir{-};
        "10";"13" **\dir{-};
        "14";"13" **\dir{-};
        "12";"13" **\dir{-};
        "12";"11" **\dir{-};
        (0,15);(0,12) **\dir{-};
        \ar@{-} "11";"16" |<(0.3)\hole
        \ar@{-} "6";"7" |\hole 
        \ar@{-}@/^3pt/ (0,15);"4" |<(0.3)\hole
    \endxy}
    ~\overset{\text{(b)}}{=}~
    \vcenter{\xy
        (0,16)*{~}; (0,-5)*{~};
        (-3, 8)="1";
        ( 3, 5)="2";
        (-3, 3)="3";
        ( 3, 5)="4";
        ( 5, 3)="5";
        ( 0, 1)="6";
        ( 0,-2)*[o]=<5pt>[Fo]{~}="7";
        (-3,-4)="8";
        ( 3,-4)="9";
        (-3,15)="10";
        ( 6,15)="11";
        ( 6,12)="12";
        ( 6, 8)="13";
        ( 6,-4)="14";
        (-3,12)="15";
        "1";"8" **\dir{-};
        "3";"6" **\dir{-};
        "4";"5" **\dir{-};
        "5";"6" **\dir{-};
        "6";"7" **\dir{-};
        "11";"14" **\dir{-};
        "10";"1" **\dir{-};
        "15";"13" **\dir{-};
        \ar@{-} "2";"9" |<(0.35){\hole} 
        \ar@{-}@/^3pt/ "2";"12" |<(0.5){\hole} 
    \endxy}
\]
show that these are morphisms in $\Q\V$.

To see that $(A,\mu,\eta)$ is a comonoid in $\Bicomod(C)$ notice that
associativity follows from the associativity of the $\mu$ viewed as a weak
bimonoid and the counit property may be seen from property (6), i.e.,
\[
    \mu(1 \ox t)\delta = 1_A
    \qquad\text{and}\qquad
    \mu(s \ox 1)c^{-1}\delta = 1_A,
\]
and so (B1) holds. (B2) follows from one application of (12), (B3) from (b),
(B4) from (c), and (B5) from (2), while the calculation
\[
    \vcenter{\xy
        (0,8)="1";
        (0,4)*[o]=<9pt>[Fo]{\scriptstyle t}="2";
        (0,0)="3";
        (-3,-3)*[o]=<9pt>[Fo]{\scriptstyle s}="4";
        (-3,-8)*[o]=<9pt>[Fo]{\scriptstyle t}="5";
        (-3,-12)="6";
        ( 3,-12)="7";
        "1";"2" **\dir{-};
        "3";"2" **\dir{-};
        "3";"4" **\dir{-};
        "4";"5" **\dir{-};
        "6";"5" **\dir{-};
        "7";"3" **\crv{(3,-5)};
    \endxy}
    ~\overset{\text{(3)}}{=}~
    \vcenter{\xy
        (0,9)*{~}; (0,-18)*{~};
        (0,8)="1";
        (0,4)*[o]=<9pt>[Fo]{\scriptstyle t}="2";
        (0,0)="3";
        (-3,-3)*[o]=<9pt>[Fo]{\scriptstyle t}="4";
        (-3,-8)*[o]=<9pt>[Fo]{\scriptstyle s}="5";
        (-3,-13)*[o]=<9pt>[Fo]{\scriptstyle t}="x";
        (-3,-17)="6";
        ( 3,-17)="7";
        "1";"2" **\dir{-};
        "3";"2" **\dir{-};
        "3";"4" **\dir{-};
        "4";"5" **\dir{-};
        "x";"5" **\dir{-};
        "6";"x" **\dir{-};
        "7";"3" **\crv{(3,-5)};
    \endxy}
    ~\overset{\text{(8)}}{=}~
    \vcenter{\xy
        (0,8)="1";
        (0,4)*[o]=<9pt>[Fo]{\scriptstyle t}="2";
        (0,0)="3";
        (-3,-3)*[o]=<9pt>[Fo]{\scriptstyle t}="4";
        (-3,-8)*[o]=<9pt>[Fo]{\scriptstyle t}="5";
        (-3,-12)="6";
        ( 3,-12)="7";
        "1";"2" **\dir{-};
        "3";"2" **\dir{-};
        "3";"4" **\dir{-};
        "4";"5" **\dir{-};
        "6";"5" **\dir{-};
        "7";"3" **\crv{(3,-5)};
    \endxy}
\]
verifies (B6). Thus, $\mathbf{A} = (A,C,s,t,\mu,\eta)$ is a quantum category
in $\Q\V$.

%========================================================================%
\subsection{Weak Hopf monoids are quantum groupoids}
%========================================================================%

Now suppose that $A = (A,1)$ is a weak Hopf monoid in $\Q\V$ with an
invertible antipode $\nu:A \ra A$, and that $\mathbf{A} = (A,C,s,t,\mu,\eta)$
is as above. The data for a quantum groupoid $(\upsilon,\nu,\theta)$ is
\begin{align*}
    \upsilon &= t \nu \nu t:C^{\o\o} \dra C \\
    \nu &= \nu: A^\o \dra A \\
    \theta &=~
    \vcenter{\xy
        (-3, 5)="1";
        ( 3, 5)="2";
        ( 3, 2)="3";
        ( 3,-2)="4";
        ( 3,-6)*[o]=<9pt>[Fo]{\scriptstyle \nu}="5";
        (-3,-10)="6";
        ( 3,-10)="7";
        "1";"4" **\crv{(-3,1)};
        "2";"5" **\dir{-};
        "5";"7" **\dir{-};
        "6";(-3,-4) **\dir{-}
        \ar@{-}@/^4pt/ (-3,-4);"3" |<(0.55)\hole
    \endxy} 
    ~:P \dra P.
\end{align*}
In the remainder of this section we will verify this claim.

The morphisms $\upsilon$ and $\nu$ are obviously morphisms in $\Q\V$, and
the two calculations
\[
    \vcenter{\xy
        (-3,15)="8";
        ( 3,15)="9";
        (-3,10)="10";
        ( 3,12)="11";
        ( 5,10)="12";
        ( 0, 8)="13";
        ( 0, 5)*[o]=<5pt>[Fo]{~}="14";
        (-3, 5)="1";
        ( 3, 5)="2";
        ( 3, 2)="3";
        ( 3,-2)="4";
        ( 3,-6)*[o]=<9pt>[Fo]{\scriptstyle \nu}="5";
        (-3,-10)="6";
        ( 3,-10)="7";
        "8";"1" **\dir{-};
        "10";"13" **\dir{-};
        "11";"12" **\dir{-};
        "13";"12" **\dir{-};
        "13";"14" **\dir{-};
        "1";"4" **\crv{(-3,1)};
        "2";"5" **\dir{-};
        "5";"7" **\dir{-};
        "6";(-3,-4) **\dir{-}
        \ar@{-}@/^4pt/ (-3,-4);"3" |<(0.55)\hole
        \ar@{-} "9";"2" |<(0.6)\hole
    \endxy}
    ~\overset{\text{(c)}}{=}~
    \vcenter{\xy
        (-3, 6)="1";
        ( 3, 6)="2";
        ( 3, 3)="3";
        (-3,-2)="4";
        ( 3, 0)="5";
        ( 7,-6)="6";
        ( 3,-9)="7";
        ( 7,-9)*[o]=<5pt>[Fo]{~}="8";
        ( 3,-12)*[o]=<9pt>[Fo]{\scriptstyle \nu}="9";
        (-6,-13)="10";
        ( 3,-16)="11";
        (-6,-19)="12";
        ( 3,-19)="13";
        "1";"4" **\dir{-};
        "6";"4" **\dir{-};
        "7";"4" **\dir{-};
        "5";"6" **\dir{-};
        "8";"6" **\dir{-};
        "9";"13" **\dir{-};
        "10";"12" **\dir{-};
        "10";"11" **\dir{-};
        "10";(-6,-2) **\dir{-};
        \ar@{-}@/^2pt/ (-6,-2);"3" |<(0.35)\hole
        \ar@{-} "2";"9" |<(0.58)\hole
    \endxy}
    ~\overset{\text{(b)}}{=}~
    \vcenter{\xy
        (-3, 5)="1";
        ( 3, 5)="2";
        ( 3, 2)="3";
        ( 3,-2)="4";
        ( 3,-6)*[o]=<9pt>[Fo]{\scriptstyle \nu}="5";
        (-3,-10)="6";
        ( 3,-10)="7";
        "1";"4" **\crv{(-3,1)};
        "2";"5" **\dir{-};
        "5";"7" **\dir{-};
        "6";(-3,-4) **\dir{-}
        \ar@{-}@/^4pt/ (-3,-4);"3" |<(0.55)\hole
    \endxy}
\]
and
\[
    \vcenter{\xy
        (-3, 5)="1"; ( 3, 5)="2"; ( 3, 2)="3"; ( 3,-2)="4";
        ( 3,-6)*[o]=<9pt>[Fo]{\scriptstyle \nu}="5";
        (-3,-12)="6"; ( 3,-10)="7"; ( 5,-12)="8"; ( 0,-14)="9";
        ( 0,-17)*[o]=<5pt>[Fo]{~}="10";
        (-3,-20)="11";
        ( 3,-20)="12";
        "1";"4" **\crv{(-3,1)}; "2";"5" **\dir{-}; "5";"7" **\dir{-};
        "6";(-3,-4) **\dir{-}; "6";"11" **\dir{-}; "6";"9" **\dir{-};
        "10";"9" **\dir{-}; "8";"9" **\dir{-}; "8";"7" **\dir{-};
        \ar@{-}@/^4pt/ (-3,-4);"3" |<(0.55)\hole
        \ar@{-} "7";"12" |<(0.3)\hole
    \endxy}
    ~\overset{\text{(17)}}{=}~
    \vcenter{\xy
        (-3,12)="1"; ( 3,12)="2"; ( 3, 9)="3"; ( 3, 5)="4";
        (-3,-2)="5"; ( 3, 2)="6";
        ( 0, 0)*[o]=<9pt>[Fo]{\scriptstyle \nu}="7";
        ( 0,-4)="8";
        ( 0,-7)*[o]=<5pt>[Fo]{~}="9";
        ( 3,-5)*[o]=<9pt>[Fo]{\scriptstyle \nu}="10";
        (-3,-10)="11"; ( 3,-10)="12";
        "1";"4" **\crv{(-3,8)}; "2";"10" **\dir{-}; "12";"10" **\dir{-};
        "6";"7" **\dir{-}; "5";"11" **\dir{-}; "5";"8" **\dir{-};
        "7";"8" **\dir{-}; "9";"8" **\dir{-};
        "5";(-3,3) **\dir{-};
        \ar@{-}@/^4pt/ (-3,3);"3" |<(0.55)\hole
    \endxy}
    ~\overset{\text{(2)}}{=}~
    \vcenter{\xy
        (-3,12)="1"; ( 3,12)="2"; ( 3, 9)="3";
        ( 3, 5)="4"; (-3,-7)="5"; ( 3, 2)="6";
        ( 0, 0)*[o]=<9pt>[Fo]{\scriptstyle \nu}="x";
        ( 0,-5)*[o]=<9pt>[Fo]{\scriptstyle r}="7";
        ( 0,-9)="8";
        ( 0,-12)*[o]=<5pt>[Fo]{~}="9";
        ( 3,-10)*[o]=<9pt>[Fo]{\scriptstyle \nu}="10";
        (-3,-15)="11";
        ( 3,-15)="12";
        "1";"4" **\crv{(-3,8)}; "2";"10" **\dir{-}; "12";"10" **\dir{-};
        "6";"x" **\dir{-}; "x";"7" **\dir{-}; "5";"11" **\dir{-};
        "5";"8" **\dir{-}; "7";"8" **\dir{-}; "9";"8" **\dir{-};
        "5";(-3,3) **\dir{-};
        \ar@{-}@/^4pt/ (-3,3);"3" |<(0.55)\hole
    \endxy}
    ~\overset{\text{(16)}}{=}~
    \vcenter{\xy
        (-3,12)="1"; ( 3,12)="2"; ( 3, 9)="3";
        ( 3, 5)="4"; (-3,-7)="5"; ( 3, 2)="6";
        ( 0, 0)*[o]=<9pt>[Fo]{\scriptstyle t}="x";
        ( 0,-5)*[o]=<9pt>[Fo]{\scriptstyle \nu}="7";
        ( 0,-9)="8";
        ( 0,-12)*[o]=<5pt>[Fo]{~}="9";
        ( 3,-10)*[o]=<9pt>[Fo]{\scriptstyle \nu}="10";
        (-3,-15)="11";
        ( 3,-15)="12";
        "1";"4" **\crv{(-3,8)}; "2";"10" **\dir{-}; "12";"10" **\dir{-};
        "6";"x" **\dir{-}; "x";"7" **\dir{-}; "5";"11" **\dir{-};
        "5";"8" **\dir{-}; "7";"8" **\dir{-}; "9";"8" **\dir{-};
        "5";(-3,3) **\dir{-};
        \ar@{-}@/^4pt/ (-3,3);"3" |<(0.55)\hole
    \endxy}
    ~\overset{\text{(4)}}{=}~
    \vcenter{\xy
        (-3,14)="1"; ( 3,14)="2"; ( 3,11)="3";
        ( 3, 7)="4"; (-3,-7)="5"; ( 3, 3)="6";
        ( 0, 0)*[o]=<9pt>[Fo]{\scriptstyle t}="x";
        ( 0,-5)*[o]=<9pt>[Fo]{\scriptstyle \nu}="7";
        ( 0,-9)="8";
        ( 0,-12)*[o]=<5pt>[Fo]{~}="9";
        ( 3,-10)*[o]=<9pt>[Fo]{\scriptstyle \nu}="10";
        (-3,-15)="11";
        ( 3,-15)="12";
        "1";"6" **\crv{(-3,9)}; "2";"10" **\dir{-}; "12";"10" **\dir{-};
        "x";"7" **\dir{-}; "5";"11" **\dir{-}; "5";"8" **\dir{-};
        "7";"8" **\dir{-}; "9";"8" **\dir{-}; "5";(-3,3) **\dir{-};
        \ar@{-}@/^4pt/ (-3,3);"3" |<(0.5)\hole
        \ar@{-}@/^4pt/ "x";"4" |<(0.55)\hole
    \endxy}
\]
\[
    ~\overset{\text{(16)}}{=}~
    \vcenter{\xy
        (-3,14)="1"; ( 3,14)="2"; ( 3,11)="3";
        ( 3, 7)="4"; (-3,-7)="5"; ( 3, 3)="6";
        ( 0, 0)*[o]=<9pt>[Fo]{\scriptstyle \nu}="x";
        ( 0,-5)*[o]=<9pt>[Fo]{\scriptstyle r}="7";
        ( 0,-9)="8";
        ( 0,-12)*[o]=<5pt>[Fo]{~}="9";
        ( 3,-10)*[o]=<9pt>[Fo]{\scriptstyle \nu}="10";
        (-3,-15)="11";
        ( 3,-15)="12";
        "1";"6" **\crv{(-3,9)}; "2";"10" **\dir{-}; "12";"10" **\dir{-};
        "x";"7" **\dir{-}; "5";"11" **\dir{-}; "5";"8" **\dir{-};
        "7";"8" **\dir{-}; "9";"8" **\dir{-}; "5";(-3,3) **\dir{-};
        \ar@{-}@/^4pt/ (-3,3);"3" |<(0.5)\hole
        \ar@{-}@/^4pt/ "x";"4" |<(0.55)\hole
    \endxy}
    ~\overset{\text{(2)}}{=}~
    \vcenter{\xy
        (-3,14)="1";
        ( 3,14)="2";
        ( 3,11)="3";
        ( 3, 7)="4";
        (-3,-2)="5";
        ( 3, 3)="6";
        ( 0, 0)*[o]=<9pt>[Fo]{\scriptstyle \nu}="7";
        ( 0,-4)="8";
        ( 0,-7)*[o]=<5pt>[Fo]{~}="9";
        ( 3,-5)*[o]=<9pt>[Fo]{\scriptstyle \nu}="10";
        (-3,-10)="11";
        ( 3,-10)="12";
        "1";"6" **\crv{(-3,9)};
        "2";"10" **\dir{-};
        "12";"10" **\dir{-};
        "5";"11" **\dir{-};
        "5";"8" **\dir{-};
        "7";"8" **\dir{-};
        "9";"8" **\dir{-};
        "5";(-3,3) **\dir{-};
        \ar@{-}@/^4pt/ (-3,3);"3" |<(0.5)\hole
        \ar@{-}@/^4pt/ "7";"4" |<(0.55)\hole
    \endxy}
    ~\overset{\text{(c)}}{=}~
    \vcenter{\xy
        (-4,14)="1";
        ( 4,14)="2";
        ( 4,11)="3";
        ( 4, 7)="4";
        ( 4, 3)="5";
        ( 0, 1)="6";
        (-2,-2)="7";
        ( 2,-2)*[o]=<9pt>[Fo]{\scriptstyle \nu}="8";
        ( 0,-5)="9";
        ( 0,-8)*[o]=<5pt>[Fo]{~}="10";
        ( 4,-6)*[o]=<9pt>[Fo]{\scriptstyle \nu}="11";
        (-4,-10)="12";
        ( 4,-10)="13";
        "1";"5" **\crv{(-3,7)};
        "2";"11" **\dir{-};
        "13";"11" **\dir{-};
        "6";"7" **\dir{-};
        "6";"8" **\dir{-};
        "9";"7" **\dir{-};
        "9";"8" **\dir{-};
        "9";"10" **\dir{-};
        "12";(-4,3) **\dir{-};
        \ar@{-}@/^5pt/ (-4,3);"3" |<(0.45)\hole
        \ar@{-}@/^3pt/ "6";"4" |<(0.5)\hole
    \endxy}
    ~\overset{\text{($\nu$)}}{=}~
    \vcenter{\xy
        (-3,14)="1";
        ( 3,14)="2";
        ( 3,11)="3";
        ( 3, 7)="4";
        ( 3, 3)="5";
        ( 0, 0)*[o]=<9pt>[Fo]{\scriptstyle r}="6";
        ( 0,-4)*[o]=<5pt>[Fo]{~}="7";
        ( 3,-5)*[o]=<9pt>[Fo]{\scriptstyle \nu}="11";
        (-3,-10)="12";
        ( 3,-10)="13";
        "1";"5" **\crv{(-3,7)};
        "2";"11" **\dir{-};
        "13";"11" **\dir{-};
        "6";"7" **\dir{-};
        "12";(-3,3) **\dir{-};
        \ar@{-}@/^4pt/ (-3,3);"3" |<(0.45)\hole
        \ar@{-}@/^3pt/ "6";"4" |<(0.5)\hole
    \endxy}
    ~\overset{\text{(2,c)}}{=}~
    \vcenter{\xy
        (-3, 5)="1";
        ( 3, 5)="2";
        ( 3, 2)="3";
        ( 3,-2)="4";
        ( 3,-6)*[o]=<9pt>[Fo]{\scriptstyle \nu}="5";
        (-3,-10)="6";
        ( 3,-10)="7";
        "1";"4" **\crv{(-3,1)};
        "2";"5" **\dir{-};
        "5";"7" **\dir{-};
        "6";(-3,-4) **\dir{-}
        \ar@{-}@/^4pt/ (-3,-4);"3" |<(0.55)\hole
    \endxy}
\]
show that $\theta$ is as well.

\begin{lemma}
An inverse for $\theta$ is given by
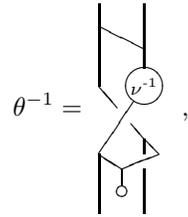
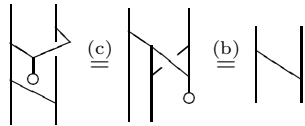
\[
    \theta^{-1} =~
    \vcenter{\xy
        (-3,14)="1";
        ( 3,14)="2";
        (-3,11)="3";
        ( 3, 8)="4";
        (-3, 3)="5";
        ( 3, 3)*[o]=<14pt>[Fo]{\scriptstyle \nu^{\text{-}1}}="6";
        (-3,-6)="7";
        ( 3,-4)="8";
        ( 5,-6)="9";
        ( 0,-8)="10";
        ( 0,-11)*[o]=<5pt>[Fo]{~}="11";
        (-3,-14)="12";
        ( 3,-14)="13";
        "1";"5" **\dir{-};
        "3";"4" **\dir{-};
        "2";"6" **\dir{-};
        "6";"7" **\dir{-};
        "10";"7" **\dir{-};
        "10";"11" **\dir{-};
        "8";"9" **\dir{-};
        "9";"10" **\dir{-};
        "7";"12" **\dir{-};
        \ar@{-} "5";"8" |<(0.56)\hole
        \ar@{-} "8";"13" |<(0.3)\hole
    \endxy}\ \ ,
\]
\end{lemma}

\begin{proof}
Since
\[
    \vcenter{\xy
        (-3,20)="1";
        ( 3,20)="2";
        (-3,15)="3";
        ( 3,17)="4";
        ( 5,15)="5";
        ( 0,13)="6";
        ( 0,10)*[o]=<5pt>[Fo]{~}="7";
        (-3,10)="8";
        ( 3,7)="9";
        (-3,4)="10";
        ( 3,4)="11";
        "1";"10" **\dir{-};
        "4";"5" **\dir{-};
        "6";"5" **\dir{-};
        "6";"3" **\dir{-};
        "6";"7" **\dir{-};
        "8";"9" **\dir{-};
        \ar@{-} "2";"11" |<(0.37)\hole
    \endxy}
    ~\overset{\text{(c)}}{=}~
    \vcenter{\xy
        (-4,20)="1";
        ( 4,20)="2";
        (-4,17)="3";
        (-1,15)="4";
        ( 4,15)="5";
        (-1,11)="6";
        ( 4,11)="7";
        ( 4,8)*[o]=<5pt>[Fo]{~}="8";
        (-4,5)="9";
        (-1,5)="10";
        "1";"9" **\dir{-};
        "2";"8" **\dir{-};
        "3";"4" **\dir{-};
        "4";"10" **\dir{-};
        "4";"7" **\dir{-};
        \ar@{-} "5";"6" |\hole
    \endxy}
    ~\overset{\text{(b)}}{=}~
    \vcenter{\xy
        (-3, 5)="1";
        ( 3, 5)="2";
        (-3, 2)="3";
        ( 3,-2)="4";
        (-3,-5)="5";
        ( 3,-5)="6";
        "1";"5" **\dir{-};
        "2";"6" **\dir{-};
        "3";"4" **\dir{-};
    \endxy}
\]
it is clear that $\theta^{-1}$ is a morphism in $\Q\V$.

That $\theta^{-1}$ is an inverse for $\theta$ may be seen in one direction
from
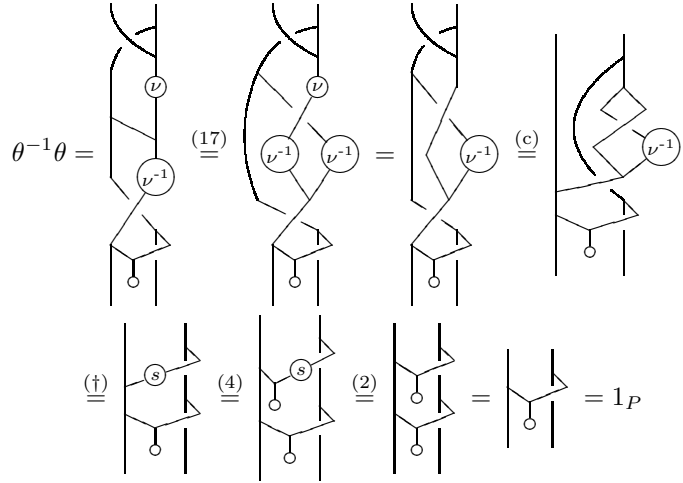
\begin{align*}
\theta^{-1} \theta &=~
    \vcenter{\xy
        (-3,25)="1";
        ( 3,25)="2";
        ( 3,22)="3";
        ( 3,18)="4";
        ( 3,14)*[o]=<9pt>[Fo]{\scriptstyle \nu}="5";
        (-3,10)="6";
        ( 3,10)="7";
        "1";"4" **\crv{(-3,21)};
        "2";"5" **\dir{-};
        "5";"7" **\dir{-};
        "6";(-3,16) **\dir{-};
        (-3,10)="a";
        ( 3, 7)="b";
        (-3, 2)="c";
        ( 3, 2)*[o]=<14pt>[Fo]{\scriptstyle \nu^{\text{-}1}}="d";
        (-3,-7)="e";
        ( 3,-5)="f";
        ( 5,-7)="g";
        ( 0,-9)="h";
        ( 0,-12)*[o]=<5pt>[Fo]{~}="i";
        (-3,-15)="j";
        ( 3,-15)="k";
        "a";"c" **\dir{-};
        "a";"b" **\dir{-};
        "b";"d" **\dir{-};
        "d";"e" **\dir{-};
        "h";"e" **\dir{-};
        "h";"i" **\dir{-};
        "f";"g" **\dir{-};
        "g";"h" **\dir{-};
        "e";"j" **\dir{-};
        "7";"b" **\dir{-};
        \ar@{-}@/^4pt/ (-3,16);"3" |<(0.55)\hole
        \ar@{-} "c";"f" |<(0.56)\hole
        \ar@{-} "f";"k" |<(0.3)\hole
    \endxy}
    ~\overset{\text{(17)}}{=}~
    \vcenter{\xy
        (-3,25)="1";
        ( 3,25)="2";
        ( 3,22)="3";
        ( 3,18)="4";
        ( 3,14)*[o]=<9pt>[Fo]{\scriptstyle \nu}="5";
        (-5,16)="6";
        "1";"4" **\crv{(-3,21)};
        "2";"5" **\dir{-};
        (-2, 5)*[o]=<14pt>[Fo]{\scriptstyle \nu^{\text{-}1}}="x";
        ( 6, 5)*[o]=<14pt>[Fo]{\scriptstyle \nu^{\text{-}1}}="y";
        (-5,-1)="c";
        ( 2,-1)="z";
        (-3,-7)="e";
        ( 3,-5)="f";
        ( 5,-7)="g";
        ( 0,-9)="h";
        ( 0,-12)*[o]=<5pt>[Fo]{~}="i";
        (-3,-15)="j";
        ( 3,-15)="k";
        "5";"x" **\dir{-};
        "z";"x" **\dir{-};
        "z";"y" **\dir{-};
        "e";"z" **\dir{-};
        "h";"e" **\dir{-};
        "h";"i" **\dir{-};
        "f";"g" **\dir{-};
        "g";"h" **\dir{-};
        "e";"j" **\dir{-};
        \ar@{-} "6";"y" |<(0.52)\hole
        \ar@{-}@/^4pt/ "6";"3" |<(0.6)\hole
        \ar@{-} "c";"f" |<(0.6)\hole
        \ar@{-}@/^5pt/ "c";"6"
        \ar@{-} "f";"k" |<(0.3)\hole
    \endxy}
    ~=~
    \vcenter{\xy
        (-3,25)="1";
        ( 3,25)="2";
        ( 3,22)="3";
        ( 3,18)="4";
        ( 3,14)="5";
        (-3,16)="6";
        "1";"4" **\crv{(-3,21)};
        "2";"5" **\dir{-};
        (-1, 5)="x";
        ( 6, 5)*[o]=<14pt>[Fo]{\scriptstyle \nu^{\text{-}1}}="y";
        (-3,-1)="c";
        ( 2,-1)="z";
        (-3,-7)="e";
        ( 3,-5)="f";
        ( 5,-7)="g";
        ( 0,-9)="h";
        ( 0,-12)*[o]=<5pt>[Fo]{~}="i";
        (-3,-15)="j";
        ( 3,-15)="k";
        "5";"x" **\dir{-};
        "z";"x" **\dir{-};
        "z";"y" **\dir{-};
        "e";"z" **\dir{-};
        "h";"e" **\dir{-};
        "h";"i" **\dir{-};
        "f";"g" **\dir{-};
        "g";"h" **\dir{-};
        "e";"j" **\dir{-};
        \ar@{-} "6";"y" |<(0.47)\hole
        \ar@{-}@/^4pt/ "6";"3" |<(0.6)\hole
        \ar@{-} "c";"f" |<(0.5)\hole
        \ar@{-} "c";"6"
        \ar@{-} "f";"k" |<(0.3)\hole
    \endxy}
    ~\overset{\text{(c)}}{=}~
    \vcenter{\xy
        (-5,25)="1";
        ( 4,25)="2";
        ( 4,22)="3";
        ( 4,18)="4";
        ( 1,15)="5";
        ( 7,15)="6";
        ( 0,10)="7";
        ( 9,10)*[o]=<14pt>[Fo]{\scriptstyle \nu^{\text{-}1}}="8";
        ( 4,6)="9";
        (-5,4)="10";
        (-5,1)="11";
        ( 4,3)="12";
        ( 6,1)="13";
        (-.5,-1)="14";
        (-.5,-4)*[o]=<5pt>[Fo]{~}="15";
        (-5,-7)="16";
        ( 4,-7)="17";
        "1";"16" **\dir{-};
        "2";"4" **\dir{-};
        "4";"5" **\dir{-};
        "4";"6" **\dir{-};
        "7";"6" **\dir{-};
        "7";"9" **\dir{-};
        "8";"9" **\dir{-};
        "10";"9" **\dir{-};
        "11";"14" **\dir{-};
        "15";"14" **\dir{-};
        "12";"13" **\dir{-};
        "14";"13" **\dir{-};
        \ar@{-} "5";"8" |<(0.4)\hole
        \ar@{-} "12";"17" |<(0.3)\hole
        \ar@{-}@/^19pt/ "12";"3" |<(0.12)\hole
    \endxy}
\\
    &~\overset{\text{($\dagger$)}}{=}~
    \vcenter{\xy
        (-4,20)="1";
        ( 4,20)="2";
        ( 4,17)="3";
        ( 6,15)="4";
        ( 0,13)*[o]=<9pt>[Fo]{\scriptstyle s}="5";
        (-4,11.5)="6";
        (-4,8)="7";
        ( 4,10)="8";
        ( 6,8)="9";
        ( 0,6)="10";
        ( 0,3)*[o]=<5pt>[Fo]{~}="11";
        (-4,0)="12";
        ( 4,0)="13";
        "1";"12" **\dir{-};
        "2";"3" **\dir{-};
        "4";"3" **\dir{-};
        "4";"5" **\dir{-};
        "5";"6" **\dir{-};
        "7";"10" **\dir{-};
        "11";"10" **\dir{-};
        "8";"9" **\dir{-};
        "9";"10" **\dir{-};
        \ar@{-} "2";"13" |<(0.3)\hole |<(0.65)\hole
    \endxy}
    ~\overset{\text{(4)}}{=}~
    \vcenter{\xy
        (-4,20)="1";
        ( 4,20)="2";
        ( 4,17)="3";
        ( 6,15)="4";
        ( 1.5,12.7)*[o]=<9pt>[Fo]{\scriptstyle s}="5";
        (-4,13)="6";
        (-2,11)="x";
        (-2,8)*[o]=<5pt>[Fo]{~}="y";
        (-4,6)="7";
        ( 4,8)="8";
        ( 6,6)="9";
        ( 0,4)="10";
        ( 0,1)*[o]=<5pt>[Fo]{~}="11";
        (-4,-2)="12";
        ( 4,-2)="13";
        "1";"12" **\dir{-};
        "2";"3" **\dir{-};
        "4";"3" **\dir{-};
        "4";"5" **\dir{-};
        "5";"x" **\dir{-};
        "6";"x" **\dir{-};
        "x";"y" **\dir{-};
        "7";"10" **\dir{-};
        "11";"10" **\dir{-};
        "8";"9" **\dir{-};
        "9";"10" **\dir{-};
        \ar@{-} "2";"13" |<(0.27)\hole |<(0.68)\hole
    \endxy}
    ~\overset{\text{(2)}}{=}~
    \vcenter{\xy
        (-3,20)="1";
        ( 3,20)="2";
        ( 3,17)="3";
        ( 5,15)="4";
        ( 0,13)="5";
        (-3,15)="6";
        ( 0,10)*[o]=<5pt>[Fo]{~}="x";
        (-3,8)="7";
        ( 3,10)="8";
        ( 5,8)="9";
        ( 0,6)="10";
        ( 0,3)*[o]=<5pt>[Fo]{~}="11";
        (-3,0)="12";
        ( 3,0)="13";
        "1";"12" **\dir{-};
        "4";"3" **\dir{-};
        "4";"5" **\dir{-};
        "6";"5" **\dir{-};
        "5";"x" **\dir{-};
        "7";"10" **\dir{-};
        "11";"10" **\dir{-};
        "8";"9" **\dir{-};
        "9";"10" **\dir{-};
        \ar@{-} "2";"13" |<(0.3)\hole |<(0.65)\hole
    \endxy}
    ~=~
    \vcenter{\xy
        (-3,20)="1";
        ( 3,20)="2";
        ( 3,17)="3";
        ( 5,15)="4";
        ( 0,13)="5";
        (-3,15)="6";
        ( 0,10)*[o]=<5pt>[Fo]{~}="7";
        (-3,7)="8";
        ( 3,7)="9";
        "1";"8" **\dir{-};
        "4";"3" **\dir{-};
        "4";"5" **\dir{-};
        "6";"5" **\dir{-};
        "5";"7" **\dir{-};
        \ar@{-} "2";"9" |<(0.45)\hole
    \endxy}
    ~= 1_P
\end{align*}
where $(\dagger)$ is given by
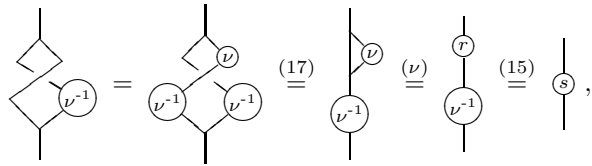
\[
    \vcenter{\xy
        ( 0,21)="1";
        ( 0,17)="2";
        (-3,14)="3";
        ( 3,14)="4";
        (-4,9)="5";
        ( 5,9)*[o]=<14pt>[Fo]{\scriptstyle \nu^{\text{-}1}}="6";
        ( 0,5)="7";
        ( 0,1)="8";
        "1";"2" **\dir{-};
        "2";"3" **\dir{-};
        "2";"4" **\dir{-};
        "4";"5" **\dir{-};
        "6";"7" **\dir{-};
        "5";"7" **\dir{-};
        "8";"7" **\dir{-};
        \ar@{-} "3";"6" |<(0.4)\hole
    \endxy}
    ~=~
    \vcenter{\xy
        ( 0,22)="1";
        ( 0,18)="2";
        (-3,15)="3";
        ( 3,15)*[o]=<9pt>[Fo]{\scriptstyle \nu}="4";
        (-5,9)*[o]=<14pt>[Fo]{\scriptstyle \nu^{\text{-}1}}="5";
        ( 5,9)*[o]=<14pt>[Fo]{\scriptstyle \nu^{\text{-}1}}="6";
        ( 0,5)="7";
        ( 0,1)="8";
        "1";"2" **\dir{-};
        "2";"3" **\dir{-};
        "2";"4" **\dir{-};
        "4";"5" **\dir{-};
        "6";"7" **\dir{-};
        "5";"7" **\dir{-};
        "8";"7" **\dir{-};
        \ar@{-} "3";"6" |<(0.4)\hole
    \endxy}
    ~\overset{\text{(17)}}{=}~
    \vcenter{\xy
        ( 0,20)="1";
        ( 0,17)="2";
        ( 3,14)*[o]=<9pt>[Fo]{\scriptstyle \nu}="3";
        ( 0,11)="4";
        ( 0, 6)*[o]=<14pt>[Fo]{\scriptstyle \nu^{\text{-}1}}="5";
        ( 0,0)="6";
        "1";"5" **\dir{-};
        "5";"6" **\dir{-};
        "2";"3" **\dir{-};
        "4";"3" **\dir{-};
    \endxy}
    ~\overset{\text{($\nu$)}}{=}~
    \vcenter{\xy
        ( 0,18)="1";
        ( 0,14)*[o]=<9pt>[Fo]{\scriptstyle r}="2";
        ( 0, 6)*[o]=<14pt>[Fo]{\scriptstyle \nu^{\text{-}1}}="3";
        ( 0,0)="4";
        "1";"2" **\dir{-};
        "2";"3" **\dir{-};
        "3";"4" **\dir{-};
    \endxy}
    ~\overset{\text{(15)}}{=}~
    \vcenter{\xy
        ( 0, 6)="1";
        ( 0, 0)*[o]=<9pt>[Fo]{\scriptstyle s}="2";
        ( 0,-6)="3";
        "1";"2" **\dir{-};
        "2";"3" **\dir{-};
    \endxy}~,
\]
and in the other direction by:
\begin{align*}
    \theta \theta^{-1} &=~
    \vcenter{\xy
        (-3,14)="1";
        ( 3,14)="2";
        (-3,11)="3";
        ( 3, 8)="4";
        (-3, 3)="5";
        ( 3, 3)*[o]=<14pt>[Fo]{\scriptstyle \nu^{\text{-}1}}="6";
        (-3,-6)="7";
        ( 3,-4)="8";
        ( 5,-6)="9";
        ( 0,-8)="10";
        ( 0,-11)*[o]=<5pt>[Fo]{~}="11";
        (-3,-12)="12";
        ( 3,-12)="13";
        ( 3,-16)="14";
        (-3,-19)="15";
        ( 3,-19)*[o]=<9pt>[Fo]{\scriptstyle \nu}="16";
        (-3,-23)="17";
        ( 3,-23)="18";
        "1";"5" **\dir{-};
        "3";"4" **\dir{-};
        "2";"6" **\dir{-};
        "6";"7" **\dir{-};
        "10";"7" **\dir{-};
        "10";"11" **\dir{-};
        "8";"9" **\dir{-};
        "9";"10" **\dir{-};
        "7";"12" **\dir{-};
        "18";"16" **\dir{-};
        "15";"17" **\dir{-};
        \ar@{-} "5";"8" |<(0.56)\hole
        \ar@{-} "8";"16" |<(0.2)\hole
        \ar@{-}@/^3pt/ "15";"13" |<(0.5)\hole
        \ar@{-}@/_3pt/ "12";"14"
    \endxy}
    ~\overset{\text{($\ddagger$)}}{=}~
    \vcenter{\xy
        (-3,14)="1";
        ( 3,14)="2";
        (-3,11)="3";
        ( 3, 8)="4";
        (-3, 3)="5";
        ( 3, 3)*[o]=<14pt>[Fo]{\scriptstyle \nu^{\text{-}1}}="6";
        (-3,-4)="12";
        ( 3,-6)="13";
        ( 3,-10)="14";
        (-3,-13)="15";
        ( 3,-13)*[o]=<9pt>[Fo]{\scriptstyle \nu}="16";
        (-3,-17)="17";
        ( 3,-17)="18";
        "1";"5" **\dir{-};
        "3";"4" **\dir{-};
        "2";"6" **\dir{-};
        "6";"12" **\dir{-};
        "18";"16" **\dir{-};
        "15";"17" **\dir{-};
        \ar@{-} "5";"13" |<(0.45)\hole
        \ar@{-} "13";"16"
        \ar@{-}@/^3pt/ "15";"13" |<(0.55)\hole
        \ar@{-}@/_3pt/ "12";"14"
    \endxy}
    ~\overset{\text{(c)}}{=}~
    \vcenter{\xy
        (-3,14)="1";
        ( 3,14)="2";
        (-3, 8)="3";
        ( 3, 6)="4";
        (-3, 0)="5";
        ( 3, 0)*[o]=<14pt>[Fo]{\scriptstyle \nu^{\text{-}1}}="6";
        (-3,-5)="7";
        ( 3,-5)="8";
        ( 0,-8)="9";
        ( 0,-12)*[o]=<9pt>[Fo]{\scriptstyle \nu}="10";
        ( 0,-17)="11";
        (-3, 11)="12";
        (-8,-17)="13";
        "1";"5" **\dir{-};
        "2";"6" **\dir{-};
        "3";"4" **\dir{-};
        "6";"7" **\dir{-};
        "7";"9" **\dir{-};
        "8";"9" **\dir{-};
        "9";"10" **\dir{-};
        "11";"10" **\dir{-};
        "13";"12" **\crv{(-8,9)};
        \ar@{-} "5";"8" |<(0.5)\hole
    \endxy}
    ~\overset{\text{(17)}}{=}~
    \vcenter{\xy
        (-3,14)="1";
        ( 3,14)="2";
        (-3, 8)="3";
        ( 3, 6)="4";
        (-3, 0)="5";
        ( 3, 0)*[o]=<14pt>[Fo]{\scriptstyle \nu^{\text{-}1}}="6";
        (-3,-7)*[o]=<9pt>[Fo]{\scriptstyle \nu}="7";
        ( 3,-7)*[o]=<9pt>[Fo]{\scriptstyle \nu}="8";
        ( 0,-11)="9";
        ( 0,-16)="11";
        (-3, 11)="12";
        (-8,-16)="13";
        "1";"5" **\dir{-};
        "2";"6" **\dir{-};
        "3";"4" **\dir{-};
        "5";"7" **\dir{-};
        "6";"8" **\dir{-};
        "7";"9" **\dir{-};
        "8";"9" **\dir{-};
        "9";"11" **\dir{-};
        "13";"12" **\crv{(-8,9)};
    \endxy}
\\
    &=~
    \vcenter{\xy
        (-3,14)="1";
        ( 3,14)="2";
        (-3, 8)="3";
        ( 3, 6)="4";
        (-3, 2)*[o]=<9pt>[Fo]{\scriptstyle \nu}="5";
        ( 3, 2)="6";
        ( 0,-1)="7";
        ( 0,-6)="8";
        (-3,11)="9";
        (-8,-6)="10";
        "1";"5" **\dir{-};
        "2";"6" **\dir{-};
        "3";"4" **\dir{-};
        "5";"7" **\dir{-};
        "6";"7" **\dir{-};
        "7";"8" **\dir{-};
        "10";"9" **\crv{(-8,8)};
    \endxy}
    ~\overset{\text{(c,$\nu$)}}{=}~
    \vcenter{\xy
        (-3, 7)="1";
        ( 3, 7)="2";
        (-3, 3)="3";
        ( 0, 0)*[o]=<9pt>[Fo]{\scriptstyle t}="4";
        ( 3,-3)="5";
        (-3,-7)="6";
        ( 3,-7)="7";
        "1";"6" **\dir{-};
        "2";"7" **\dir{-};
        "3";"4" **\dir{-};
        "4";"5" **\dir{-};
    \endxy}
    ~\overset{\text{(4)}}{=}~
    \vcenter{\xy
        (-3, 9)="1";
        ( 3, 9)="2";
        (-3, 6)="3";
        ( 3, 5)="4";
        ( 0, 4)*[o]=<9pt>[Fo]{\scriptstyle t}="5";
        ( 5, 3)="6";
        ( 0, 0)="7";
        ( 0,-3)*[o]=<5pt>[Fo]{~}="8";
        (-3,-5)="9";
        ( 3,-5)="10";
        "1";"9" **\dir{-};
        "3";"5" **\dir{-};
        "5";"8" **\dir{-};
        "4";"6" **\dir{-};
        "6";"7" **\dir{-};
        \ar@{-} "2";"10" |<(0.53)\hole
    \endxy}
    ~\overset{\text{(2)}}{=}~
    \vcenter{\xy
        (-3, 9)="1";
        ( 3, 9)="2";
        (-3, 4)="3";
        ( 3, 6)="4";
        ( 5, 4)="5";
        ( 0, 2)="6";
        ( 0,-1)*[o]=<5pt>[Fo]{~}="7";
        (-3,-4)="8";
        ( 3,-4)="9";
        "1";"8" **\dir{-};
        "3";"6" **\dir{-};
        "4";"5" **\dir{-};
        "5";"6" **\dir{-};
        "6";"7" **\dir{-};
        \ar@{-} "2";"9" |<(0.47)\hole
    \endxy}
    ~= 1_P
\end{align*}
for which the first step ($\ddagger$) holds since $\theta$ is a morphism in
$\Q\V$.
\end{proof}

That the antipode $\nu:A^\o \ra A$ is a comonoid isomorphism is our
assumption. That $\upsilon:C^{\o\o} \ra C$ is as well may be seen from the
following calculation:
\begin{align*}
    (t \ox t)\delta \upsilon &= (t \ox t)\delta t \nu \nu t \\
     &= (t \ox t)\delta \nu \nu t & \text{(3)} \\
     &= (t \ox t)c (\nu \ox \nu) \delta \nu t & \text{(17)} \\
     &= (t \ox t)c (\nu \ox \nu)c (\nu \ox \nu) \delta t & \text{(17)} \\
     &= (t \ox t)(\nu \ox \nu)(\nu \ox \nu)cc \delta t & \text{(nat)} \\
     &= (t \ox t)(t \ox t)(\nu \ox \nu)(\nu \ox \nu)cc \delta t & \text{(7)} \\
     &= (t \ox t)(\nu \ox \nu)(r \ox r)(\nu \ox \nu)cc \delta t & \text{(16)} \\
     &= (t \ox t)(\nu \ox \nu)(\nu \ox \nu)(t \ox t)cc \delta t & \text{(16)} \\
     &= (t \ox t)(\nu \ox \nu)(\nu \ox \nu)(t \ox t)cc \delta & \text{(5)} \\
     &= (\upsilon \ox \upsilon)cc \delta. & \text{(5)}
\end{align*}

An inverse for $\upsilon$ is given by the morphism
\[
    \upsilon^{-1} = t \nu^{-1} \nu^{-1} t,
\]
as may be seen in one direction by the calculation
\begin{align*}
    \upsilon^{-1}\upsilon &= t \nu^{-1} \nu^{-1} t t \nu \nu t \\
        &= t t \nu^{-1} \nu^{-1} \nu \nu t t  & \text{(16)} \\
        &= t \nu^{-1} \nu^{-1} \nu \nu t  & \text{(7)} \\
        &= t t \\
        &= t = 1_C. & \text{(7)}
\end{align*}
The other direction is similar.

Recall that the left $A \ox A$-, right $A$-coaction $\delta$ on $P$ is
defined by taking the diagonal of the commutative square:
\[
    \xygraph{{P}="1"
        [r(2.8)] {A \ox A \ox P}="2"
        [d(1.2)] {A \ox A \ox P \ox A.}="4"
     "1"[d(1.2)] {P \ox A}="3"
        "1":"2" ^-{\delta_l}
        "2":"4" ^-{1 \ox 1 \ox \delta_r}
        "1":"3" _-{\delta_r}
        "3":"4" ^-{\delta_l \ox 1}}
\]
We note that $\delta$ may be written as
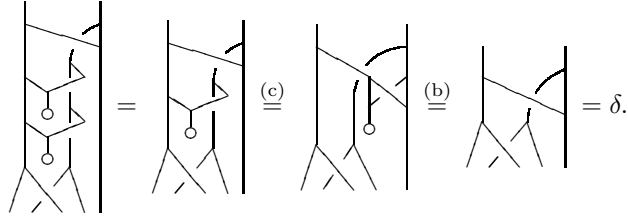
\[
    \vcenter{\xy
        (-3,20)="1"; ( 7,20)="2"; (-3,17)="3"; ( 7,17)="4"; ( 7,14)="5";
        (-3,10)="6"; ( 3,12)="7"; ( 5,10)="8"; ( 0, 8)="9";
        ( 0, 5)*[o]=<5pt>[Fo]{~}="10";
        (-3, 4)="11"; ( 3, 6)="12"; ( 5, 4)="13"; ( 0, 2)="14";
        ( 0,-1)*[o]=<5pt>[Fo]{~}="15";
        (-3,-2)="16"; ( 3,-2)="17"; (-5,-8)="18"; (-2,-8)="19";
        ( 2,-8)="20"; ( 5,-8)="21"; ( 7,-8)="22";
        "1";"16" **\dir{-}; "2";"22" **\dir{-}; "3";"5" **\dir{-};
        "6";"9" **\dir{-}; "10";"9" **\dir{-}; "11";"14" **\dir{-};
        "15";"14" **\dir{-}; "16";"18" **\dir{-}; "16";"20" **\dir{-};
        "7";"8" **\dir{-}; "9";"8" **\dir{-}; "12";"13" **\dir{-};
        "14";"13" **\dir{-}; "17";"21" **\dir{-};
        \ar@{-}@/^3pt/ "7";"4" |<(0.4)\hole
        \ar@{-} "17";"19" |<(0.6)\hole
        \ar@{-} "7";"17" |<(0.22)\hole |<(0.65)\hole
    \endxy}
    ~=~
    \vcenter{\xy
        (-3,20)="1"; ( 7,20)="2"; (-3,17)="3"; ( 7,17)="4"; ( 7,14)="5";
        (-3,10)="6"; ( 3,12)="7"; ( 5,10)="8"; ( 0, 8)="9";
        ( 0, 5)*[o]=<5pt>[Fo]{~}="10"; (-3, 3)="16"; ( 3, 3)="17";
        (-5,-3)="18"; (-2,-3)="19"; ( 2,-3)="20"; ( 5,-3)="21"; ( 7,-3)="22";
        "1";"16" **\dir{-}; "2";"22" **\dir{-}; "3";"5" **\dir{-};
        "6";"9" **\dir{-}; "10";"9" **\dir{-}; "16";"18" **\dir{-};
        "16";"20" **\dir{-}; "7";"8" **\dir{-}; "9";"8" **\dir{-};
        "17";"21" **\dir{-};
        \ar@{-}@/^3pt/ "7";"4" |<(0.4)\hole
        \ar@{-} "17";"19" |<(0.6)\hole
        \ar@{-} "7";"17" |<(0.35)\hole
    \endxy}
    ~\overset{\text{(c)}}{=}~
    \vcenter{\xy
        (-6,21)="1"; ( 6,21)="2"; (-6,18)="3"; ( 6,18)="4"; ( 1,14)="5";
        ( 6,14)="6"; ( 1,10)="7"; ( 6,10)="8"; ( 1, 7)*[o]=<5pt>[Fo]{~}="9";
        (-6, 5)="10"; (-1, 5)="11"; (-8,-1)="12"; (-5,-1)="13"; (-1,-1)="14";
        ( 2,-1)="15"; ( 6,-1)="16";
        "1";"10" **\dir{-}; "10";"12" **\dir{-}; "10";"14" **\dir{-};
        "3";"5" **\dir{-}; "5";"8" **\dir{-}; "5";"9" **\dir{-};
        "2";"16" **\dir{-}; "11";"15" **\dir{-}; "11";(-1,12) **\dir{-};
        \ar@{-}@/^5pt/ (-1,12);"4" |<(0.25)\hole
        \ar@{-} "6";"7" |\hole
        \ar@{-} "11";"13" |<(0.5)\hole
    \endxy}
    ~\overset{\text{(b)}}{=}~
    \vcenter{\xy
        (-7,15)="1"; ( 4,15)="2"; (-7,11)="3"; ( 4,12)="4"; ( 4, 6)="5";
        (-7, 5)="10"; (-1, 5)="11"; (-9,-1)="12"; (-6,-1)="13"; (-2,-1)="14";
        ( 1,-1)="15"; ( 4,-1)="16";
        "1";"10" **\dir{-}; "10";"12" **\dir{-}; "10";"14" **\dir{-};
        "3";"5" **\dir{-}; "2";"16" **\dir{-}; "11";"15" **\dir{-};
        \ar@{-}@/^4pt/ "11";"4" |<(0.3)\hole
        \ar@{-} "11";"13" |<(0.6)\hole
    \endxy}
    ~= \delta.
\]

We must show that $\theta$ is a left $A^{\ox 3}$-comodule isomorphism
$P_l \ra P_r$. That is, we must prove the commutativity of the square
\[
    \xygraph{{P_l}="1"
        [r(2.2)] {A^{\ox 3} \ox P_l}="2"
        [d(1.2)] {A^{\ox 3} \ox P_r}="4"
     "1"[d(1.2)] {P_r}="3"
        "1":"2" ^-\gamma
        "2":"4" ^-{1 \ox \theta}
        "1":"3" _-\theta
        "3":"4" ^-\gamma}
\]
where the left $A^{\ox 3}$-coactions on $P_l$ and $P_r$ were defined 
using $\delta$ (see \S\ref{sec-quangpd}).

The clockwise direction around the square is
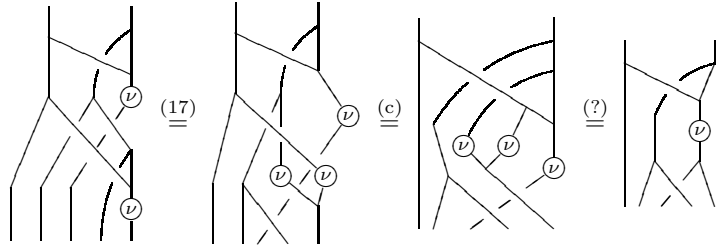
\[
    \vcenter{\xy
        (0,16)*{~}; (0,-17)*{~};
        (-7,15)="1"; ( 4,15)="2"; (-7,11)="3"; ( 4,12)="4"; ( 4, 6)="5";
        (-7, 3)="6"; (-1, 3)="7";
        ( 4, 3)*[o]=<9pt>[Fo]{\scriptstyle \nu}="8";
        (-12,-9)="9"; (-8,-9)="10"; (-4,-9)="11"; ( 4,-4)="12"; ( 4,-9)="13";
        ( 4,-12)*[o]=<9pt>[Fo]{\scriptstyle \nu}="14";
        (-12,-16)="15"; (-8,-16)="16"; (-4,-16)="17"; ( 0,-16)="18";
        ( 4,-16)="19";
        "1";"6" **\dir{-}; "3";"5" **\dir{-}; "2";"8" **\dir{-};
        "6";"9" **\dir{-}; "6";"13" **\dir{-}; "7";"12" **\dir{-};
        "15";"9" **\dir{-}; "16";"10" **\dir{-}; "11";"17" **\dir{-};
        "14";"12" **\dir{-}; "14";"19" **\dir{-};
        \ar@{-}@/^4pt/ "7";"4" |<(0.45)\hole
        \ar@{-} "7";"10" |<(0.35)\hole
        \ar@{-} "8";"11" |<(0.25)\hole |<(0.6)\hole
        \ar@{-}@/^3pt/ "18";"12" |<(0.7)\hole
    \endxy}
    ~\overset{\text{(17)}}{=}~
    \vcenter{\xy
        (-7,15)="1"; ( 4,15)="2"; (-7,11)="3"; ( 4,12)="4"; ( 4, 6)="5";
        (-7, 3)="6"; (-1, 3)="7";
        ( 8, 0)*[o]=<9pt>[Fo]{\scriptstyle \nu}="8";
        (-10,-9)="9"; (-6,-9)="10"; ( 4,-4)="12"; ( 4,-9)="13";
        (-1,-8)*[o]=<9pt>[Fo]{\scriptstyle \nu}="x";
        ( 5,-8)*[o]=<9pt>[Fo]{\scriptstyle \nu}="y";
        ( 4,-12)="z"; (-4,-17)="11"; (-10,-17)="15"; (-6,-17)="16";
        (-4,-17)="17"; ( 0,-17)="18"; ( 4,-17)="19";
        "6";"y" **\dir{-}; "z";"y" **\dir{-}; "z";"x" **\dir{-};
        "z";"19" **\dir{-}; "1";"6" **\dir{-}; "3";"5" **\dir{-};
        "2";"5" **\dir{-}; "5";"8" **\dir{-}; "6";"9" **\dir{-};
        "15";"9" **\dir{-}; "10";"16" **\dir{-}; "10";"18" **\dir{-};
        \ar@{-}@/^4pt/ "7";"4" |<(0.45)\hole
        \ar@{-} "7";"10" |<(0.3)\hole
        \ar@{-} "8";"11" |<(0.31)\hole |<(0.56)\hole |<(0.82)\hole
        \ar@{-} "7";"x" |<(0.5)\hole
    \endxy}
    ~\overset{\text{(c)}}{=}~
    \vcenter{\xy
        (-10,23)="1";
        ( 8,23)="2";
        (-10,20)="3";
        ( 8,20)="4";
        ( 8,16)="5";
        (-8,9)="6";
        (-6,2)="17";
        ( 8,9)="7";
        (-4,6)*[o]=<9pt>[Fo]{\scriptstyle \nu}="8";
        ( 2,6)*[o]=<9pt>[Fo]{\scriptstyle \nu}="9";
        (-1,3)="10";
        ( 8,3)*[o]=<9pt>[Fo]{\scriptstyle \nu}="11";
        (-10,-5)="12";
        (-8, -5)="13";
        (-3, -5)="14";
        ( 2, -5)="15";
        ( 8, -5)="16";
        "1";"12" **\dir{-};
        "2";"11" **\dir{-};
        "3";"7" **\dir{-}?(0.8)="x";
        "x";"9" **\dir{-};
        "8";"10" **\dir{-};
        "9";"10" **\dir{-};
        "17";"6" **\dir{-};
        "10";"16" **\dir{-};
        "17";"13" **\dir{-};
        "17";"15" **\dir{-};
        \ar@{-}@/^6pt/ "6";"4" |<(0.42)\hole
        \ar@{-}@/^6pt/ "8";"5" |<(0.5)\hole
        \ar@{-} "11";"14" |<(0.35)\hole |<(0.72)\hole
    \endxy}
    ~\overset{\text{(?)}}{=}~
    \vcenter{\xy
        (-6,20)="1"; ( 6,20)="2"; (-6,17)="3"; ( 6,17)="4";
        (-2,10)="5"; ( 4,12)="6";
        ( 4,8)*[o]=<9pt>[Fo]{\scriptstyle \nu}="7";
        (-2,4)="8"; ( 4,4)="9"; (-6,-2)="10"; (-4,-2)="11";
        (-1,-2)="12"; ( 3,-2)="13"; ( 6,-2)="14";
        "1";"10" **\dir{-}; "3";"6" **\dir{-}; "2";"4" **\dir{-};
        "6";"4" **\dir{-}; "6";"7" **\dir{-}; "9";"7" **\dir{-};
        "9";"14" **\dir{-}; "5";"8" **\dir{-}; "11";"8" **\dir{-};
        \ar@{-}@/^4pt/ "5";"4" |<(0.4)\hole
        \ar@{-} "12";"9" |<(0.4)\hole
        \ar@{-} "13";"8"
    \endxy}
\]
where the last step (?) is given by the following calculation

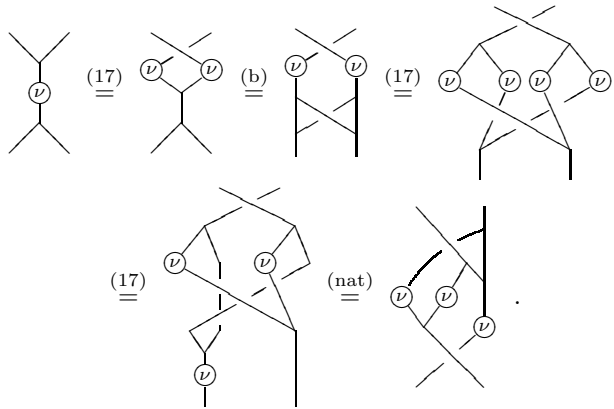
\[
    \vcenter{\xy
        (-4, 8)="1"; ( 4, 8)="2"; ( 0, 4)="3";
        ( 0, 0)*[o]=<9pt>[Fo]{\scriptstyle \nu}="4";
        ( 0,-4)="5"; (-4,-8)="6"; ( 4,-8)="7";
        "1";"3" **\dir{-}; "2";"3" **\dir{-}; "4";"3" **\dir{-};
        "4";"5" **\dir{-}; "6";"5" **\dir{-}; "7";"5" **\dir{-};
    \endxy}
    ~\overset{\text{(17)}}{=}~
    \vcenter{\xy
        (-4, 8)="1"; ( 4, 8)="2";
        (-4, 3)*[o]=<9pt>[Fo]{\scriptstyle \nu}="3";
        ( 4, 3)*[o]=<9pt>[Fo]{\scriptstyle \nu}="4";
        ( 0, 0)="5"; ( 0,-4)="6"; (-4,-8)="7"; ( 4,-8)="8";
        "1";"4" **\dir{-}; "5";"3" **\dir{-}; "4";"5" **\dir{-};
        "5";"6" **\dir{-}; "6";"7" **\dir{-}; "6";"8" **\dir{-};
        \ar@{-} "2";"3" |\hole
    \endxy}
    ~\overset{\text{(b)}}{=}~
    \vcenter{\xy
        (0,10)*{~}; (0,-9)*{~};
        (-4, 9)="1"; ( 4, 9)="2";
        (-4, 4)*[o]=<9pt>[Fo]{\scriptstyle \nu}="3";
        ( 4, 4)*[o]=<9pt>[Fo]{\scriptstyle \nu}="4";
        (-4, 0)="5"; ( 4, 0)="6"; (-4,-5)="7"; ( 4,-5)="8";
        (-4,-8)="9"; ( 4,-8)="10";
        "1";"4" **\dir{-}; "9";"3" **\dir{-}; "4";"10" **\dir{-};
        "5";"8" **\dir{-}; 
        \ar@{-} "2";"3" |\hole
        \ar@{-} "7";"6" |\hole
    \endxy}
    ~\overset{\text{(17)}}{=}~
    \vcenter{\xy
        (-4,10)="1"; ( 4,10)="2"; (-6, 5)="3"; ( 6, 5)="4";
        (-10,0)*[o]=<9pt>[Fo]{\scriptstyle \nu}="5";
        (-2, 0)*[o]=<9pt>[Fo]{\scriptstyle \nu}="6";
        ( 2, 0)*[o]=<9pt>[Fo]{\scriptstyle \nu}="7";
        (10, 0)*[o]=<9pt>[Fo]{\scriptstyle \nu}="8";
        (-6,-9)="9"; ( 6,-9)="10";
        (-6,-13)="11"; ( 6,-13)="12";
        "1";"4" **\dir{-}; "5";"3" **\dir{-}; "3";"6" **\dir{-};
        "4";"7" **\dir{-}; "4";"8" **\dir{-}; "5";"10" **\dir{-};
        "7";"10" **\dir{-}; "9";"11" **\dir{-}; "10";"12" **\dir{-};
        \ar@{-} "2";"3" |<(0.4)\hole
        \ar@{-} "6";"9" |<(0.35)\hole
        \ar@{-} "8";"9" |<(0.3)\hole |<(0.6)\hole
    \endxy}
\]
\[
    ~\overset{\text{(17)}}{=}~
    \vcenter{\xy
        (-4,10)="1"; ( 4,10)="2"; (-6, 5)="3"; ( 6, 5)="4";
        (-10,0)*[o]=<9pt>[Fo]{\scriptstyle \nu}="5";
        (-4, 0)="6";
        ( 2, 0)*[o]=<9pt>[Fo]{\scriptstyle \nu}="7";
        ( 8, 0)="8"; (-8,-9)="x"; (-4,-9)="9"; (-6,-12)="y";
        (-6,-15)*[o]=<9pt>[Fo]{\scriptstyle \nu}="z";
        ( 6,-9)="10"; (-6,-19)="11"; ( 6,-19)="12";
        "1";"4" **\dir{-}; "5";"3" **\dir{-}; "3";"6" **\dir{-};
        "4";"7" **\dir{-}; "4";"8" **\dir{-}; 
        "5";"10" **\dir{-}; "7";"10" **\dir{-}; "10";"12" **\dir{-};
        "9";"y" **\dir{-}; "x";"y" **\dir{-}; "y";"z" **\dir{-};
        "z";"11" **\dir{-};
        \ar@{-} "2";"3" |<(0.4)\hole
        \ar@{-} "6";"9" |<(0.35)\hole |<(0.8)\hole
        \ar@{-} "8";"x" |<(0.3)\hole |<(0.55)\hole
    \endxy}
    ~\overset{\text{(nat)}}{=}~
    \vcenter{\xy
        (-6,10)="1"; ( 3,10)="2"; ( 3, 7)="3";
        (-8,-2)*[o]=<9pt>[Fo]{\scriptstyle \nu}="4";
        (-2,-2)*[o]=<9pt>[Fo]{\scriptstyle \nu}="5";
        ( 3, 0)="6"; (-5,-6)="7";
        ( 3,-6)*[o]=<9pt>[Fo]{\scriptstyle \nu}="8";
        (-6,-14)="9"; ( 3,-14)="10";
        "2";"8" **\dir{-}; "1";"6" **\dir{-}?(0.75)=x; "x";"5" **\dir{-};
        "4";"7" **\dir{-}; "5";"7" **\dir{-}; "7";"10" **\dir{-};
        \ar@{-}@/_4pt/ "3";"4" |<(0.35)\hole
        \ar@{-} "8";"9" |<(0.4)\hole
    \endxy}\ \ .
\]

The counter-clockwise direction is
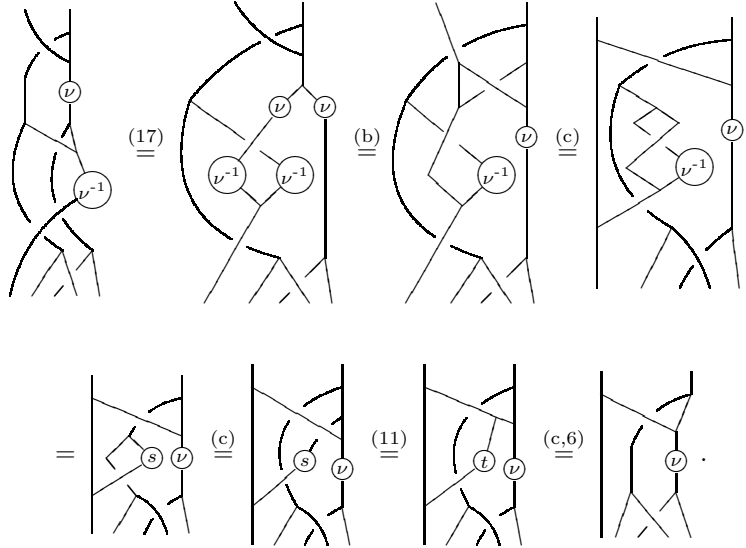
\[
    \vcenter{\xy
        (-3,20)="1"; ( 3,20)="2"; ( 3,17)="3"; ( 3,13)="4";
        ( 3, 9)*[o]=<9pt>[Fo]{\scriptstyle \nu}="5";
        (-3,11)="x"; (-3,5)="6"; ( 3,5)="7"; ( 2,-12)="8"; ( 6,-12)="9";
        ( 4,1)="10";
        ( 6,-4)*[o]=<14pt>[Fo]{\scriptstyle \nu^{\text{-}1}}="11";
        (-5,-18)="16"; (-2,-18)="17"; ( 1,-18)="18"; ( 4,-18)="19";
        ( 7,-18)="20";
        "2";"5" **\dir{-}; "7";"5" **\dir{-}; "10";"6" **\dir{-};
        "10";"7" **\dir{-}; "10";"11" **\dir{-}; "x";"6" **\dir{-};
        "8";"17" **\dir{-}; "8";"19" **\dir{-}; "9";"20" **\dir{-};
        \ar@{-}@/_6pt/ "11";"16"
        \ar@{-}@/_3pt/ "1";"4"
        \ar@{-}@/^3pt/ "x";"3" |<(0.6)\hole
        \ar@{-} "9";"18" |<(0.6)\hole
        \ar@{-}@/^11pt/ "8";"6" |<(0.2)\hole
        \ar@{-}@/^11pt/ "9";"7" |<(0.35)\hole |<(0.85)\hole
    \endxy}
    ~\overset{\text{(17)}}{=}\ \
    \vcenter{\xy
        ( 0,17)="1";
        ( 9,17)="2";
        ( 9,14)="3";
        ( 9,10)="4";
        (-6, 4)="5";
        ( 9, 6)="6";
        ( 6, 3)*[o]=<9pt>[Fo]{\scriptstyle \nu}="7";
        (12, 3)*[o]=<9pt>[Fo]{\scriptstyle \nu}="8";
        (-1,-6)*[o]=<14pt>[Fo]{\scriptstyle \nu^{\text{-}1}}="9";
        ( 8,-6)*[o]=<14pt>[Fo]{\scriptstyle \nu^{\text{-}1}}="10";
        (3.5,-10)="11";
        (6 ,-17)="x";
        (12,-17)="y";
        (-4,-23)="12";
        ( 2,-23)="13";
        ( 6,-23)="14";
        (10,-23)="15";
        (13,-23)="16";
        "2";"6" **\dir{-};
        "7";"6" **\dir{-};
        "8";"6" **\dir{-};
        "7";"9" **\dir{-};
        "11";"9" **\dir{-};
        "11";"10" **\dir{-};
        "11";"12" **\dir{-};
        "y";"8" **\dir{-};
        "x";"13" **\dir{-};
        "x";"15" **\dir{-};
        "y";"16" **\dir{-};
        \ar@{-}@/^18pt/ "x";"5" |<(0.18)\hole
        \ar@{-}@/_3pt/ "1";"4"
        \ar@{-}@/^4pt/ "5";"3" |<(0.73)\hole
        \ar@{-} "10";"5" |<(0.3)\hole
        \ar@{-} "y";"14" |<(0.6)\hole
    \endxy}
    ~\overset{\text{(b)}}{=}\ \
    \vcenter{\xy
        (0,18)*{~}; (0,-24)*{~};
        ( 0,17)="1"; (12,17)="2"; (12,14)="3"; ( 3,9)="4"; (12,9)="a";
        ( 3,3)="b"; (12,3)="c"; (-4, 4)="5"; ( 9, 6)="6"; ( 6, 3)="7";
        (12,-1)*[o]=<9pt>[Fo]{\scriptstyle \nu}="8";
        (-1,-6)="9";
        ( 8,-6)*[o]=<14pt>[Fo]{\scriptstyle \nu^{\text{-}1}}="10";
        (3.5,-10)="11"; (6 ,-17)="x"; (12,-17)="y"; (-4,-23)="12";
        ( 2,-23)="13"; ( 6,-23)="14"; (10,-23)="15"; (13,-23)="16";
        "4";"c" **\dir{-}; "4";"b" **\dir{-}; "2";"c" **\dir{-};
        "8";"c" **\dir{-}; "b";"9" **\dir{-}; "11";"9" **\dir{-};
        "11";"10" **\dir{-}; "11";"12" **\dir{-}; "y";"8" **\dir{-};
        "x";"13" **\dir{-}; "x";"15" **\dir{-}; "y";"16" **\dir{-};
        \ar@{-}@/^18pt/ "x";"5" |<(0.18)\hole
        \ar@{-} "1";"4"
        \ar@{-}@/^6pt/ "5";"3" |<(0.47)\hole
        \ar@{-} "10";"5" |<(0.4)\hole
        \ar@{-} "b";"a" |\hole
        \ar@{-} "y";"14" |<(0.6)\hole
    \endxy}
    ~\overset{\text{(c)}}{=}~
    \vcenter{\xy
        (-8,25)="1"; (10,25)="2"; (-8,22)="3"; (10,22)="4";
        (-5,16)="5"; (10,16)="6"; ( 0,13)="7"; (-3,11)="8";
        ( 3,11)="9"; (-4, 5)="10";
        ( 5, 5)*[o]=<14pt>[Fo]{\scriptstyle \nu^{\text{-}1}}="11";
        (0.5, 2)="12";
        (10,10)*[o]=<9pt>[Fo]{\scriptstyle \nu}="13";
        (-8,-3)="14"; ( 2,-3)="15"; (10,-3)="16"; (-8,-11)="17";
        (-2,-11)="18"; ( 3,-11)="19"; ( 7,-11)="20"; (11,-11)="21";
        "1";"17" **\dir{-}; "2";"13" **\dir{-}; "13";"16" **\dir{-};
        "3";"6" **\dir{-}; "5";"7" **\dir{-}; "8";"7" **\dir{-};
        "9";"7" **\dir{-}; "9";"10" **\dir{-}; "12";"10" **\dir{-};
        "12";"11" **\dir{-}; "12";"14" **\dir{-}; "15";"18" **\dir{-};
        "16";"21" **\dir{-};
        \ar@{-}@/^12pt/ "15";"5" |<(0.22)\hole
        \ar@{-}@/^4pt/ "5";"4" |<(0.35)\hole
        \ar@{-} "8";"11" |<(0.4)\hole
        \ar@{-}@/_2pt/ "16";"19" |<(0.55)\hole
        \ar@{-}@/^2pt/ "15";"20"
    \endxy}
\]
\[
    ~=~
    \vcenter{\xy
        (-6,20)="1";
        ( 6,20)="2";
        (-6,17)="3";
        ( 6,17)="4";
        (-1,12)="5";
        ( 6,12)="6";
        (-4, 9)="7";
        ( 2, 9)*[o]=<9pt>[Fo]{\scriptstyle s}="8";
        ( 6, 9)*[o]=<9pt>[Fo]{\scriptstyle \nu}="9";
        (-6, 4)="10";
        ( 0, 4)="11";
        ( 6, 4)="12";
        (-6,-1)="13";
        (-3,-1)="14";
        ( 1,-1)="15";
        ( 4,-1)="16";
        ( 7,-1)="17";
        "1";"13" **\dir{-};
        "3";"6" **\dir{-};
        "2";"9" **\dir{-};
        "5";"7" **\dir{-};
        "5";"8" **\dir{-};
        "10";"8" **\dir{-};
        "12";"9" **\dir{-};
        "11";"14" **\dir{-};
        "12";"17" **\dir{-};
        \ar@{-}@/^3pt/ "5";"4" |<(0.3)\hole
        \ar@{-} "11";"7" |\hole
        \ar@{-}@/_2pt/ "12";"15" |\hole
        \ar@{-}@/^2pt/ "11";"16"
    \endxy}
    ~\overset{\text{(c)}}{=}~
    \vcenter{\xy
        (0,24)*{~}; (0,-2)*{~};
        (-6,23)="1";
        ( 6,23)="2";
        (-6,20)="3";
        ( 6,20)="4";
        ( 6,16)="5";
        ( 6,13)="6";
        ( 1,10)*[o]=<9pt>[Fo]{\scriptstyle s}="8";
        ( 6, 9)*[o]=<9pt>[Fo]{\scriptstyle \nu}="9";
        (-6, 4)="10";
        ( 0, 4)="11";
        ( 6, 4)="12";
        (-6,-1)="13";
        (-3,-1)="14";
        ( 1,-1)="15";
        ( 4,-1)="16";
        ( 7,-1)="17";
        "1";"13" **\dir{-};
        "3";"6" **\dir{-};
        "2";"9" **\dir{-};
        "10";"8" **\dir{-};
        "12";"9" **\dir{-};
        "11";"14" **\dir{-};
        "12";"17" **\dir{-};
        \ar@{-}@/^2pt/ "8";"5" |<(0.5)\hole
        \ar@{-}@/_2pt/ "12";"15" |\hole
        \ar@{-}@/^2pt/ "11";"16"
        \ar@{-}@/_15pt/ "4";"11" |<(0.32)\hole |<(0.82)\hole
    \endxy}
    ~\overset{\text{(11)}}{=}~
    \vcenter{\xy
        (-6,23)="1"; ( 6,23)="2"; (-6,20)="3"; ( 6,20)="4"; ( 6,15)="6";
        ( 2,10)*[o]=<9pt>[Fo]{\scriptstyle t}="8";
        ( 6, 9)*[o]=<9pt>[Fo]{\scriptstyle \nu}="9";
        (-6, 4)="10"; ( 0, 4)="11"; ( 6, 4)="12"; (-6,-1)="13"; (-3,-1)="14";
        ( 1,-1)="15"; ( 4,-1)="16"; ( 7,-1)="17";
        "1";"13" **\dir{-}; "3";"6" **\dir{-}?(0.8)="x";
        "2";"9" **\dir{-}; "10";"8" **\dir{-}; "12";"9" **\dir{-};
        "11";"14" **\dir{-}; "12";"17" **\dir{-}; "x";"8" **\dir{-};
        \ar@{-}@/_2pt/ "12";"15" |\hole
        \ar@{-}@/^2pt/ "11";"16"
        \ar@{-}@/_14pt/ "4";"11" |<(0.28)\hole |<(0.84)\hole
    \endxy}
    ~\overset{\text{(c,6)}}{=}~
    \vcenter{\xy
        (-6,20)="1"; ( 6,20)="2"; (-6,17)="3"; ( 6,17)="4";
        (-2,10)="5"; ( 4,12)="6";
        ( 4,8)*[o]=<9pt>[Fo]{\scriptstyle \nu}="7";
        (-2,4)="8"; ( 4,4)="9"; (-6,-2)="10"; (-4,-2)="11";
        (-1,-2)="12"; ( 3,-2)="13"; ( 6,-2)="14";
        "1";"10" **\dir{-}; "3";"6" **\dir{-}; "2";"4" **\dir{-};
        "6";"4" **\dir{-}; "6";"7" **\dir{-}; "9";"7" **\dir{-};
        "9";"14" **\dir{-}; "5";"8" **\dir{-}; "11";"8" **\dir{-};
        \ar@{-}@/^4pt/ "5";"4" |<(0.4)\hole
        \ar@{-} "12";"9" |<(0.4)\hole
        \ar@{-} "13";"8"
    \endxy}~.
\]

Thus, $\theta$ is a left $A^{\ox 3}$-comodule morphism $P_l \ra P_r$. The
inverse of $\theta$ then is a left $A^{\ox 3}$-comodule morphism $P_r \ra
P_l$.

We now prove the properties (G1) through (G3) required of a quantum
groupoid. The calculation
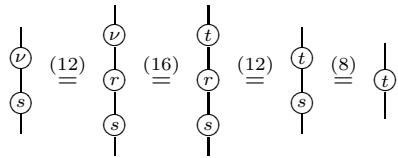
\[
    \vcenter{\xy
        (0,7)="1";
        (0,3)*[o]=<9pt>[Fo]{\scriptstyle \nu}="2";
        (0,-3)*[o]=<9pt>[Fo]{\scriptstyle s}="3";
        (0,-7)="4";
        "1";"2" **\dir{-};
        "2";"3" **\dir{-};
        "3";"4" **\dir{-};
    \endxy}
    ~\overset{\text{(12)}}{=}~
    \vcenter{\xy
        (0,7)="1";
        (0,3)*[o]=<9pt>[Fo]{\scriptstyle \nu}="2";
        (0,-3)*[o]=<9pt>[Fo]{\scriptstyle r}="3";
        (0,-9)*[o]=<9pt>[Fo]{\scriptstyle s}="4";
        (0,-13)="5";
        "1";"2" **\dir{-};
        "2";"3" **\dir{-};
        "3";"4" **\dir{-};
        "4";"5" **\dir{-};
    \endxy}
    ~\overset{\text{(16)}}{=}~
    \vcenter{\xy
        (0,7)="1";
        (0,3)*[o]=<9pt>[Fo]{\scriptstyle t}="2";
        (0,-3)*[o]=<9pt>[Fo]{\scriptstyle r}="3";
        (0,-9)*[o]=<9pt>[Fo]{\scriptstyle s}="4";
        (0,-13)="5";
        "1";"2" **\dir{-};
        "2";"3" **\dir{-};
        "3";"4" **\dir{-};
        "4";"5" **\dir{-};
    \endxy}
    ~\overset{\text{(12)}}{=}~
    \vcenter{\xy
        (0,7)="1";
        (0,3)*[o]=<9pt>[Fo]{\scriptstyle t}="2";
        (0,-3)*[o]=<9pt>[Fo]{\scriptstyle s}="3";
        (0,-7)="4";
        "1";"2" **\dir{-};
        "2";"3" **\dir{-};
        "3";"4" **\dir{-};
    \endxy}
    ~\overset{\text{(8)}}{=}~
    \vcenter{\xy
        (0,8)="1";
        (0,3)*[o]=<9pt>[Fo]{\scriptstyle t}="2";
        (0,-2)="3";
        "1";"2" **\dir{-};
        "2";"3" **\dir{-};
    \endxy}
\]
verifies (G1), and the following establishes (G2).
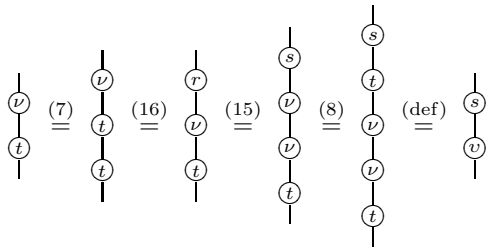
\[
    \vcenter{\xy
        (0,7)="1";
        (0,3)*[o]=<9pt>[Fo]{\scriptstyle \nu}="2";
        (0,-3)*[o]=<9pt>[Fo]{\scriptstyle t}="3";
        (0,-7)="4";
        "1";"2" **\dir{-};
        "2";"3" **\dir{-};
        "3";"4" **\dir{-};
    \endxy}
    ~\overset{\text{(7)}}{=}~
    \vcenter{\xy
        (0,7)="1";
        (0,3)*[o]=<9pt>[Fo]{\scriptstyle \nu}="2";
        (0,-3)*[o]=<9pt>[Fo]{\scriptstyle t}="3";
        (0,-9)*[o]=<9pt>[Fo]{\scriptstyle t}="4";
        (0,-13)="5";
        "1";"2" **\dir{-};
        "2";"3" **\dir{-};
        "3";"4" **\dir{-};
        "4";"5" **\dir{-};
    \endxy}
    ~\overset{\text{(16)}}{=}~
    \vcenter{\xy
        (0,7)="1";
        (0,3)*[o]=<9pt>[Fo]{\scriptstyle r}="2";
        (0,-3)*[o]=<9pt>[Fo]{\scriptstyle \nu}="3";
        (0,-9)*[o]=<9pt>[Fo]{\scriptstyle t}="4";
        (0,-13)="5";
        "1";"2" **\dir{-};
        "2";"3" **\dir{-};
        "3";"4" **\dir{-};
        "4";"5" **\dir{-};
    \endxy}
    ~\overset{\text{(15)}}{=}~
    \vcenter{\xy
        (0,7)="1";
        (0,3)*[o]=<9pt>[Fo]{\scriptstyle s}="2";
        (0,-3)*[o]=<9pt>[Fo]{\scriptstyle \nu}="3";
        (0,-9)*[o]=<9pt>[Fo]{\scriptstyle \nu}="4";
        (0,-15)*[o]=<9pt>[Fo]{\scriptstyle t}="5";
        (0,-19)="6";
        "1";"2" **\dir{-};
        "2";"3" **\dir{-};
        "3";"4" **\dir{-};
        "4";"5" **\dir{-};
        "5";"6" **\dir{-};
    \endxy}
    ~\overset{\text{(8)}}{=}~
    \vcenter{\xy
        (0,13)="0";
        (0,9)*[o]=<9pt>[Fo]{\scriptstyle s}="1";
        (0,3)*[o]=<9pt>[Fo]{\scriptstyle t}="2";
        (0,-3)*[o]=<9pt>[Fo]{\scriptstyle \nu}="3";
        (0,-9)*[o]=<9pt>[Fo]{\scriptstyle \nu}="4";
        (0,-15)*[o]=<9pt>[Fo]{\scriptstyle t}="5";
        (0,-19)="6";
        "0";"1" **\dir{-};
        "1";"2" **\dir{-};
        "2";"3" **\dir{-};
        "3";"4" **\dir{-};
        "4";"5" **\dir{-};
        "5";"6" **\dir{-};
    \endxy}
    ~\overset{\text{(def)}}{=}~
    \vcenter{\xy
        (0,7)="1";
        (0,3)*[o]=<9pt>[Fo]{\scriptstyle s}="2";
        (0,-3)*[o]=<9pt>[Fo]{\scriptstyle \upsilon}="3";
        (0,-7)="4";
        "1";"2" **\dir{-};
        "2";"3" **\dir{-};
        "3";"4" **\dir{-};
    \endxy}
\]
It remains to prove (G3), i.e., we must show that $\theta$ makes the
following square
\[
    \xygraph{{P}="1"
        [r(1.4)] {C^{\ox 3}}="2"
        [r(2)] {C^{\ox 3}}="3"
        [d(1.2)] {C^{\ox 3}}="4"
     "1"[d(1.2)] {P}="5"
        "1":"2" ^-\varsigma
        "2":"3" ^-{c_{C,C \ox C}}
        "3":"4" ^-{1 \ox 1 \ox \upsilon}
        "1":"5" _-\theta
        "5":"4" ^-\varsigma}
\]
commute.

The clockwise direction around the square is
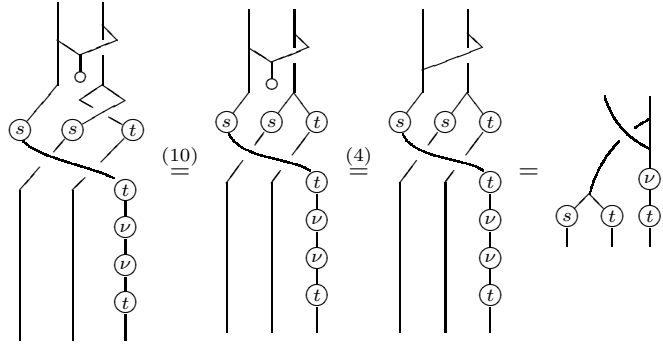
\[
    \vcenter{\xy
        (-3, 9)="1"; ( 3, 9)="2"; (-3, 4)="3";
        ( 3, 6)="4"; ( 5, 4)="5"; ( 0, 2)="6";
        ( 0,-1)*[o]=<5pt>[Fo]{~}="7";
        ( 3,-2)="8"; ( 0,-4)="9"; ( 6,-4)="10";
        (-8,-8)*[o]=<9pt>[Fo]{\scriptstyle s}="11";
        (-1,-8)*[o]=<9pt>[Fo]{\scriptstyle s}="12";
        ( 7,-8)*[o]=<9pt>[Fo]{\scriptstyle t}="13";
        (-8,-16)="14"; (-1,-16)="15";
        ( 6,-16)*[o]=<9pt>[Fo]{\scriptstyle t}="16";
        ( 6,-21)*[o]=<9pt>[Fo]{\scriptstyle \nu}="17";
        ( 6,-26)*[o]=<9pt>[Fo]{\scriptstyle \nu}="18";
        ( 6,-31)*[o]=<9pt>[Fo]{\scriptstyle t}="19";
        (-8,-36)="20"; (-1,-36)="21"; ( 6,-36)="22";
        "1";(-3,-2) **\dir{-};
        (-3,-2);"11" **\dir{-};
        "3";"6" **\dir{-}; "7";"6" **\dir{-}; "5";"6" **\dir{-};
        "5";"4" **\dir{-}; "8";"10" **\dir{-}; "8";"9" **\dir{-};
        "10";"12" **\dir{-}; "11";"16" **\crv{(-9,-12)&(6,-13)};
        "16";"17" **\dir{-}; "17";"18" **\dir{-}; "18";"19" **\dir{-};
        "19";"22" **\dir{-}; "14";"20" **\dir{-}; "15";"21" **\dir{-};
        \ar@{-} "2";"8" |<(0.55)\hole
        \ar@{-} "9";"13" |<(0.4)\hole
        \ar@{-} "12";"14" |<(0.4)\hole
        \ar@{-} "13";"15" |<(0.6)\hole
    \endxy}
    ~\overset{\text{(10)}}{=}~
    \vcenter{\xy
        (-3, 9)="1"; ( 3, 9)="2"; (-3, 4)="3"; ( 3, 6)="4";
        ( 5, 4)="5"; ( 0, 2)="6"; ( 0,-1)*[o]=<5pt>[Fo]{~}="7";
        ( 3,-2)="8"; ( 0,-4)="9"; ( 6,-4)="10";
        (-6,-6)*[o]=<9pt>[Fo]{\scriptstyle s}="11";
        ( 0,-6)*[o]=<9pt>[Fo]{\scriptstyle s}="12";
        ( 6,-6)*[o]=<9pt>[Fo]{\scriptstyle t}="13";
        (-6,-14)="14"; ( 0,-14)="15";
        ( 6,-14)*[o]=<9pt>[Fo]{\scriptstyle t}="16";
        ( 6,-19)*[o]=<9pt>[Fo]{\scriptstyle \nu}="17";
        ( 6,-24)*[o]=<9pt>[Fo]{\scriptstyle \nu}="18";
        ( 6,-29)*[o]=<9pt>[Fo]{\scriptstyle t}="19";
        (-6,-34)="20"; ( 0,-34)="21"; ( 6,-34)="22";
        "1";(-3,-2) **\dir{-}; (-3,-2);"11" **\dir{-}; "3";"6" **\dir{-};
        "7";"6" **\dir{-}; "5";"6" **\dir{-}; "5";"4" **\dir{-};
        "8";"12" **\dir{-}; "8";"13" **\dir{-};
        "11";"16" **\crv{(-6,-12)&(6,-10)};
        "16";"17" **\dir{-}; "17";"18" **\dir{-}; "18";"19" **\dir{-};
        "19";"22" **\dir{-}; "14";"20" **\dir{-}; "15";"21" **\dir{-};
        \ar@{-} "2";"8" |<(0.55)\hole
        \ar@{-} "12";"14" |<(0.4)\hole
        \ar@{-} "13";"15" |<(0.6)\hole
    \endxy}
    ~\overset{\text{(4)}}{=}~
    \vcenter{\xy
        (-3, 9)="1"; ( 3, 9)="2"; (-3, 1)="3"; ( 3, 6)="4";
        ( 5, 4)="5"; ( 3,-2)="8"; ( 0,-4)="9"; ( 6,-4)="10";
        (-6,-6)*[o]=<9pt>[Fo]{\scriptstyle s}="11";
        ( 0,-6)*[o]=<9pt>[Fo]{\scriptstyle s}="12";
        ( 6,-6)*[o]=<9pt>[Fo]{\scriptstyle t}="13";
        (-6,-14)="14";
        ( 0,-14)="15";
        ( 6,-14)*[o]=<9pt>[Fo]{\scriptstyle t}="16";
        ( 6,-19)*[o]=<9pt>[Fo]{\scriptstyle \nu}="17";
        ( 6,-24)*[o]=<9pt>[Fo]{\scriptstyle \nu}="18";
        ( 6,-29)*[o]=<9pt>[Fo]{\scriptstyle t}="19";
        (-6,-34)="20";
        ( 0,-34)="21";
        ( 6,-34)="22";
        "1";(-3,-2) **\dir{-};
        (-3,-2);"11" **\dir{-};
        "5";"3" **\dir{-};
        "5";"4" **\dir{-};
        "8";"12" **\dir{-};
        "8";"13" **\dir{-};
        "11";"16" **\crv{(-6,-12)&(6,-10)};
        "16";"17" **\dir{-};
        "17";"18" **\dir{-};
        "18";"19" **\dir{-};
        "19";"22" **\dir{-};
        "14";"20" **\dir{-};
        "15";"21" **\dir{-};
        \ar@{-} "2";"8" |<(0.55)\hole
        \ar@{-} "12";"14" |<(0.4)\hole
        \ar@{-} "13";"15" |<(0.6)\hole
    \endxy}
    ~=~
    \vcenter{\xy
        (-3, 9)="1";
        ( 3, 9)="2";
        ( 3, 6)="3";
        ( 3, 2)="4";
        (-5,-4)="5";
        ( 3,-2)*[o]=<9pt>[Fo]{\scriptstyle \nu}="6";
        (-8,-7)*[o]=<9pt>[Fo]{\scriptstyle s}="7";
        (-2,-7)*[o]=<9pt>[Fo]{\scriptstyle t}="8";
        ( 3,-7)*[o]=<9pt>[Fo]{\scriptstyle t}="9";
        (-8,-11)="10";
        (-2,-11)="11";
        ( 3,-11)="12";
        "2";"6" **\dir{-};
        "6";"9" **\dir{-};
        "9";"12" **\dir{-};
        "5";"7" **\dir{-};
        "5";"8" **\dir{-};
        "10";"7" **\dir{-};
        "11";"8" **\dir{-};
        \ar@{-}@/_3pt/ "1";"4"
        \ar@{-}@/^4pt/ "5";"3" |<(0.73)\hole
    \endxy}
\]
for which the last step holds since
\begin{align*}
t \nu \nu t s &= t \nu \nu s & \text{(8)} \\
              &= t \nu r & \text{(15)} \\
              &= t t \nu & \text{(16)} \\
              &= t \nu. & \text{(7)}
\end{align*}

The counter-clockwise direction is
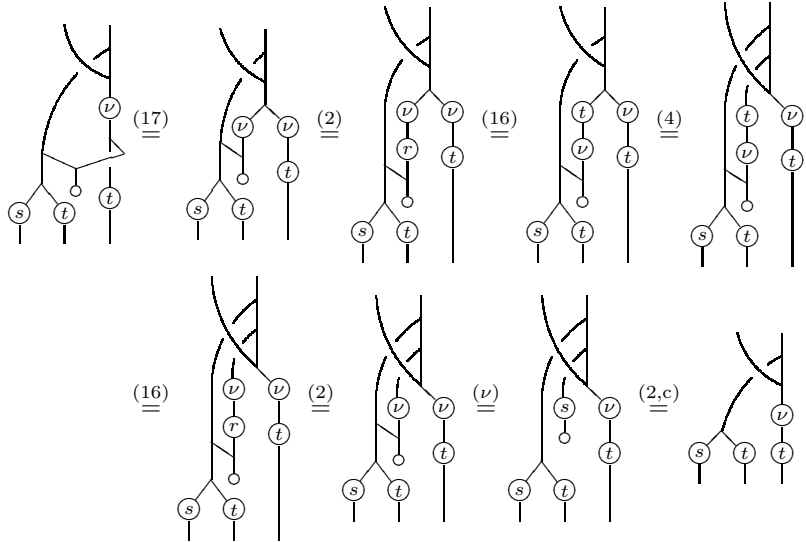
\begin{align*}
    \vcenter{\xy
        (-3, 9)="1"; ( 3, 9)="2"; ( 3, 6)="3"; ( 3, 2)="4";
        ( 3,-2)*[o]=<9pt>[Fo]{\scriptstyle \nu}="5";
        (-6,-8)="6"; ( 3,-6)="7"; ( 5,-8)="8"; (-1.5,-10)="9";
        (-1.5,-13)*[o]=<5pt>[Fo]{~}="10";
        (-6,-12)="11";
        (-9,-16)*[o]=<9pt>[Fo]{\scriptstyle s}="12";
        (-3,-16)*[o]=<9pt>[Fo]{\scriptstyle t}="13";
        ( 3,-14)*[o]=<9pt>[Fo]{\scriptstyle t}="14";
        (-9,-20)="15"; (-3,-20)="16"; ( 3,-20)="17";
        "2";"5" **\dir{-}; "14";"17" **\dir{-}; "6";"9" **\dir{-};
        "10";"9" **\dir{-}; "8";"9" **\dir{-}; "8";"7" **\dir{-};
        "6";"11" **\dir{-}; "12";"11" **\dir{-}; "13";"11" **\dir{-};
        "12";"15" **\dir{-}; "13";"16" **\dir{-};
        \ar@{-}@/_4pt/ "1";"4"
        \ar@{-}@/^6pt/ "6";"3" |<(0.78)\hole
        \ar@{-} "5";"14" |<(0.5)\hole
    \endxy}
    &\overset{\text{(17)}}{=}~
    \vcenter{\xy
        (-3, 9)="1"; ( 3, 9)="2"; ( 3, 6)="3"; ( 3, 2)="4"; ( 3,-1)="5";
        ( 0,-4)*[o]=<9pt>[Fo]{\scriptstyle \nu}="6";
        ( 6,-4)*[o]=<9pt>[Fo]{\scriptstyle \nu}="7";
        (-3,-6)="8"; ( 0,-8)="9";
        ( 0,-11)*[o]=<5pt>[Fo]{~}="10";
        ( 6,-10)*[o]=<9pt>[Fo]{\scriptstyle t}="11";
        (-3,-11)="12";
        (-6,-15)*[o]=<9pt>[Fo]{\scriptstyle s}="13";
        ( 0,-15)*[o]=<9pt>[Fo]{\scriptstyle t}="14";
        (-6,-19)="15"; ( 0,-19)="16"; ( 6,-19)="17";
        "2";"5" **\dir{-}; "5";"6" **\dir{-}; "6";"9" **\dir{-};
        "10";"9" **\dir{-}; "8";"9" **\dir{-}; "8";"12" **\dir{-};
        "13";"12" **\dir{-}; "14";"12" **\dir{-}; "13";"15" **\dir{-};
        "14";"16" **\dir{-}; "5";"7" **\dir{-}; "7";"11" **\dir{-};
        "17";"11" **\dir{-};
        \ar@{-}@/_3pt/ "1";"4"
        \ar@{-}@/^4pt/ "8";"3" |<(0.75)\hole
    \endxy}
    ~\overset{\text{(2)}}{=}
    \vcenter{\xy
        (-3,15)="1"; ( 3,15)="2"; ( 3,11)="3"; ( 3, 7)="4"; ( 3, 4)="5";
        ( 0, 1)*[o]=<9pt>[Fo]{\scriptstyle \nu}="6";
        ( 0,-4)*[o]=<9pt>[Fo]{\scriptstyle r}="x";
        ( 6, 1)*[o]=<9pt>[Fo]{\scriptstyle \nu}="7";
        (-3,-6)="8";
        ( 0,-8)="9";
        ( 0,-11)*[o]=<5pt>[Fo]{~}="10";
        ( 6,-5)*[o]=<9pt>[Fo]{\scriptstyle t}="11";
        (-3,-11)="12";
        (-6,-15)*[o]=<9pt>[Fo]{\scriptstyle s}="13";
        ( 0,-15)*[o]=<9pt>[Fo]{\scriptstyle t}="14";
        (-6,-19)="15"; ( 0,-19)="16"; ( 6,-19)="17";
        "2";"5" **\dir{-}; "5";"6" **\dir{-}; "6";"x" **\dir{-};
        "x";"9" **\dir{-}; "10";"9" **\dir{-}; "8";"9" **\dir{-};
        "8";"12" **\dir{-}; "8";(-3,2) **\dir{-}; "13";"12" **\dir{-};
        "14";"12" **\dir{-}; "13";"15" **\dir{-}; "14";"16" **\dir{-};
        "5";"7" **\dir{-}; "7";"11" **\dir{-}; "17";"11" **\dir{-};
        \ar@{-}@/_3pt/ "1";"4"
        \ar@{-}@/^4pt/ (-3,2);"3" |<(0.7)\hole
    \endxy}
    ~\overset{\text{(16)}}{=}
    \vcenter{\xy
        (-3,15)="1"; ( 3,15)="2"; ( 3,11)="3"; ( 3, 7)="4"; ( 3, 4)="5";
        ( 0, 1)*[o]=<9pt>[Fo]{\scriptstyle t}="6";
        ( 0,-4)*[o]=<9pt>[Fo]{\scriptstyle \nu}="x";
        ( 6, 1)*[o]=<9pt>[Fo]{\scriptstyle \nu}="7";
        (-3,-6)="8"; ( 0,-8)="9";
        ( 0,-11)*[o]=<5pt>[Fo]{~}="10";
        ( 6,-5)*[o]=<9pt>[Fo]{\scriptstyle t}="11";
        (-3,-11)="12";
        (-6,-15)*[o]=<9pt>[Fo]{\scriptstyle s}="13";
        ( 0,-15)*[o]=<9pt>[Fo]{\scriptstyle t}="14";
        (-6,-19)="15"; ( 0,-19)="16"; ( 6,-19)="17";
        "2";"5" **\dir{-}; "5";"6" **\dir{-}; "6";"x" **\dir{-};
        "x";"9" **\dir{-}; "10";"9" **\dir{-}; "8";"9" **\dir{-};
        "8";"12" **\dir{-}; "8";(-3,2) **\dir{-}; "13";"12" **\dir{-};
        "14";"12" **\dir{-}; "13";"15" **\dir{-}; "14";"16" **\dir{-};
        "5";"7" **\dir{-}; "7";"11" **\dir{-}; "17";"11" **\dir{-};
        \ar@{-}@/_3pt/ "1";"4"
        \ar@{-}@/^4pt/ (-3,2);"3" |<(0.7)\hole
    \endxy}
    ~\overset{\text{(4)}}{=}
    \vcenter{\xy
        (-3,16)="1"; ( 3,16)="2"; ( 3,13)="3"; ( 3, 9)="4"; ( 3, 4)="5";
        ( 0, 1)*[o]=<9pt>[Fo]{\scriptstyle t}="6";
        ( 0,-4)*[o]=<9pt>[Fo]{\scriptstyle \nu}="x";
        ( 6, 1)*[o]=<9pt>[Fo]{\scriptstyle \nu}="7";
        (-3,-6)="8"; ( 0,-8)="9";
        ( 0,-11)*[o]=<5pt>[Fo]{~}="10";
        ( 6,-5)*[o]=<9pt>[Fo]{\scriptstyle t}="11";
        (-3,-11)="12";
        (-6,-15)*[o]=<9pt>[Fo]{\scriptstyle s}="13";
        ( 0,-15)*[o]=<9pt>[Fo]{\scriptstyle t}="14";
        (-6,-19)="15"; ( 0,-19)="16"; ( 6,-19)="17";
        "2";"5" **\dir{-}; "6";"x" **\dir{-}; "x";"9" **\dir{-};
        "10";"9" **\dir{-}; "8";"9" **\dir{-}; "8";"12" **\dir{-};
        "8";(-3,2) **\dir{-}; "13";"12" **\dir{-}; "14";"12" **\dir{-};
        "13";"15" **\dir{-}; "14";"16" **\dir{-}; "5";"7" **\dir{-};
        "7";"11" **\dir{-}; "17";"11" **\dir{-};
        \ar@{-}@/^4pt/ "6";"4" |<(0.5)\hole
        \ar@{-}@/_4pt/ "1";"5"
        \ar@{-}@/^4pt/ (-3,2);"3" |<(0.55)\hole
    \endxy}
\\
    &\overset{\text{(16)}}{=}
    \vcenter{\xy
        (-3,16)="1"; ( 3,16)="2"; ( 3,13)="3"; ( 3, 9)="4"; ( 3, 4)="5";
        ( 0, 1)*[o]=<9pt>[Fo]{\scriptstyle \nu}="6";
        ( 0,-4)*[o]=<9pt>[Fo]{\scriptstyle r}="x";
        ( 6, 1)*[o]=<9pt>[Fo]{\scriptstyle \nu}="7";
        (-3,-6)="8"; ( 0,-8)="9";
        ( 0,-11)*[o]=<5pt>[Fo]{~}="10";
        ( 6,-5)*[o]=<9pt>[Fo]{\scriptstyle t}="11";
        (-3,-11)="12";
        (-6,-15)*[o]=<9pt>[Fo]{\scriptstyle s}="13";
        ( 0,-15)*[o]=<9pt>[Fo]{\scriptstyle t}="14";
        (-6,-19)="15"; ( 0,-19)="16"; ( 6,-19)="17";
        "2";"5" **\dir{-}; "6";"x" **\dir{-}; "x";"9" **\dir{-};
        "10";"9" **\dir{-}; "8";"9" **\dir{-}; "8";"12" **\dir{-};
        "8";(-3,2) **\dir{-}; "13";"12" **\dir{-}; "14";"12" **\dir{-};
        "13";"15" **\dir{-}; "14";"16" **\dir{-}; "5";"7" **\dir{-};
        "7";"11" **\dir{-};
        "17";"11" **\dir{-};
        \ar@{-}@/^4pt/ "6";"4" |<(0.5)\hole
        \ar@{-}@/_4pt/ "1";"5"
        \ar@{-}@/^4pt/ (-3,2);"3" |<(0.55)\hole
    \endxy}
    ~\overset{\text{(2)}}{=}
    \vcenter{\xy
        (-3,11)="1"; ( 3,11)="2"; ( 3, 8)="3"; ( 3, 4)="4"; ( 3,-1)="5";
        ( 0,-4)*[o]=<9pt>[Fo]{\scriptstyle \nu}="6";
        ( 6,-4)*[o]=<9pt>[Fo]{\scriptstyle \nu}="7";
        (-3,-6)="8"; ( 0,-8)="9";
        ( 0,-11)*[o]=<5pt>[Fo]{~}="10";
        ( 6,-10)*[o]=<9pt>[Fo]{\scriptstyle t}="11";
        (-3,-11)="12";
        (-6,-15)*[o]=<9pt>[Fo]{\scriptstyle s}="13";
        ( 0,-15)*[o]=<9pt>[Fo]{\scriptstyle t}="14";
        (-6,-19)="15"; ( 0,-19)="16"; ( 6,-19)="17";
        "2";"5" **\dir{-}; "6";"9" **\dir{-}; "10";"9" **\dir{-};
        "8";"9" **\dir{-}; "8";"12" **\dir{-}; "8";(-3,-2) **\dir{-};
        "13";"12" **\dir{-}; "14";"12" **\dir{-}; "13";"15" **\dir{-};
        "14";"16" **\dir{-}; "5";"7" **\dir{-}; "7";"11" **\dir{-};
        "17";"11" **\dir{-};
        \ar@{-}@/^4pt/ "6";"4" |<(0.5)\hole
        \ar@{-}@/_4pt/ "1";"5"
        \ar@{-}@/^4pt/ (-3,-2);"3" |<(0.52)\hole
    \endxy}
    ~\overset{\text{($\nu$)}}{=}
    \vcenter{\xy
        (-3,11)="1"; ( 3,11)="2"; ( 3, 8)="3"; ( 3, 4)="4"; ( 3,-1)="5";
        ( 0,-4)*[o]=<9pt>[Fo]{\scriptstyle s}="6";
        ( 6,-4)*[o]=<9pt>[Fo]{\scriptstyle \nu}="7";
        ( 0,-8)*[o]=<5pt>[Fo]{~}="10";
        ( 6,-10)*[o]=<9pt>[Fo]{\scriptstyle t}="11";
        (-3,-11)="12";
        (-6,-15)*[o]=<9pt>[Fo]{\scriptstyle s}="13";
        ( 0,-15)*[o]=<9pt>[Fo]{\scriptstyle t}="14";
        (-6,-19)="15"; ( 0,-19)="16"; ( 6,-19)="17";
        "2";"5" **\dir{-}; "6";"10" **\dir{-}; "12";(-3,-2) **\dir{-};
        "13";"12" **\dir{-}; "14";"12" **\dir{-}; "13";"15" **\dir{-};
        "14";"16" **\dir{-}; "5";"7" **\dir{-}; "7";"11" **\dir{-};
        "17";"11" **\dir{-};
        \ar@{-}@/^4pt/ "6";"4" |<(0.5)\hole
        \ar@{-}@/_4pt/ "1";"5"
        \ar@{-}@/^4pt/ (-3,-2);"3" |<(0.52)\hole
    \endxy}
    ~\overset{\text{(2,c)}}{=}
    \vcenter{\xy
        (-3, 9)="1";
        ( 3, 9)="2";
        ( 3, 6)="3";
        ( 3, 2)="4";
        (-5,-4)="5";
        ( 3,-2)*[o]=<9pt>[Fo]{\scriptstyle \nu}="6";
        (-8,-7)*[o]=<9pt>[Fo]{\scriptstyle s}="7";
        (-2,-7)*[o]=<9pt>[Fo]{\scriptstyle t}="8";
        ( 3,-7)*[o]=<9pt>[Fo]{\scriptstyle t}="9";
        (-8,-11)="10";
        (-2,-11)="11";
        ( 3,-11)="12";
        "2";"6" **\dir{-};
        "6";"9" **\dir{-};
        "9";"12" **\dir{-};
        "5";"7" **\dir{-};
        "5";"8" **\dir{-};
        "10";"7" **\dir{-};
        "11";"8" **\dir{-};
        \ar@{-}@/_3pt/ "1";"4"
        \ar@{-}@/^4pt/ "5";"3" |<(0.73)\hole
    \endxy}
\end{align*}
therefore establishing the commutativity of the square.

\begin{corollary}
Any Frobenius monoid in $\Q\V$ yields a quantum groupoid. 
\end{corollary}

By Proposition~\ref{prop-frob-wha} every Frobenius monoid $R$ in $\Q\V$
leads to a weak Hopf monoid with invertible antipode $R \ox R$. Apply
Proposition~\ref{prop-whm-qg} to this weak Hopf monoid with invertible
antipode to get a quantum groupoid.

%=========================================================================%
\appendix
%=========================================================================%

%=========================================================================%
\section{String diagrams and basic definitions}\label{sec-prelim}
%=========================================================================%

In this appendix we give a quick introduction to string
diagrams in a braided monoidal category $\V = (\V,\ox,I,c)$~\cite{JS:braided}
and use these to define monoid, module,
comonoid, comodule, and separable Frobenius monoid in $\V$. The string
calculus was shown to be rigorous in~\cite{JS:geometry}.

%=========================================================================%
\subsection{String diagrams}
%=========================================================================%

Suppose that $\V = (\V,\ox,I,c)$ is a braided (strict) monoidal category.
In a string diagram, objects label edges and morphisms label nodes. For
example, if $f:A \ox B \ra C \ox D \ox E$ is a morphism in $\V$ it is
represented as
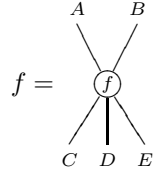
\[
    f =
    \vcenter{\xy
        (-4,10)*{\scriptstyle A};
        (4,10)*{\scriptstyle B};
        (-5,-10)*{\scriptstyle C};
        (0,-10)*{\scriptstyle D};
        (5,-10)*{\scriptstyle E};
        (-4,8)="tl";
        (4,8)="tr";
        (0,0)*[o]=<11pt>[Fo]{\scriptstyle f}="f";
        (-5,-8)="bl";
        (0,-8)="bc";
        (5,-8)="br";
        "tl";"f" **\dir{-};
        "tr";"f" **\dir{-};
        "f";"bl" **\dir{-};
        "f";"bc" **\dir{-};
        "f";"br" **\dir{-};
    \endxy}
\]
where this diagram is here meant to be read top-to-bottom. The identity
morphism on an object will be represented as the object itself as in
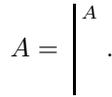
\[
    A =~
    \vcenter{\xy
        (2,5)*{\scriptstyle A};
        (0,6)="t";
        (0,-6)="b";
        **\dir{-}
    \endxy}~.
\]
A special case is the object $I \in \V$ which is represented as the empty
edge.

If, in $\V$, there are morphisms $f:A \ox B \ra C \ox D \ox E$ and
$g:D \ox E \ox F \ra G \ox H$ then they may be composed as
\[
    \xygraph{{A \ox B \ox F}
        :[r(2.6)] {C \ox D \ox E \ox F} ^-{f \ox 1}
        :[r(2.6)] {C \ox G \ox H} ^-{1 \ox g}}
\]
which may be represented as vertical concatenation
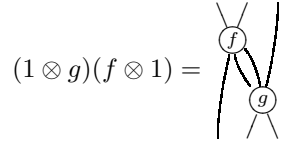
\[
    (1 \ox g)(f \ox 1) =~
    \vcenter{\xy
        (-4,9)="tl";
        (0,9)="tm";
        (4,9)="tr";
        (-2,4)*[o]=<11pt>[Fo]{\scriptstyle f}="f";
        (2,-4)*[o]=<11pt>[Fo]{\scriptstyle g}="g";
        (-4,-9)="bl";
        (0,-9)="bm";
        (4,-9)="br";
        "tl";"f" **\dir{-};
        "tm";"f" **\dir{-}; 
        "tr";"g" **\crv{(4,4)};
        "bl";"f" **\crv{(-4,-4)};
        "g";"bm" **\dir{-};
        "g";"br" **\dir{-};
        "f";"g" **\crv{(1,3)};
        "f";"g" **\crv{(-2,0)};
    \endxy}
\]
(where we have left off the objects). The tensor product of morphisms, say
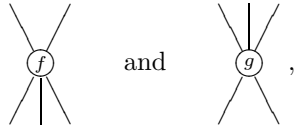
\[
    \vcenter{\xy
        (-4,8)="tl";
        (4,8)="tr";
        (0,0)*[o]=<11pt>[Fo]{\scriptstyle f}="f";
        (-4,-8)="bl";
        (0,-8)="bc";
        (4,-8)="br";
        "tl";"f" **\dir{-};
        "tr";"f" **\dir{-};
        "f";"bl" **\dir{-};
        "f";"bc" **\dir{-};
        "f";"br" **\dir{-};
    \endxy}
    \qquad \text{and} \qquad
    \vcenter{\xy
        (-4,8)="tl";
        (0,8)="tc";
        (4,8)="tr";
        (0,0)*[o]=<11pt>[Fo]{\scriptstyle g}="f";
        (-4,-8)="bl";
        (4,-8)="br";
        "tl";"f" **\dir{-};
        "tc";"f" **\dir{-};
        "tr";"f" **\dir{-};
        "f";"bl" **\dir{-};
        "f";"br" **\dir{-};
    \endxy}~,
\]
is represented as horizontal juxtaposition
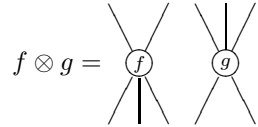
\[
    f \ox g =
    \vcenter{\xy
        (-4,8)="tl";
        (4,8)="tr";
        (0,0)*[o]=<11pt>[Fo]{\scriptstyle f}="f";
        (-4,-8)="bl";
        (0,-8)="bc";
        (4,-8)="br";
        "tl";"f" **\dir{-};
        "tr";"f" **\dir{-};
        "f";"bl" **\dir{-};
        "f";"bc" **\dir{-};
        "f";"br" **\dir{-};
    \endxy}
    \quad
    \vcenter{\xy
        (-4,8)="tl";
        (0,8)="tc";
        (4,8)="tr";
        (0,0)*[o]=<11pt>[Fo]{\scriptstyle g}="f";
        (-4,-8)="bl";
        (4,-8)="br";
        "tl";"f" **\dir{-};
        "tc";"f" **\dir{-};
        "tr";"f" **\dir{-};
        "f";"bl" **\dir{-};
        "f";"br" **\dir{-};
    \endxy}
\]
(again leaving off the objects).

The braiding $c_{A,B}:A \ox B \ra B \ox A$ is represented as a
left-over-right crossing. The inverse braiding is then represented as a
right-over-left crossing.
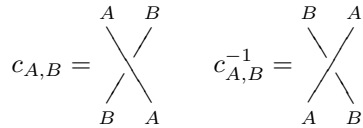
\[
    c_{A,B} =
    \vcenter{\xy
        (-3, 7)*{\scriptstyle A};
        ( 3, 7)*{\scriptstyle B};
        (-3,-7)*{\scriptstyle B};
        ( 3,-7)*{\scriptstyle A};
        (-3, 5)="tl";
        ( 3, 5)="tr";
        (-3,-5)="bl";
        ( 3,-5)="br";
        "tl";"br" **\crv{};\POS?(.5)*{\hole}="x";
        "tr";"x" **\dir{-};
        "x";"bl" **\dir{-};
    \endxy}
    \qquad
    c^{-1}_{A,B} =
    \vcenter{\xy
        (-3, 7)*{\scriptstyle B};
        ( 3, 7)*{\scriptstyle A};
        (-3,-7)*{\scriptstyle A};
        ( 3,-7)*{\scriptstyle B};
        (-3,5)="tl";
        ( 3, 5)="tr";
        (-3,-5)="bl";
        ( 3,-5)="br";
        "tr";"bl" **\crv{}; ?*{\hole}="x";
        "tl";"x" **\dir{-};
        "x";"br" **\dir{-};
    \endxy}
\]

Suppose $A \in \V$ has a left dual $A^*$, which we denote by $A^* \dashv A$
and say that $A^*$ is the left adjoint of $A$ (it is an adjunction if we were
to view $\V$ as a one object bicategory). The evaluation and coevaluation
morphisms $e_A:A^* \ox A \ra I$ and $n_A:I \ra A \ox A^*$ are represented as
\[
    e_A =~
    \vcenter{\xy
        (-3,5)*{\scriptstyle A^*};
        ( 3,5)*{\scriptstyle A};
        (-3,3);(3,3) **\crv{(-3,-5)&(3,-5)};
    \endxy}
    \quad\qquad \text{and} \quad\qquad
    n_A =~
    \vcenter{\xy
        (-3,-5)*{\scriptstyle A};
        ( 3,-5)*{\scriptstyle A^*};
        (-3,-3);(3,-3) **\crv{(-3,5)&(3,5)};
    \endxy}~.
\]
The triangle equalities become
\[
    \vcenter{\xy
        (4,10)*{\scriptstyle A};
        (-4,-10)*{\scriptstyle A};
        (2,1)*{\scriptstyle A^*};
        (4,8)="t";
        (4,0)*{}="r";
        (-0,0)*{}="m";
        (-4,0)*{}="l";
        (-4,-8)="b";
        "b";"l" **\dir{-};
        "r";"t" **\dir{-};
        "l";"m" **\crv{(-4,4)&(0,4)};
        "m";"r" **\crv{(0,-4)&(4,-4)};
    \endxy}
    \quad =~
    \vcenter{\xy
        (-2,7)*{\scriptstyle A};
        (0,8)="t";
        (0,-8)="b";
        "t";"b" **\dir{-};
    \endxy}
    \quad\qquad \text{and} \quad\qquad
    \vcenter{\xy
        (-4,10)*{\scriptstyle A^*};
        (4,-10)*{\scriptstyle A^*};
        (-1.5,1)*{\scriptstyle A};
        (-4,8)="t";
        (-4,0)="l";
        (-0,0)="m";
        (4,0)="r";
        (4,-8)="b";
        "t";"l" **\dir{-};
        "r";"b" **\dir{-};
        "l";"m" **\crv{(-4,-4)&(0,-4)};
        "m";"r" **\crv{(0,4)&(4,4)};
    \endxy}
    \quad =~
    \vcenter{\xy
        (-2,7)*{\scriptstyle A^*};
        (0,8)="t";
        (0,-8)="b";
        "t";"b" **\dir{-};
    \endxy}
\]

In what follows, in order to simplify the string diagrams, the nodes will be
omitted from certain morphisms (e.g., multiplication and comultiplication
morphisms) or simplified (e.g., unit and counit morphisms).

%=========================================================================%
\subsection{Monoids and modules}
%=========================================================================%

A \emph{monoid} $A = (A,\mu,\eta)$ in $\V$ is an object $A \in \V$ equipped
with morphisms
\[
    \mu =~
    \vcenter{\xy
        (-3,5)="l";
        (3,5)="r";
        (0,0)="m";
        (0,-5)="b";
        "l";"m" **\dir{-};
        "r";"m" **\dir{-};
        "m";"b" **\dir{-};
    \endxy}
    ~:A \ox A \dra A
    \qquad \text{and} \qquad
    \eta =~
    \vcenter{\xy
        (0,2)*[o]=<5pt>[Fo]{}="t";
        (0,-3)="b";
        "t";"b" **\dir{-};
    \endxy}
    ~:I \dra A,
\]
called the \emph{multiplication} and \emph{unit} of the monoid respectively,
satisfying
\begin{equation*}\tag{m}
    \vcenter{\xy %0;/r.16pc/: % associativity
        (0,6)*{~}; (0,-6)*{~};
        (-4,5)="l";
        (0,5)="c";
        (4,5)="r";
        (0,0)="m";
        (0,-5)="b";
        "l";"m" **\dir{-};
        "c";"m" **\dir{-};
        "r";"m" **\dir{-};
        "m";"b" **\dir{-};
    \endxy}
    ~=~
    \vcenter{\xy %0;/r.16pc/: % associativity
        (-4,5)="l";
        (0,5)="c";
        (4,5)="r";
        (0,0)="m";
        (0,-5)="b";
        "l";"m" **\dir{-};
        ?;"c" **\dir{-};
        "r";"m" **\dir{-};
        "m";"b" **\dir{-};
    \endxy}
    ~=~
    \vcenter{\xy %0;/r.16pc/:
        (-4,5)="l";
        (0,5)="c";
        (4,5)="r";
        (0,0)="m";
        (0,-5)="b";
        "l";"m" **\dir{-};
        "r";"m" **\dir{-};
        ?;"c" **\dir{-};
        "m";"b" **\dir{-};
    \endxy}
    \qquad \text{and} \qquad
    \vcenter{\xy %0;/r.16pc/: % counit
        (-2,3)*[o]=<5pt>[Fo]{~}="l";
        (4,5)="r";
        (0,0)="m";
        (0,-5)="b";
        "l";"m" **\dir{-};
        "r";"m" **\dir{-};
        "m";"b" **\dir{-};
    \endxy}
    ~=~
    \vcenter{\xy
        (0,-5)="b";
        (0,5)="t";
        "t";"b" **\dir{-};
    \endxy}
    ~=~
    \vcenter{\xy
        (-4,5)="l";
        (2,3)*[o]=<5pt>[Fo]{~}="r";
        (0,0)="m";
        (0,-5)="b";
        "l";"m" **\dir{-};
        "r";"m" **\dir{-};
        "m";"b" **\dir{-};
    \endxy}~.
\end{equation*}

If $A,B$ are monoids, a \emph{monoid morphism} $f:A \ra B$ is a morphism
in $\V$ satisfying
\[
    \vcenter{\xy
        (-3,9)*{\scriptstyle A};
        (3,9)*{\scriptstyle A};
        (-2,-7)*{\scriptstyle B};
        (-3,7)="l";
        (3,7)="r";
        (0,3)="m";
        (0,-2)*[o]=<11pt>[Fo]{\scriptstyle f}="f";
        (0,-7)="b";
        "l";"m" **\dir{-};
        "r";"m" **\dir{-};
        "m";"f" **\dir{-};
        "f";"b" **\dir{-};
    \endxy}
    ~=~
    \vcenter{\xy
        (-3,9)*{\scriptstyle A};
        (3,9)*{\scriptstyle A};
        (-2,-8)*{\scriptstyle B};
        (-3,7)="l";
        (3,7)="r";
        (-3,2)*[o]=<11pt>[Fo]{\scriptstyle f}="f1";
        (3,2)*[o]=<11pt>[Fo]{\scriptstyle f}="f2";
        (0,-3)="m";
        (0,-7)="b";
        "l";"f1" **\dir{-};
        "f1";"m" **\dir{-};
        "r";"f2" **\dir{-};
        "f2";"m" **\dir{-};
        "m";"b" **\dir{-};
    \endxy}
    \qquad \text{and} \qquad
    \vcenter{\xy % f \eta
        (-2,4)*{\scriptstyle A};
        (-2,-5)*{\scriptstyle B};
        (0,6)*[o]=<5pt>[Fo]{}="t";
        (0,0)*[o]=<11pt>[Fo]{\scriptstyle f}="f";
        (0,-5)="b";
        "t";"f" **\dir{-};
        "f";"b" **\dir{-};
    \endxy}
    ~=~
    \vcenter{\xy % \eta
        (-2,-3)*{\scriptstyle B};
        (0,2)*[o]=<5pt>[Fo]{}="t";
        (0,-3)="b";
        "t";"b" **\dir{-};
    \endxy}~.
\]

Monoids make sense in any monoidal category, however, in order that the
tensor product $A \ox B$ of monoids $A,B \in \V$ is again a monoid there
must be a ``switch'' morphism $c_{A,B}:A \ox B \ra B \ox A$ in $\V$
given by, say, a braiding. In this case $A \ox B$ becomes a monoid in $\V$ via
\[
    \mu =~
    \vcenter{\xy
        (-5,5)="l1";
        (-2,5)="l2";
        (1,5)="r1";
        (4,5)="r2";
        (-2,-1)="m1";
        (1,-1)="m2";
        (-2,-5)="b1";
        (1,-5)="b2";
        "l1";"m1" **\dir{-};
        "l2";"m2" **\dir{-};?(.5)*{\hole}="x";
        "r1";"x" **\dir{-};
        "x";"m1" **\dir{-};
        "r2";"m2" **\dir{-};
        "m1";"b1" **\dir{-};
        "m2";"b2" **\dir{-};
    \endxy}
    \qquad \text{and} \qquad
    \eta =
    \vcenter{\xy % \eta
        (-2,2)*[o]=<5pt>[Fo]{}="t";
        (-2,-3)="b";
        **\dir{-};
        (2,2)*[o]=<5pt>[Fo]{}="t";
        (2,-3)="b";
        **\dir{-};
    \endxy}~.
\]

Suppose that $A$ is a monoid in $\V$. A \emph{right $A$-module} in $\V$ is
an object $M \in \V$ equipped with a morphism
\[
    \mu =~
    \vcenter{\xy
        (-3,7)*{\scriptstyle M};
        (3,7)*{\scriptstyle A};
        (-2,-5)*{\scriptstyle M};
        (-3,5)="l";
        (3,5)="r";
        (0,0)="m";
        (0,-5)="b";
        "l";"m" **\dir{-};
        "r";"m" **\dir{-};
        "m";"b" **\dir{-};
    \endxy}
    ~:M \ox A \dra M
\]
called the \emph{action of $A$ on $M$} satisfying
\begin{equation*}\tag{m}
    \vcenter{\xy
        (-4,7)*{\scriptstyle M};
        (0,7)*{\scriptstyle A};
        (4,7)*{\scriptstyle A};
        (-2,-5)*{\scriptstyle M};
        (-4,5)="l";
        (0,5)="c";
        (4,5)="r";
        (0,0)="m";
        (0,-5)="b";
        "l";"m" **\dir{-};
        ?;"c" **\dir{-};
        "r";"m" **\dir{-};
        "m";"b" **\dir{-};
    \endxy}
    ~=~
    \vcenter{\xy
        (-4,7)*{\scriptstyle M};
        (0,7)*{\scriptstyle A};
        (4,7)*{\scriptstyle A};
        (-2,-5)*{\scriptstyle M};
        (-4,5)="l";
        (0,5)="c";
        (4,5)="r";
        (0,0)="m";
        (0,-5)="b";
        "l";"m" **\dir{-};
        "r";"m" **\dir{-};
        ?;"c" **\dir{-};
        "m";"b" **\dir{-};
    \endxy}
    \qquad \text{and} \qquad
    \vcenter{\xy
        (-2,7)*{\scriptstyle M};
        (3,1)*{\scriptstyle A};
        (-2,-4)*{\scriptstyle M};
        (-2,5)="l";
        (2,3)*[o]=<5pt>[Fo]{~}="r";
        (0,0)="m";
        (0,-5)="b";
        "l";"m" **\dir{-};
        "r";"m" **\dir{-};
        "m";"b" **\dir{-};
    \endxy}
    ~=~
    \vcenter{\xy
        (-2,4)*{\scriptstyle M};
        (0,-5)="b";
        (0,5)="t";
        "t";"b" **\dir{-};
    \endxy}~.
\end{equation*}
Notice that we use the same label ``(m)'' as the monoid axioms (and ``(c)''
below for the comodule axioms). This should not cause any confusion as the
labelling of strings disambiguates a multiplication and an action; however, the
labelling will usually be left off.

If $M,N$ are modules, a \emph{module morphism} $f:M \ra N$ is a morphism in
$\V$ satisfying
\[
    \vcenter{\xy
        (-3,9)*{\scriptstyle M};
        (3,9)*{\scriptstyle A};
        (-2,-7)*{\scriptstyle N};
        (-3,7)="l";
        (3,7)="r";
        (0,3)="m";
        (0,-2)*[o]=<11pt>[Fo]{\scriptstyle f}="f";
        (0,-7)="b";
        "l";"m" **\dir{-};
        "r";"m" **\dir{-};
        "m";"f" **\dir{-};
        "f";"b" **\dir{-};
    \endxy}
    ~=~
    \vcenter{\xy
        (-3,9)*{\scriptstyle M};
        (3,9)*{\scriptstyle A};
        (-2,-7)*{\scriptstyle N};
        (-3,7)="l";
        (3,7)="r";
        (-3,3)*[o]=<11pt>[Fo]{\scriptstyle f}="f";
        (0,-2)="m";
        (0,-7)="b";
        "l";"f" **\dir{-};
        "f";"m" **\dir{-};
        "r";"m" **\dir{-};
        "m";"b" **\dir{-};
    \endxy}~.
\]

%=========================================================================%
\subsection{Comonoids and comodules}
%=========================================================================%

Comonoids and comodules are dual to monoids and modules. Explicitly, a
\emph{comonoid} $C = (C,\delta,\epsilon)$ in $\V$ is an object $C \in \V$
equipped with morphisms
\[
    \delta =~
    \vcenter{\xy
        (0,5)="t";
        (0,0)="m";
        (-3,-5)="l";
        (3,-5)="r";
        "t";"m" **\dir{-};
        "m";"l" **\dir{-};
        "m";"r" **\dir{-};
    \endxy}
    ~:A \dra A \ox A
    \qquad \text{and} \qquad
    \epsilon =~
    \vcenter{\xy
        (0,3)="t";
        (0,-2)*[o]=<5pt>[Fo]{}="b";
        "t";"b" **\dir{-};
    \endxy}
    ~:A \dra I,
\]
called the \emph{comultiplication} and \emph{counit} of the comonoid
respectively, satisfying
\begin{equation*}\tag{c}
    \vcenter{\xy
        (0,5)="t";
        (0,0)="m";
        (-4,-5)="l";
        (0,-5)="c";
        (4,-5)="r";
        "t";"m" **\dir{-};
        "m";"l" **\dir{-};
        "m";"c" **\dir{-};
        "m";"r" **\dir{-};
    \endxy}
    ~=~
    \vcenter{\xy
        (0,5)="t";
        (0,0)="m";
        (-4,-5)="l";
        (0,-5)="c";
        (4,-5)="r";
        "t";"m" **\dir{-};
        "m";"l" **\dir{-};
        ?;"c" **\dir{-};
        "m";"r" **\dir{-};
    \endxy}
    ~=~
    \vcenter{\xy
        (0,5)="t";
        (0,0)="m";
        (-4,-5)="l";
        (0,-5)="c";
        (4,-5)="r";
        "t";"m" **\dir{-};
        "m";"l" **\dir{-};
        "m";"r" **\dir{-};
        ?;"c" **\dir{-};
    \endxy}
    \qquad \text{and} \qquad
    \vcenter{\xy
        (0,5)="t";
        (0,0)="m";
        (-2,-3)*[o]=<5pt>[Fo]{~}="l";
        (4,-5)="r";
        "t";"m" **\dir{-};
        "m";"l" **\dir{-};
        "m";"r" **\dir{-};
    \endxy}
    ~=~
    \vcenter{\xy
        (0,-5)="b";
        (0,5)="t";
        "t";"b" **\dir{-};
    \endxy}
    ~=~
    \vcenter{\xy
        (0,5)="t";
        (0,0)="m";
        (-4,-5)="l";
        (2,-3)*[o]=<5pt>[Fo]{~}="r";
        "t";"m" **\dir{-};
        "m";"l" **\dir{-};
        "m";"r" **\dir{-};
    \endxy}~.
\end{equation*}

If $C,D$ are comonoids, a \emph{comonoid morphism} $f:C \ra D$ is a
morphism in $\V$ satisfying
\[
    \vcenter{\xy
        (-3,-9)*{\scriptstyle B};
        (3,-9)*{\scriptstyle B};
        (-2,7)*{\scriptstyle A};
        (-3,-7)="l";
        (3,-7)="r";
        (0,-3)="m";
        (0,2)*[o]=<11pt>[Fo]{\scriptstyle f}="f";
        (0,7)="b";
        "l";"m" **\dir{-};
        "r";"m" **\dir{-};
        "m";"f" **\dir{-};
        "f";"b" **\dir{-};
    \endxy}
    ~=~
    \vcenter{\xy
        (-3,-9)*{\scriptstyle B};
        (3,-9)*{\scriptstyle B};
        (-2,8)*{\scriptstyle A};
        (-3,-7)="l";
        (3,-7)="r";
        (-3,-2)*[o]=<11pt>[Fo]{\scriptstyle f}="f1";
        (3,-2)*[o]=<11pt>[Fo]{\scriptstyle f}="f2";
        (0,3)="m";
        (0,7)="b";
        "l";"f1" **\dir{-};
        "f1";"m" **\dir{-};
        "r";"f2" **\dir{-};
        "f2";"m" **\dir{-};
        "m";"b" **\dir{-};
    \endxy}
    \qquad \text{and} \qquad
    \vcenter{\xy % \epsilon f
        (-2,5)*{\scriptstyle A};
        (-2,-4)*{\scriptstyle B};
        (0,5)="t";
        (0,0)*[o]=<11pt>[Fo]{\scriptstyle f}="f";
        (0,-6)*[o]=<5pt>[Fo]{}="b";
        "t";"f" **\dir{-};
        "f";"b" **\dir{-};
    \endxy}
    ~=~
    \vcenter{\xy % \epsilon
        (-2,3)*{\scriptstyle A};
        (0,3)="b";
        (0,-2)*[o]=<5pt>[Fo]{}="t";
        "t";"b" **\dir{-};
    \endxy}~.
\]

Similarly here, $\V$ must contain a switch morphism $c_{C,D}:C \ox D \ra
D \ox C$ in order that the tensor product $C \ox D$ of comonoids $C,D \in
\V$ is again a comonoid. In this case the comultiplication and counit are
given by
\[
    \delta = ~
    \vcenter{\xy % tensor product multiplication
        (-5,-5)="l1";
        (-2,-5)="l2";
        (1,-5)="r1";
        (4,-5)="r2";
        (-2,1)="m1";
        (1,1)="m2";
        (-2,5)="b1";
        (1,5)="b2";
        "m1";"l1" **\dir{-};
        "m1";"r1" **\dir{-};?(.5)*{\hole}="x";
        "m2";"r2" **\dir{-};
        "m2";"x" **\dir{-};
        "x";"l2" **\dir{-};
        "m1";"b1" **\dir{-};
        "m2";"b2" **\dir{-};
    \endxy}
    \qquad \text{and} \qquad
    \epsilon =~
    \vcenter{\xy % \eta
        (-2,-2)*[o]=<5pt>[Fo]{}="t";
        (-2,3)="b";
        **\dir{-};
        (2,-2)*[o]=<5pt>[Fo]{}="t";
        (2,3)="b";
        **\dir{-};
    \endxy}~.
\]

Suppose that $C$ is a comonoid in $\V$. A \emph{right $C$-comodule} in $\V$
is an object $M \in \V$ equipped with a morphism
\[
    \gamma =~
    \vcenter{\xy
        (-2,5)*{\scriptstyle M};
        (-3,-7)*{\scriptstyle M};
        (3,-7)*{\scriptstyle C};
        (0,5)="t";
        (0,0)="m";
        (-3,-5)="l";
        (3,-5)="r";
        "t";"m" **\dir{-};
        "m";"l" **\dir{-};
        "m";"r" **\dir{-};
    \endxy}
    ~:M \dra M \ox C
\]
called the \emph{coaction of $A$ on $M$} satisfying
\begin{equation*}\tag{c}
    \vcenter{\xy
        (-4,-7)*{\scriptstyle M};
        (0,-7)*{\scriptstyle C};
        (4,-7)*{\scriptstyle C};
        (-2,5)*{\scriptstyle M};
        (0,5)="t";
        (0,0)="m";
        (-4,-5)="l";
        (0,-5)="c";
        (4,-5)="r";
        "t";"m" **\dir{-};
        "m";"r" **\dir{-};
        ?;"c" **\dir{-};
        "m";"l" **\dir{-};
    \endxy}
    ~=~
    \vcenter{\xy
        (-4,-7)*{\scriptstyle M};
        (0,-7)*{\scriptstyle C};
        (4,-7)*{\scriptstyle C};
        (-2,5)*{\scriptstyle M};
        (0,5)="t";
        (0,0)="m";
        (-4,-5)="l";
        (0,-5)="c";
        (4,-5)="r";
        "t";"m" **\dir{-};
        "m";"l" **\dir{-};
        ?;"c" **\dir{-};
        "m";"r" **\dir{-};
    \endxy}
    \qquad \text{and} \qquad
    \vcenter{\xy
        (-2,-7)*{\scriptstyle M};
        (3,-1)*{\scriptstyle C};
        (-2,5)*{\scriptstyle M};
        (0,5)="t";
        (0,0)="m";
        (-2,-5)="l";
        (2,-3)*[o]=<5pt>[Fo]{~}="r";
        "t";"m" **\dir{-};
        "m";"l" **\dir{-};
        "m";"r" **\dir{-};
    \endxy}
    ~=~
    \vcenter{\xy
        (-2,5)*{\scriptstyle M};
        (0,5)="b";
        (0,-5)="t";
        "t";"b" **\dir{-};
    \endxy}~.
\end{equation*}

If $M,N$ are $C$-comodules, a \emph{comodule morphism} $f:M \ra N$ is a
morphism in $\V$ satisfying
\[
    \vcenter{\xy
        (-3,-9)*{\scriptstyle N};
        (3,-9)*{\scriptstyle C};
        (-2,7)*{\scriptstyle M};
        (-3,-7)="l";
        (3,-7)="r";
        (0,-3)="m";
        (0,2)*[o]=<11pt>[Fo]{\scriptstyle f}="f";
        (0,7)="b";
        "l";"m" **\dir{-};
        "r";"m" **\dir{-};
        "m";"f" **\dir{-};
        "f";"b" **\dir{-};
    \endxy}
    ~=~
    \vcenter{\xy
        (-3,-9)*{\scriptstyle N};
        (3,-9)*{\scriptstyle C};
        (-2,7)*{\scriptstyle M};
        (-3,-7)="l";
        (3,-7)="r";
        (-3,-3)*[o]=<11pt>[Fo]{\scriptstyle f}="f";
        (0,2)="m";
        (0,7)="b";
        "l";"f" **\dir{-};
        "f";"m" **\dir{-};
        "r";"m" **\dir{-};
        "m";"b" **\dir{-};
    \endxy}~.
\]

In this paper we also make use of $C$-bicomodules. Suppose that $M$ is
both a left $C$-comodule and a right $C$-comodule with coactions
\begin{align*}
    \gamma_l &: M \dra C \ox M \\
    \gamma_r &: M \dra M \ox C.
\end{align*}
If the square
\[
    \xygraph{{M}="1"
        [r(2.4)] {C \ox M}="2"
        "1"[d(1.2)] {M \ox C}="3"
        "2"[d(1.2)] {C \ox M \ox C}="4"
        "1":"2" ^-{\gamma_l}
        "1":"3" _-{\gamma_r}
        "2":"4" ^-{1 \ox \gamma_r}
        "3":"4" ^-{\gamma_l \ox 1}}
\]
commutes, meaning
\[
    \vcenter{\xy
        (0,-9)*{~}; (0, 7)*{~};
        (-4,-7)*{\scriptstyle M};
        (0,-7)*{\scriptstyle C};
        (4,-7)*{\scriptstyle M};
        (-2,5)*{\scriptstyle C};
        (0,5)="t";
        (0,0)="m";
        (-4,-5)="l";
        (0,-5)="c";
        (4,-5)="r";
        "t";"m" **\dir{-};
        "m";"r" **\dir{-};
        ?;"c" **\dir{-};
        "m";"l" **\dir{-};
    \endxy}
    ~=~
    \vcenter{\xy
        (-4,-7)*{\scriptstyle M};
        (0,-7)*{\scriptstyle C};
        (4,-7)*{\scriptstyle M};
        (-2,5)*{\scriptstyle C};
        (0,5)="t";
        (0,0)="m";
        (-4,-5)="l";
        (0,-5)="c";
        (4,-5)="r";
        "t";"m" **\dir{-};
        "m";"l" **\dir{-};
        ?;"c" **\dir{-};
        "m";"r" **\dir{-};
    \endxy}
\]
in string diagrams, then $M$ is called a \emph{$C$-bicomodule}. The diagonal
of the square will be denoted by
\[
    \gamma:M \dra C \ox M \ox C.
\]

%=========================================================================%
\subsection{Frobenius monoids}
%=========================================================================%

A \emph{Frobenius monoid} $R$ in $\V$ is both a monoid and a comonoid
in $\V$ which additionally satisfies the ``Frobenius condition'':
\[
    \vcenter{\xygraph{{R \ox R}="0"
        :[r(2.2)] {R \ox R \ox R} ^-{\delta \ox 1}
        :[d(1.2)] {R \ox R.}="1" ^-{1 \ox \mu}
    "0":[d(1.2)] {R \ox R \ox R} _-{1 \ox \delta}
        :"1" ^-{\mu \ox 1}}}
\]
In strings the Frobenius condition is displayed as
\[
    \vcenter{\xy
        (0,-7)*{~}; (0,-7)*{~};
        (-3,-5)="bl";
        (-3,5)="tl";
        (-3,3)="s";
        (3,-5)="br";
        (3,5)="tr";
        (3,-3)="e";
        "bl";"tl" **\dir{-};
        "br";"tr" **\dir{-};
        "s";"e" **\dir{-};
    \endxy}
    ~=~
    \vcenter{\xy
        (-3,-5)="bl";
        (-3,5)="tl";
        (-3,-3)="e";
        (3,-5)="br";
        (3,5)="tr";
        (3,3)="s";
        "bl";"tl" **\dir{-};
        "br";"tr" **\dir{-};
        "s";"e" **\dir{-};
    \endxy}~.
\]

We will now review some basic facts about Frobenius monoids.

\begin{lemma}
    $(1 \ox \mu)(\delta \ox 1) = \delta \mu = (\mu \ox 1)(1 \ox \delta):
    R \ox R \ra R \ox R.$
\end{lemma}

\begin{proof}
The left-hand identity is proved by the following string calculation.
\[
    \vcenter{\xy
        (0,-7)*{~};
        (0,-7)*{~};
        (-3,5)="tl";
        (-3,-5)="bl";
        (3,5)="tr";
        (3,-5)="br";
        (-3,3)="s";
        (3,-3)="e";
        "bl";"tl" **\dir{-};
        "br";"tr" **\dir{-};
        "s";"e" **\dir{-};
    \endxy}
    ~=~
    \vcenter{\xy
        (-3,5)="tl";
        (-3,-5)="bl";
        (3,5)="tr";
        (3,-5)="br";
        (-3,3)="s";
        (3,-3)="e";
        (3,3)="u1";
        (6,0)*[o]=<5pt>[Fo]{~}="u2";
        "bl";"tl" **\dir{-};
        "br";"tr" **\dir{-};
        "s";"e" **\dir{-};
        "u1";"u2" **\dir{-};
    \endxy}
    ~=~
    \vcenter{\xy
        (-3,5)="tl";
        (-3,-5)="bl";
        (3,5)="tr";
        (3,-3)*[o]=<5pt>[Fo]{~}="br";
        (-3,3)="s";
        (3,-1)="e";
        (0,-5)="m";
        "bl";"tl" **\dir{-};
        "br";"tr" **\dir{-};
        "s";"e" **\dir{-};
        ?;"m" **\dir{-};
    \endxy}
    ~=~
    \vcenter{\xy
        (-3,5)="tl";
        (-3,-5)="bl";
        (3,5)="tr";
        (3,-3)*[o]=<5pt>[Fo]{~}="br";
        (-3,3)="s";
        (3,-1)="e";
        (-3,-2)="m0";
        (0,-5)="m";
        "bl";"tl" **\dir{-};
        "br";"tr" **\dir{-};
        "s";"e" **\dir{-};
        "m0";"m" **\dir{-};
    \endxy}
    ~=~
    \vcenter{\xy
        (-3,5)="tl";
        (-3,-5)="bl";
        (3,5)="tr";
        (3,-3)*[o]=<5pt>[Fo]{~}="br";
        (-3,0)="s";
        (3,3)="e";
        (-3,-2)="m0";
        (0,-5)="m";
        "bl";"tl" **\dir{-};
        "br";"tr" **\dir{-};
        "s";"e" **\dir{-};
        "m0";"m" **\dir{-};
    \endxy}
    ~=~
    \vcenter{\xy
        (3,5)="tl";
        (-3,5)="tr";
        (0,2)="mt";
        (0,-2)="mb";
        (3,-5)="bl";
        (-3,-5)="br";
        "tl";"mt" **\dir{-};
        "tr";"mt" **\dir{-};
        "mt";"mb" **\dir{-};
        "mb";"bl" **\dir{-};
        "mb";"br" **\dir{-};
    \endxy}
\]
The right-hand identity follows from a similar calculation.
\end{proof}

Define morphisms $\rho$ and $\sigma$ by
\begin{align*}
    & \rho = \big(\xygraph{{I}
        :[r(1)] {R} ^-\eta
        :[r(1.4)] {R \ox R} ^-\delta}\big) = 
    \vcenter{\xy
        (0,5)*{~};
        (0,-5)*{~};
        (0,4)*[o]=<5pt>[Fo]{~}="t";
        (0,0)="m";
        (-3,-4)="l";
        (3,-4)="r";
        "t";"m" **\dir{-};
        "m";"l" **\dir{-};
        "m";"r" **\dir{-};
    \endxy} \\
    & \sigma = \big(\xygraph{{R \ox R}
        :[r(1.4)] {R} ^-\mu
        :[r(1)] {I} ^-\epsilon}\big) =
    \vcenter{\xy
        (0,5)*{~};
        (0,-5)*{~};
        (-3,4)="l";
        (3,4)="r";
        (0,0)="m";
        (0,-4)*[o]=<5pt>[Fo]{~}="b";
        "l";"m" **\dir{-};
        "r";"m" **\dir{-};
        "m";"b" **\dir{-};
    \endxy}~.
\end{align*}

\begin{proposition}
The morphisms $\rho$ and $\sigma$ form the unit and counit of an adjunction
$R \dashv R$.
\end{proposition}

\begin{proof}
One of the triangle identities is given as
\[
    \vcenter{\xy
        (-3,8)="tl";
        (3,5)*[o]=<5pt>[Fo]{~}="tr";
        (-3,-5)*[o]=<5pt>[Fo]{~}="bl";
        (3,-8)="br";
        (3,2)="s";
        (-3,-2)="e";
        "tl";"bl" **\dir{-};
        "tr";"br" **\dir{-};
        "s";"e" **\dir{-};
    \endxy}
    ~\overset{(\text{f})}{=}~
    \vcenter{\xy
        (-3,8)="tl";
        (3,3)*[o]=<5pt>[Fo]{~}="tr";
        (-3,-3)*[o]=<5pt>[Fo]{~}="bl";
        (3,-8)="br";
        (-3,2)="s";
        (3,-2)="e";
        "tl";"bl" **\dir{-};
        "tr";"br" **\dir{-};
        "s";"e" **\dir{-};
    \endxy}
    ~\overset{(\text{m,c})}{=}~
    \vcenter{\xy
        (0,8)="t";
        (0,-8)="b";
        "t";"b" **\dir{-};
    \endxy}
\]
and the other should now be clear.
\end{proof}

A \emph{morphism of Frobenius monoids} $f:R \ra S$ is a morphism in $\V$
which is both a monoid and comonoid morphism.

\begin{proposition}\label{prop-frobiso}
Any morphism of Frobenius monoids $f:R \ra S$ is an isomorphism.
\end{proposition}

\begin{proof}
Given $f:R \ra S$ define $f^{-1}:S \ra R$ by
\[
    f^{-1} = \xygraph{{S}
        :[r(1.7)] {R \ox R \ox S} ^-{\rho \ox 1}
        :[r(2.5)] {R \ox S \ox S} ^-{1 \ox f \ox 1}
        :[r(1.7)] {R} ^-{1 \ox \sigma}}
    =~
    \vcenter{\xy % f inverse
        (-4,6)*[o]=<5pt>[Fo]{~}="1";
        (-4,3)="2";
        (-4,-9)="3";
        (0,0)*[o]=<11pt>[Fo]{\scriptstyle f}="4";
        (4,9)="5";
        (4,-3)="6";
        (4,-6)*[o]=<5pt>[Fo]{~}="7";
        "1";"3" **\dir{-};
        "2";"4" **\dir{-};
        "4";"6" **\dir{-};
        "5";"7" **\dir{-};
    \endxy}~.
\]
It is then an easy calculation to show that $f^{-1}$ is the inverse of $f$,
namely,
\[
    \vcenter{\xy
        (0, 9)="1";
        (0, 4)*[o]=<11pt>[Fo]{\scriptstyle f}="2";
        (0,-4)*[o]=<16pt>[Fo]{\scriptstyle f^{-1}}="3";
        (0,-10)="4";
        "1";"2" **\dir{-};
        "2";"3" **\dir{-};
        "3";"4" **\dir{-};
    \endxy}
    ~=~
    \vcenter{\xy
        (-4,6)*[o]=<5pt>[Fo]{~}="1";
        (-4,3)="2";
        (-4,-9)="3";
        (0,0)*[o]=<11pt>[Fo]{\scriptstyle f}="4";
        (6,0)*[o]=<11pt>[Fo]{\scriptstyle f}="6";
        (6,9)="8";
        (3,-3)="5";
        (3,-6)*[o]=<5pt>[Fo]{~}="7";
        "1";"3" **\dir{-};
        "2";"4" **\dir{-};
        "4";"5" **\dir{-};
        "5";"7" **\dir{-};
        "5";"6" **\dir{-};
        "8";"6" **\dir{-};
    \endxy}
    ~=~
    \vcenter{\xy
        (-3,6)*[o]=<5pt>[Fo]{~}="1";
        (-3,3)="2";
        (-3,-12)="3";
        (3,9)="4";
        (3,-1)="5";
        (3,-5)*[o]=<11pt>[Fo]{\scriptstyle f}="6";
        (3,-10)*[o]=<5pt>[Fo]{~}="7";
        "1";"3" **\dir{-};
        "2";"5" **\dir{-};
        "4";"6" **\dir{-};
        "6";"7" **\dir{-};
    \endxy}
    ~=~
    \vcenter{\xy
        (-3, 6)*[o]=<5pt>[Fo]{~}="1";
        (-3, 2)="2";
        (-3,-8)="3";
        ( 3, 8)="4";
        ( 3,-2)="5";
        ( 3,-6)*[o]=<5pt>[Fo]{~}="7";
        "1";"3" **\dir{-};
        "2";"5" **\dir{-};
        "4";"7" **\dir{-};
    \endxy}
    ~\overset{\text{(f)}}{=}~
    \vcenter{\xy
        (-3, 4)*[o]=<5pt>[Fo]{~}="1";
        (-3,-2)="2";
        (-3,-8)="3";
        ( 3, 8)="4";
        ( 3, 2)="5";
        ( 3,-4)*[o]=<5pt>[Fo]{~}="7";
        "1";"3" **\dir{-};
        "2";"5" **\dir{-};
        "4";"7" **\dir{-};
    \endxy}
    ~\overset{\text{(m,c)}}{=}~
    \vcenter{\xy
        (0, 7); (0,-7); **\dir{-};
    \endxy}
\]
and the same viewed upside down.
\end{proof}

A similar calculation shows that
\[
    \xygraph{{S}
        :[r(1.7)] {S \ox R \ox R} ^-{1 \ox \rho}
        :[r(2.5)] {S \ox S \ox R} ^-{1 \ox f \ox 1}
        :[r(1.7)] {R} ^-{\sigma \ox 1}}
    =~
    \vcenter{\xy % f inverse
        ( 4,6)*[o]=<5pt>[Fo]{~}="1";
        ( 4,3)="2";
        ( 4,-9)="3";
        (0,0)*[o]=<11pt>[Fo]{\scriptstyle f}="4";
        (-4,9)="5";
        (-4,-3)="6";
        (-4,-6)*[o]=<5pt>[Fo]{~}="7";
        "1";"3" **\dir{-};
        "2";"4" **\dir{-};
        "4";"6" **\dir{-};
        "5";"7" **\dir{-};
    \endxy}~.
\]
is also an inverse of $f$. Therefore,

\begin{corollary}
For any morphism of Frobenius monoids $f:R \ra S$ we have
\[
    \vcenter{\xy % f inverse
        (-4,6)*[o]=<5pt>[Fo]{~}="1";
        (-4,3)="2";
        (-4,-9)="3";
        (0,0)*[o]=<11pt>[Fo]{\scriptstyle f}="4";
        (4,9)="5";
        (4,-3)="6";
        (4,-6)*[o]=<5pt>[Fo]{~}="7";
        "1";"3" **\dir{-};
        "2";"4" **\dir{-};
        "4";"6" **\dir{-};
        "5";"7" **\dir{-};
    \endxy}
    ~=~
    \vcenter{\xy % f inverse
        ( 4,6)*[o]=<5pt>[Fo]{~}="1";
        ( 4,3)="2";
        ( 4,-9)="3";
        (0,0)*[o]=<11pt>[Fo]{\scriptstyle f}="4";
        (-4,9)="5";
        (-4,-3)="6";
        (-4,-6)*[o]=<5pt>[Fo]{~}="7";
        "1";"3" **\dir{-};
        "2";"4" **\dir{-};
        "4";"6" **\dir{-};
        "5";"7" **\dir{-};
    \endxy}\ \ .
\]
\end{corollary}

\begin{definition}\label{def-sepfrob}
A Frobenius monoid $R$ is said to be \emph{separable} if and only if $\mu
\delta = 1$, i.e.,
\[
    \vcenter{\xy
        (0,6)="1";
        (0,3)="2";
        (-3,0)="a";
        (3,0)="b";
        (0,-3)="3";
        (0,-6)="4";
        "1";"2" **\dir{-};
        "2";"a" **\dir{-};
        "2";"b" **\dir{-};
        "3";"a" **\dir{-};
        "3";"b" **\dir{-};
        "3";"4" **\dir{-};
    \endxy}
    ~=~
    \vcenter{\xy
        (0,6); (0,-6); **\dir{-};
    \endxy}~.
\]
\end{definition}

%=========================================================================%
\section{Proofs of the properties of $s$, $t$, and $r$}\label{sec-proofst}
%=========================================================================%

As we have noted in \S\ref{sec-wb}, $s:A \ra A$ is invariant under rotation
by $\pi$, $t:A \ra A$ is invariant under horizontal reflection, and $r$ is
$t$ rotated by $\pi$. This reduces the number of proofs we present as the
others are derivable.

\begin{enumerate}
\item % This is 1.a
    $\vcenter{\xy
        (0,5)="1";
        (0,2)*[o]=<9pt>[Fo]{\scriptstyle s}="2";
        (0,-1)="3";
        (-3,-5)="4";
        ( 3,-5)="5";
        "1";"2" **\dir{-};
        "2";"3" **\dir{-};
        "3";"4" **\dir{-};
        "3";"5" **\dir{-};
    \endxy}
    ~\overset{\text{($s$)}}{=}~
    \vcenter{\xy
        (3,8)="1";
        (0,5)*[o]=<5pt>[Fo]{~}="2";
        (0,3)="3";
        (0,0)="4";
        (0,-2)*[o]=<5pt>[Fo]{~}="5";
        (-3,-3)="6";
        (-6,-6)="7";
        ( 0,-6)="8";
        "1";"4" **\crv{(3,4)};
        "2";"5" **\dir{-};
        "6";"3" **\crv{(-3,0)};
        "6";"7" **\dir{-};
        "6";"8" **\dir{-};
    \endxy}
    ~\overset{\text{(c)}}{=}~
    \vcenter{\xy
        (3,9)="1";
        (0,6)*[o]=<5pt>[Fo]{~}="2";
        (0,3)="3";
        (-3,-3)="4";
        ( 0,-3)="5";
        ( 3, 0)*[o]=<5pt>[Fo]{~}="6";
        "2";"5" **\dir{-};
        "3";"4" **\dir{-};
        "3";"6" **\dir{-}?="x";
        "1";"x" **\crv{(3,5)};
    \endxy}
    ~\overset{\text{(w)}}{=}~
    \vcenter{\xy
        (6,8.5)="1";
        (-3,5.5)*[o]=<5pt>[Fo]{~}="2";
        ( 3,5.5)*[o]=<5pt>[Fo]{~}="3";
        ( 3,3)="4";
        (-3,3)="5";
        ( 0,1)="6";
        ( 3,0)="7";
        ( 3,-2.5)*[o]=<5pt>[Fo]{~}="8";
        (-3,-5.5)="9";
        ( 0,-5.5)="10";
        "2";"9" **\dir{-};
        "3";"8" **\dir{-};
        "6";"10" **\dir{-};
        "4";"6" **\dir{-};
        "5";"6" **\dir{-};
        "1";"7" **\crv{(6,4)};
    \endxy}
    ~\overset{\text{($s$)}}{=}~
    \vcenter{\xy
        (0,7)*{~};
        (0,-6)*{~};
        (-3,3)*[o]=<5pt>[Fo]{~}="1";
        (-3,0)="2";
        (-3,-5)="3";
        (3,6)="4";
        (1.5,2)*[o]=<9pt>[Fo]{\scriptstyle s}="5";
        (0,-2)="6";
        (0,-5)="7";
        "1";"3" **\dir{-};
        "2";"6" **\dir{-};
        "6";"7" **\dir{-};
        "5";"6" **\dir{-};
        "4";"5" **\dir{-};
    \endxy}$
\\ 
    $\vcenter{\xy
        (0,5)="1";
        (0,2)*[o]=<9pt>[Fo]{\scriptstyle s}="2";
        (0,-1)="3";
        (-3,-5)="4";
        ( 3,-5)="5";
        "1";"2" **\dir{-};
        "2";"3" **\dir{-};
        "3";"4" **\dir{-};
        "3";"5" **\dir{-};
    \endxy}
    ~\overset{\text{($s$)}}{=}~
    \vcenter{\xy
        (3,8)="1";
        (0,5)*[o]=<5pt>[Fo]{~}="2";
        (0,3)="3";
        (0,0)="4";
        (0,-2)*[o]=<5pt>[Fo]{~}="5";
        (-3,-3)="6";
        (-6,-6)="7";
        ( 0,-6)="8";
        "1";"4" **\crv{(3,4)};
        "2";"5" **\dir{-};
        "6";"3" **\crv{(-3,0)};
        "6";"7" **\dir{-};
        "6";"8" **\dir{-};
    \endxy}
    ~\overset{\text{(c)}}{=}~
    \vcenter{\xy
        (3,9)="1";
        (0,6)*[o]=<5pt>[Fo]{~}="2";
        (0,3)="3";
        (-3,-3)="4";
        ( 0,-3)="5";
        ( 3, 0)*[o]=<5pt>[Fo]{~}="6";
        "2";"5" **\dir{-};
        "3";"4" **\dir{-};
        "3";"6" **\dir{-}?="x";
        "1";"x" **\crv{(3,5)};
    \endxy}
    \overset{\text{(w)}}{=}~
    \vcenter{\xy
        (8,9)="1";
        (-4,5.5)*[o]=<5pt>[Fo]{~}="2";
        ( 4,5.5)*[o]=<5pt>[Fo]{~}="3";
        (-4,3)="4";
        ( 4,3)="5";
        (-2.5,0)="x";
        ( 2.5,0)="y";
        ( 0,-2)="6";
        ( 4,0)="7";
        ( 4,-2.5)*[o]=<5pt>[Fo]{~}="8";
        (-4,-5.5)="9";
        ( 0,-5.5)="10";
        "2";"9" **\dir{-};
        "3";"8" **\dir{-};
        "6";"10" **\dir{-};
        "6";"x" **\dir{-};
        "6";"y" **\dir{-};
        "5";"x" **\dir{-};
        "1";"7" **\crv{(8,4)};
        \ar@{-} "4";"y" |<(0.57)\hole
    \endxy}
    ~\overset{\text{($s$)}}{=}~
    \vcenter{\xy
        (0,8)*{~};
        (0,-5)*{~};
        (-5,4)*[o]=<5pt>[Fo]{~}="1";
        (-5,1)="2";
        (-5,-4)="3";
        (3,7)="4";
        (0,3)*[o]=<9pt>[Fo]{\scriptstyle s}="5";
        (0,-2)="6";
        (0,-4)="7";
        "1";"3" **\dir{-};
        "6";"7" **\dir{-};
        "5";"6" **\crv{(-3,0)}?(0.5)*{\hole}="x";
        "2";"x" **\dir{-};
        "6";"x" **\crv{(3,0)};
        "4";"5" **\dir{-};
    \endxy}$ \\
%
% This is 1.b
%
    $\vcenter{\xy
        (0,5)="t";
        (0,2)*[o]=<9pt>[Fo]{\scriptstyle t}="f";
        (0,-1)="m";
        (-3,-5)="l";
        ( 3,-5)="r";
        "t";"f" **\dir{-};
        "f";"m" **\dir{-};
        "m";"l" **\dir{-};
        "m";"r" **\dir{-};
    \endxy}
    ~\overset{\text{($t$)}}{=}
    \vcenter{\xy
        (-3,8)="1";
        (0,5)*[o]=<5pt>[Fo]{~}="2";
        (0,3)="3";
        (0,0)="4";
        (0,-2)*[o]=<5pt>[Fo]{~}="5";
        (-3,-3)="6";
        (-6,-6)="7";
        ( 0,-6)="8";
        "1";"4" **\crv{(-3,1)};
        "2";"5" **\dir{-};
        "6";"7" **\dir{-};
        "6";"8" **\dir{-};
        \ar@{-}@/^3pt/ "6";"3" |<(0.6)\hole
    \endxy}
    ~\overset{\text{(c)}}{=}~
    \vcenter{\xy
        (-4,9)="1";
        (0,6)*[o]=<5pt>[Fo]{~}="2";
        (0,3)="3";
        (-4,-5)="4";
        ( 0,-5)="5";
        ( 4,-2)*[o]=<5pt>[Fo]{~}="6";
        "3";"6" **\dir{-}?(0.75)="x";
        "1";"x" **\crv{(-4,-1)};
        \ar@{-} "2";"5" |<(0.55)\hole
        \ar@{-} "3";"4" |<(0.34)\hole
    \endxy}
    ~\overset{\text{(w)}}{=}~
    \vcenter{\xy
        ( 2,8.5)="1";
        (-4,5.5)*[o]=<5pt>[Fo]{~}="2";
        ( 6,5.5)*[o]=<5pt>[Fo]{~}="3";
        (-4,3)="4";
        ( 6,3)="5";
        ( 0,1)="6";
        ( 6,0)="7";
        ( 6,-2.5)*[o]=<5pt>[Fo]{~}="8";
        (-4,-5.5)="9";
        ( 0,-5.5)="10";
        "2";"9" **\dir{-};
        "3";"8" **\dir{-};
        "6";"10" **\dir{-};
        "4";"6" **\dir{-};
        "1";"7" **\crv{(2,3)};
        \ar@{-} "5";"6" |<(0.46)\hole
    \endxy}
    ~\overset{\text{($t$)}}{=}~
    \vcenter{\xy
        (0,7)*{~};
        (0,-6)*{~};
        (-3,3)*[o]=<5pt>[Fo]{~}="1";
        (-3,0)="2";
        (-3,-5)="3";
        (3,6)="4";
        (1.5,2)*[o]=<9pt>[Fo]{\scriptstyle t}="5";
        (0,-2)="6";
        (0,-5)="7";
        "1";"3" **\dir{-};
        "2";"6" **\dir{-};
        "6";"7" **\dir{-};
        "5";"6" **\dir{-};
        "4";"5" **\dir{-};
    \endxy}$
\\ 
    $\vcenter{\xy
        (0,5)="t";
        (0,2)*[o]=<9pt>[Fo]{\scriptstyle t}="f";
        (0,-1)="m";
        (-3,-5)="l";
        ( 3,-5)="r";
        "t";"f" **\dir{-};
        "f";"m" **\dir{-};
        "m";"l" **\dir{-};
        "m";"r" **\dir{-};
    \endxy}
    ~\overset{\text{($t$)}}{=}
    \vcenter{\xy
        (-3,8)="1";
        (0,5)*[o]=<5pt>[Fo]{~}="2";
        (0,3)="3";
        (0,0)="4";
        (0,-2)*[o]=<5pt>[Fo]{~}="5";
        (-3,-3)="6";
        (-6,-6)="7";
        ( 0,-6)="8";
        "1";"4" **\crv{(-3,1)};
        "2";"5" **\dir{-};
        "6";"7" **\dir{-};
        "6";"8" **\dir{-};
        \ar@{-}@/^3pt/ "6";"3" |<(0.6)\hole
    \endxy}
    ~\overset{\text{(c)}}{=}~
    \vcenter{\xy
        (-4,9)="1";
        (0,6)*[o]=<5pt>[Fo]{~}="2";
        (0,3)="3";
        (-4,-5)="4";
        ( 0,-5)="5";
        ( 4,-2)*[o]=<5pt>[Fo]{~}="6";
        "3";"6" **\dir{-}?(0.75)="x";
        "1";"x" **\crv{(-4,-1)};
        \ar@{-} "2";"5" |<(0.55)\hole
        \ar@{-} "3";"4" |<(0.34)\hole
    \endxy}
    \overset{\text{(w)}}{=}~
    \vcenter{\xy
        (1,9)="1";
        (-4,5.5)*[o]=<5pt>[Fo]{~}="2";
        ( 6,5.5)*[o]=<5pt>[Fo]{~}="3";
        (-4,3)="4";
        ( 6,3)="5";
        (-2,0)="x";
        ( 2,-1)="y";
        ( 0,-2.5)="6";
        ( 6,0)="7";
        ( 6,-2.5)*[o]=<5pt>[Fo]{~}="8";
        (-4,-5.5)="9";
        ( 0,-5.5)="10";
        "2";"9" **\dir{-};
        "3";"8" **\dir{-};
        "6";"10" **\dir{-};
        "6";"x" **\dir{-};
        "6";"y" **\dir{-};
        "1";"7" **\crv{(1,3)};
        \ar@{-} "4";"y" |<(0.57)\hole
        \ar@{-} "5";"x" |<(0.34)\hole
    \endxy}
    ~\overset{\text{($t$)}}{=}~
    \vcenter{\xy
        (0,8)*{~};
        (0,-5)*{~};
        (-5,4)*[o]=<5pt>[Fo]{~}="1";
        (-5,1)="2";
        (-5,-4)="3";
        (3,7)="4";
        (0,3)*[o]=<9pt>[Fo]{\scriptstyle t}="5";
        (0,-2)="6";
        (0,-4)="7";
        "1";"3" **\dir{-};
        "6";"7" **\dir{-};
        "5";"6" **\crv{(-3,0)}?(0.5)*{\hole}="x";
        "2";"x" **\dir{-};
        "6";"x" **\crv{(3,0)};
        "4";"5" **\dir{-};
    \endxy}$

\item % (2) This is 2.a
    $\vcenter{\xy
        (0,4)*[o]=<5pt>[Fo]{~}="t";
        (0,1)="m";
        (-3,-2)*[o]=<9pt>[Fo]{\scriptstyle s}="g";
        (-3,-5)="l";
        ( 3,-5)="r";
        "t";"m" **\dir{-};
        "r";"m" **\crv{(3,-2)};
        "m";"g" **\dir{-};
        "g";"l" **\dir{-};
    \endxy}
    ~\overset{\text{($s$)}}{=}~
    \vcenter{\xy
        (0,6)*[o]=<5pt>[Fo]{~}="1";
        (0,3)="2";
        (-3, 4)*[o]=<5pt>[Fo]{~}="3";
        (-3, 2)="4";
        (-3,-1)="5";
        (-3,-3)*[o]=<5pt>[Fo]{~}="6";
        (-6,-6)="7";
        ( 3,-6)="8";
        "1";"2" **\dir{-};
        "2";"5" **\dir{-};
        "3";"6" **\dir{-};
        "8";"2" **\crv{(3,-2)};
        "7";"4" **\crv{(-6,-2)};
    \endxy}
    ~\overset{\text{(nat)}}{=}~
    \vcenter{\xy
        (-3,5)*[o]=<5pt>[Fo]{~}="1";
        ( 3,5)*[o]=<5pt>[Fo]{~}="2";
        (-3,2)="3";
        ( 3,2)="4";
        (0,0)="5";
        (0,-3)*[o]=<5pt>[Fo]{~}="6";
        (-3,-5)="7";
        ( 3,-5)="8";
        "1";"7" **\dir{-};
        "2";"8" **\dir{-};
        "5";"3" **\dir{-};
        "5";"4" **\dir{-};
        "5";"6" **\dir{-};
    \endxy}
    ~\overset{\text{(w)}}{=}~
    \vcenter{\xy
        (0,3)*[o]=<5pt>[Fo]{~}="1";
        (0,0)="2";
        (0,-3)*[o]=<5pt>[Fo]{~}="n";
        (-4,-5)="3";
        ( 4,-5)="4";
        "1";"2" **\dir{-};
        "2";"3" **\dir{-};
        "2";"4" **\dir{-};
        "2";"n" **\dir{-};
    \endxy}
    ~\overset{\text{(c)}}{=}~
    \vcenter{\xy
        (0,2.5)*[o]=<5pt>[Fo]{~}="1";
        (0,0)="2";
        (-2.5,-3)="3";
        ( 2.5,-3)="4";
        "1";"2" **\dir{-};
        "2";"3" **\dir{-};
        "2";"4" **\dir{-};
    \endxy}$ \\
%
% This is 2.b
%
    $\vcenter{\xy
        (0,-4.5)="t";
        (0,0)*[o]=<9pt>[Fo]{\scriptstyle s}="g";
        (0,4)*[o]=<5pt>[Fo]{}="b";
        "t";"g" **\dir{-};
        "g";"b" **\dir{-};
    \endxy}
    ~\overset{\text{($s$)}}{=}~
    \vcenter{\xy
        (0, 5)*[o]=<5pt>[Fo]{~}="1";
        (0, 3)="2";
        (0, 0)="3";
        (0,-2)*[o]=<5pt>[Fo]{~}="4";
        (3,3)*[o]=<5pt>[Fo]{~}="5";
        (-3,-5)="6";
        "1";"4" **\dir{-};
        "5";"3" **\dir{-};
        "6";"2" **\crv{(-3,-1)};
    \endxy}
    ~\overset{\text{(m)}}{=}~
    \vcenter{\xy 
        (0,2)*[o]=<5pt>[Fo]{~}="t";
        (0,0)="w1";
        (2,-2)*[o]=<5pt>[Fo]{~}="w2";
        (0,-5)="b";
        "w1";"w2" **\dir{-};
        "t";"b" **\dir{-};
    \endxy}
    ~\overset{\text{(c)}}{=}~
    \vcenter{\xy 
        (0,2)*[o]=<5pt>[Fo]{}="t";
        (0,-3)="b";
        "t";"b" **\dir{-};
    \endxy}$ \\
%
% This is 2.c
%
    $\vcenter{\xy
        (0,4)*[o]=<5pt>[Fo]{~}="t";
        (0,1)="m";
        (-3,-2)*[o]=<9pt>[Fo]{\scriptstyle t}="g";
        (-3,-5)="l";
        ( 3,-5)="r";
        "t";"m" **\dir{-};
        "m";"g" **\dir{-};
        "g";"l" **\dir{-};
        "r";"m" **\crv{(3,-2)};
    \endxy}
    ~\overset{\text{(t)}}{=}~
    \vcenter{\xy
        ( 0,10)*{~}; ( 0,-9)*{~};
        ( 0, 9)*[o]=<5pt>[Fo]{~}="1";
        ( 0, 6)="2";
        (-1, 3)*[o]=<5pt>[Fo]{~}="3";
        (-1, 1)="4";
        (-1,-2)="5";
        (-1,-5)*[o]=<5pt>[Fo]{~}="6";
        (-5,-8)="7";
        ( 5,-8)="8";
        "1";"2" **\dir{-};
        "3";"6" **\dir{-};
        "2";"5" **\crv{(-6,4)&(-6,-2)};
        "8";"2" **\crv{(4,3)};
        \ar@{-}@/^3pt/ "7";"4" |<(0.65)\hole
    \endxy}
    ~\overset{\text{(nat)}}{=}~
    \vcenter{\xy
        (-4,5.5)*[o]=<5pt>[Fo]{~}="2";
        ( 4,5.5)*[o]=<5pt>[Fo]{~}="3";
        (-4,3)="4";
        ( 4,3)="5";
        (-2.5,0)="x";
        ( 2.5,0)="y";
        ( 0,-2)="6";
        ( 4,0)="7";
        (-4,-8)="9";
        ( 4,-8)="8";
        ( 0,-5)*[o]=<5pt>[Fo]{~}="10";
        "2";"9" **\dir{-};
        "3";"8" **\dir{-};
        "6";"10" **\dir{-};
        "6";"x" **\dir{-};
        "6";"y" **\dir{-};
        "5";"x" **\dir{-};
        \ar@{-} "4";"y" |<(0.57)\hole
    \endxy}
    ~\overset{\text{(w)}}{=}~
    \vcenter{\xy
        (0,3)*[o]=<5pt>[Fo]{~}="1";
        (0,0)="2";
        (0,-3)*[o]=<5pt>[Fo]{~}="n";
        (-4,-5)="3";
        ( 4,-5)="4";
        "1";"2" **\dir{-};
        "2";"3" **\dir{-};
        "2";"4" **\dir{-};
        "2";"n" **\dir{-};
    \endxy}
    ~\overset{\text{(c)}}{=}~
    \vcenter{\xy
        (0,2.5)*[o]=<5pt>[Fo]{~}="1";
        (0,0)="2";
        (-2.5,-3)="3";
        ( 2.5,-3)="4";
        "1";"2" **\dir{-};
        "2";"3" **\dir{-};
        "2";"4" **\dir{-};
    \endxy}$ \\
%
% This is 2.d
%
    $\vcenter{\xy
        (0,-4.5)="t";
        (0,0)*[o]=<9pt>[Fo]{\scriptstyle t}="g";
        (0,4)*[o]=<5pt>[Fo]{}="b";
        "t";"g" **\dir{-};
        "g";"b" **\dir{-};
    \endxy}
    ~\overset{\text{(t)}}{=}~
    \vcenter{\xy
        (0,7)*{~}; ( 0,-6)*{~};
        (0, 5)*[o]=<5pt>[Fo]{~}="1";
        (0, 3)="2";
        (0, 0)="3";
        (0,-2)*[o]=<5pt>[Fo]{~}="4";
        (-3,6)*[o]=<5pt>[Fo]{~}="5";
        (-3,-5)="6";
        "1";"4" **\dir{-};
        \ar@{-}@/^3pt/ "6";"2" |<(0.69)\hole
        \ar@{-}@/_3pt/ "5";"3"
    \endxy}
    ~\overset{\text{(m)}}{=}~
    \vcenter{\xy 
        (0,2)*[o]=<5pt>[Fo]{~}="t";
        (0,0)="w1";
        (2,-2)*[o]=<5pt>[Fo]{~}="w2";
        (0,-5)="b";
        "w1";"w2" **\dir{-};
        "t";"b" **\dir{-};
    \endxy}
    ~\overset{\text{(c)}}{=}~
    \vcenter{\xy 
        (0,2)*[o]=<5pt>[Fo]{}="t";
        (0,-3)="b";
        "t";"b" **\dir{-};
    \endxy}$

\item % This is 3.a
    $\vcenter{\xy
        (0,7)*{~};
        (0,-7)*{~};
        (0,6)="t";
        (0,3)*[o]=<9pt>[Fo]{\scriptstyle s}="f";
        (0,-0)="m";
        (-3,-3)*[o]=<9pt>[Fo]{\scriptstyle s}="g";
        (-3,-6)="l";
        ( 3,-6)="r";
        "t";"f" **\dir{-};
        "f";"m" **\dir{-};
        "r";"m" **\crv{(3,-3)};
        "m";"g" **\dir{-};
        "g";"l" **\dir{-};
    \endxy}
    ~\overset{\text{(1)}}{=}~
    \vcenter{\xy
        (-3,4)*[o]=<5pt>[Fo]{~}="1";
        (-3,1)="2";
        (-3,-2)*[o]=<9pt>[Fo]{\scriptstyle s}="x";
        (-3,-5)="3";
        (3,6)="4";
        (1.5,3)*[o]=<9pt>[Fo]{\scriptstyle s}="5";
        (0, 0)="6";
        (0,-5)="7";
        "1";"x" **\dir{-};
        "x";"3" **\dir{-};
        "2";"6" **\dir{-};
        "6";"7" **\dir{-};
        "5";"6" **\dir{-};
        "4";"5" **\dir{-};
    \endxy}
    ~\overset{\text{(2)}}{=}~
    \vcenter{\xy
        (-3,3)*[o]=<5pt>[Fo]{~}="1";
        (-3,0)="2";
        (-3,-5)="3";
        (3,6)="4";
        (1.5,2)*[o]=<9pt>[Fo]{\scriptstyle s}="5";
        (0,-2)="6";
        (0,-5)="7";
        "1";"3" **\dir{-};
        "2";"6" **\dir{-};
        "6";"7" **\dir{-};
        "5";"6" **\dir{-};
        "4";"5" **\dir{-};
    \endxy}
    ~\overset{\text{(1)}}{=}
    \vcenter{\xy
        (0,5)="t";
        (0,2)*[o]=<9pt>[Fo]{\scriptstyle s}="f";
        (0,-1)="m";
        (-3,-5)="l";
        ( 3,-5)="r";
        "t";"f" **\dir{-};
        "f";"m" **\dir{-};
        "m";"l" **\dir{-};
        "m";"r" **\dir{-};
    \endxy}$ \\
%
% This is 3.b
%
    $\vcenter{\xy
        (0,6)="t";
        (0,3)*[o]=<9pt>[Fo]{\scriptstyle t}="f";
        (0,-0)="m";
        (-3,-3)*[o]=<9pt>[Fo]{\scriptstyle t}="g";
        (-3,-6)="l";
        ( 3,-6)="r";
        "t";"f" **\dir{-};
        "f";"m" **\dir{-};
        "r";"m" **\crv{(3,-3)};
        "m";"g" **\dir{-};
        "g";"l" **\dir{-};
    \endxy}
    ~\overset{\text{(1)}}{=}~
    \vcenter{\xy
        (-3,4)*[o]=<5pt>[Fo]{~}="1";
        (-3,1)="2";
        (-3,-2)*[o]=<9pt>[Fo]{\scriptstyle t}="x";
        (-3,-5)="3";
        (3,6)="4";
        (1.5,3)*[o]=<9pt>[Fo]{\scriptstyle t}="5";
        (0, 0)="6";
        (0,-5)="7";
        "1";"x" **\dir{-};
        "x";"3" **\dir{-};
        "2";"6" **\dir{-};
        "6";"7" **\dir{-};
        "5";"6" **\dir{-};
        "4";"5" **\dir{-};
    \endxy}
    ~\overset{\text{(2)}}{=}~
    \vcenter{\xy
        (-3,3)*[o]=<5pt>[Fo]{~}="1";
        (-3,0)="2";
        (-3,-5)="3";
        (3,6)="4";
        (1.5,2)*[o]=<9pt>[Fo]{\scriptstyle t}="5";
        (0,-2)="6";
        (0,-5)="7";
        "1";"3" **\dir{-};
        "2";"6" **\dir{-};
        "6";"7" **\dir{-};
        "5";"6" **\dir{-};
        "4";"5" **\dir{-};
    \endxy}
    ~\overset{\text{(1)}}{=}~
    \vcenter{\xy
        (0,5)="t";
        (0,2)*[o]=<9pt>[Fo]{\scriptstyle t}="f";
        (0,-1)="m";
        (-3,-5)="l";
        ( 3,-5)="r";
        "t";"f" **\dir{-};
        "f";"m" **\dir{-};
        "m";"l" **\dir{-};
        "m";"r" **\dir{-};
    \endxy}$

\item % This is 4 sw a.
    $\vcenter{\xy
        (-3,8)="1";
        ( 3,8)="2";
        ( 3,5)*[o]=<9pt>[Fo]{\scriptstyle s}="3";
        ( 0,2)="4";
        ( 0,-1)="5";
        (-3,-5)="6";
        ( 3,-5)="7";
        "1";"4" **\crv{(-3,5)};
        "2";"3" **\dir{-};
        "3";"4" **\dir{-};
        "4";"5" **\dir{-};
        "5";"6" **\dir{-};
        "5";"7" **\dir{-};
    \endxy}
    ~\overset{\text{(b)}}{=}~
    \vcenter{\xy
        (-3,6)="1";
        ( 3,6)="2";
        ( 3,3)*[o]=<9pt>[Fo]{\scriptstyle s}="3";
        (-3,0)="4";
        ( 3,0)="5";
        (-3,-4)="6";
        ( 3,-4)="7";
        (-3,-7)="8";
        ( 3,-7)="9";
        "1";"8" **\dir{-};
        "2";"3" **\dir{-};
        "3";"9" **\dir{-};
        "4";"7" **\dir{-};
        \ar@{-} "5";"6" |\hole
    \endxy}
    ~\overset{\text{(1)}}{=}~
    \vcenter{\xy
        (-3,6)="1";
        ( 3,6)="2";
        ( 3,3)*[o]=<5pt>[Fo]{~}="3";
        (-3,0)="4";
        ( 3,0)="5";
        (-3,-5)="6";
        ( 3,-5)="7";
        (-3,-8)="8";
        ( 3,-8)="9";
        ( 3,-2.5)="x";
        ( 7,1)*[o]=<9pt>[Fo]{\scriptstyle s}="y";
        ( 7,6)="z";
        "1";"8" **\dir{-};
        "3";"9" **\dir{-};
        "4";"7" **\dir{-};
        "x";"y" **\dir{-};
        "y";"z" **\dir{-};
        \ar@{-} "5";"6" |\hole
    \endxy}
    ~\overset{\text{(m)}}{=}~
    \vcenter{\xy
        (0,7)*{~}; (0,-11)*{~};
        (-3,6)="1";
        ( 3,6)="2";
        ( 3,3)*[o]=<5pt>[Fo]{~}="3";
        (-3,0)="4";
        ( 3,0)="5";
        (-3,-5)="6";
        ( 3,-5)="7";
        (-3,-10)="8";
        ( 3,-10)="9";
        ( 3,-7.5)="x";
        ( 7,-4)*[o]=<9pt>[Fo]{\scriptstyle s}="y";
        ( 7,6)="z";
        "1";"8" **\dir{-};
        "3";"9" **\dir{-};
        "4";"7" **\dir{-};
        "x";"y" **\dir{-};
        "y";"z" **\dir{-};
        \ar@{-} "5";"6" |\hole
    \endxy}
    ~\overset{\text{(b)}}{=}~
    \vcenter{\xy
        (-2,8)="1";
        ( 2,8)="2";
        ( 2,4)*[o]=<9pt>[Fo]{\scriptstyle s}="3";
        (-2,3)="4";
        ( 2,-1)="5";
        (-2,-5)="6";
        ( 2,-5)="7";
        "1";"6" **\dir{-};
        "2";"3" **\dir{-};
        "3";"7" **\dir{-};
        "4";"5" **\dir{-};
    \endxy}$ \\
%
% This is 4 sw b.
%
    $\vcenter{\xy
        ( 3,8)="1";
        (-3,8)="2";
        (-3,5)*[o]=<9pt>[Fo]{\scriptstyle s}="3";
        ( 0,2)="4";
        ( 0,-1)="5";
        (3,-5)="6";
        (-3,-5)="7";
        "1";"4" **\crv{(3,5)};
        "2";"3" **\dir{-};
        "3";"4" **\dir{-};
        "4";"5" **\dir{-};
        "5";"6" **\dir{-};
        "5";"7" **\dir{-};
    \endxy}
    ~\overset{\text{(b)}}{=}~
    \vcenter{\xy
        (-3,6)="1";
        ( 3,6)="2";
        (-3,3)*[o]=<9pt>[Fo]{\scriptstyle s}="3";
        (-3,0)="4";
        ( 3,0)="5";
        (-3,-4)="6";
        ( 3,-4)="7";
        (-3,-7)="8";
        ( 3,-7)="9";
        "1";"3" **\dir{-};
        "3";"8" **\dir{-};
        "2";"9" **\dir{-};
        "4";"7" **\dir{-};
        \ar@{-} "5";"6" |\hole
    \endxy}
    ~\overset{\text{(1)}}{=}~
    \vcenter{\xy
        (-1, 9)="1";
        ( 6, 9)="2";
        (-6, 5)*[o]=<5pt>[Fo]{~}="3";
        (-1, 5)*[o]=<9pt>[Fo]{\scriptstyle s}="4";
        (-6, 2)="5";
        ( 1,-2)="6";
        ( 6,-3)="7";
        (-6,-7)="8";
        ( 6,-7)="9";
        (-6,-10)="10";
        ( 6,-10)="11";
        ( 3,0)="x";
        "1";"4" **\dir{-};
        "2";"11" **\dir{-};
        "3";"10" **\dir{-};
        "6";"9" **\dir{-};
        "6";"x" **\dir{-};
        "4";"6" **\crv{(-1,0)};
        \ar@{-} "8";"7" |<(0.75)\hole
        \ar@{-} "5";"x" |<(0.6)\hole
    \endxy}
    ~\overset{\text{(m)}}{=}~
    \vcenter{\xy
        (0,9)*{~}; (0,-11)*{~};
        (-6,8)="1";
        ( 5,8)="2";
        ( 5,3)*[o]=<5pt>[Fo]{~}="3";
        (-6,-2)="4";
        ( 5,0)="5";
        (-6,-5)="6";
        ( 5,-5)="7";
        (-6,-10)="8";
        ( 5,-10)="9";
        ( 5,-7.5)="x";
        ( 2,4)*[o]=<9pt>[Fo]{\scriptstyle s}="y";
        ( 2,8)="z";
        "1";"8" **\dir{-};
        "3";"9" **\dir{-};
        "y";"x" **\crv{(1,-5)};
        "y";"z" **\dir{-};
        \ar@{-} "5";"6" |<(0.25)\hole |<(0.6)\hole
        \ar@{-} "4";"7" |<(0.75)\hole
    \endxy}
    ~\overset{\text{(b)}}{=}~
    \vcenter{\xy
        (0,9)*{~};
        (0,-6)*{~};
        (-2,8)="1";
        ( 2,8)="2";
        ( 2,5)="3";
        ( 2,-2)="4";
        (-2,-5)="5";
        ( 2,-5)="6";
        (-2,5)*[o]=<9pt>[Fo]{\scriptstyle s}="7";
        "2";"6" **\dir{-};
        "1";"7" **\dir{-};
        "7";"4" **\crv{(-2,2)}?(0.5)*{\hole}="x";
        "3";"x" **\crv{(1,5)};
        "5";"x" **\crv{(-2,-2)};
    \endxy}$ \\
%
% This is 4 tw a.
%
    $\vcenter{\xy
        (-3,8)="1";
        ( 3,8)="2";
        ( 3,5)*[o]=<9pt>[Fo]{\scriptstyle t}="3";
        ( 0,2)="4";
        ( 0,-1)="5";
        (-3,-5)="6";
        ( 3,-5)="7";
        "1";"4" **\crv{(-3,5)};
        "2";"3" **\dir{-};
        "3";"4" **\dir{-};
        "4";"5" **\dir{-};
        "5";"6" **\dir{-};
        "5";"7" **\dir{-};
    \endxy}
    ~\overset{\text{(b)}}{=}~
    \vcenter{\xy
        (-3,6)="1";
        ( 3,6)="2";
        ( 3,3)*[o]=<9pt>[Fo]{\scriptstyle t}="3";
        (-3,0)="4";
        ( 3,0)="5";
        (-3,-4)="6";
        ( 3,-4)="7";
        (-3,-7)="8";
        ( 3,-7)="9";
        "1";"8" **\dir{-};
        "2";"3" **\dir{-};
        "3";"9" **\dir{-};
        "4";"7" **\dir{-};
        \ar@{-} "5";"6" |\hole
    \endxy}
    ~\overset{\text{(1)}}{=}~
    \vcenter{\xy
        (-3,6)="1";
        ( 3,6)="2";
        ( 3,3)*[o]=<5pt>[Fo]{~}="3";
        (-3,0)="4";
        ( 3,0)="5";
        (-3,-5)="6";
        ( 3,-5)="7";
        (-3,-8)="8";
        ( 3,-8)="9";
        ( 3,-2.5)="x";
        ( 7,1)*[o]=<9pt>[Fo]{\scriptstyle t}="y";
        ( 7,6)="z";
        "1";"8" **\dir{-};
        "3";"9" **\dir{-};
        "4";"7" **\dir{-};
        "x";"y" **\dir{-};
        "y";"z" **\dir{-};
        \ar@{-} "5";"6" |\hole
    \endxy}
    ~\overset{\text{(m)}}{=}~
    \vcenter{\xy
        (0,7)*{~}; (0,-11)*{~};
        (-3,6)="1";
        ( 3,6)="2";
        ( 3,3)*[o]=<5pt>[Fo]{~}="3";
        (-3,0)="4";
        ( 3,0)="5";
        (-3,-5)="6";
        ( 3,-5)="7";
        (-3,-10)="8";
        ( 3,-10)="9";
        ( 3,-7.5)="x";
        ( 7,-4)*[o]=<9pt>[Fo]{\scriptstyle t}="y";
        ( 7,6)="z";
        "1";"8" **\dir{-};
        "3";"9" **\dir{-};
        "4";"7" **\dir{-};
        "x";"y" **\dir{-};
        "y";"z" **\dir{-};
        \ar@{-} "5";"6" |\hole
    \endxy}
    ~\overset{\text{(b)}}{=}~
    \vcenter{\xy
        (-2,8)="1";
        ( 2,8)="2";
        ( 2,4)*[o]=<9pt>[Fo]{\scriptstyle t}="3";
        (-2,3)="4";
        ( 2,-1)="5";
        (-2,-5)="6";
        ( 2,-5)="7";
        "1";"6" **\dir{-};
        "2";"3" **\dir{-};
        "3";"7" **\dir{-};
        "4";"5" **\dir{-};
    \endxy}$ \\
%
% This is 4 tw b. 
%
    $\vcenter{\xy
        (0,10)*{~};
        (0,-7)*{~};
        ( 3,8)="1";
        (-3,8)="2";
        (-3,5)*[o]=<9pt>[Fo]{\scriptstyle t}="3";
        ( 0,2)="4";
        ( 0,-1)="5";
        ( 3,-5)="6";
        (-3,-5)="7";
        "1";"4" **\crv{(3,5)};
        "2";"3" **\dir{-};
        "3";"4" **\dir{-};
        "4";"5" **\dir{-};
        "5";"6" **\dir{-};
        "5";"7" **\dir{-};
    \endxy}
    ~\overset{\text{(b)}}{=}~
    \vcenter{\xy
        (-3,6)="1";
        ( 3,6)="2";
        (-3,3)*[o]=<9pt>[Fo]{\scriptstyle t}="3";
        (-3,0)="4";
        ( 3,0)="5";
        (-3,-4)="6";
        ( 3,-4)="7";
        (-3,-7)="8";
        ( 3,-7)="9";
        "1";"3" **\dir{-};
        "3";"8" **\dir{-};
        "2";"9" **\dir{-};
        "4";"7" **\dir{-};
        \ar@{-} "5";"6" |\hole
    \endxy}
    ~\overset{\text{(1)}}{=}~
    \vcenter{\xy
        (-1, 9)="1";
        ( 6, 9)="2";
        (-6, 5)*[o]=<5pt>[Fo]{~}="3";
        (-1, 5)*[o]=<9pt>[Fo]{\scriptstyle t}="4";
        (-6, 2)="5";
        ( 1,-2)="6";
        ( 6,-3)="7";
        (-6,-7)="8";
        ( 6,-7)="9";
        (-6,-10)="10";
        ( 6,-10)="11";
        ( 3,0)="x";
        "1";"4" **\dir{-};
        "2";"11" **\dir{-};
        "3";"10" **\dir{-};
        "6";"9" **\dir{-};
        "6";"x" **\dir{-};
        "4";"6" **\crv{(-1,0)};
        \ar@{-} "8";"7" |<(0.75)\hole
        \ar@{-} "5";"x" |<(0.6)\hole
    \endxy}
    ~\overset{\text{(m)}}{=}~
    \vcenter{\xy
        (0,9)*{~}; (0,-11)*{~};
        (-6,8)="1";
        ( 5,8)="2";
        ( 5,3)*[o]=<5pt>[Fo]{~}="3";
        (-6,-2)="4";
        ( 5,0)="5";
        (-6,-5)="6";
        ( 5,-5)="7";
        (-6,-10)="8";
        ( 5,-10)="9";
        ( 5,-7.5)="x";
        ( 2,4)*[o]=<9pt>[Fo]{\scriptstyle t}="y";
        ( 2,8)="z";
        "1";"8" **\dir{-};
        "3";"9" **\dir{-};
        "y";"x" **\crv{(1,-5)};
        "y";"z" **\dir{-};
        \ar@{-} "5";"6" |<(0.25)\hole |<(0.6)\hole
        \ar@{-} "4";"7" |<(0.75)\hole
    \endxy}
    ~\overset{\text{(b)}}{=}~
    \vcenter{\xy
        (0,9)*{~};
        (0,-6)*{~};
        (-2,8)="1";
        ( 2,8)="2";
        ( 2,5)="3";
        ( 2,-2)="4";
        (-2,-5)="5";
        ( 2,-5)="6";
        (-2,5)*[o]=<9pt>[Fo]{\scriptstyle t}="7";
        "2";"6" **\dir{-};
        "1";"7" **\dir{-};
        "7";"4" **\crv{(-2,2)}?(0.5)*{\hole}="x";
        "3";"x" **\crv{(1,5)};
        "5";"x" **\crv{(-2,-2)};
    \endxy}$

\item % This is 5 sw
    $\vcenter{\xy
        ( 0,-6)="1";
        ( 0,-3)*[o]=<9pt>[Fo]{\scriptstyle s}="2";
        ( 0,-0)="3";
        (-3, 3)*[o]=<9pt>[Fo]{\scriptstyle s}="4";
        ( 3, 3)*[o]=<9pt>[Fo]{\scriptstyle s}="5";
        (-3, 7)="6";
        ( 3, 7)="7";
        "1";"2" **\dir{-};
        "2";"3" **\dir{-};
        "3";"4" **\dir{-};
        "4";"6" **\dir{-};
        "3";"5" **\dir{-};
        "5";"7" **\dir{-};
    \endxy}
    ~\overset{\text{(4)}}{=}~
    \vcenter{\xy
        (0,7)*{~};
        (0,-10)*{~};
        (-2, 6)="1";
        ( 2, 6)="2";
        (-2, 3)*[o]=<9pt>[Fo]{\scriptstyle s}="3";
        ( 2, 3)*[o]=<9pt>[Fo]{\scriptstyle s}="4";
        (-2, 0)="5";
        ( 2,-3)="6";
        (-2,-6)*[o]=<9pt>[Fo]{\scriptstyle s}="7";
        ( 2,-6)*[o]=<5pt>[Fo]{~}="8";
        (-2,-9)="9";
        "1";"3" **\dir{-};
        "3";"7" **\dir{-};
        "7";"9" **\dir{-};
        "2";"4" **\dir{-};
        "4";"8" **\dir{-};
        "5";"6" **\dir{-};
    \endxy}
    ~\overset{\text{(3)}}{=}~
    \vcenter{\xy
        (-2, 6)="1";
        ( 2, 6)="2";
        (-2, 3)*[o]=<9pt>[Fo]{\scriptstyle s}="3";
        ( 2, 3)*[o]=<9pt>[Fo]{\scriptstyle s}="4";
        (-2, 0)="5";
        ( 2,-3)="6";
        ( 2,-6)*[o]=<5pt>[Fo]{~}="8";
        (-2,-8)="9";
        "1";"3" **\dir{-};
        "3";"9" **\dir{-};
        "2";"4" **\dir{-};
        "4";"8" **\dir{-};
        "5";"6" **\dir{-};
    \endxy}
    ~\overset{\text{(4)}}{=}~
    \vcenter{\xy
        ( 0,-3)="2";
        ( 0, 0)="3";
        (-3, 3)*[o]=<9pt>[Fo]{\scriptstyle s}="4";
        ( 3, 3)*[o]=<9pt>[Fo]{\scriptstyle s}="5";
        (-3, 7)="6";
        ( 3, 7)="7";
        "2";"3" **\dir{-};
        "3";"4" **\dir{-};
        "4";"6" **\dir{-};
        "3";"5" **\dir{-};
        "5";"7" **\dir{-};
    \endxy}$ \\
%
% This is 5 tw
%
    $\vcenter{\xy
        ( 0,-6)="1";
        ( 0,-3)*[o]=<9pt>[Fo]{\scriptstyle t}="2";
        ( 0,-0)="3";
        (-3, 3)*[o]=<9pt>[Fo]{\scriptstyle t}="4";
        ( 3, 3)*[o]=<9pt>[Fo]{\scriptstyle t}="5";
        (-3, 7)="6";
        ( 3, 7)="7";
        "1";"2" **\dir{-};
        "2";"3" **\dir{-};
        "3";"4" **\dir{-};
        "4";"6" **\dir{-};
        "3";"5" **\dir{-};
        "5";"7" **\dir{-};
    \endxy}
    ~\overset{\text{(4)}}{=}~
    \vcenter{\xy
        (0,7)*{~};
        (0,-10)*{~};
        (-2, 6)="1";
        ( 2, 6)="2";
        (-2, 3)*[o]=<9pt>[Fo]{\scriptstyle t}="3";
        ( 2, 3)*[o]=<9pt>[Fo]{\scriptstyle t}="4";
        (-2, 0)="5";
        ( 2,-3)="6";
        (-2,-6)*[o]=<9pt>[Fo]{\scriptstyle t}="7";
        ( 2,-6)*[o]=<5pt>[Fo]{~}="8";
        (-2,-9)="9";
        "1";"3" **\dir{-};
        "3";"7" **\dir{-};
        "7";"9" **\dir{-};
        "2";"4" **\dir{-};
        "4";"8" **\dir{-};
        "5";"6" **\dir{-};
    \endxy}
    ~\overset{\text{(3)}}{=}~
    \vcenter{\xy
        (-2, 6)="1";
        ( 2, 6)="2";
        (-2, 3)*[o]=<9pt>[Fo]{\scriptstyle t}="3";
        ( 2, 3)*[o]=<9pt>[Fo]{\scriptstyle t}="4";
        (-2, 0)="5";
        ( 2,-3)="6";
        ( 2,-6)*[o]=<5pt>[Fo]{~}="8";
        (-2,-8)="9";
        "1";"3" **\dir{-};
        "3";"9" **\dir{-};
        "2";"4" **\dir{-};
        "4";"8" **\dir{-};
        "5";"6" **\dir{-};
    \endxy}
    ~\overset{\text{(4)}}{=}~
    \vcenter{\xy
        ( 0,-3)="2";
        ( 0, 0)="3";
        (-3, 3)*[o]=<9pt>[Fo]{\scriptstyle t}="4";
        ( 3, 3)*[o]=<9pt>[Fo]{\scriptstyle t}="5";
        (-3, 7)="6";
        ( 3, 7)="7";
        "2";"3" **\dir{-};
        "3";"4" **\dir{-};
        "4";"6" **\dir{-};
        "3";"5" **\dir{-};
        "5";"7" **\dir{-};
    \endxy}$

\item % This is 6 a
    $\vcenter{\xy
        (0,8)="1";
        (0,5)="2";
        (-3,3)="3";
        ( 3,3)="4";
        (-3,-1)*[o]=<9pt>[Fo]{\scriptstyle s}="5";
        ( 3,-1)="6";
        ( 0,-4)="7";
        ( 0,-7)="8";
        "1";"2" **\dir{-}; 
        "2";"3" **\dir{-}; 
        "2";"4" **\dir{-}; 
        "4";"5" **\dir{-}; 
        {\ar@{-} |\hole "3";"6"}; 
        "5";"7" **\dir{-}; 
        "6";"7" **\dir{-}; 
        "8";"7" **\dir{-}; 
    \endxy}
    ~\overset{\text{(s)}}{=}~
    \vcenter{\xy
        (0,10)*{~}; (0,-13)*{~};
        ( 0, 9)="1";
        ( 0, 6)="2";
        (-4, 3)*[o]=<5pt>[Fo]{~}="3";
        (-4, 1)="4";
        (-4,-2)="5";
        (-4,-4)*[o]=<5pt>[Fo]{~}="6";
        ( 0,-9)="7";
        ( 0,-12)="8";
        (-2, 4)="x";
        ( 2,-7)="y";
        "1";"2" **\dir{-};
        "7";"8" **\dir{-};
        "3";"6" **\dir{-};
        "5";"2" **\crv{(5,3)};
        "4";"7" **\crv{(-11,-5)};
        "7";"y" **\dir{-};
        "2";"x" **\dir{-};
        \ar@{-} "x";"y" |<(0.35)\hole
    \endxy}
    ~\overset{\text{(nat)}}{=}~
    \vcenter{\xy
        ( 3, 7)="1";
        (-3, 5)*[o]=<5pt>[Fo]{~}="2";
        (-3, 2)="3";
        ( 3, 2)="4";
        (-3,-2)="5";
        ( 3,-2)="6";
        ( 3,-5)*[o]=<5pt>[Fo]{~}="7";
        (-3,-7)="8";
        "2";"8" **\dir{-};
        "1";"7" **\dir{-};
        "3";"6" **\dir{-};
        \ar@{-} "5";"4" |\hole
    \endxy}
    ~\overset{\text{(b)}}{=}~
    \vcenter{\xy
        (-2, 4)*[o]=<5pt>[Fo]{~}="1";
        ( 3, 6)="2";
        ( 0, 2)="3";
        ( 0,-1)="4";
        (-3,-5)="5";
        ( 2,-3)*[o]=<5pt>[Fo]{~}="6";
        "1";"3" **\dir{-};
        "2";"3" **\dir{-};
        "3";"4" **\dir{-};
        "4";"5" **\dir{-};
        "4";"6" **\dir{-};
    \endxy}
    ~\overset{\text{(m,c)}}{=}~
    \vcenter{\xy
        (0,6);
        (0,-5);
        **\dir{-}
    \endxy}$\\
% This is 6 b
    $\vcenter{\xy
        (0,7)="1";
        (0,5)="2";
        (3,0)*[o]=<9pt>[Fo]{\scriptstyle t}="3";
        (-3,0)="4";
        (0,-5)="5";
        (0,-7)="6";
        "1";"2" **\dir{-};
        "3";"2" **\crv{(3,2)};
        "3";"5" **\crv{(3,-2)};
        "4";"2" **\crv{(-3,2)};
        "4";"5" **\crv{(-3,-2)};
        "5";"6" **\dir{-};
    \endxy}
    ~\overset{\text{(t)}}{=}~
    \vcenter{\xy
        (0,6)*{~}; (0,-8)*{~};
        ( 3, 5)*[o]=<5pt>[Fo]{~}="1";
        (-3, 7)="2";
        (-3, 2)="3";
        ( 3, 2)="4";
        (-3,-2)="5";
        ( 3,-2)="6";
        ( 3,-5)*[o]=<5pt>[Fo]{~}="7";
        (-3,-7)="8";
        "2";"8" **\dir{-};
        "1";"7" **\dir{-};
        "3";"6" **\dir{-};
        \ar@{-} "5";"4" |\hole
    \endxy}
    ~\overset{\text{(b)}}{=}~
    \vcenter{\xy
        (-3, 6)="1";
        ( 2, 4)*[o]=<5pt>[Fo]{~}="2";
        ( 0, 2)="3";
        ( 0,-1)="4";
        (-3,-5)="5";
        ( 2,-3)*[o]=<5pt>[Fo]{~}="6";
        "1";"3" **\dir{-};
        "2";"3" **\dir{-};
        "3";"4" **\dir{-};
        "4";"5" **\dir{-};
        "4";"6" **\dir{-};
    \endxy}
    ~\overset{\text{(m,c)}}{=}~
    \vcenter{\xy
        (0,6);
        (0,-5);
        **\dir{-}
    \endxy}$

\item % This is 7s
    $\vcenter{\xy
        (0,7)="t";
        (0,3)*[o]=<9pt>[Fo]{\scriptstyle s}="m1";
        (0,-3)*[o]=<9pt>[Fo]{\scriptstyle s}="m2";
        (0,-7)="b";
        "t";"m1" **\dir{-};
        "m1";"m2" **\dir{-};
        "m2";"b" **\dir{-};
    \endxy}
    ~\overset{\text{($s$)}}{=}~
    \vcenter{\xy
        (0,10)*{~}; (0,-7)*{~};
        ( 3, 9)="1";
        ( 0, 6)*[o]=<5pt>[Fo]{~}="2";
        ( 0, 4)="3";
        ( 0, 1)="4";
        ( 0,-1)*[o]=<5pt>[Fo]{~}="5";
        (-3,-2)*[o]=<9pt>[Fo]{\scriptstyle s}="6";
        (-3,-6)="7";
        "1";"4" **\crv{(3,4)};
        "2";"5" **\dir{-};
        "6";"3" **\crv{(-3,1)};
        "6";"7" **\dir{-};
    \endxy}
    ~\overset{\text{(2)}}{=}~
    \vcenter{\xy
        ( 3, 9)="1";
        ( 0, 6)*[o]=<5pt>[Fo]{~}="2";
        ( 0, 4)="3";
        ( 0, 1)="4";
        ( 0,-1)*[o]=<5pt>[Fo]{~}="5";
        (-3,-4)="6";
        "1";"4" **\crv{(3,4)};
        "2";"5" **\dir{-};
        "6";"3" **\crv{(-3,1)};
    \endxy}
    ~\overset{\text{($s$)}}{=}~
    \vcenter{\xy
        (0,4.5)="1";
        (0,0)*[o]=<9pt>[Fo]{\scriptstyle s}="2";
        (0,-4.5)="3";
        "1";"2" **\dir{-};
        "2";"3" **\dir{-};
    \endxy}$ \\
%
% This is 7t
%
    $\vcenter{\xy
        (0,7)="t";
        (0,3)*[o]=<9pt>[Fo]{\scriptstyle t}="m1";
        (0,-3)*[o]=<9pt>[Fo]{\scriptstyle t}="m2";
        (0,-7)="b";
        "t";"m1" **\dir{-};
        "m1";"m2" **\dir{-};
        "m2";"b" **\dir{-};
    \endxy}
    ~\overset{\text{($t$)}}{=}~
    \vcenter{\xy
        (0,10)*{~}; (0,-7)*{~};
        (-3, 9)="1";
        ( 0, 6)*[o]=<5pt>[Fo]{~}="2";
        ( 0, 4)="3";
        ( 0, 1)="4";
        ( 0,-1)*[o]=<5pt>[Fo]{~}="5";
        (-3,-2)*[o]=<9pt>[Fo]{\scriptstyle t}="6";
        (-3,-6)="7";
        "2";"5" **\dir{-};
        "6";"7" **\dir{-};
        \ar@{-}@/^3pt/ "6";"3" |<(0.55)\hole
        \ar@{-}@/_3pt/ "1";"4"
    \endxy}
    ~\overset{\text{(2)}}{=}~
    \vcenter{\xy
        (-3, 9)="1";
        ( 0, 6)*[o]=<5pt>[Fo]{~}="2";
        ( 0, 4)="3";
        ( 0, 1)="4";
        ( 0,-1)*[o]=<5pt>[Fo]{~}="5";
        (-3,-4)="6";
        "2";"5" **\dir{-};
        \ar@{-}@/_3pt/ "1";"4"
        \ar@{-}@/^3pt/ "6";"3" |<(0.72)\hole
    \endxy}
    ~\overset{\text{($t$)}}{=}~
    \vcenter{\xy
        (0,4.5)="1";
        (0,0)*[o]=<9pt>[Fo]{\scriptstyle t}="2";
        (0,-4.5)="3";
        "1";"2" **\dir{-};
        "2";"3" **\dir{-};
    \endxy}$
    
%
% Now on to the interactions.
%

\item % This is 8s
    $\vcenter{\xy
        (0,7)="t";
        (0,3)*[o]=<9pt>[Fo]{\scriptstyle s}="m1";
        (0,-3)*[o]=<9pt>[Fo]{\scriptstyle t}="m2";
        (0,-7)="b";
        "t";"m1" **\dir{-};
        "m1";"m2" **\dir{-};
        "m2";"b" **\dir{-};
    \endxy}
    ~\overset{\text{($s$)}}{=}~
    \vcenter{\xy
        (0,10)*{~}; (0,-7)*{~};
        ( 3, 9)="1";
        ( 0, 6)*[o]=<5pt>[Fo]{~}="2";
        ( 0, 4)="3";
        ( 0, 1)="4";
        ( 0,-1)*[o]=<5pt>[Fo]{~}="5";
        (-3,-2)*[o]=<9pt>[Fo]{\scriptstyle t}="6";
        (-3,-6)="7";
        "1";"4" **\crv{(3,4)};
        "2";"5" **\dir{-};
        "6";"3" **\crv{(-3,1)};
        "6";"7" **\dir{-};
    \endxy}
    ~\overset{\text{(2)}}{=}~
    \vcenter{\xy
        ( 3, 9)="1";
        ( 0, 6)*[o]=<5pt>[Fo]{~}="2";
        ( 0, 4)="3";
        ( 0, 1)="4";
        ( 0,-1)*[o]=<5pt>[Fo]{~}="5";
        (-3,-4)="6";
        "1";"4" **\crv{(3,4)};
        "2";"5" **\dir{-};
        "6";"3" **\crv{(-3,1)};
    \endxy}
    ~\overset{\text{($s$)}}{=}~
    \vcenter{\xy
        (0,4.5)="1";
        (0,0)*[o]=<9pt>[Fo]{\scriptstyle s}="2";
        (0,-4.5)="3";
        "1";"2" **\dir{-};
        "2";"3" **\dir{-};
    \endxy}$ \\
%
% This is 8t
%
    $\vcenter{\xy
        (0,7)="t";
        (0,3)*[o]=<9pt>[Fo]{\scriptstyle t}="m1";
        (0,-3)*[o]=<9pt>[Fo]{\scriptstyle s}="m2";
        (0,-7)="b";
        "t";"m1" **\dir{-};
        "m1";"m2" **\dir{-};
        "m2";"b" **\dir{-};
    \endxy}
    ~\overset{\text{($t$)}}{=}~
    \vcenter{\xy
        (0,10)*{~}; (0,-7)*{~};
        (-3, 9)="1";
        ( 0, 6)*[o]=<5pt>[Fo]{~}="2";
        ( 0, 4)="3";
        ( 0, 1)="4";
        ( 0,-1)*[o]=<5pt>[Fo]{~}="5";
        (-3,-2)*[o]=<9pt>[Fo]{\scriptstyle s}="6";
        (-3,-6)="7";
        "2";"5" **\dir{-};
        "6";"7" **\dir{-};
        \ar@{-}@/^3pt/ "6";"3" |<(0.55)\hole
        \ar@{-}@/_3pt/ "1";"4"
    \endxy}
    ~\overset{\text{(2)}}{=}~
    \vcenter{\xy
        (-3, 9)="1";
        ( 0, 6)*[o]=<5pt>[Fo]{~}="2";
        ( 0, 4)="3";
        ( 0, 1)="4";
        ( 0,-1)*[o]=<5pt>[Fo]{~}="5";
        (-3,-4)="6";
        "2";"5" **\dir{-};
        \ar@{-}@/_3pt/ "1";"4"
        \ar@{-}@/^3pt/ "6";"3" |<(0.72)\hole
    \endxy}
    ~\overset{\text{($t$)}}{=}~
    \vcenter{\xy
        (0,4.5)="1";
        (0,0)*[o]=<9pt>[Fo]{\scriptstyle t}="2";
        (0,-4.5)="3";
        "1";"2" **\dir{-};
        "2";"3" **\dir{-};
    \endxy}$

\item % This is 9a
    $\vcenter{\xy
        (0,6)="t";
        (0,3)*[o]=<9pt>[Fo]{\scriptstyle t}="f";
        (0, 0)="m";
        (-3,-3)*[o]=<9pt>[Fo]{\scriptstyle s}="g";
        (-3,-6)="l";
        ( 3,-6)="r";
        "t";"f" **\dir{-};
        "f";"m" **\dir{-};
        "r";"m" **\crv{(3,-3)};
        "m";"g" **\dir{-};
        "g";"l" **\dir{-};
    \endxy}
    ~\overset{\text{(3)}}{=}~
    \vcenter{\xy
        (0,7)*{~};
        (0,-12)*{~};
        (0,6)="t";
        (0,3)*[o]=<9pt>[Fo]{\scriptstyle t}="f";
        (0, 0)="m";
        (-3,-3)*[o]=<9pt>[Fo]{\scriptstyle t}="g";
        (-3,-8)*[o]=<9pt>[Fo]{\scriptstyle s}="a";
        (-3,-11)="l";
        ( 3,-11)="r";
        "t";"f" **\dir{-};
        "f";"m" **\dir{-};
        "r";"m" **\crv{(3,-5)};
        "m";"g" **\dir{-};
        "g";"a" **\dir{-};
        "a";"l" **\dir{-};
    \endxy}
    ~\overset{\text{(8)}}{=}~
    \vcenter{\xy
        (0,6)="t";
        (0,3)*[o]=<9pt>[Fo]{\scriptstyle t}="f";
        (0,-0)="m";
        (-3,-3)*[o]=<9pt>[Fo]{\scriptstyle t}="g";
        (-3,-6)="l";
        ( 3,-6)="r";
        "t";"f" **\dir{-};
        "f";"m" **\dir{-};
        "r";"m" **\crv{(3,-3)};
        "m";"g" **\dir{-};
        "g";"l" **\dir{-};
    \endxy}
    ~\overset{\text{(3)}}{=}~
    \vcenter{\xy
        (0,5)="t";
        (0,2)*[o]=<9pt>[Fo]{\scriptstyle t}="f";
        (0,-1)="m";
        (-3,-5)="l";
        ( 3,-5)="r";
        "t";"f" **\dir{-};
        "f";"m" **\dir{-};
        "m";"l" **\dir{-};
        "m";"r" **\dir{-};
    \endxy}$ \\
%
% This is 9b
%
    $\vcenter{\xy
        (0,6)="t";
        (0,3)*[o]=<9pt>[Fo]{\scriptstyle s}="f";
        (0,-0)="m";
        (-3,-3)*[o]=<9pt>[Fo]{\scriptstyle t}="g";
        (-3,-6)="l";
        ( 3,-6)="r";
        "t";"f" **\dir{-};
        "f";"m" **\dir{-};
        "r";"m" **\crv{(3,-3)};
        "m";"g" **\dir{-};
        "g";"l" **\dir{-};
    \endxy}
    ~\overset{\text{(3)}}{=}~
    \vcenter{\xy
        (0,7)*{~};
        (0,-12)*{~};
        (0,6)="t";
        (0,3)*[o]=<9pt>[Fo]{\scriptstyle s}="f";
        (0, 0)="m";
        (-3,-3)*[o]=<9pt>[Fo]{\scriptstyle s}="g";
        (-3,-8)*[o]=<9pt>[Fo]{\scriptstyle t}="a";
        (-3,-11)="l";
        ( 3,-11)="r";
        "t";"f" **\dir{-};
        "f";"m" **\dir{-};
        "r";"m" **\crv{(3,-5)};
        "m";"g" **\dir{-};
        "g";"a" **\dir{-};
        "a";"l" **\dir{-};
    \endxy}
    ~\overset{\text{(8)}}{=}~
    \vcenter{\xy
        (0,6)="t";
        (0,3)*[o]=<9pt>[Fo]{\scriptstyle s}="f";
        (0,-0)="m";
        (-3,-3)*[o]=<9pt>[Fo]{\scriptstyle s}="g";
        (-3,-6)="l";
        ( 3,-6)="r";
        "t";"f" **\dir{-};
        "f";"m" **\dir{-};
        "r";"m" **\crv{(3,-3)};
        "m";"g" **\dir{-};
        "g";"l" **\dir{-};
    \endxy}
    ~\overset{\text{(3)}}{=}~
    \vcenter{\xy
        (0,5)="t";
        (0,2)*[o]=<9pt>[Fo]{\scriptstyle s}="f";
        (0,-1)="m";
        (-3,-5)="l";
        ( 3,-5)="r";
        "t";"f" **\dir{-};
        "f";"m" **\dir{-};
        "m";"l" **\dir{-};
        "m";"r" **\dir{-};
    \endxy}$

\item % This is 10
    $\vcenter{\xy
        ( 0, 4)="2";
        ( 0, 0)="3";
        (-3,-3)*[o]=<9pt>[Fo]{\scriptstyle t}="4";
        ( 3,-3)*[o]=<9pt>[Fo]{\scriptstyle s}="5";
        (-3,-7)="6";
        ( 3,-7)="7";
        "2";"3" **\dir{-};
        "3";"4" **\dir{-};
        "4";"6" **\dir{-};
        "3";"5" **\dir{-};
        "5";"7" **\dir{-};
    \endxy}
    ~\overset{\text{($s,t$)}}{=}~
    \vcenter{\xy
        (0,11)="1";
        (0,8)="2";
        (-3, 4)*[o]=<5pt>[Fo]{~}="3";
        (-3, 2)="4";
        (-3,-2)="5";
        (-3,-4)*[o]=<5pt>[Fo]{~}="6";
        ( 3, 4)*[o]=<5pt>[Fo]{~}="7";
        ( 3, 2)="8";
        ( 3,-2)="9";
        ( 3,-4)*[o]=<5pt>[Fo]{~}="10";
        (-7,-7)="11";
        ( 0,-7)="12";
        "3";"6" **\dir{-};
        "7";"10" **\dir{-};
        "1";"2" **\dir{-};
        \ar@{-}@/^3pt/ "12";"8"
        \ar@{-}@/^3pt/ "11";"4" |<(0.73){\hole} 
        \ar@{-}@/^13pt/ "5";"2"
        \ar@{-}@/_13pt/ "9";"2"
    \endxy}
    ~\overset{\text{(nat)}}{=}~
    \vcenter{\xy
        ( 0,11)="1";
        ( 0,8)="2";
        (-3,5)="x";
        ( 3,5)="y";
        (-6, 4)*[o]=<5pt>[Fo]{~}="3";
        (-6, 2)="4";
        (-6,-2)="5";
        (-6,-4)*[o]=<5pt>[Fo]{~}="6";
        ( 6, 4)*[o]=<5pt>[Fo]{~}="7";
        ( 6, 2)="8";
        ( 6,-2)="9";
        ( 6,-4)*[o]=<5pt>[Fo]{~}="10";
        (-4,-12)="11";
        ( 4,-12)="12";
        "3";"6" **\dir{-};
        "7";"10" **\dir{-};
        "1";"2" **\dir{-};
        "2";"x" **\dir{-};
        "2";"y" **\dir{-};
        "5";"y" **\dir{-};
        "4";"12" **\crv{(-10,0)&(-10,-8)&(2,-8)};
        \ar@{-}@/^1pt/ "11";"8" |<(0.27){\hole}  |<(0.8){\hole} 
        \ar@{-} "x";"9" |<(0.35){\hole}
    \endxy}
    ~\overset{\text{(v)}}{=}~
    \vcenter{\xy
        ( 0,7)="1";
        ( 0,2)="2";
        (-4, 4)*[o]=<5pt>[Fo]{~}="3";
        (-4, 2)="4";
        (-4,-2)="5";
        (-4,-4)*[o]=<5pt>[Fo]{~}="6";
        ( 4, 4)*[o]=<5pt>[Fo]{~}="7";
        ( 4, 2)="8";
        ( 4,-2)="9";
        ( 4,-4)*[o]=<5pt>[Fo]{~}="10";
        (-4,-12)="11";
        ( 4,-12)="12";
        "3";"6" **\dir{-};
        "7";"10" **\dir{-};
        "1";"2" **\dir{-};
        "5";"2" **\dir{-};
        "9";"2" **\dir{-};
        "4";"12" **\crv{(-8,0)&(-8,-8)&(2,-8)};
        \ar@{-}@/^1pt/ "11";"8" |<(0.3){\hole}  |<(0.8){\hole} 
    \endxy}
    ~\overset{\text{($s,t$)}}{=}~
    \vcenter{\xy
        (0,5)*{~};
        (0,-8)*{~};
        ( 0, 3)="2";
        ( 0, 1)="3";
        (-4,-2)*[o]=<9pt>[Fo]{\scriptstyle s}="4";
        ( 4,-2)*[o]=<9pt>[Fo]{\scriptstyle t}="5";
        (-4,-7)="6";
        ( 4,-7)="7";
        "2";"3" **\dir{-};
        "3";"4" **\dir{-};
        "3";"5" **\dir{-};
        "4";"7" **\dir{-}?(0.5)*{\hole}="x";
        "5";"x" **\dir{-};
        "x";"6" **\dir{-};
    \endxy}$

\item % This is 11
    $\vcenter{\xy
        (-2,8)="1";
        ( 2,8)="2";
        (-2,4)="3";
        ( 2,0)="4";
        (-2,-1)*[o]=<9pt>[Fo]{\scriptstyle t}="5";
        (-2,-5)="6";
        ( 2,-5)="7";
        "1";"5" **\dir{-};
        "5";"6" **\dir{-};
        "2";"7" **\dir{-};
        "3";"4" **\dir{-};
    \endxy}
    ~\overset{\text{(m)}}{=}~
    \vcenter{\xy
        (0,12)*{~};
        (0,-6)*{~};
        (-2,11)="1";
        ( 3,11)="2";
        (-2,6)="x";
        ( 0,8)*[o]=<5pt>[Fo]{~}="y";
        (-2,4)="3";
        ( 3,0)="4";
        (-2,-1)*[o]=<9pt>[Fo]{\scriptstyle t}="5";
        (-2,-5)="6";
        ( 3,-5)="7";
        "1";"5" **\dir{-};
        "5";"6" **\dir{-};
        "2";"7" **\dir{-};
        "3";"4" **\dir{-};
        "x";"y" **\dir{-};
    \endxy}
    ~\overset{\text{(4)}}{=}~
    \vcenter{\xy
        (0,10)*{~};
        (0,-7)*{~};
        (-3, 9)="1";
        ( 3, 9)="2";
        ( 1, 6)*[o]=<5pt>[Fo]{~}="3";
        ( 1, 4)="4";
        ( 1, 1)="5";
        (-3,-2)*[o]=<9pt>[Fo]{\scriptstyle t}="6";
        ( 3,-2)="7";
        (-3,-6)="8";
        ( 3,-6)="9";
        "1";"5" **\crv{(-3,4)};
        "3";"5" **\dir{-};
        "2";"9" **\dir{-};
        "5";"7" **\dir{-};
        "6";"8" **\dir{-};
        \ar@{-}@/^3pt/ "6";"4" |<(0.56)\hole
    \endxy}
    ~\overset{\text{(m)}}{=}~
    \vcenter{\xy
        (-3, 9)="1";
        ( 3, 9)="2";
        ( 1, 6)*[o]=<5pt>[Fo]{~}="3";
        ( 1, 4)="4";
        ( 3, 1)="5";
        (-3,-2)*[o]=<9pt>[Fo]{\scriptstyle t}="6";
        ( 3,-2)="7";
        (-3,-6)="8";
        ( 3,-6)="9";
        "1";"7" **\crv{(-3,2)};
        "3";"4" **\dir{-};
        "4";"5" **\dir{-};
        "2";"9" **\dir{-};
        "6";"8" **\dir{-};
        \ar@{-}@/^3pt/ "6";"4" |<(0.53)\hole
    \endxy}
    ~\overset{\text{(2)}}{=}~
    \vcenter{\xy
        (-3, 9)="1";
        ( 3, 9)="2";
        ( 1, 6)*[o]=<5pt>[Fo]{~}="3";
        ( 1, 4)="4";
        ( 3, 1)="5";
        ( 3,-2)="7";
        (-3,-6)="8";
        ( 3,-6)="9";
        "1";"7" **\crv{(-3,2)};
        "3";"4" **\dir{-};
        "4";"5" **\dir{-};
        "2";"9" **\dir{-};
        \ar@{-}@/^3pt/ "8";"4" |<(0.7)\hole
    \endxy}
    ~\overset{\text{(2)}}{=}~
    \vcenter{\xy
        (-3, 9)="1";
        ( 3, 9)="2";
        ( 1, 6)*[o]=<5pt>[Fo]{~}="3";
        ( 1, 4)="4";
        ( 3, 1)="5";
        (-3,-2)*[o]=<9pt>[Fo]{\scriptstyle s}="6";
        ( 3,-2)="7";
        (-3,-6)="8";
        ( 3,-6)="9";
        "1";"7" **\crv{(-3,2)};
        "3";"4" **\dir{-};
        "4";"5" **\dir{-};
        "2";"9" **\dir{-};
        "6";"8" **\dir{-};
        \ar@{-}@/^3pt/ "6";"4" |<(0.53)\hole
    \endxy}
    ~\overset{\text{(4)}}{=}~
    \vcenter{\xy
        (-2,8)="1";
        ( 2,8)="2";
        ( 2,5)="3";
        ( 2,-2)="4";
        (-2,-5)="5";
        ( 2,-5)="6";
        (-2,-2)*[o]=<9pt>[Fo]{\scriptstyle s}="7";
        "2";"6" **\dir{-};
        "1";"4" **\crv{(-2,2)}?(0.5)*{\hole}="x";
        "3";"x" **\crv{(1,5)};
        "5";"7" **\dir{-};
        "7";"x" **\crv{(-2,0)};
    \endxy}$

\item % This is 12
    $\vcenter{\xy
        (0,8)="1";
        (0,3)*[o]=<9pt>[Fo]{\scriptstyle r}="2";
        (0,-3)*[o]=<9pt>[Fo]{\scriptstyle s}="3";
        (0,-8)="4";
        "1";"2" **\dir{-};
        "2";"3" **\dir{-};
        "3";"4" **\dir{-};
    \endxy}
    ~\overset{\text{(r,s)}}{=}~
    \vcenter{\xy
        (4,13)="1";
        (0,10)*[o]=<5pt>[Fo]{~}="2";
        (0, 8)="3";
        (0, 4)="4";
        (0, 2)*[o]=<5pt>[Fo]{~}="5";
        (0,-2)*[o]=<5pt>[Fo]{~}="6";
        (0,-4)="7";
        (0,-8)="8";
        (0,-10)*[o]=<5pt>[Fo]{~}="9";
        (-4,-13)="10";
        "3";"8" **\crv{(5,5)&(5,-5)};
        "10";"7" **\crv{(-4,-9)};
        "2";"5" **\dir{-};
        "6";"9" **\dir{-};
        \ar@{-}@/^3pt/ "1";"4" |<(0.7){\hole} 
    \endxy}
    ~\overset{\text{(nat)}}{=}~
    \vcenter{\xy
        (-4,5)="1";
        (-8,-9)="2";
        (-4,8)="bl";
        (0,5)*[o]=<5pt>[Fo]{~}="bc";
        (4,8)="br";
        (-4,-5)*[o]=<5pt>[Fo]{~}="tl";
        (-4,-2)="ml";
        (0,2)="mc";
        (4,-2)="mr";
        (4,-5)*[o]=<5pt>[Fo]{~}="tr";
        "bl";"tl" **\dir{-};
        "br";"tr" **\dir{-};
        "bc";"mc" **\dir{-};
        "mc";"ml" **\crv{(4,0)}?(0.6)*{\hole}="x";
        "mc";"x" **\crv{(-3,0.5)};
        "mr";"x" **\crv{(1,-1.5)};
        "2";"1" **\crv{(-8,0)};
    \endxy}
    ~\overset{\text{(v)}}{=}~
    \vcenter{\xy
        (-4,5)*[o]=<5pt>[Fo]{~}="1";
        ( 0,5)*[o]=<5pt>[Fo]{~}="2";
        ( 4,8)="3";
        ( 0,0)="4";
        (0,-3)*[o]=<5pt>[Fo]{~}="5";
        (-4,-6)="7";
        "2";"5" **\dir{-};
        "3";"4" **\crv{(4,4)};
        "1";"4" **\dir{-}?="6";
        "7";"6" **\crv{(-4,-2)};
    \endxy}
    ~\overset{\text{(c)}}{=}~
    \vcenter{\xy
        ( 3, 7)="1";
        ( 0, 4)*[o]=<5pt>[Fo]{~}="2";
        ( 0, 2)="3";
        ( 0,-2)="4";
        ( 0,-4)*[o]=<5pt>[Fo]{~}="5";
        (-3,-7)="6";
        "2";"5" **\dir{-};
        "1";"4" **\crv{(3,2)};
        "6";"3" **\crv{(-3,-2)};
    \endxy}
    ~\overset{\text{(s)}}{=}~
    \vcenter{\xy
        (0, 6)="1";
        (0, 0)*[o]=<9pt>[Fo]{\scriptstyle s}="2";
        (0,-6)="3";
        "1";"2" **\dir{-};
        "2";"3" **\dir{-};
    \endxy}$

\item % This is 13
    $\vcenter{\xy
        ( 0, 4)="2";
        ( 0, 0)="3";
        (-3,-3)*[o]=<9pt>[Fo]{\scriptstyle t}="4";
        ( 3,-3)*[o]=<9pt>[Fo]{\scriptstyle r}="5";
        (-3,-7)="6";
        ( 3,-7)="7";
        "2";"3" **\dir{-};
        "3";"4" **\dir{-};
        "4";"6" **\dir{-};
        "3";"5" **\dir{-};
        "5";"7" **\dir{-};
    \endxy}
    ~\overset{\text{($t,r$)}}{=}~
    \vcenter{\xy
        (0,11)="1";
        (0,8)="2";
        (-3, 4)*[o]=<5pt>[Fo]{~}="3";
        (-3, 2)="4";
        (-3,-2)="5";
        (-3,-4)*[o]=<5pt>[Fo]{~}="6";
        ( 3, 4)*[o]=<5pt>[Fo]{~}="7";
        ( 3, 2)="8";
        ( 3,-2)="9";
        ( 3,-4)*[o]=<5pt>[Fo]{~}="10";
        (-7,-7)="11";
        ( 7,-7)="12";
        "3";"6" **\dir{-};
        "7";"10" **\dir{-};
        "1";"2" **\dir{-};
        \ar@{-}@/_3pt/ "12";"8"
        \ar@{-}@/^3pt/ "11";"4" |<(0.73){\hole} 
        \ar@{-}@/^13pt/ "5";"2"
        \ar@{-}@/_13pt/ "9";"2" |<(0.17){\hole} 
    \endxy}
    ~\overset{\text{(nat)}}{=}~
    \vcenter{\xy
        ( 0,11)="1";
        ( 0,8)="2";
        (-3,5)="x";
        ( 3,5)="y";
        (-6, 4)*[o]=<5pt>[Fo]{~}="3";
        (-6, 2)="4";
        (-6,-2)="5";
        (-6,-4)*[o]=<5pt>[Fo]{~}="6";
        ( 6, 4)*[o]=<5pt>[Fo]{~}="7";
        ( 6, 2)="8";
        ( 6,-2)="9";
        ( 6,-4)*[o]=<5pt>[Fo]{~}="10";
        (-4,-9)="11";
        ( 4,-9)="12";
        "3";"6" **\dir{-};
        "7";"10" **\dir{-};
        "1";"2" **\dir{-};
        "2";"x" **\dir{-};
        "2";"y" **\dir{-};
        \ar@{-}@/^1pt/ "4";"12"
        \ar@{-}@/^1pt/ "11";"8" |<(0.4){\hole}  |<(0.8){\hole} 
        \ar@{-} "x";"9" |<(0.35){\hole}
        \ar@{-} "y";"5" |<(0.7){\hole} 
    \endxy}
    ~\overset{\text{(v)}}{=}~
    \vcenter{\xy
        ( 0,8)="1";
        ( 0,4)="2";
        (-6, 4)*[o]=<5pt>[Fo]{~}="3";
        (-6, 2)="4";
        (-6,-2)="5";
        (-6,-4)*[o]=<5pt>[Fo]{~}="6";
        ( 6, 4)*[o]=<5pt>[Fo]{~}="7";
        ( 6, 2)="8";
        ( 6,-2)="9";
        ( 6,-4)*[o]=<5pt>[Fo]{~}="10";
        (-4,-9)="11";
        ( 4,-9)="12";
        "3";"6" **\dir{-};
        "7";"10" **\dir{-};
        "1";"2" **\dir{-};
        "2";"9" **\dir{-};
        \ar@{-}@/^1pt/ "4";"12"
        \ar@{-}@/^1pt/ "11";"8" |<(0.4){\hole}  |<(0.8){\hole} 
        \ar@{-} "2";"5" |<(0.6){\hole} 
    \endxy}
    ~\overset{\text{($t,r$)}}{=}~
    \vcenter{\xy
        ( 0, 4)="2";
        ( 0, 1)="3";
        (-3.5,-2)*[o]=<9pt>[Fo]{\scriptstyle r}="4";
        ( 3.5,-2)*[o]=<9pt>[Fo]{\scriptstyle t}="5";
        (-3,-8)="6";
        ( 3,-8)="7";
        "2";"3" **\dir{-};
        "3";"4" **\dir{-};
        "3";"5" **\dir{-};
        "4";"7" **\dir{-}?(0.53)*{\hole}="x";
        "5";"x" **\dir{-};
        "x";"6" **\dir{-};
    \endxy}$

\item % This is 14
    $\vcenter{\xy
        (-3, 7)="1";
        ( 3, 7)="2";
        ( 3, 3)*[o]=<9pt>[Fo]{\scriptstyle r}="3";
        ( 0, 0)="4";
        ( 0,-3)*[o]=<9pt>[Fo]{\scriptstyle s}="5";
         (0,-7)="6";
        "1";"4" **\crv{(-3,3)};
        "2";"3" **\dir{-};
        "3";"4" **\dir{-};
        "4";"5" **\dir{-};
        "5";"6" **\dir{-};
    \endxy}
    ~\overset{\text{(3)}}{=}~
    \vcenter{\xy
        (-3, 7)="1";
        ( 3, 7)="2";
        ( 3, 3)*[o]=<9pt>[Fo]{\scriptstyle r}="3";
        ( 3,-2)*[o]=<9pt>[Fo]{\scriptstyle s}="4";
        ( 0,-5)="5";
        ( 0,-8)*[o]=<9pt>[Fo]{\scriptstyle s}="6";
         (0,-13)="7";
        "1";"5" **\crv{(-3,-2)};
        "2";"3" **\dir{-};
        "3";"4" **\dir{-};
        "4";"5" **\dir{-};
        "5";"6" **\dir{-};
        "6";"7" **\dir{-};
    \endxy}
    ~\overset{\text{(12)}}{=}~
    \vcenter{\xy
        (-3, 7)="1";
        ( 3, 7)="2";
        ( 3, 3)*[o]=<9pt>[Fo]{\scriptstyle s}="3";
        ( 0, 0)="4";
        ( 0,-3)*[o]=<9pt>[Fo]{\scriptstyle s}="5";
         (0,-7)="6";
        "1";"4" **\crv{(-3,3)};
        "2";"3" **\dir{-};
        "3";"4" **\dir{-};
        "4";"5" **\dir{-};
        "5";"6" **\dir{-};
    \endxy}
    ~\overset{\text{(3)}}{=}~
    \vcenter{\xy
        (-3, 7)="1";
        ( 3, 7)="2";
        ( 0, 3)="3";
        ( 0, 0)*[o]=<9pt>[Fo]{\scriptstyle s}="4";
        ( 0,-4)="5";
        "1";"3" **\dir{-};
        "2";"3" **\dir{-};
        "3";"4" **\dir{-};
        "4";"5" **\dir{-};
    \endxy}$ \\
    $\vcenter{\xy
        (-3, 7)="1";
        ( 3, 7)="2";
        ( 3, 3)*[o]=<9pt>[Fo]{\scriptstyle s}="3";
        ( 0, 0)="4";
        ( 0,-3)*[o]=<9pt>[Fo]{\scriptstyle r}="5";
         (0,-7)="6";
        "1";"4" **\crv{(-3,3)};
        "2";"3" **\dir{-};
        "3";"4" **\dir{-};
        "4";"5" **\dir{-};
        "5";"6" **\dir{-};
    \endxy}
    ~\overset{\text{(3)}}{=}~
    \vcenter{\xy
        (-3, 7)="1";
        ( 3, 7)="2";
        ( 3, 3)*[o]=<9pt>[Fo]{\scriptstyle s}="3";
        ( 3,-2)*[o]=<9pt>[Fo]{\scriptstyle r}="4";
        ( 0,-5)="5";
        ( 0,-8)*[o]=<9pt>[Fo]{\scriptstyle r}="6";
         (0,-13)="7";
        "1";"5" **\crv{(-3,-2)};
        "2";"3" **\dir{-};
        "3";"4" **\dir{-};
        "4";"5" **\dir{-};
        "5";"6" **\dir{-};
        "6";"7" **\dir{-};
    \endxy}
    ~\overset{\text{(12)}}{=}~
    \vcenter{\xy
        (-3, 7)="1";
        ( 3, 7)="2";
        ( 3, 3)*[o]=<9pt>[Fo]{\scriptstyle r}="3";
        ( 0, 0)="4";
        ( 0,-3)*[o]=<9pt>[Fo]{\scriptstyle r}="5";
         (0,-7)="6";
        "1";"4" **\crv{(-3,3)};
        "2";"3" **\dir{-};
        "3";"4" **\dir{-};
        "4";"5" **\dir{-};
        "5";"6" **\dir{-};
    \endxy}
    ~\overset{\text{(3)}}{=}~
    \vcenter{\xy
        (-3, 7)="1";
        ( 3, 7)="2";
        ( 0, 3)="3";
        ( 0, 0)*[o]=<9pt>[Fo]{\scriptstyle r}="4";
        ( 0,-4)="5";
        "1";"3" **\dir{-};
        "2";"3" **\dir{-};
        "3";"4" **\dir{-};
        "4";"5" **\dir{-};
    \endxy}$ \\
\end{enumerate}

\newpage
%=========================================================================%

%=========================================================================%

\end{document}